\documentclass[a4paper]{article}
\pdfoutput=1
\usepackage{amsmath}
\usepackage{amssymb}
\usepackage{layout}
\usepackage{ifthen}
\usepackage{url}
\usepackage{pgf}
\usepackage{tikz}
\usepackage[retainorgcmds]{IEEEtrantools}
\usepackage{lscape}
\usetikzlibrary{calc,intersections,through,backgrounds}
\author{Roger Sewell\\
\href{mailto:roger.sewell@cantab.net}{\scriptsize{roger.sewell@cantab.net}}
}
\usepackage{graphicx}
\usepackage{xr}
\usepackage{float}
\usepackage[pdftex,colorlinks=true,linkcolor=blue]{hyperref}

\newboolean{isfalse}
\setboolean{isfalse}{false}

\title{Hypothesis testing and confidence sets: why Bayesian not
  frequentist, and how to set a prior with a regulatory
  authority\footnote{Seventh version submitted to arxiv.org; RFS version
3.79.1.1
; added appendix \ref{hypoconf} on exact relationship between complete
  families of critical regions and confidence sets, and anecdote on
  the adverse frequentist effects of hard work.}
  }

\DeclareMathOperator*{\mean}{mean}
\DeclareMathOperator*{\E}{E}

\begin{document}
\maketitle

Acknowledgements: I am grateful to: Lalita Ramakrishan for eliciting
the comments which provoked me to start writing; Pawe{\l} Piwek for
constructive criticism and inspiring sections
\ref{admissiblesolutions} and \ref{admmeth}; Filip Tokarski for
inspiring section \ref{cannotpostfix}; Mark Owen for detailed comments
and in particular for telling me to be less combative; Paul Edelstein
for encouragement and for inspiring section \ref{example3}; Mark Troll
and Luke Sewell for encouragement and for reflecting the viewpoint of
scientists less inclined to abstraction in pushing for additional
examples; and Brendan Kelly for suggesting possible further points to
include. Nonetheless responsibility for the text and any errors
etc. remains my own -- reports of typos and errata welcome at the
above email address.

\begin{center}
\textbf{Abstract}
\end{center}

We aim to marshall in one place the reasons why Bayesian methods are
preferable to both frequentist \mbox{hypothesis} tests (FHT) and
frequentist confidence intervals/sets (FCS), and how while FHT and FCS
may be considered a reasonable approximation in some circumstances,
they cannot be relied upon to always give correct answers. We argue
that therefore both FHT and FCS should be phased out.

We first specify what we mean by an inference problem. We define an
admissible solution as one which gives a probability distribution on
the desired unknowns consistent with the likelihood and the observed
data, and an admissible method as one which for any inference problem
yields an admissible solution, noting that for any prior the Bayesian
method is admissible. We then consider seven weaker common-sense
criteria satisfied by any admissible method and find that frequentist
hypothesis testing violates all of them, even on problems with
one-dimensional data possessing uniformly optimal families of critical
regions. As an \textit{aide-memoire} we name these criteria
Complementarity, Inclusion, Intention, Conjunction, Disjunction,
Multiplicity, and Sequential Optimality (to which we could add
Admissibility and Information Optimality). In passing we note that
pseudo-Bayesian methods, constructed from Bayesian methods by
handicapping them to satisfy constraints on type I error probability,
which are sometimes made out to be Bayesian, are in fact frequentist
in nature.

To convey the issues intuitively, we then consider five example
hypothesis-testing problems with two-dimensional data (or sufficient
statistics) in which inference can be visualised (plus a trivial
discrete example with no nuisance variables for those who haven't
encountered these issues before). The first is designed to require
almost no calculus in its solution; the second is an abstract problem
designed to dramatically illustrate the defects in both FHT and FCS,
but with only a single data point; the third uses many data points,
and has the advantage that both true and apparent Shannon information
content of the different solutions can be easily calculated; the
fourth is an everyday problem showing that the same issues arise also
in common situations, albeit to a less dramatic extent; and the fifth
illustrates that on some real-life problems fixed sample-size
(resp. pseudo-Bayesian) versions of FHT may require more data than
Bayesian methods by a factor of more than 3000 (resp. 300) without the
latter making any recourse to ``informative'' priors. This last
example also illustrates the difficulty of passing some frequentist
regulatory tests without cheating, and the resulting temptations to
cheat with consequent invalidation of product safety measures.

To address the issue of different parties with opposing interests
disagreeing about a prior on which they need to reach a common mind,
we illustrate the beneficial effects of a Bayesian ``Let the data
decide'' policy both on results under a wide variety of conditions and
on motivation to reach a common prior by consent.

We respond to suggestions that frequentists can use FHT to reach a
conclusion just as good as the Bayesian one by intuitively
post-processing the frequentist confidence level, showing that with
many choices of critical regions this is not possible because the
frequentist confidence contains strictly less Shannon information than
the Bayesian posterior about the desired unknowns. Indeed the
frequentist confidence without further post-processing is often
misleading in the sense that it often provides negative
\textit{apparent} Shannon information (ASI). We give an example where
no deterministic choice of critical regions yields any Shannon
information at all about the desired unknowns even though the Bayesian
posterior may contain up to 1 bit of information on the two possible
values -- the maximum one could hope for. On the other hand we show
that if we make use of the Bayesian prior we can always construct a
non-deterministic set of critical regions from whose resulting
frequentist confidence suitable post-processing can recover the
Bayesian posterior (although we can see no point in doing this).

Finally we discuss some other counter-arguments before concluding that
Bayesian methods are fundamentally correct and that FHT and FCS are
wrong in principle as methods to solve inference problems. We
therefore suggest that we should stop teaching FHT and FCS to
high-school and undergraduate students and teach Bayesian methods
first, informing them only at a later stage about frequentist methods
and their defects as being of historical interest. Further, academics
should admit that frequentist methods are inappropriate ways of
solving inference problems rather than running a truce between
frequentists and Bayesians as is currently in vogue, and regulators
should stop insisting on frequentist criteria (e.g. type I error rate)
in approvals testing.

\tableofcontents

\section{Introduction}

\subsection{Background}
\label{introbackground}

When testing which of two hypotheses are likely to be true in the
light of some data, two different methods are commonly used: Bayesian
(\cite{Fienberg}) and frequentist (\cite{Neyman}) hypothesis testing
(FHT). When estimating unknown parameters (including as a special case
which hypotheses may hold), both methods can also be used to provide
confidence sets (also known as credible sets in the Bayesian case) (we
abbreviate \textit{frequentist} confidence sets by FCS). In the
author's experience, during the period 1980 - 2020, frequentist
methods have been used much more frequently than Bayesian ones in
reported analyses, though this now is gradually changing. In
particular, regulators such as the American Food and Drug
Administration (FDA) and the European Medicines Agency (EMA) currently
still insist on using frequentist methods for medical drug and device
approvals testing (\cite{FDApseudoBayes, Regulators}).

However, the two methods give different answers -- they therefore
cannot both be correct answers to the same question; indeed neither
answers the exact question usually asked by the user, namely:
\begin{description}
\item[User's question:] Given the data, which of $H_0$ and $H_1$ is
  true ?
\end{description}
In the problem settings of interest, it is usually not possible to say
with complete certainty which hypothesis is true. Applied to
hypothesis testing, with the hypotheses partitioning parameter space
of unknowns $H=H_0\cup H_1$ with $H_0\cap H_1=\emptyset$, $H_0\neq
\emptyset$, and $H_1\neq \emptyset$, and the data consisting of $x$
lying in a space $X$, the questions correctly answered by the two
methods are:
\begin{description}

\item[Bayesian question:] What are the probabilities $P(H_0|x)$ and
  $P(H_1|x)$ that $H_0$ or $H_1$ respectively is true given the data
  $x$ ?

\item[Frequentist question:] Did we, before collecting the data $x$,
  choose a nested set of critical regions $(C_\eta\subseteq
  X)_{\eta\in [0,1]}$, such that $\eta_1>\eta_2\implies
  C_{\eta_1}\subseteq C_{\eta_2}$, and such that for all
  $\eta\in[0,1]$, and all $h\in H_0$, $P(x\in C_\eta|h)\leq 1-\eta$;
  and if so what is, after collecting the data, the frequentist
  confidence $c=\sup(\{\eta\in[0,1]:x\in C_\eta\}\cup\{0\})$ that
  $H_1$ is true ?

\end{description}
As a user, one therefore needs to decide which of these questions one
wants to answer instead of the original question asked. One may
\textit{not} answer the frequentist question and pretend (either to
oneself or to others) that it is the answer to the Bayesian
question. However, in order to use either method of attempting to
answer the user's question, additional input is required beyond what
is given in the question, namely:
\begin{description}
\item[Bayesian input:] The prior probability $P(H_1)$ that $H_1$
  holds; and often also $P(h)$, the probability distribution over the
  whole of hypothesis space $H$, not just on the subsets $H_0$ and
  $H_1$.
\item[Frequentist input:] The specific choice of $(C_\eta)_{\eta\in[0,1]}$,
  which must be made \textit{before} collecting the data.
\end{description}

The user therefore has to make a choice between these two methods of
answering something different from his original question. Researchers
have argued during most of the last hundred years about which is more
appropriate. More recently we observe that a type of truce has
emerged, where emphasis has been placed more on what the two methods
have in common (for example in their behaviour as the amount of data
approaches infinity under nice conditions), and both sides have
largely ceased to argue their case over the other\cite{LeeChu}.

However, it is our case in this paper that the difference between
these two approaches is of vital importance to the future of science,
to the regulation of medical drugs and devices, and to safety in
technical industries. Without wanting to claim originality for what
follows, we aim to marshall in one place the various arguments why we
believe that the Bayesian approach is the correct way of tackling such
questions, and the frequentist approach is not. We are not alone in
thinking thus (e.g. \cite{RMatthews}, \cite{DJCMBook}), but previous
authors' efforts have so far not had sufficient effect on the
scientific community, hence our present desire to press the points
further. 

To avoid misunderstanding, we are not aiming to discuss heuristic
methods of signal processing, which aim to work most of the time (in
some sense), but do not calculate any explicit levels of probability,
confidence, etc.

In the following we therefore start in section \ref{infprobdef} by
defining what we mean by an inference problem, what we mean by an
admissible solution to such a problem, and criteria that correct
inference methods should meet. To be clear what the different
inference methods in question are, in section \ref{threemethods} we
then specify three inference methods (Bayesian, Frequentist, and
Pseudo-Bayesian), which really amount to just two methods, as the
Pseudo-Bayesian method is just a particular type of frequentist
method. In section \ref{critsection} we discuss their theoretical
merits against a much weaker list of criteria than posited in section
\ref{infprobdef}, pointing out where the criteria are and are not
met. In section \ref{example1} we then provide an example of a
somewhat fanciful inference problem that with only a single data point
is easily visualisable, illustrating the differences between the
various methods, which also gives an intuitive and dramatic view of
why we believe that the frequentist approach is wrong in principle. In
section \ref{shells} we address a more realistic problem with many
data points, in which both true Shannon information (TSI) and apparent
Shannon information (ASI) can be easily calculated, to illustrate the
extent to which frequentist methods can mislead. In section
\ref{example3} we treat similarly a less dramatic and more everyday
example. In section \ref{example2} we give some examples involving
sequences of independent identically distributed random variables such
as are often treated in statistics courses, this time for the purpose
of showing how the Bayesian method produces more accurate results from
less data than frequentist methods. In section \ref{regulatorprior} we
discuss the issue of how one should set prior distributions in the
specific circumstance that two parties disagree about which prior
should be used, for example where one is an industry regulator and the
other is a manufacturer of a product seeking approval. In section
\ref{cannotpostfix} we discuss whether it is possible for a
researcher, given the output of a frequentist method, to apply
intuitive or precise reasoning to deduce what the output of a Bayesian
analysis would have been. We consider some of the counterarguments
often used in section \ref{counterarguments} before concluding in
section \ref{conclusion}, and finally discuss what we should do about
it in section \ref{discussion}.

\subsection{If you find this hard to read at any point...}

Any reader who while reading this finds themselves struggling to
understand abstract mathematics is invited to consider reading about
one or more of the example problems in sections \ref{example1} to
\ref{example3} before coming back to the start. Indeed if you have
never encountered the distinction between Bayesian and frequentist
methods before, we suggest you start with the trivial example in
appendix \ref{dice}.

\section{Inference problems, their solutions, and how to find them}
\label{infprobdef}

\subsection{Inference problems}
\label{infprobdefsub}

By an \textbf{inference problem} we mean a problem of the following
form:

\begin{itemize}

\item We want to know the value, or possible values, of $\theta$,
  which takes values in some set $\Theta$; we also want to know how
  likely each is to be the true value.

\item We observe the value of $x$ which takes values in some set
  $X$.

\item We are not interested in knowing $\phi$ which takes values in
  some set $\Phi$, but it affects the problem because of the next
  point.

\item For each value of $\theta \in \Theta$ and each value of $\phi\in
  \Phi$, we know the probability distribution of $x$ that results,
  which we denote $P(x|\theta,\phi)$, and call ``the likelihood''.

\item In many cases we also know that actually-occurring combinations
  of values $(\theta, \phi)$ can only lie in some subset $H$ of
  $\Theta\times\Phi$; we will often refer to a particular combination
  $(\theta,\phi)$ as $h$ which then takes values in $H$.

\end{itemize}

In other words, we want to answer the question ``How likely is each
possible value of $\theta$ given that we have observed $x$ ?''.

Note that although each of $x$, $\theta$, and $\phi$ may be a scalar
variable, each may also be a vector of many variables. In addition $x$
is assumed to contain not just the observed values, but any decisions
that were reached during the data collection, e.g. on how many data
points to collect, or on which data points to collect.

The sets $\Theta, \Phi, X, H$ often carry some additional structure,
and may in that case (and anyway) be referred to as e.g. ``the space
$H$''. For example $X$ might be the real line $\mathbb{R}$, which
carries algebraic, topological, order, and measure structure (in other
words we can add up real numbers, know which ones are how near which
others, know which are larger and smaller, and know how to assign a
``length'' to subsets (or at least to the Borel sets)). In practice
these spaces always carry a $\sigma$-field of subsets to which a measure
can be applied, and they are usually subsets of $n$-dimensional real
space $\mathbb{R}^n$ for some $n$.

Note that when applying an inference problem (and any solution to it)
to real life it is likely that the real-life probabilities are not
exactly as given by the likelihood in the formal specification of the
problem. In principle this raises the same issue for both frequentist
and Bayesian approaches, although in practice our experience is that
Bayesians are likely to expend more effort than frequentists in
choosing distribution families that model the real-life data well. For
example, many (less academic) frequentists rely greatly on Gaussian
likelihoods, despite it being almost universally the case that the
relevant real-life distributions are much better modelled as
e.g. Student, log-Gaussian, log-Student, Gamma, etc. It is also
arguably easier to deal with non-Gaussian distributions Bayesianly
than frequentistly; indeed we can make the likelihood a mixture of a
variety of distributions selected by a variable in $\phi$. However,
for the purposes of this paper, we will consider methods for solving
the abstract inference problem, rather than the inaccuracies in its
specification, and so lay this issue to one side.

\subsection{Examples of inference problems}

Table \ref{infprobtable} gives a range of fairly typical inference
problems to serve as examples. It is immediately clear that although
some of them are hypothesis tests, others are of a rather more general
form. 

\clearpage

\begin{landscape}
\begin{table}[ht]
\begin{tabular}{p{3cm}p{6cm}p{6cm}p{6cm}p{2cm}}
\textbf{Name} & \textbf{We want to know $\theta$} & \textbf{We observe $x$} & \textbf{$\phi$ (nuisance variables)} & \textbf{Why ?} \\ \hline

Beta-binomial & Probability that biscuit lands chocolate side down &
Number of times dropped, number of times chocolate down & & Curiosity,
Pedagogy \\ \hline

Hypothesis test & Whether probability $p$ that product breaks when dropped $< 0.025$ &
Number of times dropped, number of times broke & $p$ & Regulatory
approval \\ \hline

Multiparallel \mbox{hypothesis} test & Whether each of 250 similar parameters
below threshold &
Number of times tested, number of times failed, for each & $(p_k)_{k=1,...,250}$ & Regulatory
approval \\ \hline

Polynomial fitting stage 1 & Degree $N$ of polynomial $f$ & $x_k=f(z_k)+n_k$
where $n_k$ are unknown noise samples &
Coefficients of polynomial & Polynomial fitting\\ \hline

Proofreading & Number of typos left in book & Which typos found by which proofreader
& Poisson rate of typos in this publisher's books, probability each
proofreader notices a typo & Quality \mbox{control}\\ \hline

Bomb-finding & Where to dig to find an atom bomb buried in your school
grounds & Locations of $\gamma$-ray photons observed at ground level &
Depth of bomb, locations of unobserved photons, fraction of observed
photons that are from bomb & Save your school \\ \hline

Weather part 1 & Exact state of planet & Observed meteorological data from finitely
many sensors & & Forecast tomorrow's weather \\ \hline

\end{tabular}
\caption{A range of inference problems}
\label{infprobtable}
\end{table}
\end{landscape}

\subsection{Admissible solutions}
\label{admissiblesolutions}

Having defined what we mean by an inference problem, we now need to
consider what sort of solutions are acceptable. 

For some observed $x$ and some subset $\Theta_1$ of $\Theta$, we will denote the
output of a solution to an inference problem by $R_x(\Theta_1)$, which tells
us how sure we are if we have observed $x$ that $\theta \in \Theta_1$.

For example, if $\Phi=\{\phi_0\}$ is a single-point space (so that
there are no nuisance variables), $\Theta=[0,1]$, and
$H=\Theta\times\Phi = [0,1]\times\{\phi_0\}\cong[0,1]$, i.e. $\theta$ could
effectively be any number between 0 and 1, then we suggest that any of
the following would be \underline{\textbf{un}}acceptable as
conclusions having observed $x$:

\begin{itemize}

\item $R_x(\{0\}) = 0.95$ and at the same time $R_x(\{1\}) = 0.95$;
  i.e. it is unreasonable to suppose that we are both 95\% sure that
  $\theta=0$ and 95\% sure that $\theta=1$.

\item $R_x(\{1\}) = -3$; i.e. it is unreasonable to be less than 0\%
  sure that $\theta=1$.

\item $R_x([\frac{1}{2},1])=1$ and at the same time
  $R_x([\frac{1}{4},1])=\frac{1}{2}$; i.e. it is unreasonable to be
  completely sure that $\theta\geq \frac{1}{2}$ but only 50\% sure that
  $\theta\geq \frac{1}{4}$.

\end{itemize}

These, and the unacceptability of other similar conclusions, can be
summarised by saying that for any observed $x\in X$, $R_x$ must be a
\textit{probability measure} on $\Theta$ (the reader may like to think
of a probability measure as a probability distribution; alternatively
see appendix \ref{probmeas}, for a precise definition of a probability
measure, related matters, and for more on why it is appropriate here,
or see \cite{Chung} or \cite{probmeaswiki} for other accounts.)

But just being a probability measure is not sufficient to make a
solution acceptable. For example, if $R_x(\{0\})=1$ for all $x$ that
might be observed, but the likelihood is $P(x|\theta) = 1$ for $x=\theta$ and
$P(x|\theta)=0$ otherwise, then we would be saying that even though $x$ and
$\theta$ are almost surely equal, we will conclude that $\theta=0$ whatever
value of $x$ we observe -- which would obviously be nonsense. So we
need some consistency conditions: consistency with the likelihood, but
also consistency with some wider view of the world that includes all
of $\theta, \phi$, and $x$. This brings us to the following definition.

By an \textbf{admissible solution} to an inference problem we mean a
family of probability measures $R_x$ on $\Theta$ satisfying the following
condition (the interpretation is that, for each observed $x$, $R_x$
tells us how likely each value of $\theta$ is):

There exists a joint probability measure $P$ on $\Theta \times \Phi
\times X$ such that\footnote{In all these conditions we allow that
  there may be a subset $N$ of $\Theta\times\Phi\times X$ with
  $P(N)=0$ on which the condition doesn't hold; i.e. we only require
  the conditions to hold ``almost surely''.}:

\begin{description}

\item[Answer consistent:] for all $x\in X$ and all $\theta\in\Theta$,
  the induced\footnote{See appendix section \ref{induced} if the
    meaning of this is not obvious to you.} conditional distribution
  on $\Theta$ given $x$ satisfies $P(\theta | x) = R_x(\theta)$; and

\item[Likelihood consistent:] for all $x\in X$, $\theta\in\Theta$, and
  all $\phi\in\Phi$, the induced conditional distribution
  $P(x|\theta,\phi)$ on $X$ is equal to the likelihood given in the
  problem; and

\item[Restriction consistent:] the induced marginal distribution
  $P(\theta, \phi)$ on $\Theta\times\Phi$ gives
  $P((\Theta\times\Phi)\setminus H) = 0$.

\item[Reality consistent:] if there is further information known which
  tells us how likely each combination of values of $h=(\theta, \phi)$
  is, then the induced distribution $P(\theta, \phi)$ is consistent
  with that information.

\end{description}

To give intuition to this requirement, asking that $R_x$ is a
probability measure on $\Theta$ implies among other things that:

\begin{enumerate}

\item $R_x(\{\theta\})$ cannot be negative or greater than 1;

\item if $\theta_1 \neq \theta_2$ then we cannot have both
  $R_x(\{\theta_1\})>\frac{1}{2}$ and
  $R_x(\{\theta_2\})>\frac{1}{2}$;

\item if $\Theta_1 \subseteq \Theta_2 \subseteq \Theta$ then we cannot
  have $R_x(\Theta_1) > R_x(\Theta_2)$;

\item for all $x\in X$, $R_x(\Theta)=1$ and $R_x(\emptyset)=0$.

\end{enumerate}

Similarly the existence of $P$ satisfying the consistency conditions
means that the answer is at least consistent with some view of the
world that complies with the laws of probability and implies the given
likelihood.

It is our view that any solution failing to meet any of these points
should not be considered as a serious contender to be a solution to
the given problem. However, for those who may at this point want to
object that these requirements are too strict and mean that only
Bayesian solutions stand a chance of being admissible, we will
consider some other, weaker, criteria in section \ref{critsection}
below.

\subsection{Admissible methods}
\label{admmeth}

So now that we have defined an admissible solution, we need a way or
ways of getting to one. To that end we make the following definition:

An \textbf{admissible method} is a method which when applied to any
inference problem results in an admissible solution; to clarify, it
must work on \textit{any} inference problem, although it is only
required to work in principle -- we do not require that the
computation required should be possible on any particular computing
platform or within any particular timescale. Specifically it should
work on problems where $\theta$ and $x$ are independent, including the cases
where $X$ is a single-point space (e.g. when no data is collected),
even though these cases are of no practical interest, and \textit{a
  fortiori} it should work e.g. when $X=\mathbb{R}$, i.e. when only a
single data value has been collected.

\section{The three (or just two) methods}
\label{threemethods}

\subsection{The problem scenario being addressed}

Because of the limitations of frequentist hypothesis testing, we now
concentrate on a subclass of inference problems, where we have two
hypotheses $H_0$ and $H_1$, partioning the space of unknowns $H$, and
we plan to collect some data $x\in X$. In terms of our general
framework for inference problems from section \ref{infprobdefsub}
above, we are considering problems
where $$\Theta=\{0,1\},$$ $$H_0=(\{0\}\times \Phi) \cap
H\subset\Theta\times\Phi,$$ and $$H_1=(\{1\}\times \Phi) \cap
H\subset\Theta\times\Phi,$$ so that $\theta$ is either $0$ or $1$
and tells us whether $H_0$ or $H_1$ holds respectively; wanting to
know $\theta$ is then equivalent to wanting to know whether $h\in H_0$
or $h\in H_1$. Note that each of $H_0$ and $H_1$ may contain either a
single point or many, depending on the nature and relationship of $H$
and $\Phi$.

Indeed as far as the frequentist methods go, we concentrate in this
paper specifically on (a) frequentist hypothesis testing and (b) the
closely related frequentist confidence sets (or intervals), rather
than on any other frequentist ideas.

For all these methods except for the last one in section
\ref{pseudoBayes} we assume that the data collection plan has been
fixed in advance.

\subsection{The Bayesian method}
\label{Bayesmeth}

For this method we choose:

\begin{itemize}

\item a prior probability distribution $P(h)$ on $H$ which expresses
  what we knew\footnote{Note that it is not possible to ``know
    nothing'', although in some cases we can say that each possible
    value of $h$ is equally likely.} before collecting the data about
  the different possible values of $h$; if appropriate to the context
  one may choose one which is in some sense ``uninformative''.

\end{itemize}

Then we collect the data $x$.

As output we calculate\footnote{See appendix \ref{basicBayes} for how.}:

\begin{itemize} 

\item the posterior probability $P(h\in H_1|x)$ given by Bayes
  theorem. 

\end{itemize}

Finally we report $P(h\in H_1|x)$, which tells us how likely it is
that $H_1$ is true given the data.

This is a \textit{purely} Bayesian method. We will consider a
pseudo-Bayesian method below in section \ref{pseudoBayes}.

Note that there is nothing to stop us considering a range of priors
after collecting the data, and calculating the posterior probabilities
for each one; indeed this is encouraged as a way of determining to
what extent the posterior is determined by the prior and to what
extent by the data. When doing that we report each prior together with
its consequent posterior.

\subsection{The frequentist method}
\label{freqmeth}

For this method we first choose:

\begin{itemize}

\item a nested set of critical regions $(C_\eta\subseteq
  X)_{\eta\in [0,1]}$ such that for all $\eta\in [0,1]$ and all
  $h\in H_0$, $P(x\in C_\eta|h)\leq 1 - \eta$, and such that
  $\eta_1 \leq \eta_2$ implies that $C_{\eta_1} \supseteq
  C_{\eta_2}$.

\end{itemize}

Then we collect the data $x$.

As output we calculate:

\begin{itemize} 

\item the frequentist confidence that $h\in H_1$, given by $c=
  \sup{(\{\eta\in [0,1] : x\in C_\eta\}\cup\{0\})}$.

\end{itemize}

Finally we report $c$.

For readers who are unfamiliar with the concept of ``critical
region'', a nested set of critical regions corresponds to a
frequentist ``test'' (e.g. $t$-test, Mann-Whitney $U$-test, log-rank
test, etc). The term ``nested'' refers to the fact that for $\eta_1
\leq \eta_2$, the critical region $C_{\eta_1}$ for the lower
confidence $\eta_1$ must completely contain that for the higher
confidence, like a set of Russian matryoshka dolls. For example when
the likelihood is Gaussian of known unit variance and mean $h$,
$H_0=\{0\}$, and we observe a single sample $x$ of data from the
likelihood, we might set $$C_\eta=\left\{x:G(x)\leq
\frac{1-\eta}{2}\right\}\cup\left\{x:G(x)\geq
1-\frac{1-\eta}{2}\right\},$$ where $G$ is the Gaussian
cdf\footnote{Often denoted $\Phi$, which we avoid here in order not to
  cause confusion.}, and call this a ``simple Gaussian test''.

Note that in this case, in contrast to the Bayesian case, it is
essential that the nested set of critical regions (or frequentist
test) is chosen \textit{before} collecting the data, as otherwise one
can usually choose critical regions deliberately to contain the
observed data and not much else and thereby obtain arbitrarily high
frequentist confidence that $H_1$ holds.

In this context frequentist opinion comes in two flavours:
\begin{description}
\item[Strict frequentist:] If $c$ is greater than some conventional
  value (usually $0.95$), we conclude that $H_1$ is true, but
  otherwise we draw no conclusion, so that in particular we never
  conclude that $H_0$ is true.
\item[Non-strict frequentist:] If $c$ is greater than some
  conventional value (usually $0.95$), we conclude that $H_1$ is true,
  and otherwise we assume that $H_0$ is true.
\end{description}
Most users of frequentist methods (e.g. \cite{pivot}) are non-strict
in their interpretation, while most academic statisticians would
probably adopt the strict interpretation.

\subsection{The pseudo-Bayesian method}
\label{pseudoBayes}

\subsubsection{Preamble}

The pseudo-Bayesian approach essentially aims to take a Bayesian
method and handicap it in such a way that it becomes a frequentist
method. It is unfortunate that we need to discuss it here at all --
indeed the only reasons we do are that it is often described as a
Bayesian method (e.g. \cite{FDApseudoBayes}), which it isn't (as
agreed by the same source in \cite{FDA2018}), and in order to use it
to provide examples of the misbehaviour of frequentist methods.

We therefore urge the reader not to get too bogged down in the details
of how it works, but rather to remember that it is not Bayesian in
nature, even though it results from handicapping a Bayesian method.

It comes in two versions, a basic deterministic version and a full
non-deterministic version. To save verbiage we will often refer to the
basic version simply as the pseudo-Bayesian method and specify
``full'' when needed.

\subsubsection{Intuitive description of basic pseudo-Bayesian method}
\label{pseudoBayesintuitive}

We carry out this frequentist method by first choosing a desired
frequentist confidence level $c_0\in (0,1)$ (which will remain fixed),
and trial values of a prior $P(h)$, a data collection plan, a possibly
enlarged null hypothesis $H_0'$ such that $H_0\subseteq H_0'\subseteq
H$, and a threshold probability $p_0$. We let $H_1'$ be the complement
of $H_0'$ in $H$.

We then determine (either analytically or by extensive simulation) the
region $C\subseteq X$ such that if it later turns out that $x\in C$ we
will calculate a posterior probability $P(h\in H_1'|x) > p_0$. We
can then calculate the frequentist confidence $c$ that would be
achieved by observing $x\in C$.

We then note whether $c\geq c_0$ or not. If we are lucky enough that
$c=c_0$, we proceed to collect data according to the given plan, and
conclude whether or not we have frequentist confidence $c_0$ that
$h\in H_1$ according to whether $x\in C$ or not respectively.

However, if $c < c_0$, we choose some adjustment of some of $P(h)$,
$p_0$, $H_0'$, or the data collection plan, with an aim of increasing
$c$; while if $c > c_0$ by an important margin we choose some
adjustments in the opposite direction.

Eventually we reach a prior $P(h)$, a threshold posterior probability
$p_0$, an $H_0'$, and a data collection plan, such that $c\geq c_0$
with $c - c_0$ being small enough for us not to worry about it.

Finally we collect data, and if $x\in C$ we report frequentist
confidence $c$ that $h\in H_1$.

\subsubsection{Formal description of basic pseudo-Bayesian method}
\label{concisedescription}

Alternatively and more formally: We choose and fix a data collection
plan, an enlarged null hypothesis $H_0'$ with $H_0\subseteq H_0'
\subseteq H$, and a prior $P(h)$ (\textit{perhaps} using the approach
of section \ref{pseudoBayesintuitive} and perhaps without any
consideration of what we actually think about the likely values of $h$
before collecting the data). We define $H_1'=H\setminus H_0'$. We
define non-strictly increasing functions $f,f_1,g:[0,1]\to [0,1]$ and
$B,B_1:[0,1] \to \mathbb{S}(X)$ (the set of subsets of X) by $$B(p) =
\{x\in X:P(h\in H_1'|x)> p\}$$ $$B_1(p) = \{x\in X:P(h\in H_1'|x)\geq
p\}$$ $$f(p) = 1 - \sup_{h\in H_0}P(x\in B(p)|h)$$ $$f_1(p) = 1 -
\sup_{h\in H_0}P(x\in B_1(p)|h)$$ $$ g(\eta) = \inf\{p\in
    [0,1]:f(p)\geq \eta\},$$ and define $$C_\eta =
    \left\{\begin{matrix} B_1(g(\eta)) & (f_1(g(\eta)) \geq \eta)\\ 
B(g(\eta)) & (f_1(g(\eta)) < \eta)
\end{matrix}\right.$$

We then collect the data and report the frequentist confidence $c$
that $h\in H_1$ given by $$c = \sup{(\{\eta \in [0,1]:x\in
  C_\eta\}\cup\{0\})}.$$ (See appendix \ref{conciseproofs} for proofs
that the claimed properties hold and that $(C_\eta)_{\eta \in [0,1]}$
are a valid nested set of critical regions.)

\subsubsection{The full pseudo-Bayesian method}
\label{fullpseudoBayes}

However, we can also construct a uniformly more powerful
pseudo-Bayesian (and still frequentist) method by making use of added
randomness. Keeping the notation of section \ref{concisedescription},
we use that randomness to fill in the gaps between values of
$f(g(\eta))$ where this function is discontinuous. The main motivation
is that doing so results, in the case that $H_0$ and $H_1$ are sets
with only one element each (and therefore $H_0'=H_0$), in a uniformly
optimal frequentist solution.

However, the details are fairly complicated. 

First, building on the notation of section \ref{concisedescription},
we expand the data space to $$X'=X\times [0,1]$$ with the
likelihood $$P(x,u|\theta,\phi)=P(x|\theta,\phi)[u\in[0,1]]$$ so that
$u$ is uniformly distributed in $[0,1]$ independently of all the other
variables. We can then ``observe'' $x'=(x,u)$ by actually observing
$x$ then drawing an independent uniform random $u$ from $[0,1]$.

Then we set about constructing a new family of critical
regions $(C'_\eta\subseteq X')_{\eta\in[0,1]}$ by making the following
definitions: 
 $$S(\eta) = \{\eta' < \eta : C_{\eta'} \supsetneq C_\eta\}$$ $$D_\eta
= \left(X\cap \bigcap_{\eta'\in S(\eta)}{C_{\eta'}}\right) \supseteq
C_\eta$$ $$\zeta_1(\eta) = 1 - \sup_{h\in H_0}{P(x\in
  C_\eta|h)}$$ $$\zeta_0(\eta) = \sup_{\eta'\in S(\eta)}{\left(1 -
  \sup_{h\in H_0}{P(x\in C_{\eta'}|h)}\right)} = 1 - \inf_{\eta'\in
  S(\eta)}{\sup_{h\in H_0}{P(x\in C_{\eta'}|h)}}\leq \zeta_1(\eta),$$
where if $S(\eta)=\emptyset$ the supremum is taken to be zero and the
infimum one. Here the idea is that $D_\eta$ is something like the next
larger $C_\eta$ if there is one (otherwise $D_\eta=C_\eta$ or
$D_\eta=X$) and that $\zeta_0$ and $\zeta_1$ give the lower and upper
ends of the discontinuity in frequentist confidence between the two.

Then we can define $$C'_\eta=\left\{\begin{matrix} (C_\eta\times
[0,1]) \cup \left(D_\eta\times\left[0,\frac{
    \zeta_1(\eta)-\eta}{\zeta_1(\eta)-\zeta_0(\eta)}\right]\right) &
(\zeta_0(\eta)<\zeta_1(\eta))\\ C_\eta\times [0,1] &
(\zeta_0(\eta)=\zeta_1(\eta)),
\end{matrix}\right.$$ and as usual report the frequentist
confidence $$c = \sup{(\{\eta \in [0,1]:x'\in C'_\eta\}\cup\{0\})}.$$

See appendix \ref{fullpseudoBayesproofs} for proofs that this gives a
valid set of critical regions.

Noting (without proof) that in the case of both $H_0$ and $H_1$ both
being single-point sets the full pseudo-Bayesian method is uniformly
optimal but non-deterministic, it is in contrast \textit{not} the case
that then the basic pseudo-Bayesian method is always uniformly optimal
among deterministic methods -- for a counter-example see sections
\ref{freqpseudoBayes} and \ref{nearpseudoBayes}, which apply despite
the Neyman-Pearson lemma which this does not contradict.

\subsubsection{Comments on the pseudo-Bayesian method}

This method meets all the criteria desired of a frequentist test, and
it involves Bayesian calculations in setting it up. But it is
\textit{not} a Bayesian method, because in general:

\begin{itemize}

\item the prior may have been changed from the originally chosen prior
  or chosen arbitrarily rather than according to ones prior beliefs;
  and

\item the effective null hypothesis $H_0'$ may have been changed from
  the original one; and

\item it is the frequentist confidence $c$ that is the primary output
  and not the posterior probability $P(h\in H_1|x)$ -- and decisions
  about whether to approve the method or to publish, submit, or
  approve the result are taken on the basis of the value of $c$.

\end{itemize}

Note also that when using this method it is obligatory to make the
choices of prior and data collection plan \textit{before} collecting
any data, as otherwise one could potentially achieve arbitrarily high
frequentist confidence that $h\in H_1$ by adjusting the data
collection plan and prior to suit the actual outcome.

Thus the pseudo-Bayesian method is just a special case of a
frequentist method -- and we are really discussing only two types of
method, the Bayesian and the frequentist.

\section{Criteria for judging between inference methods}
\label{critsection}

\subsection{Preamble}

Comparing these methods, it is clear that the Bayesian and frequentist
methods are \textit{not} equivalent to each other, whatever prior is
used -- they give different answers, so they cannot both be the right
answer to the same question. Nor is the Bayesian method equivalent to
the pseudo-Bayesian method, which latter is emphatically frequentist
and not Bayesian.

The question remains as to which of these methods (if any) is better
than the others as an inference method. In order to answer that, we
give in the following sections some desirable criteria for inference
methods in general, without necessarily assuming that we are dealing
with an admissible method as defined in section \ref{admmeth}
above. We take it as read that the frequentist methods report
frequentist confidence (and therefore imply a ``type I'' error
probability, the fraction of false positives among the actual
negatives, but with in general no control of the fraction of false
negatives among the actual positives\footnote{Although sometimes
  information is given controlling the fraction of false negatives
  among a particular subset of the true positives.}) while the
Bayesian ones report posterior probability (and therefore imply
probabilities for false positives as a fraction of apparent positives
and for false negatives among the apparent negatives)\footnote{This
  sentence repays more careful reading than it is likely to get at
  first glance.}. Whether it is a
good idea to base an inference method on controlling the type I error
probability is a different matter, which we bear in mind as we
consider first some criteria for judging inference methods, and later
some specific example problems.

In the following we are \textit{not} addressing a number of peripheral
questions, namely:
\begin{itemize}

\item Computational effort required to conduct any particular
  analysis;

\item Ease of understanding any method (we believe that this depends a
  great deal on which method one met first in life);

\item Whether typical results are in analytical closed form or require
  a computer to numerically evaluate them.

\end{itemize}
Rather, we concentrate on the fundamental question of which method
produces correct answers.

We will give some appropriate criteria in section \ref{criteria}, and
consider which methods satisfy them in section \ref{whichsatisfy},
before also looking at some somewhat more nebulous but nonetheless
important considerations in section \ref{further}.

\subsection{Criteria for distinguishing inference methods}
\label{criteria}

In the following we use the term ``$\eta$-sure'' to cover whatever
method of quantitating residual uncertainty that is relevant to the
method, whether that be frequentist confidence or posterior
probability, and we assume $\eta\in [0,1]$.

To formulate our criteria we continue to assume a hypothesis-testing
framework with $\Theta$ a two-point set. 

\begin{enumerate}

\item \label{complementarity} \textbf{``Complementarity''}: If given
  some data $x$ we are $\eta$-sure that $H_1$ holds, and $\eta$ is
  high (i.e. near 1), then if we set $$H_0'=H_1,$$ $$H_1'=H_0,$$ then
  applying the same inference method to $H_0', H_1'$ with the same
  data $x$ we should be $\eta'$-sure that $H_1'$ holds only for values
  of $\eta'$ that are low (i.e. near 0); indeed ideally we should have
  $\eta'=1-\eta$ at most.

  Intuitively, this says that what we conclude should not depend on
  the names we give to the different possibilities.

\item \label{containment} \textbf{``Inclusion''}: If $H_1 \subseteq H_1'$
  then it should not be possible to conclude that we are more sure
  that $h\in H_1$ than that $h\in H_1'$.

  Example: If I am $90\%$ sure that the Queen is in Buckingham Palace
  then I should not be able to be less than $90\%$ sure that she is in
  London.

\item \label{whetherproceed} \textbf{``Intention''}: If we collect data
  $x_1$, then the inference resulting should not depend on whether or
  not we intended to perhaps, depending only on the value of $x_1$,
  collect more data $x_2$ on the same system.

  Example: If I toss a coin five times getting five heads, then the
  inference \textit{at that time} on whether the coin is fair or not
  should not depend on whether, had I not got five heads, I would have
  tossed it a further 1000 times in the hope of getting a sequence of
  five consecutive heads at some point.

\item \label{productprob} \textbf{``Conjunction''}: If we have $N > 1$
  independent and identical systems with null hypotheses $H_{0,n}$ for
  $n=1,...,N$ and complementary alternative hypotheses $H_{1,n}$, and
  we collect data on each that for some $\eta \in (0,1)$ and for all
  $n$ makes us exactly $\eta$-sure that $h_n\in H_{1,n}$, then we
  should always be strictly less than $\eta$-sure that for all $n$,
  $h_n\in H_{1,n}$.

  Example: If I toss 100 coins five times each, and by some
  amazing luck happen to get five heads from each coin, and my
  inference method makes me exactly 93.75\% sure that the first coin
  is biased, and exactly 93.75\% sure that the second coin is biased,
  and so on, then it should not make me 93.75\% sure that every single
  coin is biased. (Just as my being 50\% sure that if I toss an
  unbiased coin once I will get heads should not make me 50\% sure
  that if I toss 100 unbiased coins once each I will get heads from
  all 100 coins.)

\item \label{sumprob} \textbf{``Disjunction''}: If we again have $N >
  1$ independent systems, identical apart from the values of $h_n$,
  with null hypotheses as before, and we collect data on each such
  that for some $\eta\in (0,1)$ and for all $n$ makes us exactly
  $\eta$-sure that $h_n\in H_{1,n}$, then we should always be strictly
  more than $\eta$-sure that for at least one $n$, $h_n\in H_{1,n}$.

  Example: If I toss 100 coins five times each, and by some
  amazing luck happen to get five heads from each coin, and my
  inference method makes me exactly 93.75\% sure that the first coin
  is biased, and exactly 93.75\% sure that the second coin is biased,
  and so on, then it should make me strictly more than 93.75\% sure
  that at least one coin is biased. (Just as my being 50\% sure that
  if I toss an unbiased coin once I will get heads should make me
  more than 50\% sure that if I toss 100 unbiased coins once each I will get
  at least one head among the 100 coins.)

  Note: A weaker version of this criterion (with ``at least'' instead
  of ``strictly more than'') is automatic for any method that
  satisfies criterion \ref{containment}. Moreover for any method that
  treats $H_0$ and $H_1$ identically and is admissible, this criterion
  is equivalent to criterion \ref{productprob} by de Morgan's
  laws. However frequentist hypothesis testing does not treat $H_0$
  and $H_1$ identically, so we treat this as a separate criterion.

\item \label{multisystem} \textbf{``Multiplicity''}: If we again have
  $N$ independent systems, identical apart from the values of $h_n$,
  with null hypotheses as before, and for each $n$ we have collected
  data $x_n$, then the conclusion about whether or not $h_1\in
  H_{1,1}$ should not depend on the value of $N$ if the value of $x_1$
  remains constant.

  Example: If I toss a coin five times, getting five heads,
  then the inference on whether that coin is fair or not should not
  depend on whether I also toss twenty other unrelated coins or what
  results if and when I do toss them.

\item \label{unifpow} \textbf{``Sequential Optimality''}: Suppose we
  have a single system for which we collect data $x_1,x_2,x_3,...$
  sequentially (all the $x_k$ being conditionally independent from
  each other given $\theta$ and $\phi$), and we want to show that
  $h\in H_1$. Then it is desirable that we should be able to decide on
  an optimal data collection and analysis plan without having first to
  guess the value of $h$.

  More precisely, it is desirable that for each $\eta\in [0,1]$ there
  exists a uniformly most powerful data collection and analysis plan
  $D$. In other words for each $\eta\in [0,1]$ there should exist a
  data collection and analysis plan $D$ such that for all alternative
  plans $D'$ and all $h\in H_1$ and all $k \geq 0$, $p_k(D,h) \geq
  p_k(D',h)$, where $p_k(D,h)$ denotes the probability for that
  particular true value of $h$ that under plan $D$ we conclude that we
  are $\eta$-sure that $h\in H_1$ at some point before the collection
  of $x_{k+1}$.

  By a ``data collection and analysis plan for $\eta$'' we mean a
  sequence of measurable quit-probability functions $q_0, q_1(x_1),
  q_2(x_1,x_2), q_3(x_1,x_2,x_3), ...$ ($q_0$ just being a constant)
  taking values in $[0,1]$ such that we stop data collection
  immediately after collecting $x_k$ if an independent uniform random
  number $u_k$ from $(0,1)$ is less than $q_k(x_1,...,x_k)$, and such
  that whenever we so stop collecting data we can conclude that we are
  $\eta$-sure that $h\in H_1$ (otherwise we go on collecting data for
  ever). (We can encode the idea that we never continue collecting
  after $x_K$ by setting $q_k$ to be identically zero for all $k>K$.)
  (Note that $q_0$ denotes the probability that we conclude without
  collecting any data at all.)

  Example: I have a factory making torches. My quality
  control criterion is that when I turn the switch on the torch should
  light with probability at least 0.9 . I do not know exactly what the
  probability $p$ of a torch from this factory working is, and I want
  to test whether the factory satisfies the quality control criterion,
  and if it does, become 95\% sure that it does. It is desirable that
  in deciding how many torches to test I should not need to know the
  actual value of $p$, and it is desirable that I don't test more
  torches than necessary.

\end{enumerate}

To these we could add \textbf{``Admissibility''} (that the method is
admissible, section \ref{admissiblesolutions} above) and
\textbf{``Information Optimality''} (that the solution contains at
least as much information about $\theta$ as any other solution,
sections \ref{Shannon} and \ref{basicSi} below).

\subsection{Which inference methods satisfy these criteria ?}
\label{whichsatisfy}

We note that for any prior probability measure $P(h)$ on $H$, its
associated Bayesian solution is admissible, by setting $$P(\theta,
\phi, x) = P(\theta, \phi)P(x|\theta, \phi) = P(h)P(x|\theta, \phi).$$
Moreover, any admissible solution corresponds to a Bayesian solution
using a prior $P(h)=P(\theta,\phi)$ induced by a $P(\theta, \phi, x)$
specified in the definition of an admissible solution\footnote{Note,
  however, that there may be many such $P$ and many different such
  priors.}. Showing that any Bayesian solution satisfies these criteria
is therefore equivalent to showing that any admissible solution
satisfies these criteria (and we have thus in passing proved the
``Admissibility'' criterion).

But it is easy to show that any (pure) Bayesian solution satisfies all
of these criteria; proofs are in appendix \ref{bayesproofs}
below. Moreover in relation to criterion \ref{unifpow} it is shown in
theorem B1 of \cite{CTM} that under weak additional conditions there
even exists a uniformly most powerful data collection plan $D$ such
that for all $h\in H_1$, $p_k(D,h)\to 1$ as $k\to \infty$. Further, it
is clear from the calculations in appendix \ref{basicBayes} that we
can apply the Bayesian method to \textit{any} inference problem,
making it an admissible method according to the definition of section
\ref{admmeth} above.

In contrast, for each of these criteria, examples exist in which the
frequentist method violates the criterion -- specific examples are
given in appendix \ref{freqexamples} below. It follows that the
frequentist method is \textit{not} admissible. In the case of
criterion \ref{unifpow} the example given is even one where a
uniformly most powerful nested set of critical regions for a fixed
number of data samples \textit{does} exist -- but even then there is
no uniformly optimal critical region for a variable number of data
samples and no uniformly most powerful data collection and analysis
plan -- and with only a fixed number of data samples we cannot hope
that $p_k(D,h)\to 1$ as $k\to \infty$.

We note that the combination of violations of criteria
\ref{whetherproceed}, \ref{multisystem}, and \ref{unifpow} has an
adverse effect on the ability of industry to perform frequentist
testing of factories producing safety-critical devices, of such a
severity that it becomes near-impossible to pass some such tests
without cheating \cite{CTM}. In the same setting Bayesian methods
instead produce efficient testing that can be passed with probability
approaching 1 as the amount of data approaches infinity whenever the
factory meets the specification and without having to guess the true
margins by which it does so. So the differences are important in the
real world, not just in theory.

Moreover it is apparent from the definition of the frequentist
hypothesis-testing method in section \ref{freqmeth} that it is not
readily adaptable to general inference problems, unlike the Bayesian
method. Note however that (many different) frequentist confidence set
functions can be found for most inference problems - this is explored
more in section \ref{fconfsets} below. Indeed for any given inference
problem there is a bijective relationship between \textit{complete}
families of critical regions and frequentist confidence set functions,
as defined and shown in appendix \ref{hypoconf}. Thus each
counter-example in appendix \ref{freqexamples} showing that
frequentist hypothesis testing fails one of the criteria can be
immediately translated into a counter-example showing that frequentist
confidence set functions also fail that criterion -- and we see that
frequentist confidence sets are only any better than frequentist
hypothesis testing in that some (but not all) frequentist confidence
set functions are symmetric.

\subsection{Further considerations}
\label{further}

\subsubsection{Which variables should follow the ``given bar'' $|$ in
  probabilities of interest ?}
\label{whenwhich}

By the ``given bar'' we here mean the symbol $|$ in expressions such
as $P(A|B)$ or $P(h|x_1,x_2)$.

Suppose that we want to consider ``predictions'' in the broadest
sense: \begin{itemize}\item we may wish to predict the value of a
  variable which is not yet determined, such as the value of something
  on the stock market one week from now;
\item we may wish to predict the number that has already been rolled
  on a die that is currently out of sight;
\item we may wish to predict the (currently unknown but already
  determined) mean of a distribution from which we have seen some
  samples.
\end{itemize}

In all these cases probabilities can be used to express our current
state of knowledge and uncertainty, and in all these cases our
knowledge may increase gradually or in several steps: the value of the
stock market item tomorrow and the day after will usually contain
information about its value in one week's time; somebody who can see
the die may tell us that the number rolled is even; and we may gather
some more samples from the distribution whose mean is wanted.

Let us take a specific example: we want to know the sum $S_{10}$ of
the first ten rolls of a fair 6-sided die. Before rolling the die at
all, we know that the expectation of $S_{10}$ is $35$, and the
distribution of $S_{10}$ is that of the sum of ten independent
variables each of which is uniform on the integers from $1$ to
$6$. Before the die has been rolled at all, if we want to make
predictions or bets about $S_{10}$, this is the right distribution to
use: it would be reasonable at this point to offer even odds on
$S_{10}$ being greater or less than $35$.

Now suppose that the first eight rolls of the die give us $(2, 3, 1,
4, 2, 2, 5, 1)$ with a total of $S_8 = 20$. Is it still reasonable to
offer the same odds to new bets on whether $S_{10}$ will be greater or
less than $35$ ? Of course it isn't -- we now know that $S_{10}$
cannot be greater than $32$, and that the probability of getting even
$32$ is only $\frac{1}{36}$, while the expectation of $S_{10}$ is now
just $27$. In symbols we write $P(S_{10}=32|S_8)=\frac{1}{36}$ and
$\E(S_{10}|S_8)=27$, where the LHS of each of these equations is a
function of the value of $S_8$ actually observed, which happens to be
$20$. 

The same would apply if the die had already been rolled all ten times
but we had been blindfolded after the eighth roll.

\textit{A fortiori}, when making predictions (including those of
things already determined), the probabilities that matter are those
where all the variables whose values we already know are listed on the
right side of the ``given bar'' $|$ in expressions such as
$P(S_{10}|S_8)$. The probabilities we had to start with, before the
values of these variables became known, are no longer useful for
making predictions.

Thus when (perish the thought) we consider buying shares in Facebook,
Apple, or Microsoft with a view to selling them in two weeks, the
price we are prepared to pay is based on a distribution around their
price today, not on the far lower price that these shares fetched when
first floated.

(Obviously where variables become known that are independent of those
which we are predicting, they make no difference and can safely be
omitted.) 

Notice then that the frequentist method (section \ref{freqmeth})
considers only probabilities which have $h$ (which we don't know) to
the right of the given bar, contrary to this principle, while the
Bayesian method (section \ref{Bayesmeth}) calculates probabilities in
which $x$ (the only variable whose value we do know) appears to the
right of the given bar, consistent with this principle.

\subsubsection{Which points of $H\times X$ are relevant ?}

A further point we might consider is this: in the space $H\times X$,
each method considers the probabilities of certain particular
combinations $(h,x)$ occurring. Which of these combinations should be
considered, and which does each method actually consider ?

Pictorially, it may look as shown in figure \ref{htimesxplot1} (note that the $H$-axis goes from
left to right while the $X$-axis goes from bottom to top):

\begin{figure}[h]
\begin{center}
\includegraphics{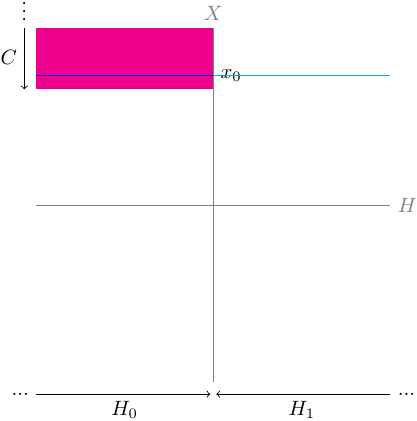}

\caption{The $H\times X$ plane for a problem where $x$ is distributed
  Gaussianly with mean $h$ and variance $1$. In calculating whether
  the critical region $C$ is valid, the frequentist approach considers
  all the points in the magenta region, while the Bayesian approach
  considers only the points on the cyan line. In this case the
  frequentist rejects $H_0$.
\label{htimesxplot1}
}
\end{center}
\end{figure}

Here we might be considering a problem where $x$ is distributed
Gaussianly with mean $h$ and variance $1$; if $H_0=\{h\in
\mathbb{R}:h\leq 0\}$, then a typical choice of critical region
$C=C_\eta$ for $\eta=0.975$ would be $\{x\in \mathbb{R}:x \geq
1.96\}$.

Each point on the plot represents a combination $(h,x)$ and has a
probability density, given by $$P(h,x)=P(h)P(x|h),$$ the product of
the prior and the likelihood.

In evaluating $\eta$ for this critical region $C$ the frequentist
method considers the likelihood at all the points coloured magenta in
the diagram and ignores the prior. On the other hand if we observe the
specific data value $x_0$, then in calculating $P(h\in H_1|x=x_0)$ the
Bayesian method considers both the likelihood and the prior at all the
points on the cyan line in the diagram -- it integrates $P(h|x=x_0)$
over $H_1$, but in calculating $P(h|x=x_0)$ in the first place it
considers all the points on the entire cyan line from $(-\infty,x_0)$
to $(+\infty,x_0)$,
as $$P(h|x)=\frac{P(h)P(x|h)}{\int_{-\infty}^\infty{P(h)P(x|h)\,dh}}.$$

Which of these two is appropriate ? The Bayesian method considers all
the combinations of $h$ and $x$ that could actually have occurred --
any not on the line are irrelevant as we know they haven't
occurred. On the other hand the frequentist method not only considers
a pile of points not on the cyan line (that couldn't be relevant as we
now know that the value of $x$ is $x_0$), but it also fails to
consider the points $(h,x_0)$ for $h\in H_1$ -- and while in this
problem they are more likely than those in the magenta region, in
other problems they might not be, and concluding that $h\in H_1$
simply because it is unlikely that we will see $x_0$ when $h\in H_0$ is
not sensible if it could be even less likely when $h\in H_1$.

\begin{figure}[h]
\begin{center}
\includegraphics{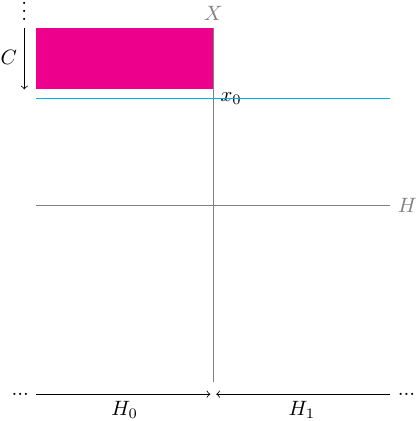}

\caption{The $H\times X$ plane for a problem where $x$ is distributed
  Gaussianly with mean $h$ and variance $1$. In calculating whether
  the critical region $C$ is valid, the frequentist approach considers
  all the points in the magenta region, while the Bayesian approach
  considers only the points on the cyan line. In this case the
  frequentist does not reject $H_0$.
\label{htimesxplot2}
}
\end{center}
\end{figure}

Moreover the picture might instead look like Figure
\ref{htimesxplot2}, in which case the frequentist method won't even
consider \textit{any} of the $(h,x)$ combinations that actually could
have occurred -- but nonetheless most \textit{users} of frequentist
methods will conclude that $h\in H_0$ (e.g. \cite{pivot}), even if
stict frequentists would say that that isn't appropriate.

\subsubsection{What if there are more than two hypotheses ?}
\label{morethantwo}

Yet another consideration is what to do if there are more than two
hypotheses. The Bayesian method applies straightforwardly,
consistently, and symmetrically to any number of competing
hypotheses. On the other hand, for frequentist hypothesis testing, it
is already the case that $H_0$ is treated differently from
$H_1$. Suppose $H_0$ is rejected -- should we then treat $H_1$ as the
new null hypothesis for a further test against an alternative
hypothesis of $H_2$ ? Or the other way round ? And if we reject $H_0$
at the $96\%$ level and then go on to make $H_1$ the null hypothesis
and reject it at the $98\%$ level, with what frequentist confidence
can we claim to have shown that $H_2$ holds ?

We remark in passing that for the frequentist ``conclusive''
confidence set functions of section \ref{conclusive} below, dealing
with more than two hypotheses is not a problem as there is an obvious,
and symmetric, extension of this approach available -- not that it
helps much, as we will see.

\section{An example problem comparing the different methods}

\label{example1}

We now give an example problem which illustrates the various methods
visually. In order to make the issues clear, this is a 2-dimensional
problem, but nonetheless only slightly more complicated than problems
with which most readers will be familiar.

We realise that many frequentists would avoid using frequentist
hypothesis testing (or indeed frequentist confidence sets) to address
this problem, and would instead opt for a Bayesian approach (or
something that genuinely is equivalent to a Bayesian approach with a
flat prior). Nonetheless, it \textit{is} an inference problem; and the
point of using it is to apply \textit{reductio ad absurdum} to the
frequentist approach -- the principles of a good inference method,
applied to any inference problem, should not yield nonsensical
answers.

We note further that frequentists might object to us considering a
problem where we are given only a single data point. However, there is
no difference in principle between (a) being given $N$ data points,
each of which is a vector of $K$ dimensions, (b) being given a single
data point which is a vector of $NK$ dimensions, and (c) being given a
single sufficient statistic (which may or may not come from $NK$
original data dimensions). Here we will take a single data point which
is a vector of 2 dimensions: again, a good inference method should
work well on any number of data points, but illustration is more
difficult when the total number of dimensions exceeds 2. It is our
belief that insistence on having data of e.g. 20 dimensions or more
would serve only to obscure what is actually going on geometrically,
and that this 2-dimensional problem is ideally suited to illustrate
the issues involved.

For those, however, who think that using a problem with only a single
data point is just providing a frequentist straw man to be easily
knocked over, we do also provide examples of a more statistical
flavour (and with many identically distributed data points) in
sections \ref{shells}, \ref{example3}, and \ref{example2} below,
though the last of these is more designed to show the advantages of
the Bayesian method than why the frequentist method is wrong. We also
provide an extremely simple and fully worked example in appendix
\ref{simpleexample} that requires almost no calculus for those who are
totally unfamiliar with an abstract approach to inference problems.

We ask the reader's indulgence for us presenting this problem in
fanciful language, which we believe helps us to remember what is
happening as we analyse the situation, while avoiding bias towards any
particular scientific discipline.

\subsection{The problem}

On a large flat sandy desert plain there are two towers, each of
height $z = 1 \text{ km}$, and difficult to climb. One is painted with
a large digit zero and is located at coordinates $(-2z, 0)$ on the
plain (i.e. it is 2 km to the West of the origin), while the other is
at coordinates $(+2z, 0)$ (i.e. it is 2 km to the East of the origin)
and is painted with a large digit one.

Princess Arabella has been captured and imprisoned by placing her on
the top of one of these towers, but we don't know which
one\footnote{\textit{Lord of the Rings} fans might think of Gandalf
  captured by Saruman and stuck on the top of Orthanc in
  Isengard.}. Unusually, her captor accidentally left his gun behind
on the top of the tower with her, loaded with just one bullet. The
innocent princess has never seen a gun before, and there isn't much to
do on top of the tower on which she finds herself, so she starts
investigating the gun, and, without intending to shoot anybody, pulls
the trigger at a time when the gun is pointing in a random direction
around the sphere. The bullet - fortunately missing her and unaffected
by gravity - travels in a straight line, either off into space, or
hitting the sandy plain at some point $(x,y)$. The wind blows
overnight, leaving the bullet (if it lands) in the same place and
exposed (in a bed of particles larger grains move to the top when
agitated), but obscuring any clues as to its direction of arrival.

A prince knows that she has been so imprisoned, and wants to rescue
her. Knowing both princess and captor well, he knows what is likely to
have happened with the gun and the bullet. There are two cases to
consider: the bullet lands somewhere, or it travels off into
space. Either way he wonders which tower to climb to rescue his
princess. He decides to call the hypothesis that the princess is in
the tower painted with a zero $H_0$ and the other possibility he calls
$H_1$.

But which of the three methods will he employ to work this out ?

\subsection{The parts of the solution common to all three approaches}
\label{commonpartsall}

Let $h=0$ if $H_0$ holds and $h=1$ if $H_1$ holds, so that we identify
$H_0$ with $\{0\}$ and $H_1$ with $\{1\}$ and $H$ with $\{0,1\}$.

\begin{figure}[p]
\begin{center}

\includegraphics[scale=0.5]{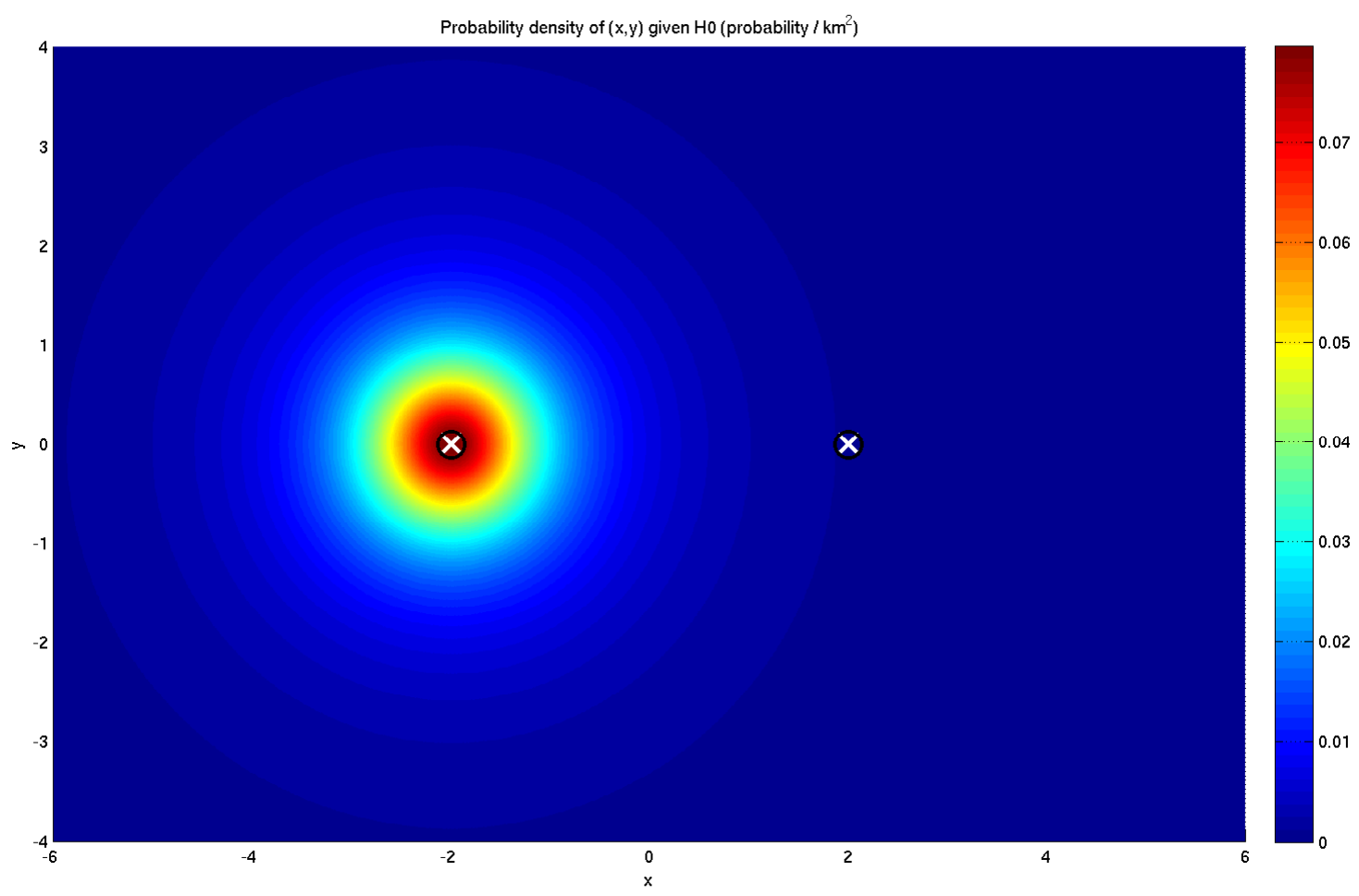}

\caption{The probability density of $(x,y)$ given $H_0$ (probability
  per square kilometre). The positions of the two towers are also
  shown.
\label{likelihood0}
}

\end{center}
\end{figure}

\begin{figure}[hp]
\begin{center}

\includegraphics[scale=0.5]{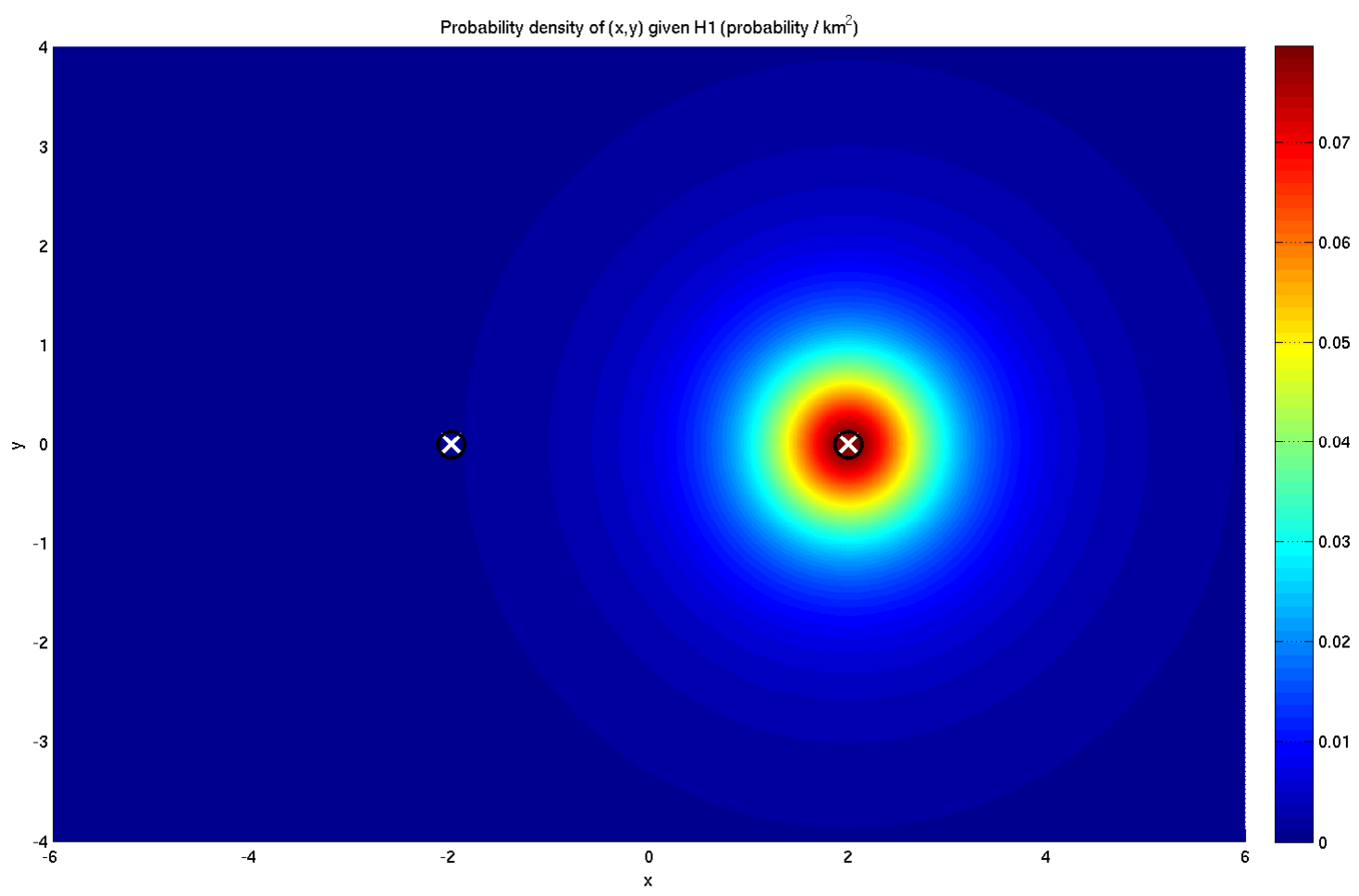}

\caption{The probability density of $(x,y)$ given $H_1$ (probability
  per square kilometre). The positions of the two towers are also
  shown.
\label{likelihood1}
}

\end{center}
\end{figure}

Whichever inference method is employed, if $(a,b)$ denotes the
coordinates of the base of the tower holding the princess, i.e.
$(-2z,0)$ or $(+2z,0)$ as appropriate, we first calculate\footnote{We
  recommend just taking this on trust as it is agreed by both
  frequentists and Bayesians, but the calculation can be found if
  wanted in appendix \ref{densitycalcs}.} that the probability of the
bullet heading off into space is $\frac{1}{2}$, and the probability
density of its landing point $(x,y)$ is given by $$P(x,y|a,b) =
\frac{z}{4 \pi (z^2 + (x-a)^2 + (y-b)^2)^{\frac{3}{2}}}.$$ This
distribution is shown in figure \ref{likelihood0} for the case that
$h=0$ and in figure \ref{likelihood1} for the case that $h=1$.

(In a similar Gaussian problem it is only the value of $x$ that
carries information about $h$; note that for this problem, however,
the value of $y$ also carries information about $h$ if $x$ is already
known -- specifically keeping $x$ alone carries on average
$\approx0.233$ bits of information about $h$ for a single fired
bullet, while keeping $y$ as well carries a further $\approx0.025$
bits -- see section \ref{Shannon} for the details of how this is
calculated.)

Similarly, if $R$ denotes the distance of $(x,y)$ from $(a,b)$, we can
also in all cases calculate that $$P(R) = \frac{zR}{2(z^2 +
  R^2)^{\frac{3}{2}}},$$ and hence also that $$P(R > R_0) =
\frac{z}{2\sqrt{z^2 + R_0^2}},$$ so that $$P(R > 0) = \frac{1}{2},$$
with the remaining probability of $\frac{1}{2}$ being taken up by the
bullet going into space rather than landing.

But whichever method we decide to employ now, we next face a
dilemma. We will deal with the various methods in order of simplicity,
dealing with the simplest first.

\subsection{The Bayesian approach}

\subsubsection{Recap of the Bayesian approach}

Let us denote by $D$ the data obtained, which is either that the
bullet has not landed anywhere, or is a pair of coordinates $(x,y)$
where the bullet landed.

Then by Bayes' theorem, we have $$P(h|D) =
\frac{P(h)P(D|h)}{\sum_{h=0}^1{P(h)P(D|h)}},$$ or in
particular $$P(H_1|D) = \frac{P(H_1)P(D|H_1)}{P(H_1)P(D|H_1) +
  P(H_0)P(D|H_0)},$$ where $P(H_1)$ is the prior probability that the
princess in in tower $1$ and $P(D|H_1)$ is the likelihood, given
by $$P(D|H_1) = \left\{\begin{matrix}\frac{1}{2} & \text{(bullet did
  not land)}\\ \\ \frac{z}{4\pi (z^2 + (x-2z)^2 + y^2)^{\frac{3}{2}}} &
\text{(bullet landed at }(x,y))\end{matrix}\right.$$ and similarly for
$H_0$ but with $(x+2z)$ instead of $(x-2z)$.

\subsubsection{The easy case}

First let us deal with the easy case, namely that the bullet has not
landed. Then the princess is more likely to be in whichever tower the
prince's prior favours, or equally likely to be in either if the
prince has equal prior probability of $\frac{1}{2}$ on each,
since $$P(H_1|D) =
\frac{P(H_1)\frac{1}{2}}{P(H_0)\frac{1}{2}+P(H_1)\frac{1}{2}} =
\frac{P(H_1)}{P(H_0)+P(H_1)} = P(H_1).$$ (We have learned nothing, so
the posterior is equal to the prior.)

\subsubsection{Using the obvious prior}

But let us now assume that the bullet has landed at $(x,y)$. The
dilemma is again that we have to put a prior on the variable $h$,
which is $0$ in the case that $(a,b) = (-2z,0)$ and $1$ in case that
$(a,b) = (+2z,0)$.

Let us start by assuming that we choose
the prior that $P(h=0)=\frac{1}{2}=P(h=1)$ (we will discuss what
happens with other priors later).

\begin{figure}[p]
\begin{center}

\includegraphics[scale=0.5]{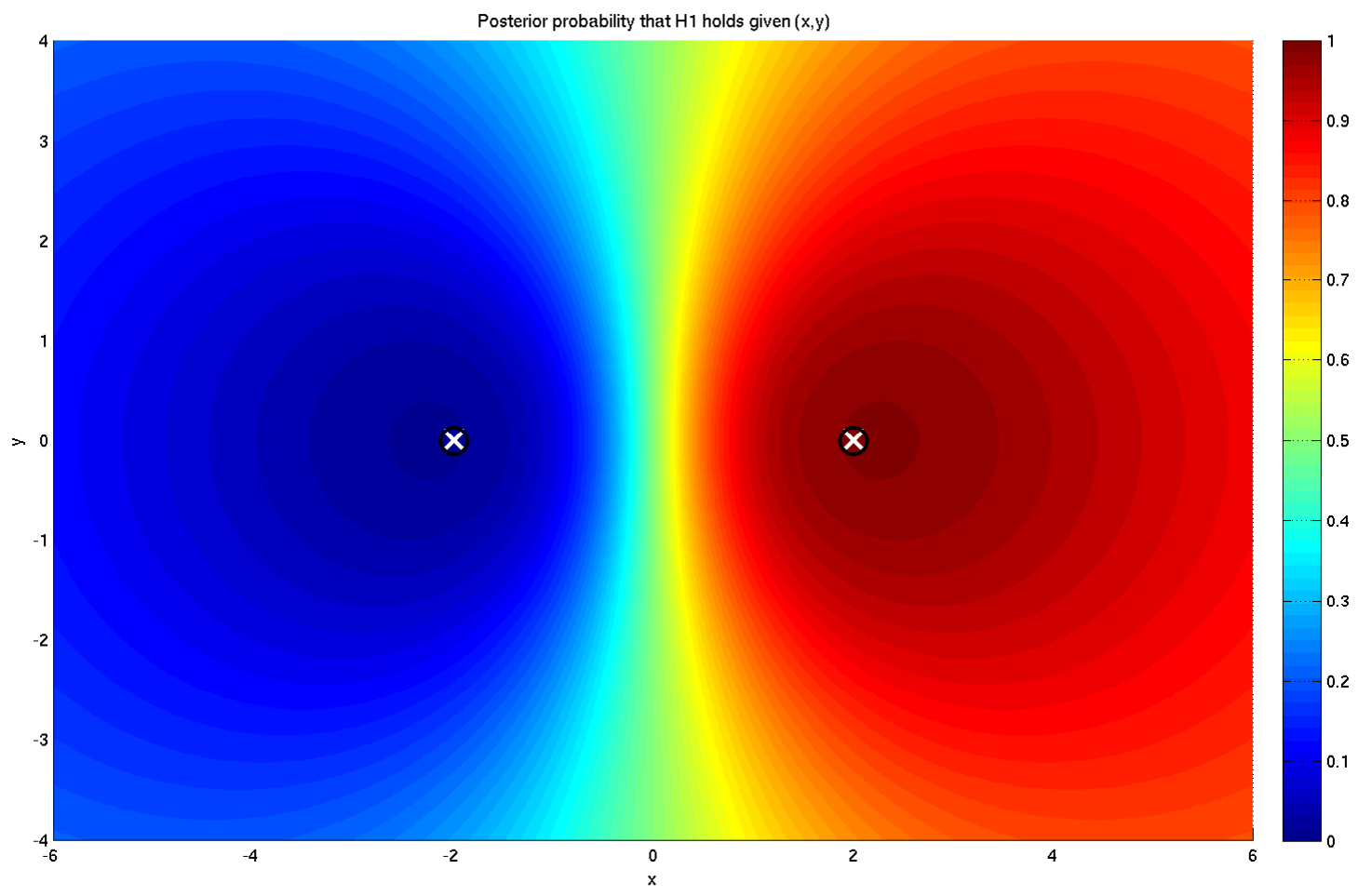}

\caption{The probability of $H_1$ given $(x,y)$, i.e. the posterior
  probability that $H_1$ is true. The positions of the two towers are
  also shown.
\label{postprobeven}
}

\end{center}
\end{figure}

Then application of Bayes' theorem tells us that $$P(h=1|x,y) =
\frac{\frac{1}{2}\frac{z}{4 \pi (z^2 + (x-2z)^2 +
    y^2)^{\frac{3}{2}}}}{\frac{1}{2}\frac{z}{4 \pi (z^2 + (x-2z)^2 +
    y^2)^{\frac{3}{2}}} + \frac{1}{2}\frac{z}{4 \pi (z^2 + (x+2z)^2 +
    y^2)^{\frac{3}{2}}}}.$$ In other words, the posterior probability
of $H_1$ given $(x,y)$ is as shown in figure \ref{postprobeven}.

\begin{figure}[htp]
\begin{center}

\includegraphics[scale=0.5]{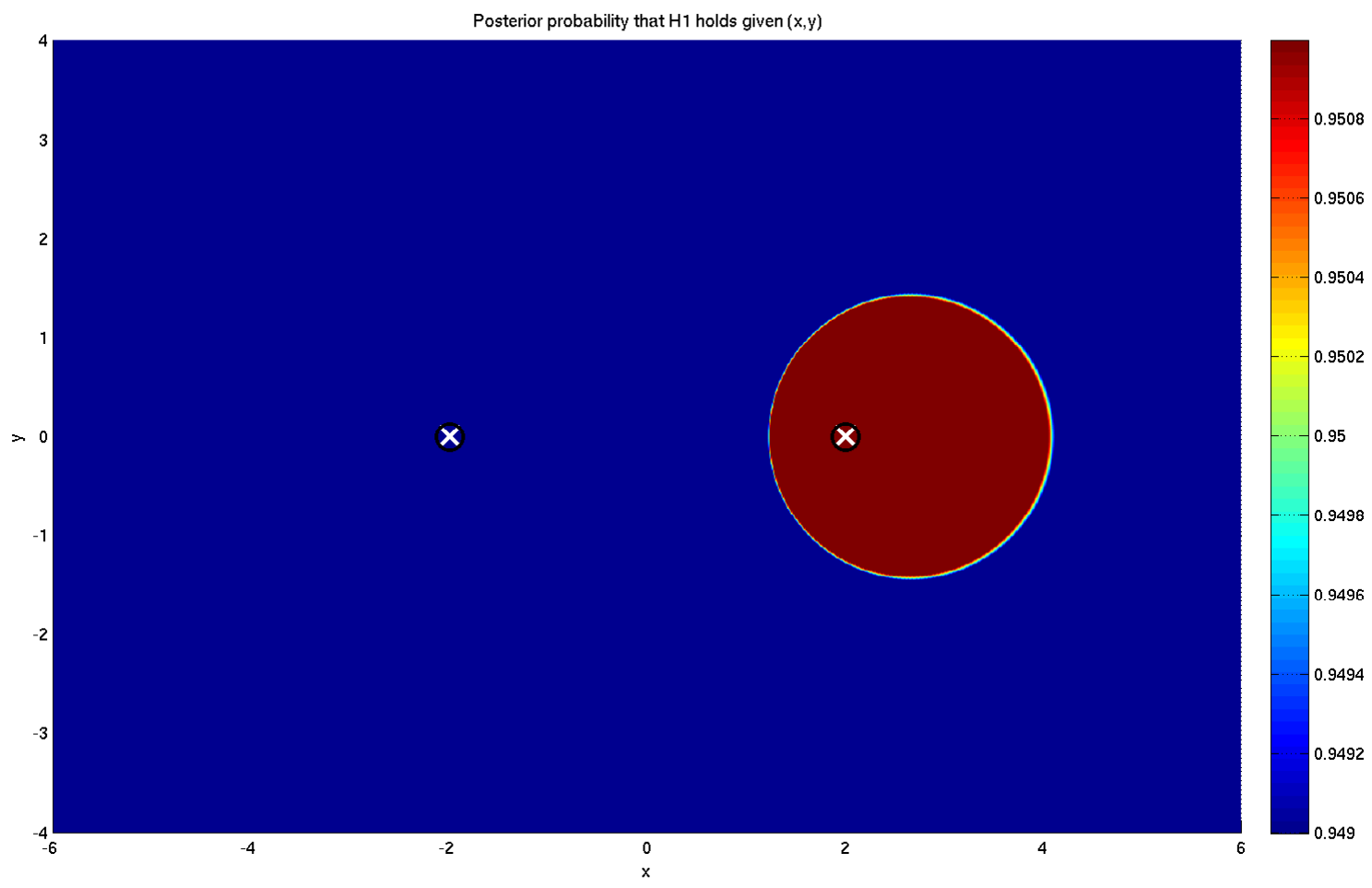}

\caption{The region of the $(x,y)$ plane in which the posterior
  probability that $H_1$ holds is at least $0.95$ is shown in
  brown. The positions of the two towers are also shown.
\label{postprobeven95}
}

\end{center}
\end{figure}

Now, this plot tells us that if $x>0$ then $H_1$ is more likely to
hold than $H_0$ and if $x<0$ then $H_0$ is more likely than $H_1$. So
if $x>0$ then the prince would be well advised to climb tower $1$.

However, the prince may not be somebody who wants to risk his life on
something that is less than 95\% certain. In that case, before he goes
climbing tower $1$, he would like to see $(x,y)$ lie in the dark brown
area of figure \ref{postprobeven95} (or in a similar area on the left
of the plot for tower $0$); and otherwise this less-dashing prince
will give up his attempt to rescue the princess.

\subsubsection{Using another prior}

\begin{figure}[p]
\begin{center}

\includegraphics[scale=0.5]{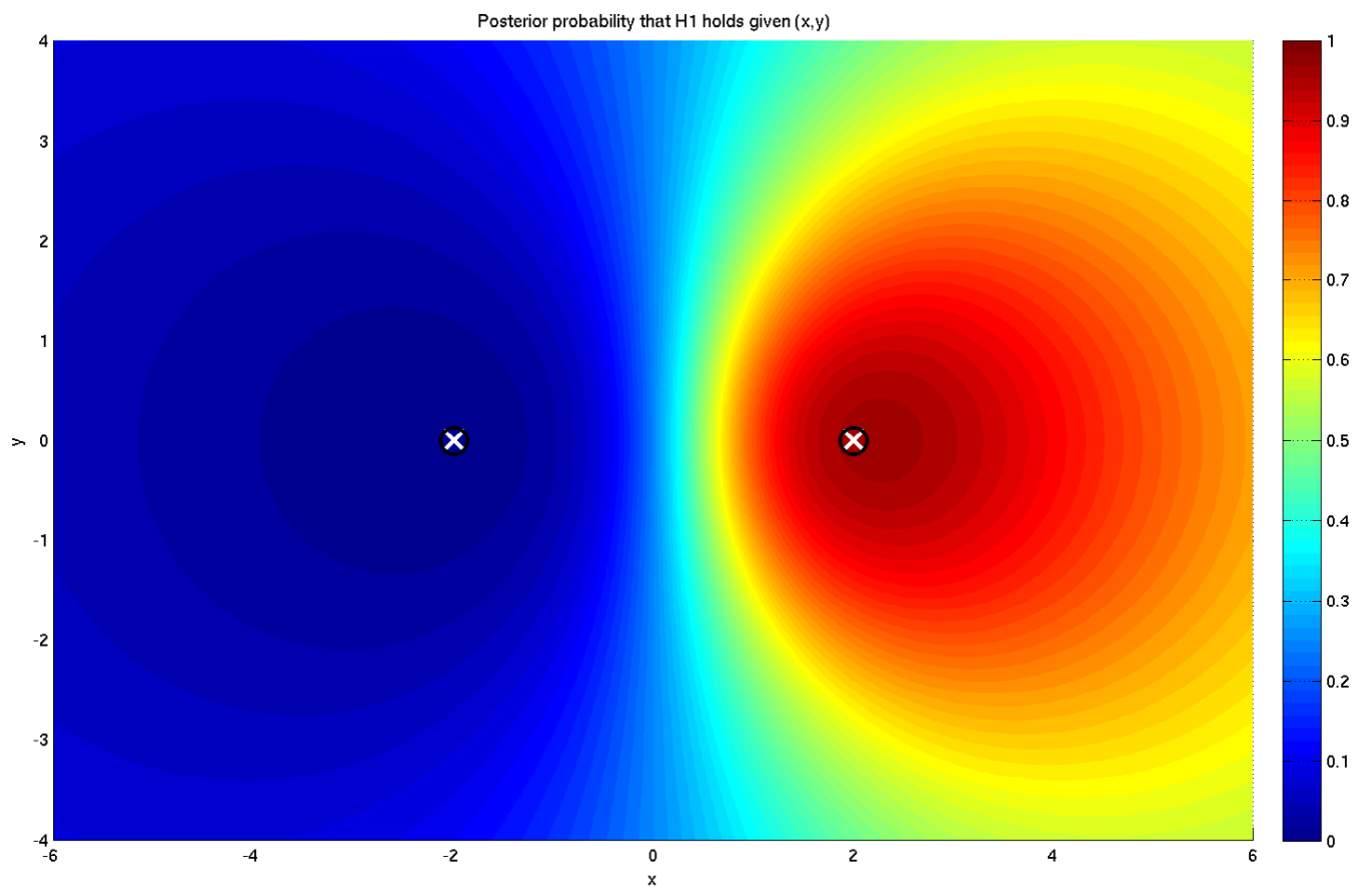}

\caption{The posterior probability that $H_1$ holds as a function of
  $(x,y)$, working from a prior that puts three-quarters of the prior
  probability on $H_0$. The positions of the two towers are also shown.
\label{postprobbias}
}

\end{center}
\end{figure}

\begin{figure}[p]
\begin{center}

\includegraphics[scale=0.5]{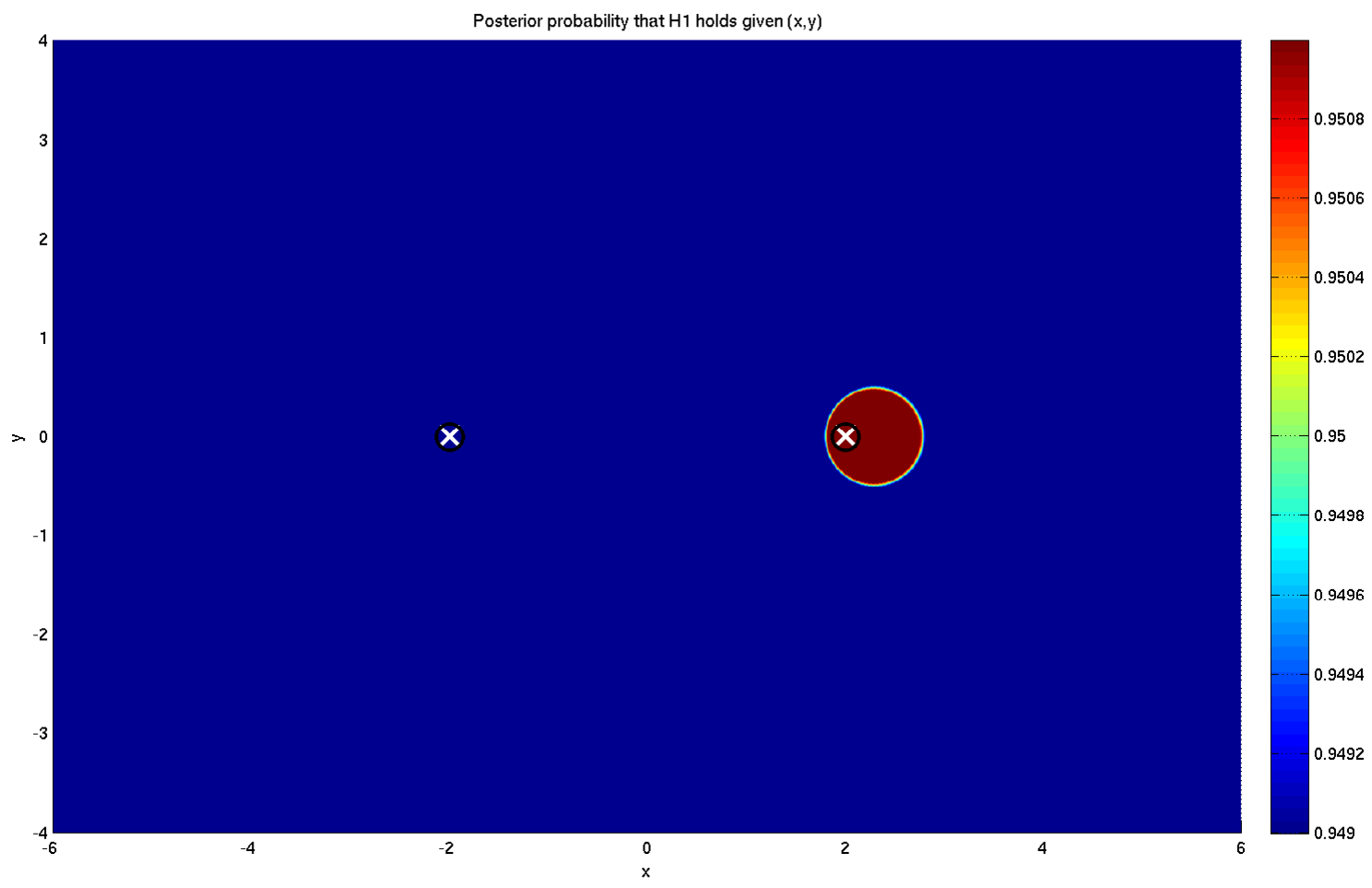}

\caption{The region of the $(x,y)$ plane in which, using the prior
  biased in favour of $H_0$, the posterior probability that $H_1$
  holds is at least $0.95$ is shown in brown. The positions of the two
  towers are also shown.
\label{postprobbias95}
}

\end{center}
\end{figure}

Of course, our prince may have a reason to suppose that the kidnapper
is more likely to put the princess in tower $0$, in which case he
might choose to set $P(h=0)=\frac{3}{4}$ and $P(h=1)=\frac{1}{4}$. In
that case instead of figures \ref{postprobeven} and
\ref{postprobeven95} we get figures \ref{postprobbias} and
\ref{postprobbias95}, using the formula $$P(h=1|x,y) =
\frac{\frac{1}{4}\frac{z}{4 \pi (z^2 + (x-2z)^2 +
    y^2)^{\frac{3}{2}}}}{\frac{1}{4}\frac{z}{4 \pi (z^2 + (x-2z)^2 +
    y^2)^{\frac{3}{2}}} + \frac{3}{4}\frac{z}{4 \pi (z^2 + (x+2z)^2 +
    y^2)^{\frac{3}{2}}}}.$$ 

Notice that the actual sets in the $x,y$-plane giving each level of
posterior probability are the same set of sets as with the even prior,
but they have different posterior probability levels attached to them.

\clearpage

\subsection{Frequentist hypothesis testing solution(s) - including the pseudo-Bayesian approach}

\subsubsection{Recap of frequentist basic principles}

Turning to frequentist approaches, we first remind ourselves of the
basics of frequentist hypothesis testing, as applied to the simple
case of two hypotheses that each contain just one possibility. Given
two alternative hypotheses, we first name one of them $H_0$ (our
``null hypothesis'') and its alternative $H_1$.

The key priority of a frequentist is then to ensure that we ``control
the type I error probability'', or in other words that we conduct our
inference in such a way that we ensure that the probability of
concluding that $H_1$ holds when actually $H_0$ is true is less than
or equal to $1-c$, where $c$ is the degree of frequentist confidence
we wish to claim. In order to achieve this, we need to define a
``critical region'' $C_c$, a subset of data space $X$, for values of
$(x,y)$ in which we will conclude that $H_1$ holds. For values of
$(x,y)$ not in $C_c$ we will either conclude that $H_0$ holds, or (if
we are a strict flavour of frequentist) that we cannot tell which
hypothesis holds (in which case this setup gives no way of ever
concluding that $H_0$ holds). In order to control the type I error
probability, we require that $P((x,y)\in C_c|H_0) \leq 1-c$.

Now, in the case (as in this example) that the likelihood $P((x,y)|h)$
is a continuous distribution, the probability, given either $h$, of
getting a single particular value of $(x,y)$ is zero. So if we are
allowed to define $C_c$ \textit{after} collecting the data, we can
always set $C_c=\{(x,y)\}$, easily satisfying the requirement, and
conclude that $H_1$ holds, at least so long as the bullet did land
somewhere. We therefore have to insist that $C_c$ be defined
\textit{before} collecting the data (although this requirement is
often ignored by frequentists, and where it is not ignored, it can be
difficult to prove that it hasn't been).

This can be generalised slightly in order to avoid having to say what
value of $c$ we are interested in before collecting the data. If
instead we define a nested family of critical regions $(C_\eta)_{ \eta
  \in [0,1]}$ such that $$\eta_1\leq \eta_2\implies
C_{\eta_2}\subseteq C_{\eta_1}$$ and such that for all $\eta\in [0,1]$
and all $h\in H_0$, $$P((x,y)\in C_\eta|h)\leq 1-\eta,$$ then we can
uniquely determine the frequentist confidence that $H_1$ holds that is
achieved by observing any particular $(x,y)$ by concluding
that $$c=\sup{(\{\eta\in [0,1]:(x,y)\in C_\eta\}\cup\{0\})}.$$ In this
particular case $H_0$ consists of only a single value of $h$, so the
condition can be simplified to read $$P((x,y)\in C_\eta | H_0) \leq
1-\eta.$$

\subsubsection{The part common to all frequentist hypothesis testing
  approaches}
\label{commonparts}

The prince's first job, then, is to decide which hypothesis is $H_0$
and which is $H_1$. Since one tower is painted with a large zero, he
sets $H_0$ to be the hypothesis that the prince is in tower $0$. Note
that this breaks the symmetry of the problem, as $H_0$ plays a
different role in the frequentist approach than $H_1$; and as we will
see, this symmetry is (in this problem) never restored, at least in
the hypothesis testing formulation.

He then needs to decide whether he is a strict frequentist (who can
never conclude that the princess is in tower $0$) or a non-strict
one. Since in the former case the princess is doomed if she is in
tower $0$, we will assume he is a non-strict frequentist.

The next issue to consider is how we are going to deal with the case
that the bullet has gone upwards into space and has not landed. There
are two classes of approach:

\begin{description}

\item[Randomised]: The prince could, for example, toss an eicosahedral
  (20-sided) die, and if it comes down $20$ decide that the princess
  is in tower $1$, otherwise that she is in tower $0$. That uses up
  probability of $\frac{1}{2}\times \frac{1}{20}=0.025$ of making a
  type I error. So if he now ensures that what he does if the bullet
  lands somewhere can't accumulate any more than a further $0.025$ of
  probability of type I error, he can stay within a total type I error
  rate of $0.05$ and therefore potentially get a $95\%$ frequentist
  confident conclusion.

\item[Deterministic]: \begin{description} 
\item[Favouring $H_0$]: He could just decide that if the bullet goes
  upwards into space then he will conclude that the princess is in
  tower $0$, which doesn't contribute anything to the type I error
  rate.
\item[Favouring $H_1$]: He could just decide that if the bullet goes
  upwards into space then he will conclude that the princess is in
  tower $1$, which contributes $\frac{1}{2}$ to the type I error rate,
  so that whatever happens under other circumstances he can never be
  more than $50\%$ frequentist confident that she is in tower 1.
\end{description}

\end{description}

We will start by assuming that the prince deals with this case in the
Deterministic way favouring $H_0$. That leaves all of the type I error
rate to be handled in the case that the bullet lands at $(x,y)$. Later
we will consider two non-deterministic approaches.

Just as the Bayesian approach next had to consider what prior
probabilities to set on $H_0$ (and hence $H_1$), the frequentist
approach now also faces a dilemma, namely which nested family of
critical regions to pick.

There are an infinite range of possible choices; we will explore a few
of them.

\subsubsection{Minimum excluded area}

\begin{figure}[ht]
\begin{center}

\includegraphics[scale=0.5]{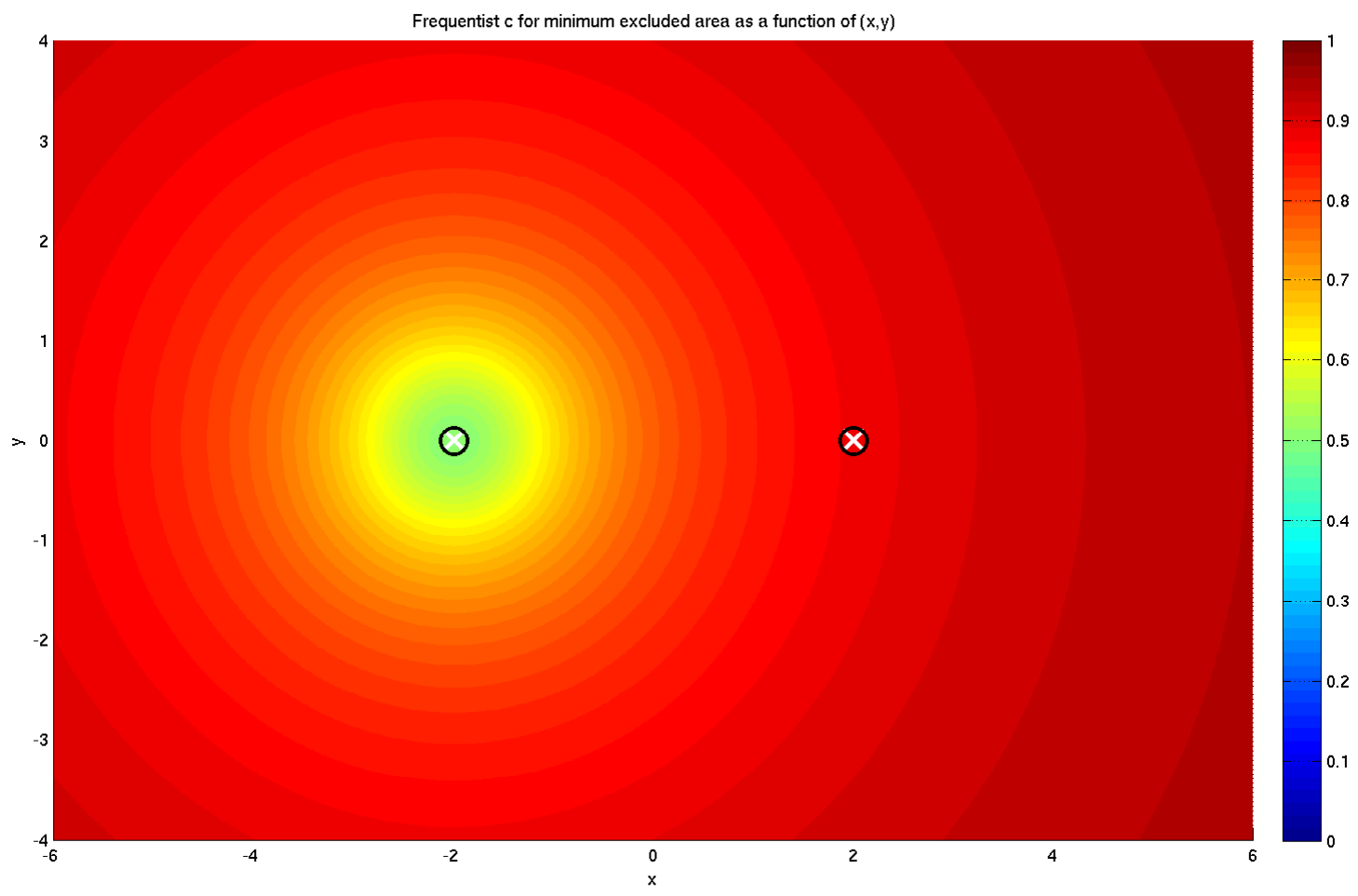}

\caption{Frequentist confidence $c$ that $H_1$ holds as a function of
  $(x,y)$ resulting from using the critical regions that exclude the
  minimum possible area of the plane. The positions of the two
  towers are also shown.
\label{freqmaxarea}
}

\end{center}
\end{figure}

First, for each $\eta\in [0,1]$, we could set $C_\eta$ such that the
area of the plane excluded from $C_\eta$ is as small as possible. That
gets us the nested set of critical regions that result in the
frequentist confidence $c$ that $H_1$ holds, as a function of $(x,y)$,
being as shown in figure \ref{freqmaxarea}; here $$C_\eta = \{(x,y)\in
\mathbb{R}^2: (x+2z)^2 + y^2 > R_\eta^2\}$$ where $$R_\eta =
\left\{\begin{matrix} z \sqrt{\frac{1}{4(1-\eta)^2}-1} &
(\eta\geq\frac{1}{2})\\ 0 & (\eta < \frac{1}{2})\end{matrix}\right.;$$
in other words $C_\eta$ is the region outside a circle of radius
$R_\eta$ centred on tower $0$.

This figure is shown on the same scale as the previous plots -- yet
we notice that for \textit{no} value of $(x,y)$ in the region of the
plot can we conclude that we are 95\% frequentist-confident that $H_1$
holds; the maximum frequentist-confidence we can get while staying
within the plot is at the top and bottom right-hand corners, where we
get about $94.4\%$ frequentist-confidence that $H_1$ holds. On the
other hand we \textit{could} conclude that we had 95\% frequentist
confidence that $H_1$ holds if we observed $(x,y)=(-22,0)$, a point
way off the \textit{left}-hand edge of the plot (note that $H_1$ is
that the true value of $(a,b)$ is $(+2,0)$).

I.e. even if the bullet were found right at the foot of tower $1$,
this particular frequentist approach would not be able to reach 95\%
confidence that the princess was in tower $1$, but would if it were
found 22 km on the far side of tower $0$ (which is nonsensical).

\clearpage

\subsubsection{Circular critical regions around tower $1$}

\begin{figure}[p]
\begin{center}

\includegraphics[scale=0.5]{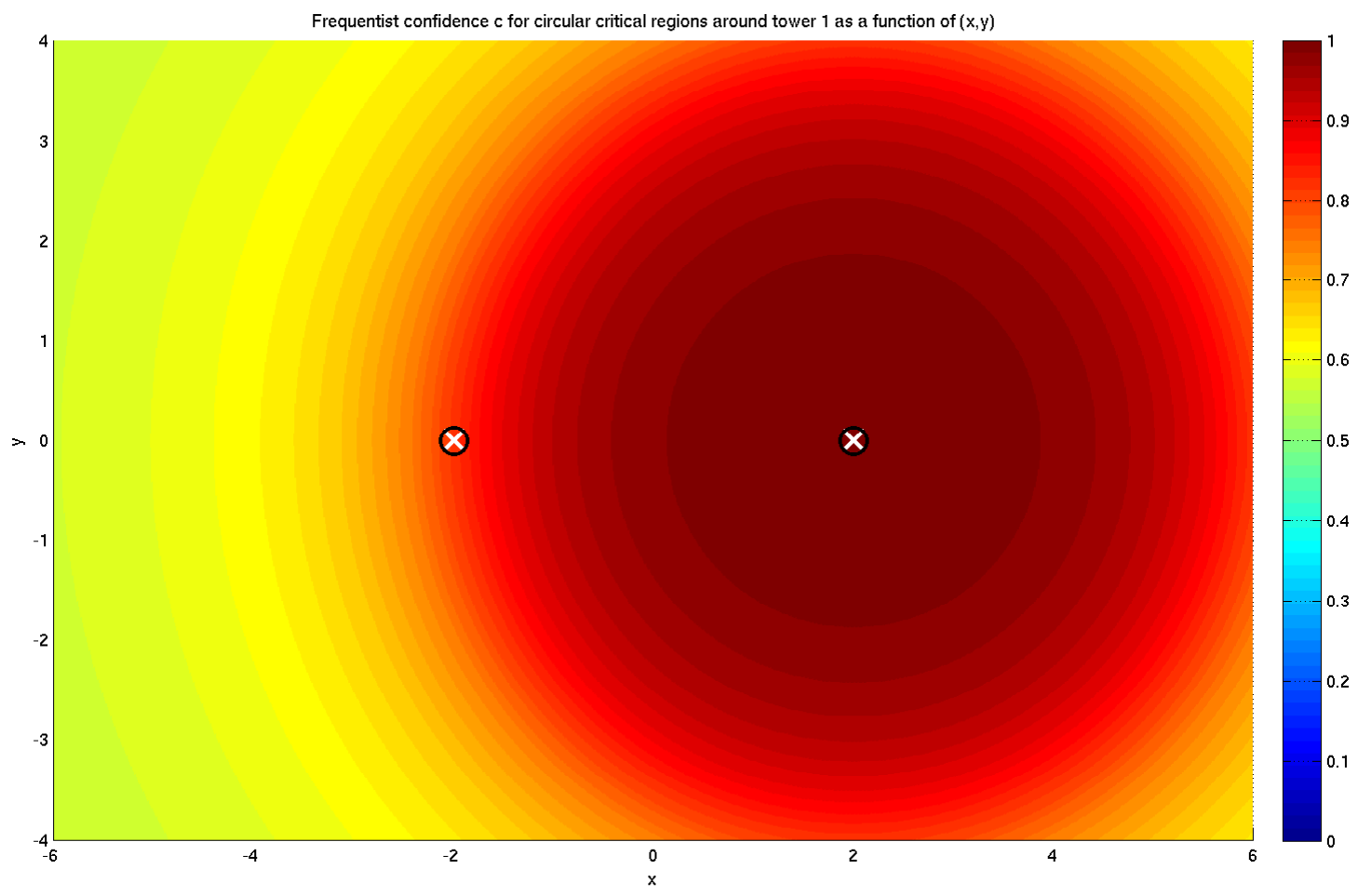}

\caption{Frequentist confidence $c$ that $H_1$ holds as a function of
  $(x,y)$ resulting from using the circular critical regions centred
  on tower $1$. The positions of the two towers are also shown.
\label{freqct1}
}

\end{center}
\end{figure}

\begin{figure}[p]
\begin{center}

\includegraphics[scale=0.5]{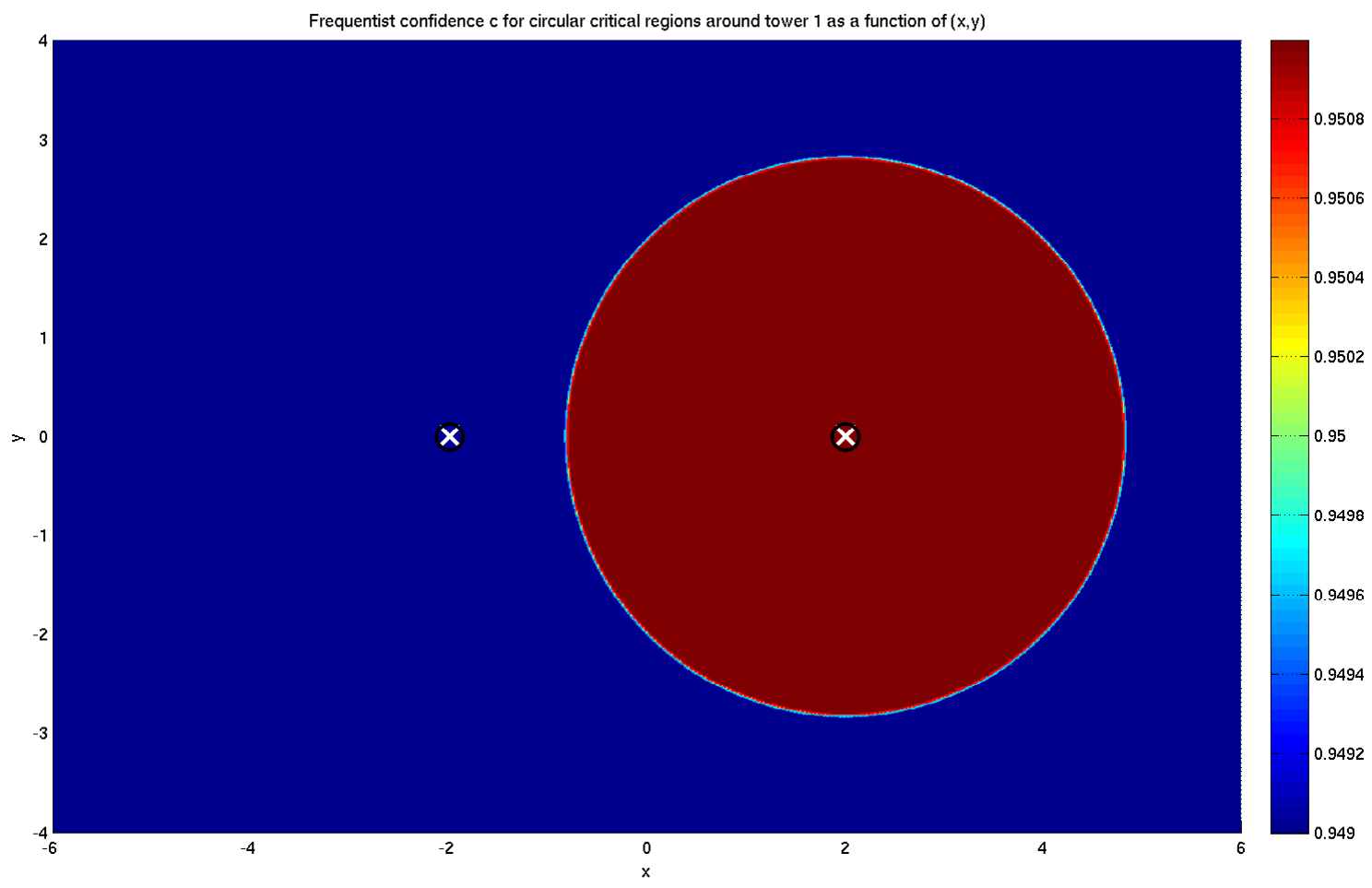}

\caption{The 95\% critical region $C_{0.95}$ from the set of critical
  regions that are circular and centred on tower $1$ is shown in
  brown.  The positions of the two towers are also shown.
\label{freqct195}
}

\end{center}
\end{figure}

Noting that the critical regions of minimum excluded area are the
complements of discs centred on tower $0$, we might consider critical
regions that are discs centred on tower $1$. This idea yields the
frequentist confidence plot shown in figure \ref{freqct1}. Notice that
this plot is \textit{not} the complement of that of figure
\ref{freqmaxarea} obtained by subtracting the values of the latter
from one (the base of tower $0$ gives $50\%$ frequentist confidence
that the princess is in tower $1$ in figure \ref{freqmaxarea} but
$82\%$ in figure \ref{freqct1}). The corresponding region giving
$95\%$ frequentist confidence that the princess is in tower $1$ is
shown in figure \ref{freqct195} -- and unlike in figure
\ref{freqmaxarea} there is a large region where we become $95\%$
frequentist confident of this. Indeed, should the bullet be found
exactly half way between the two towers at $(0,0)$ we are now about
$98\%$ frequentist confident that the princess is in tower $1$,
despite the fact that frequentist solutions above all want to avoid
concluding that the princess is in tower $1$ when she is actually in
tower $0$ ! Even worse, if the bullet is found there we become $>95\%$
frequentist confident that the princess is in tower $1$; yet with the
same data but opposite $H_0$ we become $>95\%$ frequentist confident
that she is in tower $0$ -- so we are either $>95\%$ confident that
she is in tower $0$ or $>95\%$ confident that she is in tower $1$,
depending only on the arbitrary decision of which tower to make the
null hypothesis.

\clearpage

\subsubsection{Critical regions that depend only on $x$-coordinate}

\begin{figure}[p]
\begin{center}

\includegraphics[scale=0.5]{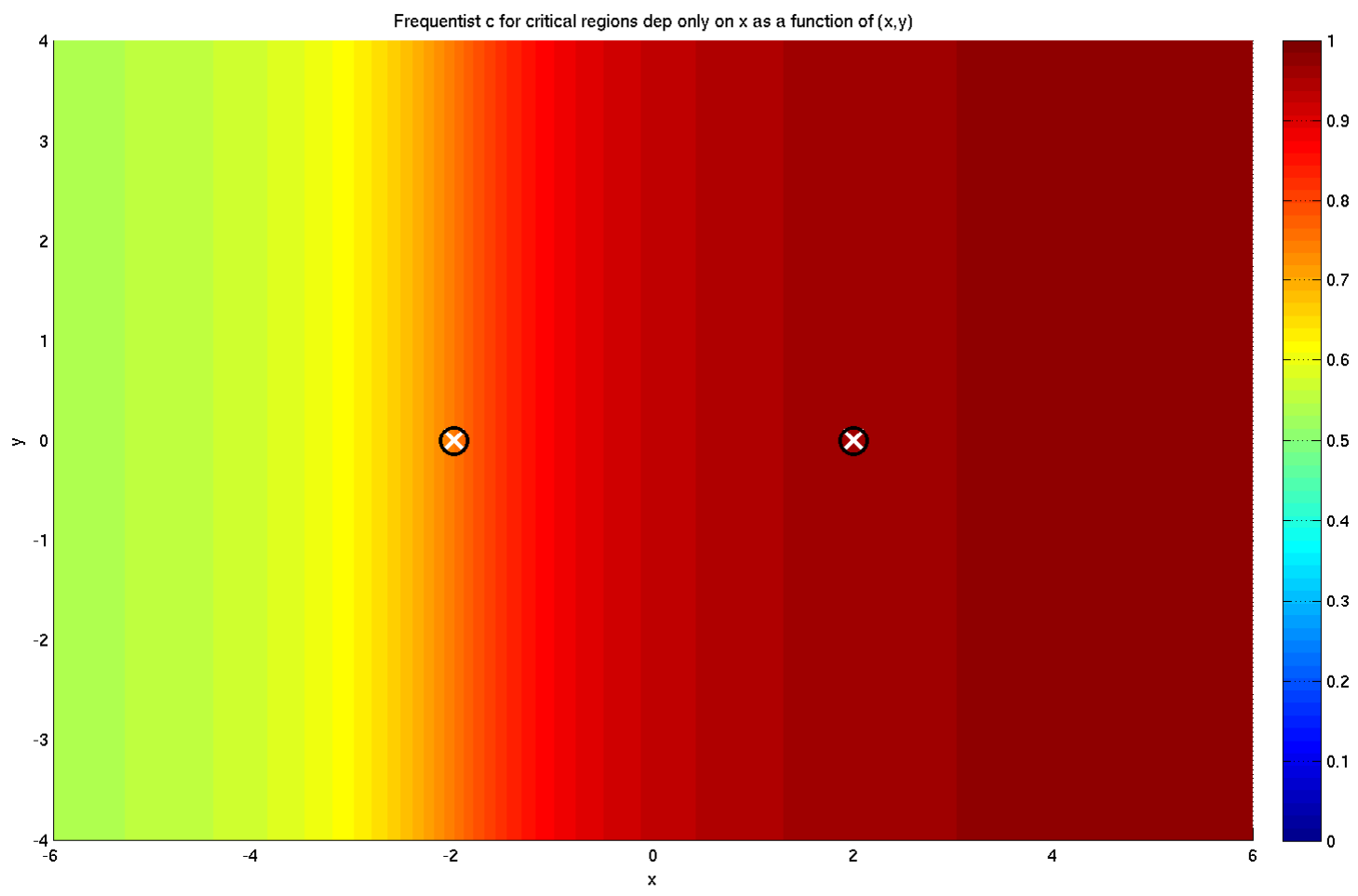}

\caption{Frequentist confidence $c$ that $H_1$ holds as a function of
  $(x,y)$ resulting from using critical regions that depend only on
  the $x$-coordinate and which shrink to the right. The positions of
  the two towers are also shown.
\label{freqxonly}
}

\end{center}
\end{figure}

Another idea for making the critical region $C_\eta$ might be to
choose an $x_\eta$ such that $$C_\eta= \{(x,y)\in \mathbb{R}^2:x >
x_\eta\}$$ for some appropriate $x_\eta$ satisfying $P(x >
x_\eta|H_0)=1 - \eta$. This gets us the nested set of critical regions
that result in the frequentist confidence $c$ that $H_1$ holds, as a
function of $(x,y)$, being given as shown in figure \ref{freqxonly}.

\begin{figure}[htp]
\begin{center}

\includegraphics[scale=0.5]{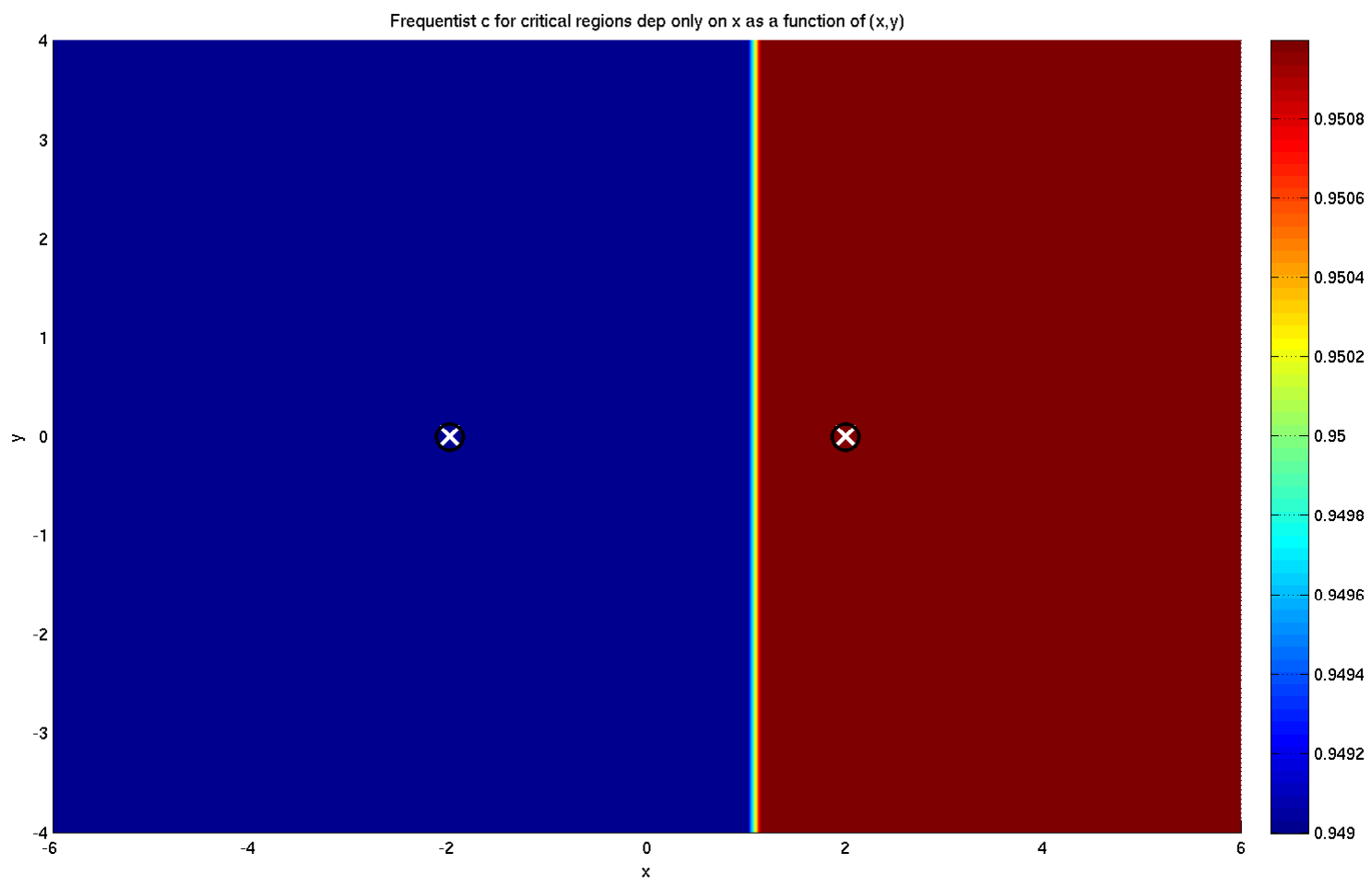}

\caption{The 95\% critical region $C_{0.95}$ from the set of critical
  regions that depend only on $x$ and which shrink to the right is
  shown in brown. The positions of the two towers are also shown.
\label{freqxonly95}
}

\end{center}
\end{figure}

In this case we find that we can become $95\%$ frequentist-confident
that $H_1$ holds if $x$ is greater than about $1.08$, as shown in
figure \ref{freqxonly95}. On the other hand even if the bullet is
found at the base of tower $0$ we become $75\%$ frequentist confident
that the princess is in tower $1$ (though the Bayesian solution with
even prior tells us that e.g. for $(x,y)=(-2,0)$ we are about $98.7\%$
sure that $H_0$ holds).

(Moreover we could be really perverse, and pick a set of critical
regions made by taking those just discussed, and reflecting them about
the line $x=-2$: this makes an equally valid set of critical regions
which now become 95\% frequentist confident that $H_1$ holds only when
$x < -5.08$ approx., i.e. when we are a long way on the opposite side
of tower $0$.)

\begin{figure}[p]
\begin{center}

\includegraphics[scale=0.5]{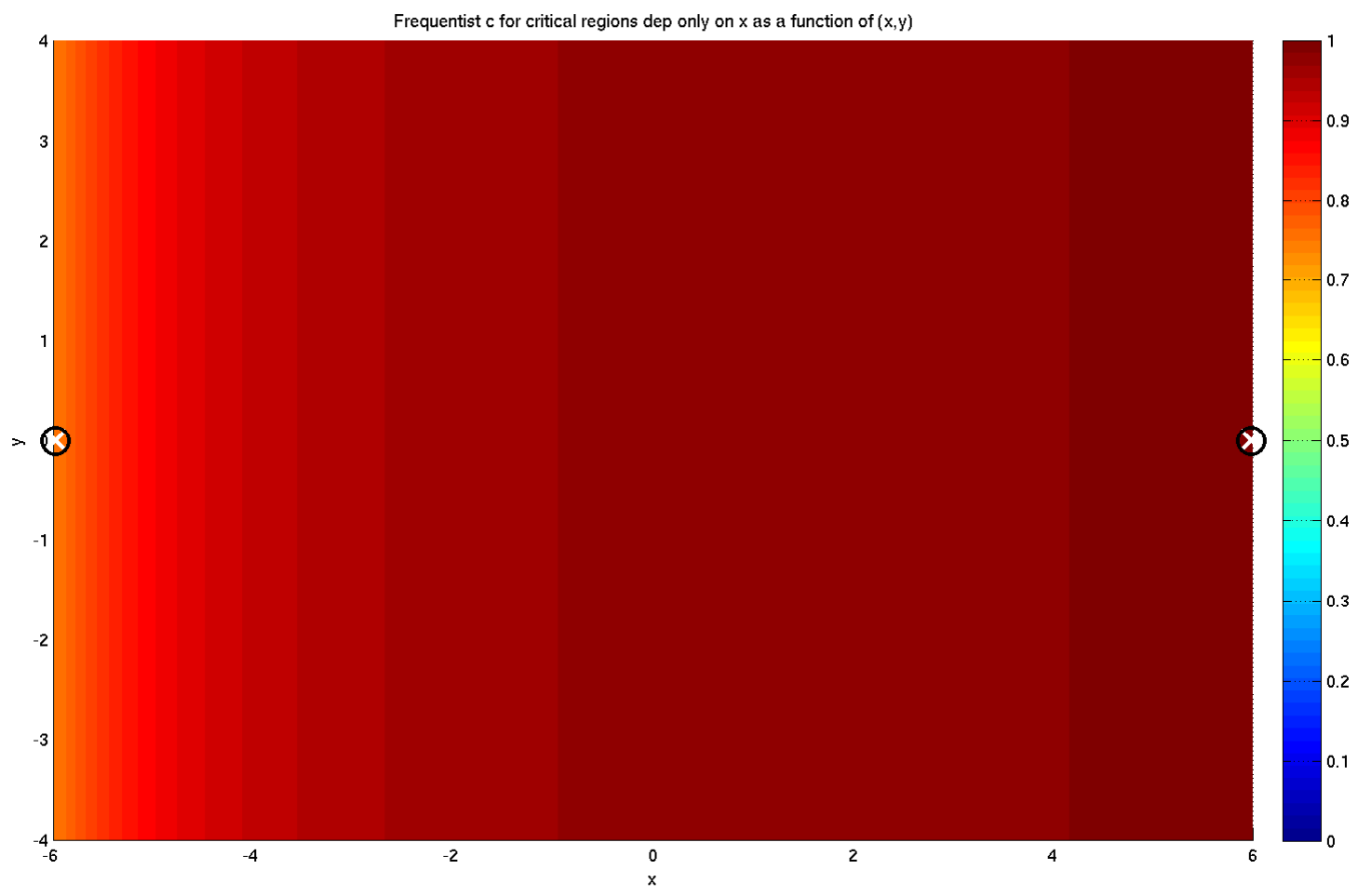}

\caption{Frequentist confidence $c$ that $H_1$ holds as a function of
  $(x,y)$ resulting from using critical regions that depend only on
  the $x$-coordinate and which shrink to the right. The positions of
  the two towers are also shown - in this case they are 12 km apart
  instead of the usual 4 km apart.
\label{freqxaonly}
}

\end{center}
\end{figure}

\begin{figure}[htp]
\begin{center}

\includegraphics[scale=0.5]{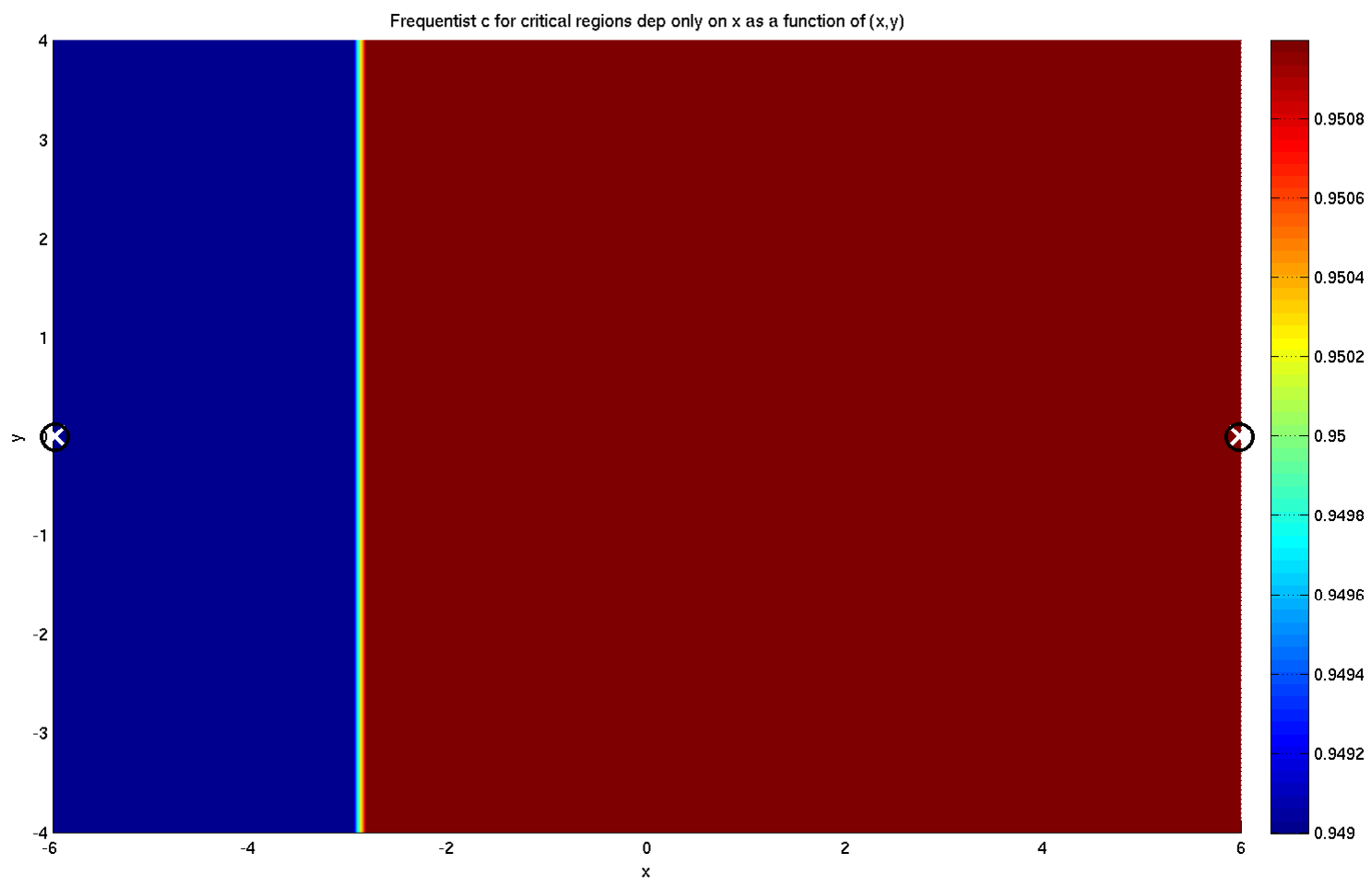}

\caption{The 95\% critical region $C_{0.95}$ from the set of critical
  regions that depend only on $x$ and which shrink to the right is
  shown in brown. The positions of the two towers are also shown - in
  this case they are 12 km apart instead of the usual 4 km apart.
\label{freqxaonly95}
}

\end{center}
\end{figure}

Some may be tempted to think that the solution of figure
\ref{freqxonly} is the ``obvious'' correct solution to this
problem. To see that it really isn't, let us consider what happens if
we temporarily move the two towers to being at $(-6,0)$ and $(+6,0)$
instead of their usual positions. Then we get figures \ref{freqxaonly}
and \ref{freqxaonly95}. 

Now, if the bullet lands at an undisclosed point on the line $x=-2.92$
we are $95\%$ frequentist confident that the princess is on tower $1$,
despite this event being nearly $8$ times more likely if the princess
is on tower $0$ than on tower $1$. Indeed we also so conclude if the
bullet lands specifically at $(x,y)=(-2.92,0)$ despite this event
being nearly $21$ times more likely if the princess is on tower $0$
than on tower $1$. At the same time if the bullet lands at any point
on the line $x=0$ we are over $97\%$ frequentist confident that the
princess is on tower $1$, despite the obvious symmetry of the
situation. Hopefully these facts are persuasive that this choice of
critical regions doesn't remove the difficulties.

\clearpage

\subsubsection{Critical regions based on direction from origin}

\begin{figure}[ht]
\begin{center}

\includegraphics[scale=0.5]{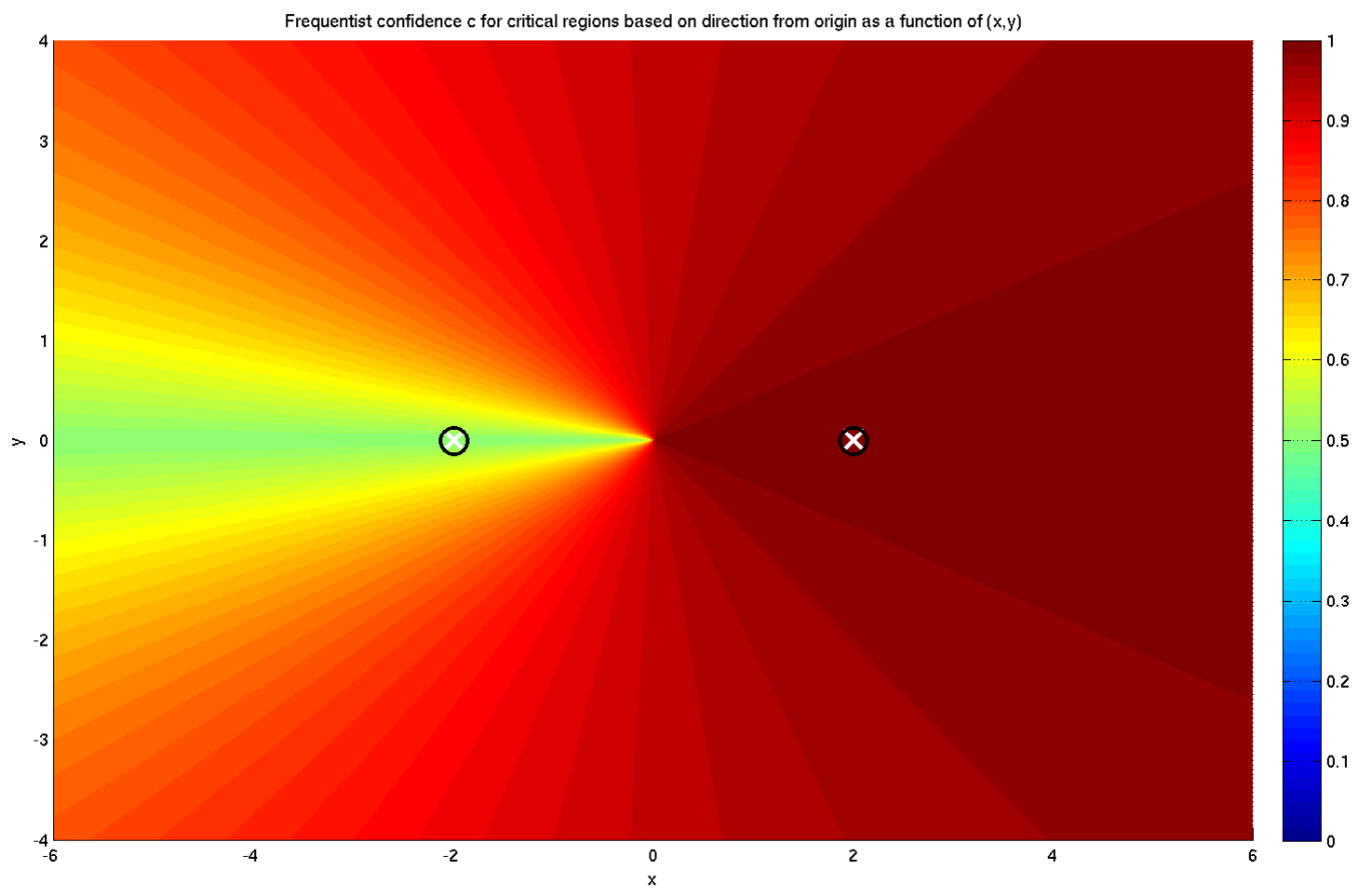}

\caption{Frequentist confidence $c$ that $H_1$ holds as a function of
  $(x,y)$ resulting from using critical regions that depend only on
  the direction from the origin. The positions of the two towers are
  also shown.
\label{freqdirection}
}

\end{center}
\end{figure}

Returning to the usual positions of the towers at $(-2,0)$ and
$(+2,0)$, another way to pick critical regions would be to base them
on direction from origin, making $C_\eta$ be a sector symmetric about
the $x$-axis; this gives the frequentist confidence shown in figure
\ref{freqdirection}. Here the origin is on the boundary of $C_\eta$
for all $\eta \in [0,1]$, so if we make the $C_\eta$s be closed
sectors, the bullet being found at the origin gives us $100\%$
frequentist confidence that the princess is in tower 1, while if we
make the $C_\eta$ be open sectors (i.e. not including the boundary) we
instead get $50\%$.

Either way, if the bullet is found 1 mm to the East of the origin, we
become $100\%$ frequentist confident that the princess is in tower 1.

\subsubsection{The basic pseudo-Bayesian solution}
\label{freqpseudoBayes}

Assuming that our prince doesn't try deciding which tower to climb
before the gun has been fired, there is then only one possible basic
pseudo-Bayesian solution for this problem: there are no options to
change the data collection plan (we are given only one data point),
$H_0$ cannot be enlarged without making it the whole of $H$, and
changing the prior makes no difference to the $C_\eta$s (as it turns
out that although it changes the value of $p=g(\eta)$ attached to each
$C_\eta$, the $C_\eta$s themselves remain unchanged).

\begin{figure}[htp]
\begin{center}

\includegraphics[scale=0.5]{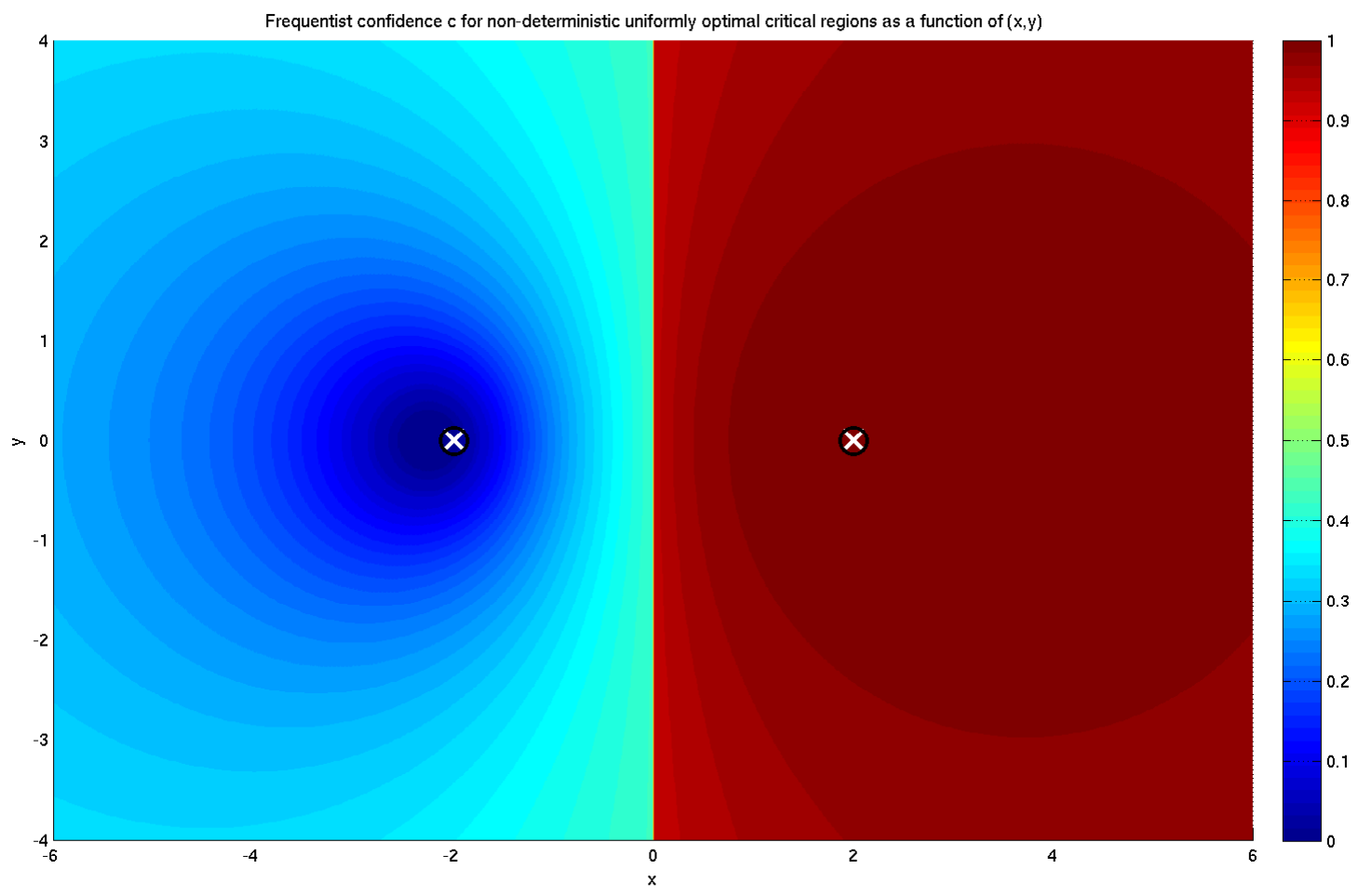}

\caption{The frequentist confidence that $H_1$ holds as a function of
  $(x,y)$ when using the basic pseudo-Bayesian set of critical
  regions. The positions of the two towers are also shown. Note that
  with this approach if the bullet is not found, the frequentist confidence value
  $\frac{\tan^{-1}2}{2 \pi}+\frac{1}{4}\approx 0.426$ is reported.
\label{freqpseudoBayesfig}
}

\end{center}
\end{figure}

Applying then the method described in section
\ref{concisedescription}, we get figure \ref{freqpseudoBayesfig}. It
is immediately obvious that something odd is going on, as there is a
big jump in the frequentist confidence as we cross the y-axis from
about 0.426 to about 0.926. The reason is that this solution has decided
to include the event that the bullet hasn't been found, which has
probability 0.5, in all $C_\eta$ for $\eta\leq\frac{\tan^{-1}2}{2
  \pi}+\frac{1}{4}\approx 0.426$.

The critical regions are in fact the same \textit{set} of regions as
the regions of various levels of posterior probability, but the
frequentist confidence values attached to them are very different from
the posterior probability values attached to them.

Note, for comparison with section \ref{nearpseudoBayes},
that $$P((x,y)\in
C_{\frac{1}{2}}|H_1)=\frac{1}{4}-\frac{\tan^{-1}2}{2\pi} \approx
0.074\ .$$

\subsubsection{A nearly pseudo-Bayesian solution uniformly optimal
  deterministic for $\eta \geq 0.5$}
\label{nearpseudoBayes}

Now let us aim to remove the discontinuity from the solution of section
\ref{freqpseudoBayes} by modifying the nested set of critical regions
as follows. First remove the event that the bullet isn't found from
all the $C_\eta$; then reevaluate the appropriate frequentist
confidence for each of the remaining $C_\eta$, relabelling them (they
now all have $\eta\geq 0.5$); then set $C_\eta=C_{\frac{1}{2}}$ for
all $\eta \in (0,0.5)$ and $C_0=X$, so that we again have zero
frequentist confidence that the princess is in tower 1 if the bullet
isn't found.

\begin{figure}[p]
\begin{center}

\includegraphics[scale=0.5]{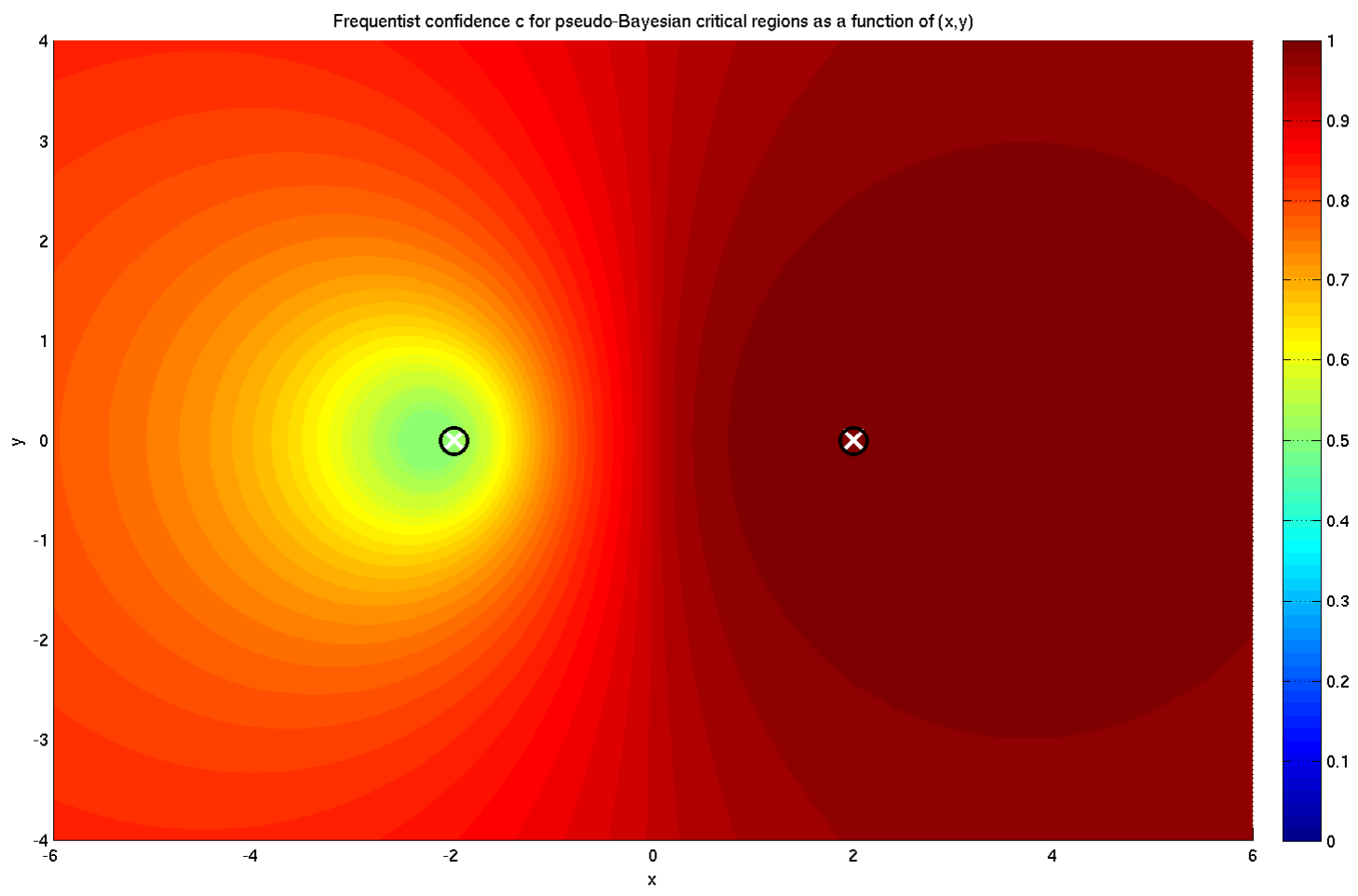}

\caption{The frequentist confidence that $H_1$ holds as a function of
  $(x,y)$ when using the nearly pseudo-Bayesian, uniformly optimal
  deterministic, set of critical regions. The positions of the two
  towers are also shown.
\label{freqpseudoB}
}

\end{center}
\end{figure}

The resulting solution is in fact what frequentists call ``uniformly
optimal'' for this problem, at least among deterministic approaches
and for frequentist confidence levels $\eta \geq 0.5$. In other words,
for each possible level $\eta\geq 0.5$ of frequentist confidence, and
for any other valid critical region $C'_\eta$, under $H_1$ we are at
least as likely to get the data falling into $C_\eta$ as into
$C'_\eta$. I.e. this is the set of critical regions which makes it
most likely that we will conclude that $H_1$ holds when it
does\footnote{Frequentists are not usually interested in values of
  $\eta<0.5$. However, for completeness, the critical regions shown in
  figure \ref{unifopt} below are overall uniformly optimal, and are
  deterministic except for $\eta$ in the open set
  $\frac{\tan^{-1}2}{2\pi}+\left(\frac{1}{4},\frac{3}{4}\right)\approx
  (0.426,0.926)$. Above $0.926$ they coincide with those in figure
  \ref{freqpseudoB}, and below $0.426$ they are uniformly optimal
  deterministic. In the range
  $\left(\frac{\tan^{-1}2}{2\pi}+\frac{1}{4},\frac{1}{2}\right)$ there
  is a third nested family of critical regions which is uniformly
  optimal deterministic, which the interested reader may like to
  find. However no two of these three nested sets of critical regions
  can be combined into a nested set.}. (Note that for most problems no
such uniformly optimal family of critical regions exists.) In
particular, $P(x\in C_{\frac{1}{2}}|H_1)=\frac{1}{2}$, showing that
for $\eta=\frac{1}{2}$ the basic pseudo-Bayesian solution is not
uniformly optimal.

\begin{figure}[p]
\begin{center}

\includegraphics[scale=0.5]{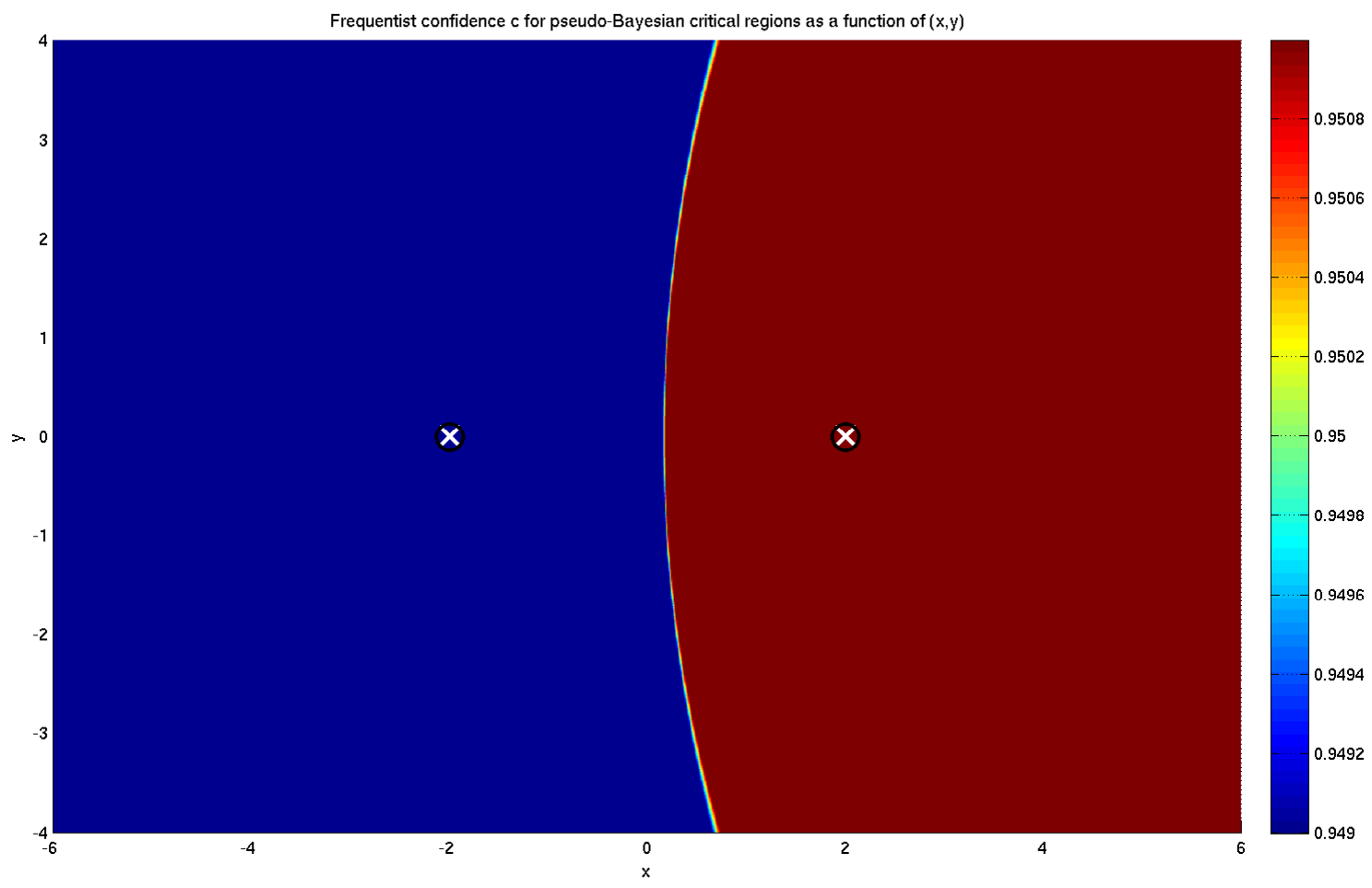}

\caption{The 95\% critical region $C_{0.95}$ from the nearly
  pseudo-Bayesian, uniformly optimal deterministic, set of critical
  regions is shown in brown. The positions of the two towers are also
  shown.
\label{freqpseudoB95}
}

\end{center}
\end{figure}

The critical regions are again the same \textit{set} of regions as
the regions of various levels of posterior probability, except that
the event that the bullet doesn't land is only included for $\eta=0$;
but the frequentist confidence values attached to them are again very
different from the posterior probability values attached to them. The
frequentist confidence $c$ that $H_1$ holds for each possible value of
$(x,y)$ is shown in figure \ref{freqpseudoB}, and the specific $95\%$
region in figure \ref{freqpseudoB95}.

Comparing figures \ref{postprobeven} and \ref{freqpseudoB}, we see
that \textit{this} frequentist solution is far more likely to conclude
that $H_1$ holds than is justified according to the Bayesian solution
-- don't let anybody ever tell you that Bayesian methods are always
more likely to get you a significant result ! Indeed, if
$(x,y)=(0,0)$, the point mid-way between the two towers, this
frequentist solution says it is $92.6\%$-frequentist-confident that the
princess is in tower $1$, despite both the obvious symmetry of the
problem, and its supposed concern to avoid thinking that she is in
tower $1$ when actually she is in tower $0$.

\begin{figure}[htp]
\begin{center}

\includegraphics[scale=0.5]{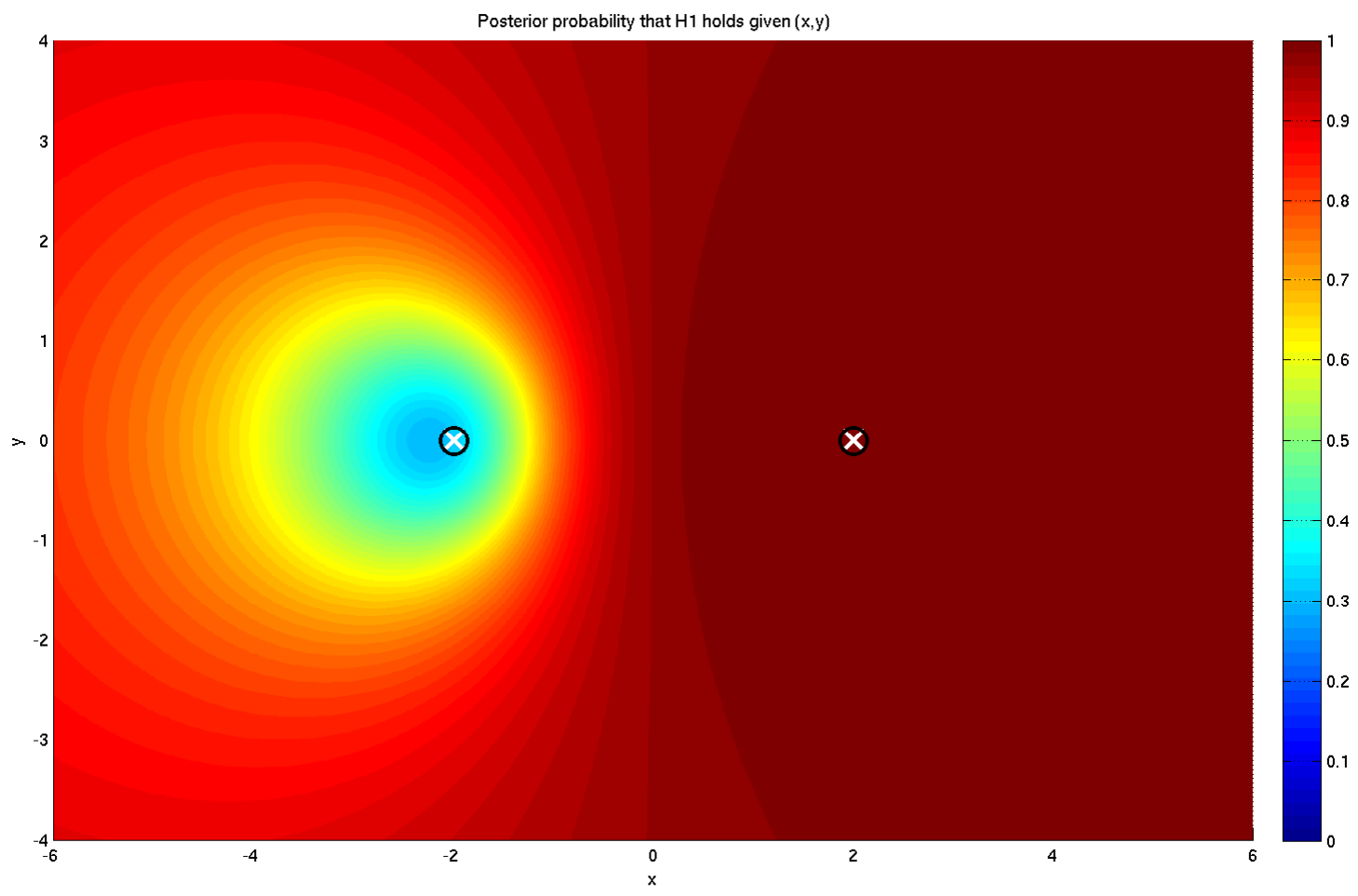}

\caption{The posterior probability that $H_1$ holds as a function of
  $(x,y)$ using $P(h=1)=0.97$ . The positions of the two towers are
  also shown.
\label{postprobupbias}
}

\end{center}
\end{figure}

One might justifiably wonder whether the Bayesian and nearly 
pseudo-Bayesian solutions could be made to match up by setting an
appropriately high value for the prior probability $P(h=1)$. It turns
out one can not do this: as $P(h=1)$ is gradually increased from $0.5$
to $0.9$, the dark red area on the right of figure \ref{postprobeven}
expands, while the blue one on the left shrinks, but at $0.97$ one
gets figure \ref{postprobupbias}, and on comparing this with figure
\ref{freqpseudoB}, we see that the right hand side of the plot is
already too dark (red), while the left hand side has still got too
much blue in it. But it is also noteworthy that a solution whose primary
concern is to control the type I error rate is closest to a Bayesian
solution which favours $H_1$ with a prior probability of about $0.97$.

\subsubsection{Solution based on random numbers}
\label{freqrandom}

Turning now to consider non-deterministic critical regions, let us
draw a random number $u$ uniformly distributed on $[0,1]$ in addition
to knowing $(x,y)$ (or that the bullet didn't land anywhere). Our data
is now either $(x,y,u)$ or $(\text{not found}, u)$. Let us set $C_\eta
= \{(x,y,u):u\geq \eta\} \cup \{(\text{not found},u):u\geq
\eta\}$. Then the probability, if $H_0$ holds, of landing in $C_\eta$
is $1-\eta$, showing that these are a valid nested set of critical
regions.

But using them is equivalent to ignoring the position (if any) of the
bullet and simply drawing a random number $u$ uniformly distributed in
$[0,1]$ and reporting it as the frequentist confidence.

\subsubsection{The full pseudo-Bayesian solution}

\begin{figure}[htp]
\begin{center}

\includegraphics[scale=0.5]{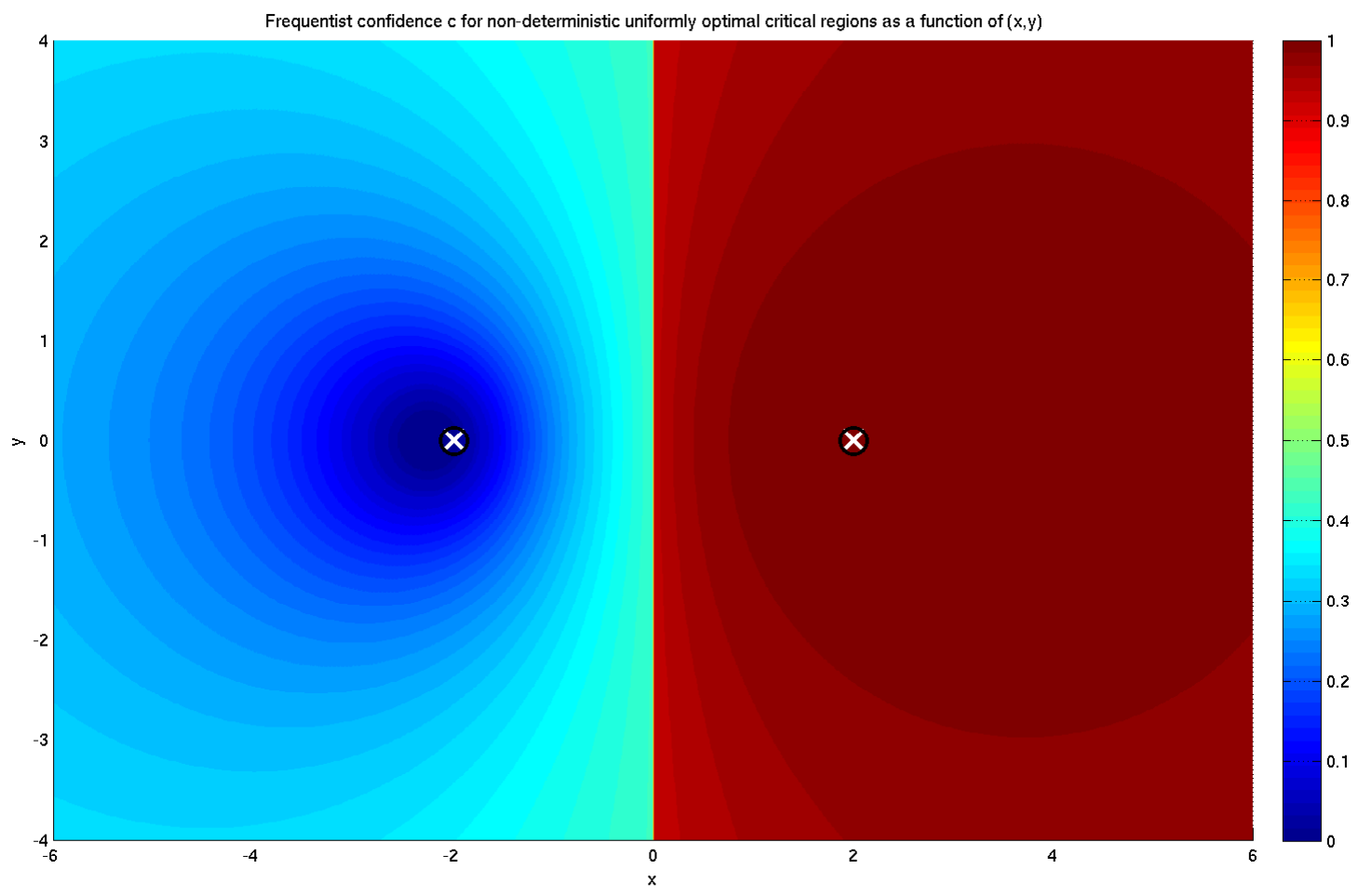}

\caption{The frequentist confidence that $H_1$ holds as a function of
  $(x,y)$ when using the overall uniformly optimal (and
  non-deterministic) set of critical regions. The positions of the two
  towers are also shown. Note that with this approach if the bullet is
  not found, an independent random number $u$ drawn from the uniform
  distribution on $[0,1]$ is taken and the frequentist confidence
  value $\frac{u}{2}+\frac{\tan^{-1}2}{2 \pi}+\frac{1}{4}$ is reported.
\label{unifopt}
}

\end{center}
\end{figure}

Continuing the non-deterministic line of thinking, we can also make a
nested set of critical regions that is uniformly optimal not just
among deterministic sets, but among non-deterministic sets as well,
and for all $\eta\in [0,1]$. In this case we determine frequentist
confidence as in figure \ref{unifopt} if the bullet is found, but if
the bullet is not found we instead report the value
$\frac{u}{2}+\frac{\tan^{-1}2}{2\pi}+\frac{1}{4} \approx
\frac{u}{2}+0.426$. This coincides with the full pseudo-Bayesian
solution. Note that in this case there is a sudden discontinuity in the
frequentist confidence reported as the $x$-coordinate of the place the
bullet is found moves from $x\geq 0$ to $x<0$, indeed frequentist
confidence that $H_1$ holds falls by exactly $\frac{1}{2}$; moreover
there is still no symmetry between how the two towers are treated.

\subsection{Methods based on frequentist confidence sets}
\label{fconfsets}

However, some frequentist reader may say that instead of using
frequentist hypothesis testing, they would rather consider frequentist
confidence intervals/sets. The two are, of course, closely
related. For completeness, we now turn to consider the same example
from the point of view of frequentist confidence sets.

\subsubsection{Frequentist confidence sets - a reminder}
\label{confsetreminder}

We first recall the definition of a frequentist confidence set
function for a problem. We consider only the simple case that there
are no nuisance variables (as is the case for the present example). We
have a space of hypotheses $H$: in this case $H=\{0,1\}$. We have a
space $X$ of possible observed variables: in this case $X =
\mathbb{R}^2\cup \{\text{not found}\}$. We will denote the set of
subsets of a set $A$ (also known as the power set of $A$) by
$\mathbb{S}(A)$. Finally we have the set $I=[0,1]$ of possible levels
of confidence.

Then a frequentist confidence set function is defined to be a
function $$f:X\times I \to \mathbb{S}(H)$$ such that the following
three conditions hold:

\begin{enumerate}

\item $f$ is chosen \textit{before} any data is collected (as is of
  course also required of a (nested family of) critical region(s));

\item \label{nested} for each $x\in X$ and $\eta_1 \leq \eta_2 \in
  I$, $$f(x,\eta_1)\subseteq f(x,\eta_2);$$

\item \label{bigger} for each $h\in H$ and $\eta \in I$, $$P(h\in
  f(x,\eta) | h) \geq \eta.$$

\end{enumerate}

\subsubsection{Asymmetric ``inconclusive'' confidence sets}

\label{inconclusive} 
The simplest way of constructing confidence set functions from the
families of critical regions considered above proceeds as follows --
we will consider more sensible symmetric ones in a later section.

We choose one of figures \ref{freqmaxarea}, \ref{freqct1},
\ref{freqxonly}, \ref{freqdirection}, or \ref{freqpseudoB} (or some
other similar result) whose critical regions we will use.

For each $(x,y)\in \mathbb{R}^2$ and each $\eta\in I$, we set
$f((x,y),\eta)$ to be $\{0,1\}$ if the frequentist confidence shown in
the relevant plot of figures \ref{freqmaxarea}, \ref{freqct1},
\ref{freqxonly}, \ref{freqdirection}, or \ref{freqpseudoB} at $(x,y)$
is less than $\eta$, and otherwise we set $f((x,y),\eta)$ to be
$\{1\}$. In particular we set $f(\text{not found},\eta)=\{0,1\}$.

In other words we
set $$f((x,y),\eta)=\left\{\begin{matrix}\{1\}&((x,y)\in
C_\eta)\\\{0,1\}&((x,y)\notin C_\eta).\end{matrix}\right.$$

We name this the ``inconclusive'' version because the only possible
confidence sets that arise are $\{0,1\}$ and $\{1\}$ - we can never
conclude that $h=0$, i.e. that $H_0$ holds. It therefore leaves one
saying either ``We are $\eta$-confident that the princess is in tower
1'', or ``We are not $\eta$-confident about which tower the princess
is in''.

It is easy to see that these various functions $f$ define valid
frequentist confidence set functions: for condition \ref{nested} of
section \ref{confsetreminder}, if $\eta_1\leq\eta_2$, then
\begin{IEEEeqnarray*}{lrCl}
& 0 & \in & f((x,y),\eta_1)\\
\implies & (x,y) & \notin & C_{\eta_1} \supseteq C_{\eta_2}\\
\implies & (x,y) & \notin & C_{\eta_2}\\
\implies & 0 & \in & f((x,y), \eta_2),
\end{IEEEeqnarray*}
so that (since we always have $1\in f((x,y),\eta)$)
$$f((x,y),\eta_1)\subseteq f((x,y),\eta_2).$$ For condition
\ref{bigger}, if $h=1$ we always have $h\in f((x,y),\eta)$ so that
$P(h\in f((x,y),\eta)|h)=1\geq \eta$, while if $h=0$ then
\begin{IEEEeqnarray*}{rCl}
P(h\in f((x,y),\eta)|h) &=& P((x,y)\notin C_\eta | h)\\
&=& 1-P((x,y)\in C_\eta|h)\\
&\geq& 1-(1-\eta)\\
&=&\eta.
\end{IEEEeqnarray*}

But the different choices of nested sets of critical regions give
different confidence set functions, the choice between which is
arbitrary -- and these are just a few of many such functions that
could be chosen.

\subsubsection{Symmetric ``conclusive'' confidence sets}
\label{conclusive}

But we can also be slightly more inventive in constructing confidence
set functions, in a way that allows us to also conclude when
appropriate that $H_0$ holds, as follows.

Having chosen one of the five (or other) plots and constructed a
confidence set function $f_1$ according to the recipe in
\ref{inconclusive} above, taking possible values $\{0,1\}$ and
$\{1\}$, we then reverse the roles of the two towers and construct a
second confidence set function $f_2$ taking the possible values
$\{0,1\}$ and $\{0\}$. Finally we set $$f((x,y),\eta) =
f_1((x,y),\eta)\cap f_2((x,y),\eta).$$

In this way we get a frequentist confidence set function which can
take any of the four values $\emptyset$, $\{0\}$, $\{1\}$, or
$\{0,1\}$. To see that $f$ is then a valid frequentist confidence set
function, note first that since $f_1$ and $f_2$ are both increasing
with $\eta$, so is their intersection. For condition \ref{bigger},
note that for $h=1$,
\begin{IEEEeqnarray*}{rClr}
P(h\in f((x,y),\eta)|h) &=& P(h\in f_2((x,y),\eta)|h) & \text{ (since
  we always have }1\in f_1((x,y),\eta))\\ &\geq &\eta,$$
\end{IEEEeqnarray*}
with an exactly similar argument for $h=0$.

\begin{figure}[h]
\begin{center}

\includegraphics[scale=0.5]{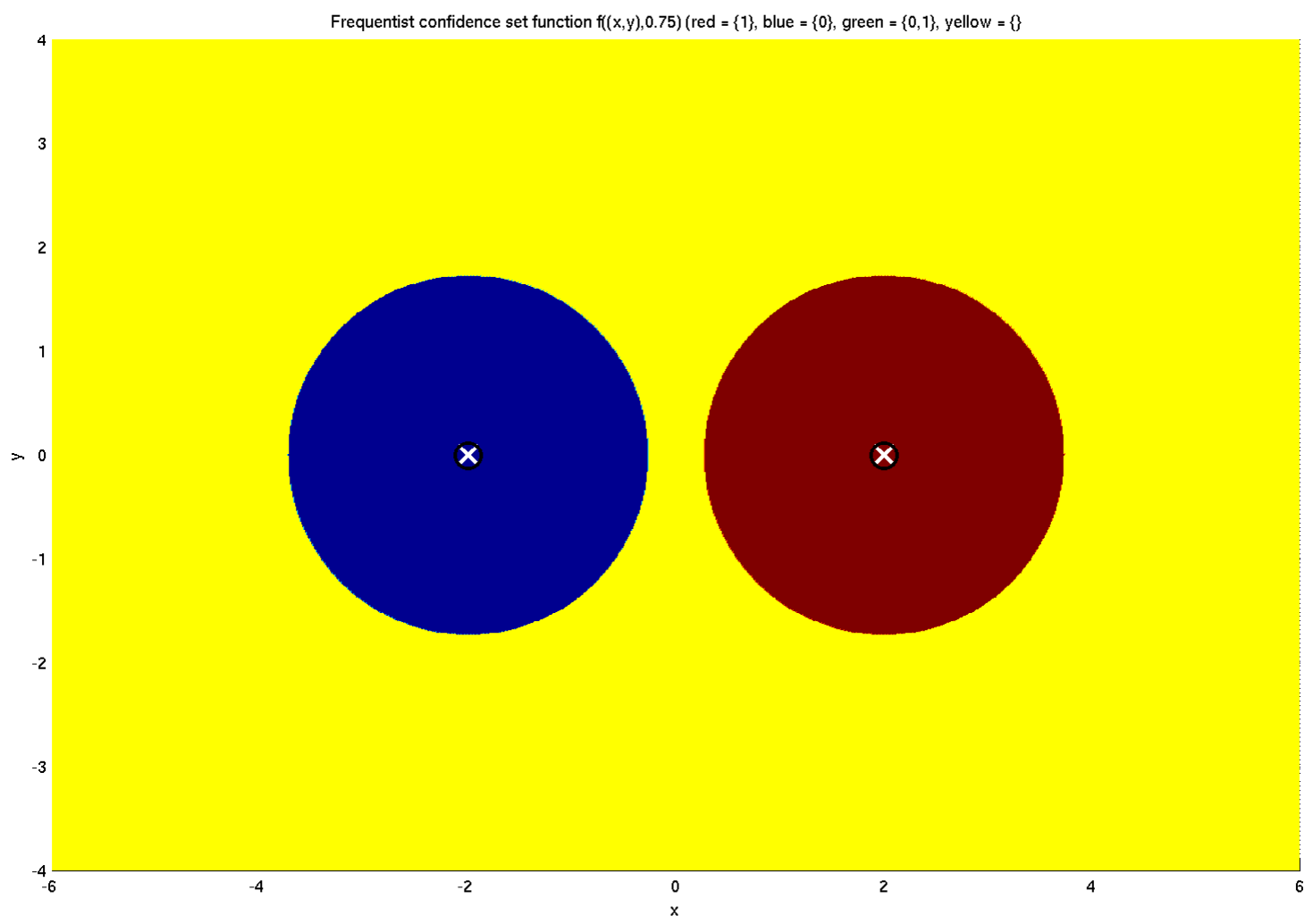}

\caption{Frequentist confidence set function $f((x,y),0.75)$ derived
  from the critical regions that gave figure
  \ref{freqmaxarea}. Red-brown indicates where the function takes the
  value $\{1\}$, blue where it is $\{0\}$, green where it is
  $\{0,1\}$, and yellow where it is $\emptyset$. The positions of the
  two towers are also shown.
\label{freqmaxareaf75}
}

\end{center}
\end{figure}

\begin{figure}[p]
\begin{center}

\includegraphics[scale=0.5]{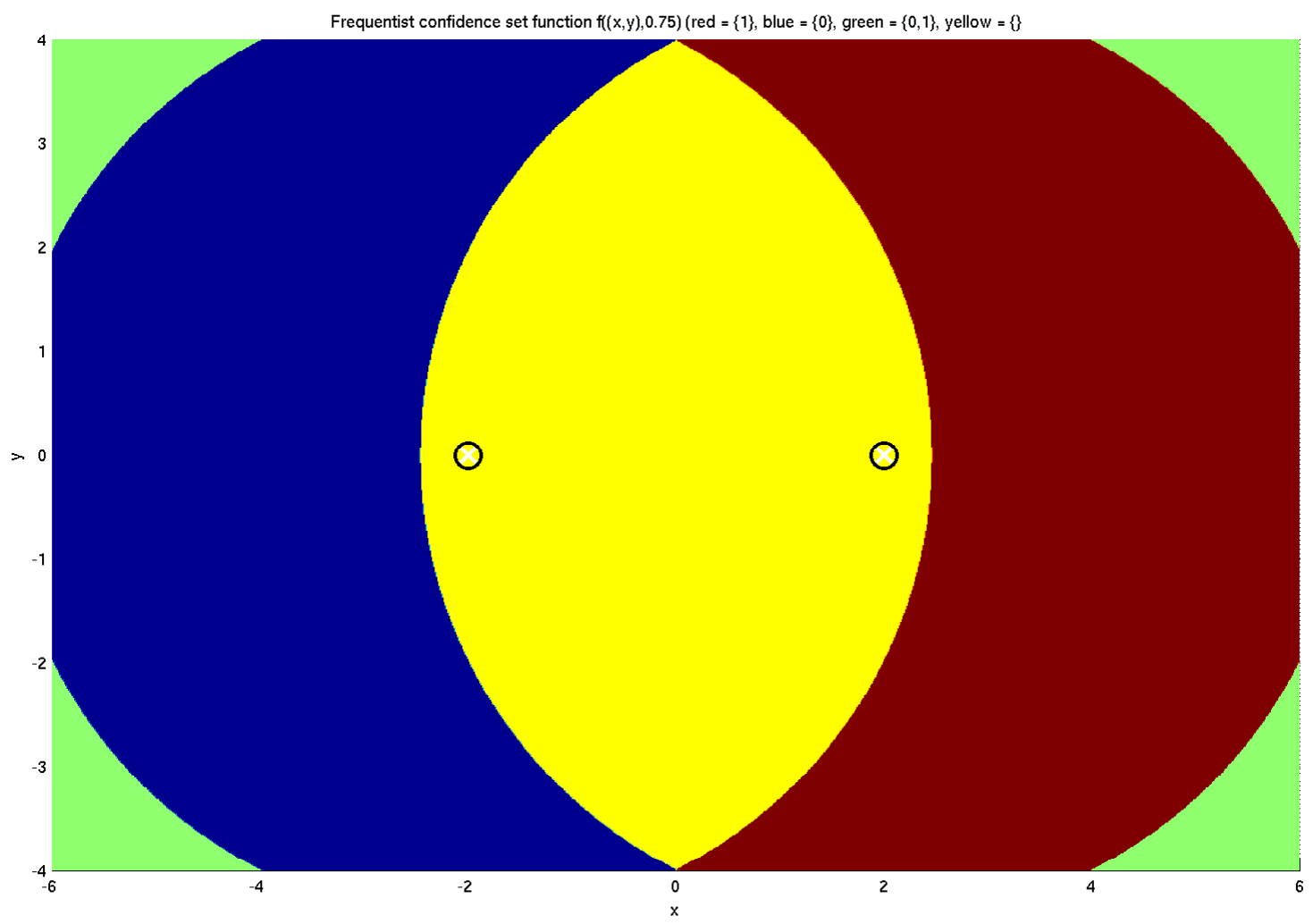}

\caption{Frequentist confidence set function $f((x,y),0.75)$ derived
  from the critical regions that gave figure \ref{freqct1}. Red-brown
  indicates where the function takes the value $\{1\}$, blue where it
  is $\{0\}$, green where it is $\{0,1\}$, and yellow where it is
  $\emptyset$. The positions of the two towers are also shown.
\label{freqct1f75}
}

\end{center}
\end{figure}

\begin{figure}[p]
\begin{center}

\includegraphics[scale=0.5]{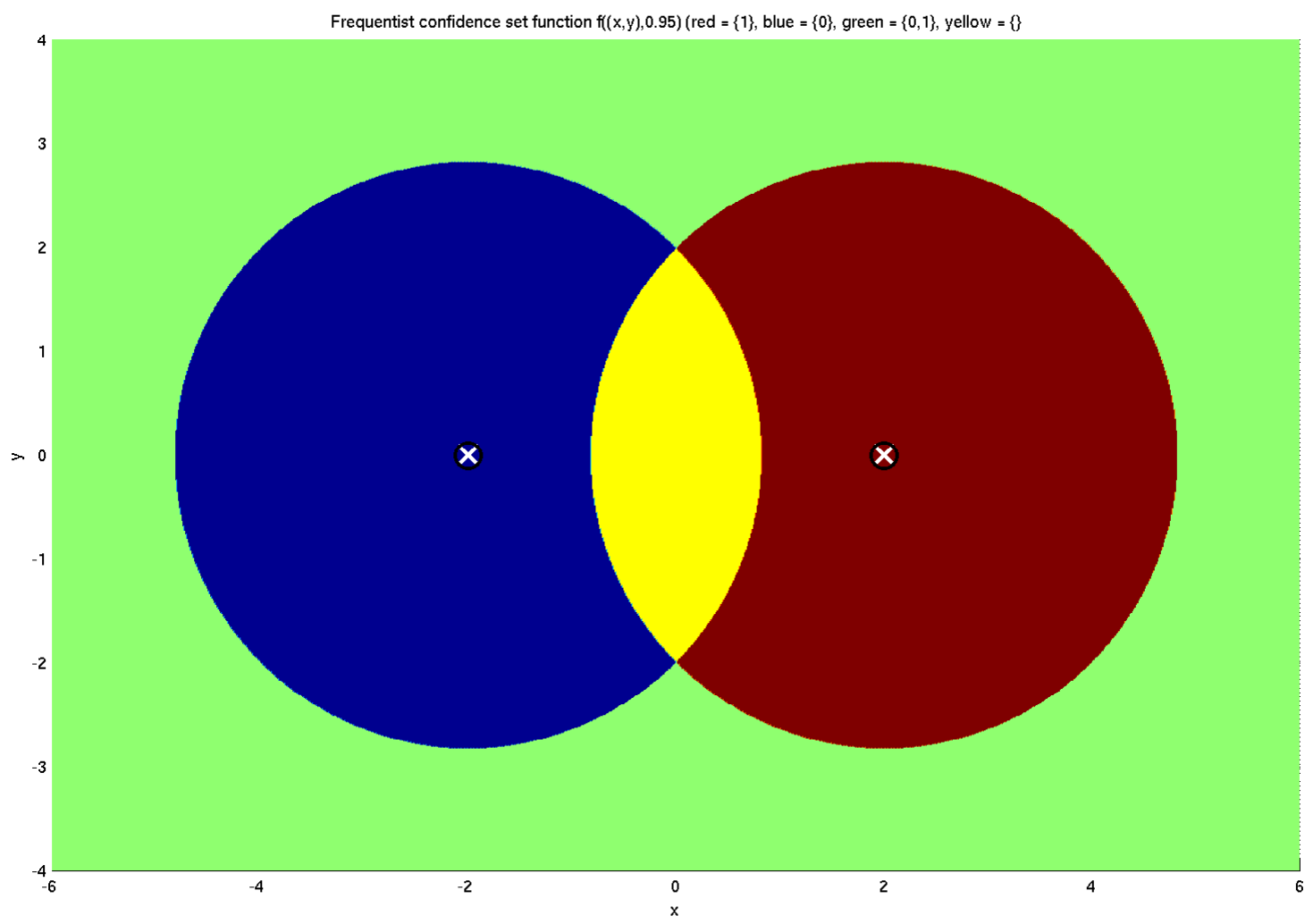}

\caption{Frequentist confidence set function $f((x,y),0.95)$ derived
  from the critical regions that gave figure \ref{freqct1}. Red-brown
  indicates where the function takes the value $\{1\}$, blue where it
  is $\{0\}$, green where it is $\{0,1\}$, and yellow where it is
  $\emptyset$. The positions of the two towers are also shown.
\label{freqct1f95}
}

\end{center}
\end{figure}

\begin{figure}[p]
\begin{center}

\includegraphics[scale=0.5]{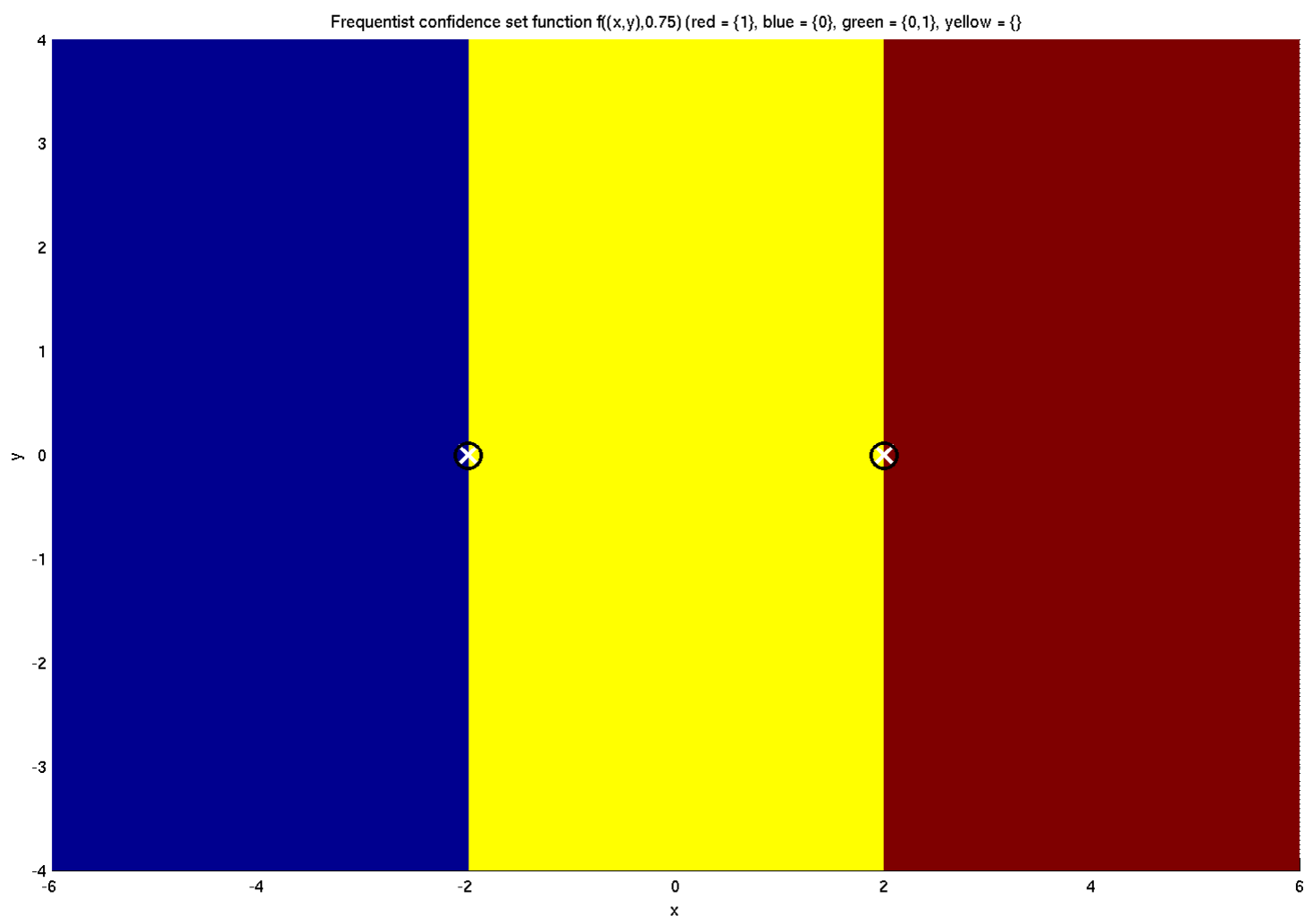}

\caption{Frequentist confidence set function $f((x,y),0.75)$ derived
  from the critical regions that gave figure
  \ref{freqxonly}. Red-brown indicates where the function takes the
  value $\{1\}$, blue where it is $\{0\}$, green where it is
  $\{0,1\}$, and yellow where it is $\emptyset$. The positions of the
  two towers are also shown.
\label{freqxonlyf75}
}

\end{center}
\end{figure}

\begin{figure}[p]
\begin{center}

\includegraphics[scale=0.5]{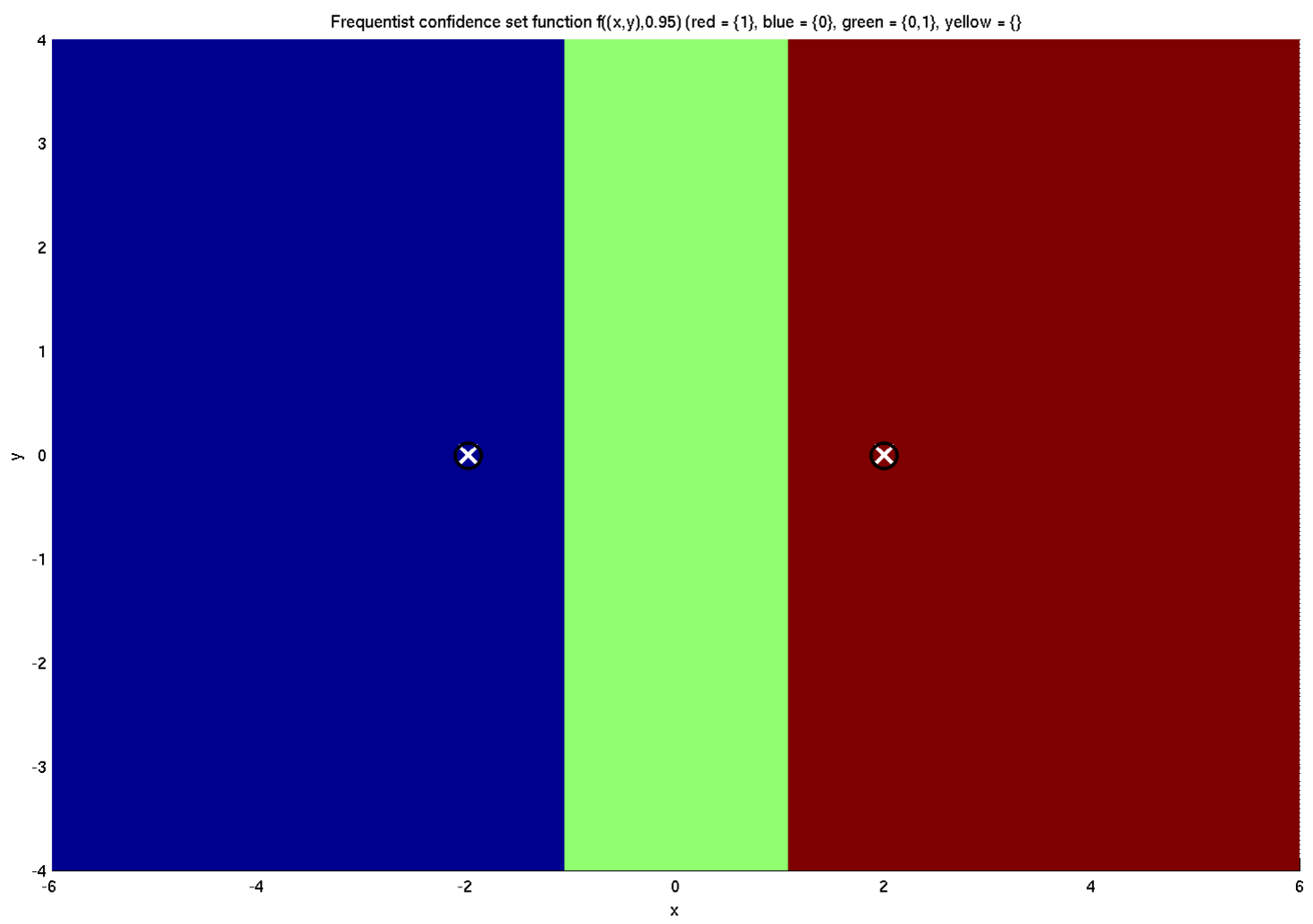}

\caption{Frequentist confidence set function $f((x,y),0.95)$ derived
  from the critical regions that gave figure
  \ref{freqxonly}. Red-brown indicates where the function takes the
  value $\{1\}$, blue where it is $\{0\}$, green where it is
  $\{0,1\}$, and yellow where it is $\emptyset$. The positions of the
  two towers are also shown.
\label{freqxonlyf95}
}

\end{center}
\end{figure}

\begin{figure}[p]
\begin{center}

\includegraphics[scale=0.5]{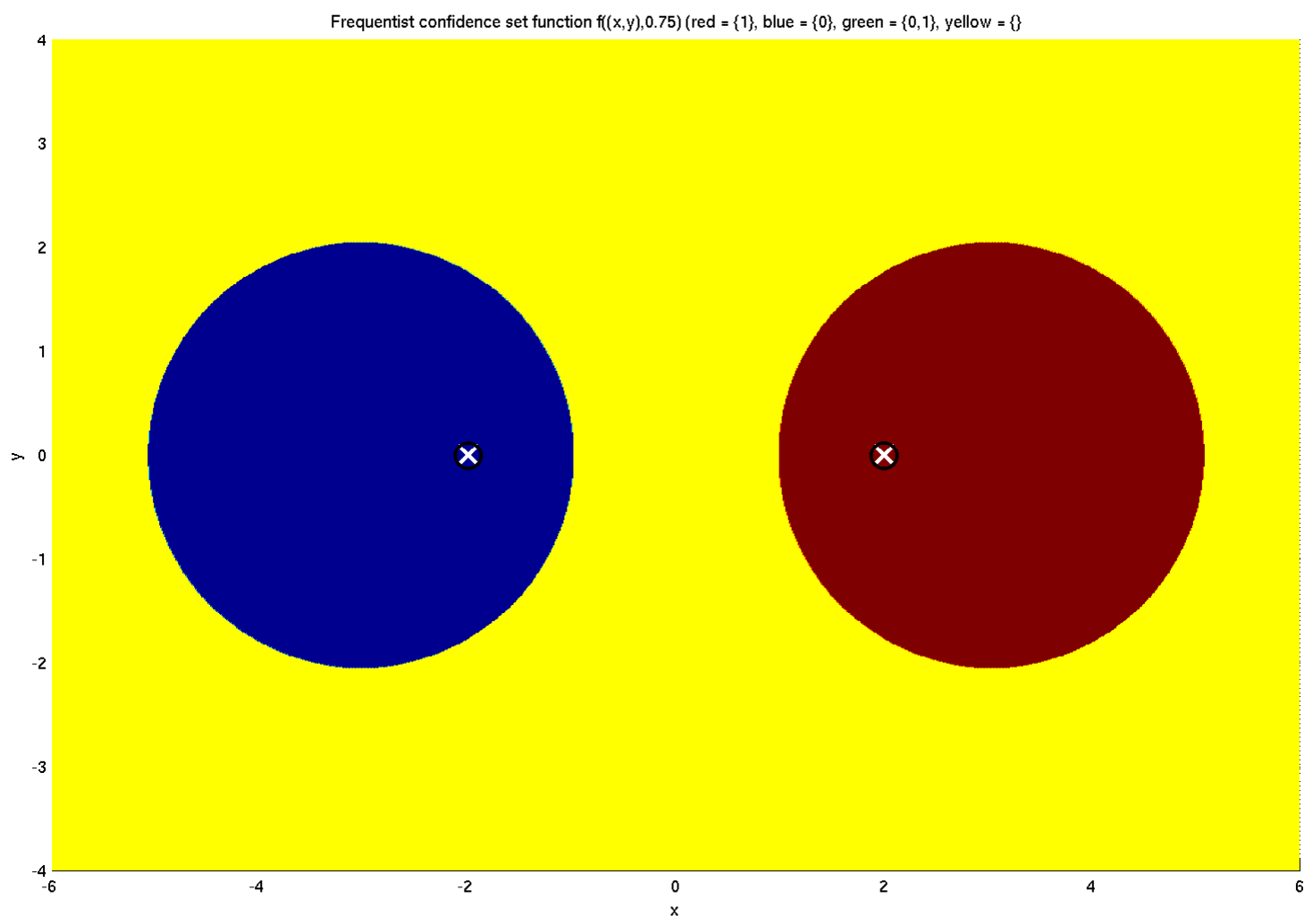}

\caption{Frequentist confidence set function $f((x,y),0.75)$ derived
  from the critical regions that gave figure
  \ref{freqpseudoB}. Red-brown indicates where the function takes the
  value $\{1\}$, blue where it is $\{0\}$, green where it is
  $\{0,1\}$, and yellow where it is $\emptyset$. The positions of the
  two towers are also shown.
\label{freqpseudoBf75}
}

\end{center}
\end{figure}

\begin{figure}[p]
\begin{center}

\includegraphics[scale=0.5]{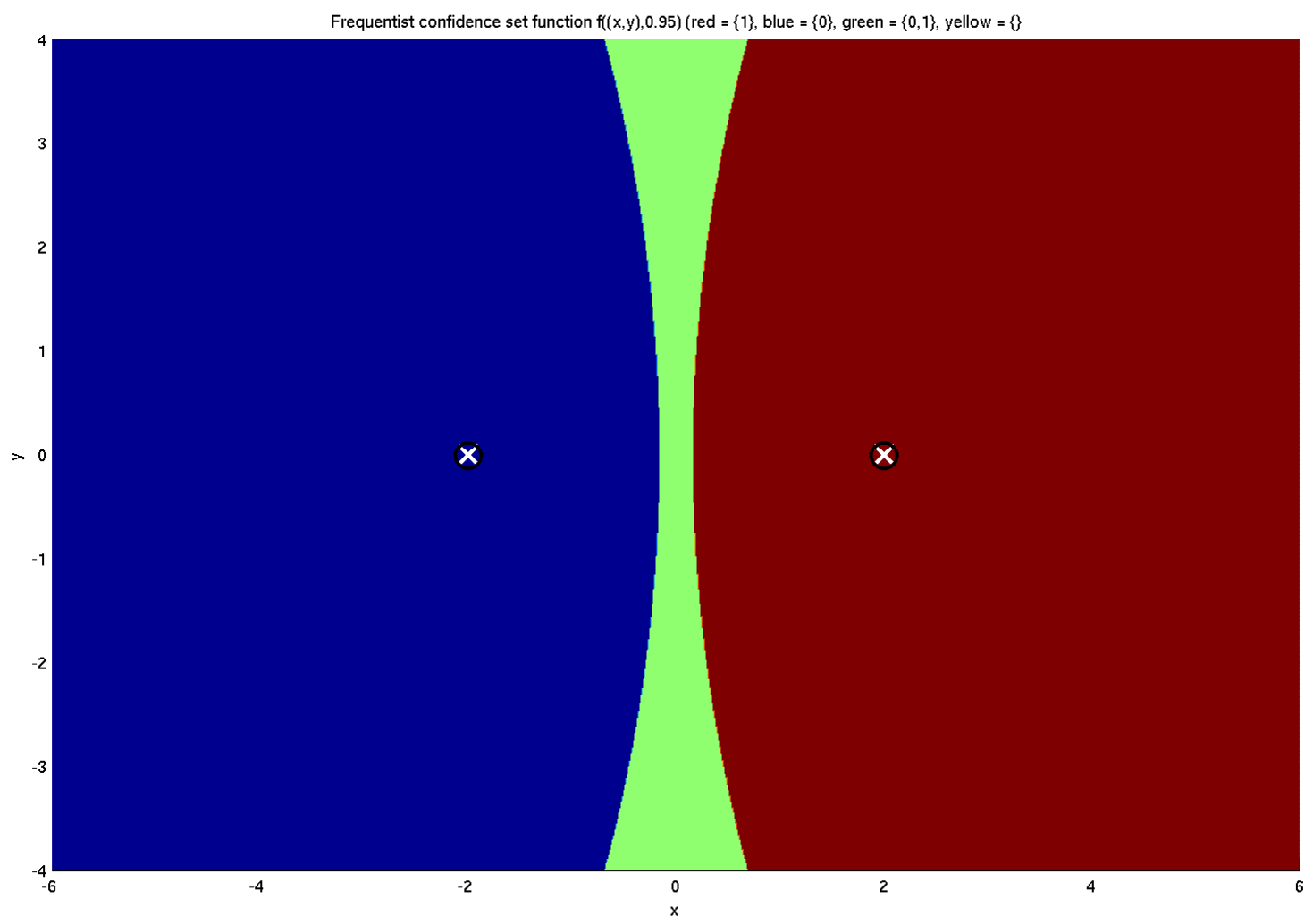}

\caption{Frequentist confidence set function $f((x,y),0.95)$ derived
  from the critical regions that gave figure
  \ref{freqpseudoB}. Red-brown indicates where the function takes the
  value $\{1\}$, blue where it is $\{0\}$, green where it is
  $\{0,1\}$, and yellow where it is $\emptyset$. The positions of the
  two towers are also shown.
\label{freqpseudoBf95}
}

\end{center}
\end{figure}

We illustrate the various such functions that result in figures
\ref{freqmaxareaf75} to \ref{freqpseudoBf95} for various values of
$\eta$.

We see from figure \ref{freqmaxareaf75} that at the $75\%$ level,
depending on $(x,y)$, we may be $75\%$ frequentist confident that
$H_0$ holds (i.e. that $h\in H_0$), that $H_1$ holds, or that
$h\in\emptyset$ (which is actually impossible). (In principle we might
also have had the option of concluding that one or the other holds,
which we knew already - see e.g. figure \ref{freqct1f95}.)

An alternative interpretation illustrating the perversity of the
frequentist approach comes from noting that the yellow region of
figure \ref{freqmaxareaf75} arises from concluding for $(x,y)$
in this region first that $H_1$ holds with 75\% frequentist
confidence, then after switching the roles of the two towers, that
$H_0$ holds with 75\% frequentist confidence - thus violating
the first bullet point of section \ref{admissiblesolutions}, albeit
with 75 in place of 95.

Similar observations may be made from the other figures,
none of which corresponds to any of the possible Bayesian
solutions. 

For example from figure \ref{freqct1f95} we see that using the
``conclusive'' frequentist confidence set function corresponding to
figure \ref{freqct1}, if $(x,y)=(0,0)$ (i.e. if the bullet lands
mid-way between the two towers) we are now 95\% frequentist confident
that the princess is in neither tower - so where did the bullet come
\mbox{from ?} Even if the bullet lands at the base of one or other tower, we
are $75\%$ frequentist confident that the princess is in neither
tower (figure \ref{freqct1f75}).

On the other hand using figure \ref{freqct1f95} again, if
$(x,y)=(1000, 1000)$ (i.e. if the bullet lands 1414 km North-East of
the origin), we do not know with 95\% frequentist confidence which
tower the princess is in, but have not excluded either.

\begin{figure}[p]
\begin{center}

\includegraphics[scale=0.5]{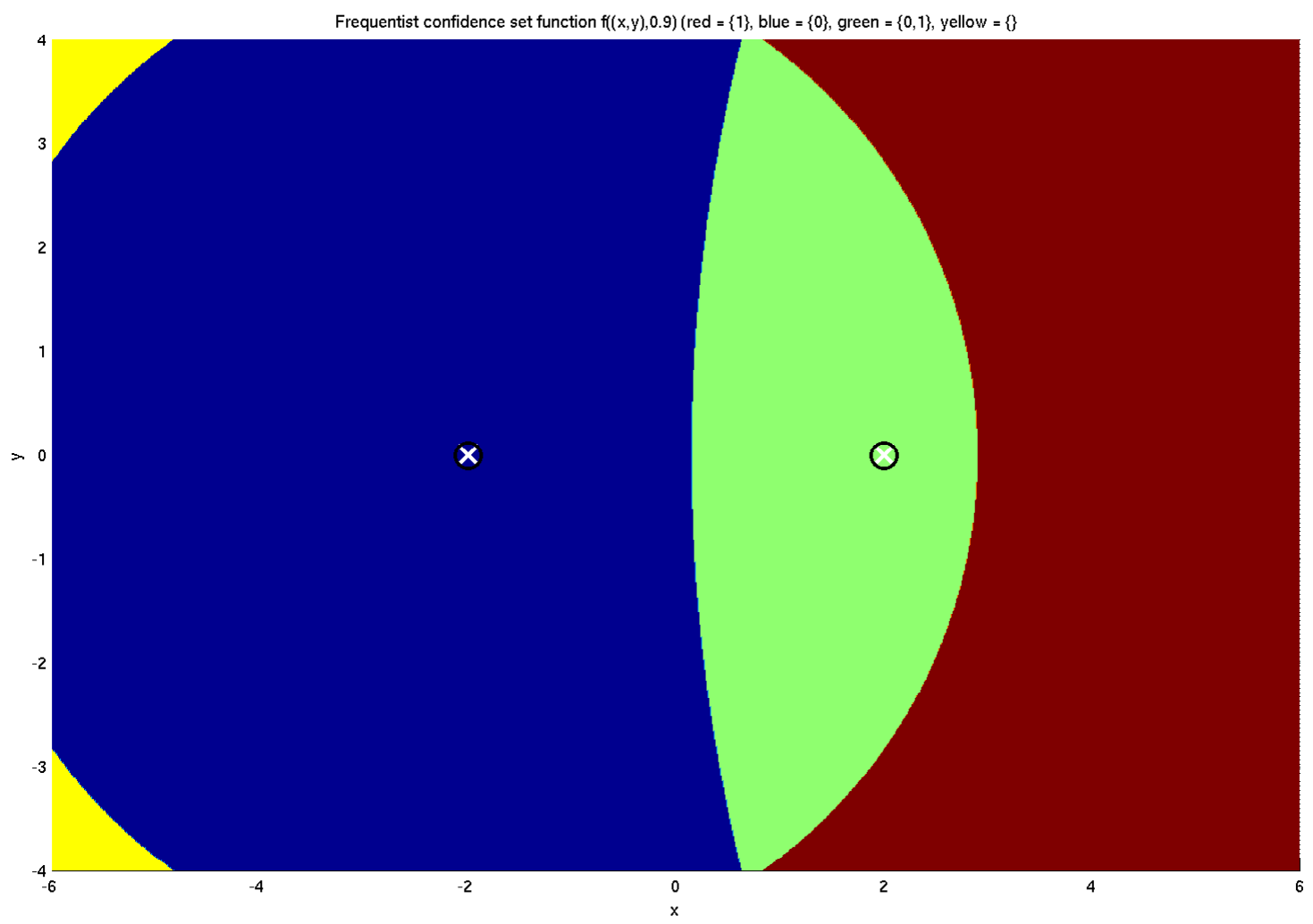}

\caption{Frequentist confidence set function $f((x,y),0.9)$ derived
  from the critical regions that gave figures \ref{freqmaxarea} (as
  shown) and \ref{freqpseudoB} (with the roles of the two towers
  reversed). Red-brown indicates where the function takes the value
  $\{1\}$, blue where it is $\{0\}$, green where it is $\{0,1\}$, and
  yellow where it is $\emptyset$. The positions of the two towers are
  also shown.
\label{freqmixedf}
}

\end{center}
\end{figure}

\begin{figure}[htp]
\begin{center}

\includegraphics[scale=0.5]{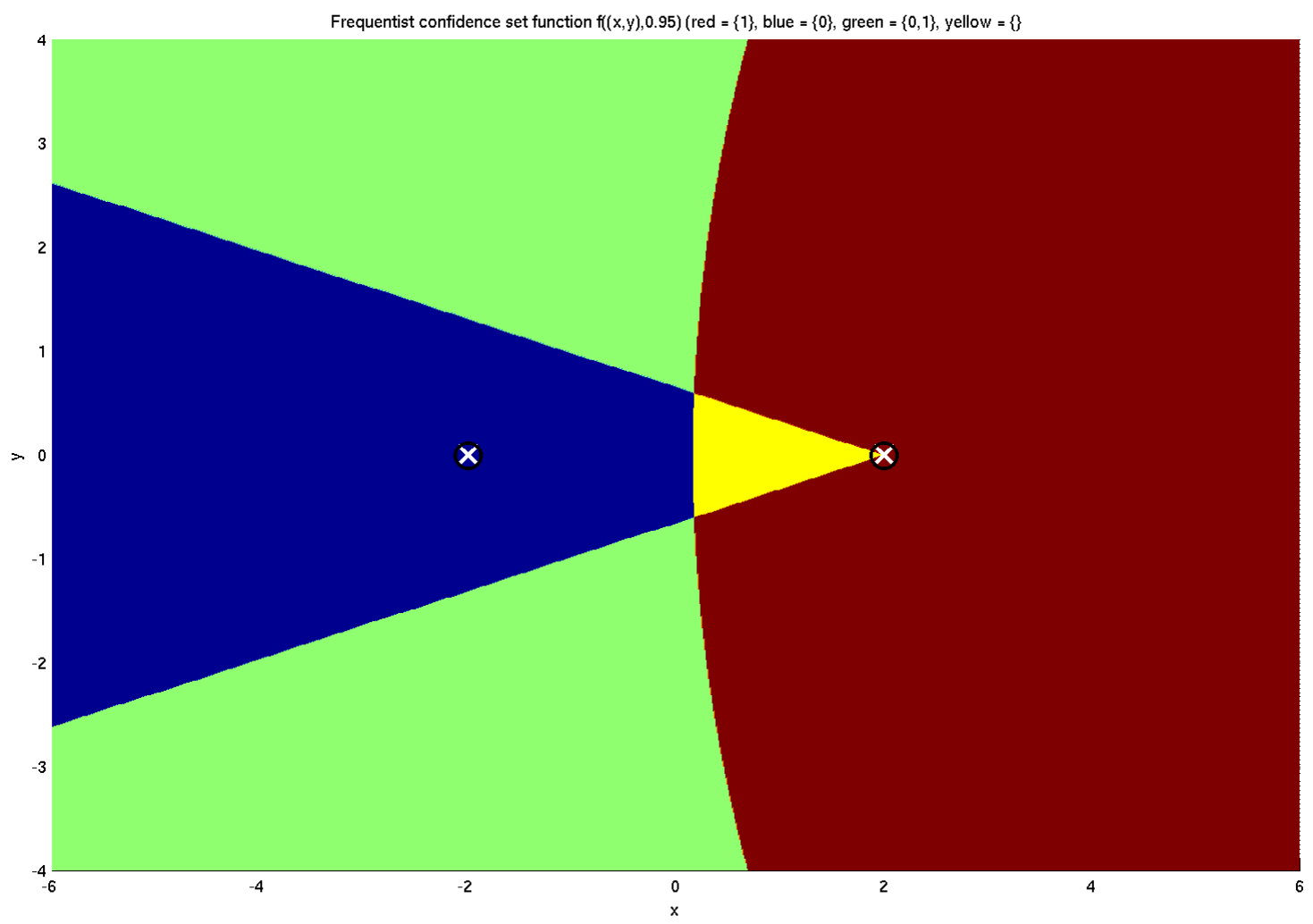}

\caption{Frequentist confidence set function $f((x,y),0.95)$ derived
  from the set of critical regions that gave figure \ref{freqpseudoB}
  (as shown) and one based on direction from tower 0 (with the roles
  of the two towers reversed). Red-brown indicates where the function
  takes the value $\{1\}$, blue where it is $\{0\}$, green where it is
  $\{0,1\}$, and yellow where it is $\emptyset$. The positions of the
  two towers are also shown.
\label{freqmixedf2}
}

\end{center}
\end{figure}

We can even mix the two halves of such frequentist confidence set
functions, giving us functions such as those illustrated in figures
\ref{freqmixedf} and \ref{freqmixedf2}. In the latter, if the bullet
lands at the origin, we are $95\%$ frequentist confident that the
princess is on tower $0$, but between there and the base of tower $1$
is a roughly triangular region, if the bullet lands in which we are
$95\%$ frequentist confident that the princess is on neither tower and
that $h\in\emptyset$.

So in the frequentist confidence-set paradigm we are likewise faced
with arbitrary choices which disagree in their conclusions, some of
which may make us fairly confident that impossible things have
happened -- which is indeed \textit{reductio ad absurdum}.

\subsubsection{Confidence sets based on random numbers}

As in section \ref{freqrandom} above, a frequentist confidence set
function based on random numbers is possible; for the symmetric one we
draw $u_0,u_1$ each independently uniformly distributed in $[0,1]$,
and report a frequentist confidence set of confidence $c$ that
excludes $0$ iff $u_0\geq c$ and excludes $1$ iff $u_1\geq c$.

\subsection{Conclusion on the example}

We make a number of points. 

First, some frequentists might say that they wouldn't use hypothesis
testing methods for this particular problem. To this we respond that
the point of examples is to show the defects in a method: if an
inference method (in this case frequentist hypothesis testing)
produces nonsensical answers for any single inference problem to which
it can be applied, then it is flawed -- this is the point of the
mathematical concept of a ``counter-example''.

Second, we note that it is not just Bayesians who have to invoke
choices not mentioned in the original problem. Bayesians choose
priors, but frequentists choose (nested sets of) critical regions --
and while it is entirely appropriate that what we conclude should
depend on what we knew beforehand, it is not appropriate that it
should depend on arbitrary choices of critical regions -- and above
we have barely started to consider the possibilities (others include
e.g. complements of diamond shapes centred on tower 0, diamonds
centred on tower 1, direction from tower 0, sets of rectangles of
particular aspect ratio, ellipses, etc., with the range of
possibilities increasing exponentially with the number of data
dimensions). (In the case of confidence sets/intervals, frequentists
instead choose equally arbitrary confidence set functions, even though
this may not be obvious to many users -- see section \ref{fconfsets}
above.)

For those not used to the concept of critical regions, we could also
describe the various nested sets of critical regions considered as
frequentist ``tests'', according to the correspondence in table
\ref{tests}. For example, the set of critical regions determined by
minimum excluded area might be termed ``the distance-from-tower-0
test''; or in \cite{FDApseudoBayes} the FDA describes a test of the
pseudo-Bayesian type as a ``Bayesian test'', even though of course it
is a frequentist test and not Bayesian at all (as they have later
agreed \cite{FDA2018}). The frequentist then has an arbitrary choice
of test to make.  

\begin{table}[hp]
\begin{tabular}{ll}
\textbf{Critical region description} & \textbf{Name of corresponding
  frequentist test}\\
Minimal excluded area & Distance from tower 0 test\\
Discs centre tower 1 & Distance from tower 1 test\\
$x$-coordinate only & $x$-coordinate test\\
Direction from origin & Direction test\\
Nearly  pseudo-Bayesian & ``Bayesian'' test\\
Full pseudo-Bayesian & Uniformly optimal test
\end{tabular}
\caption{
\label{tests}
How the various nested sets of critical regions could instead be
described as frequentist tests.
}
\end{table}

In particular we offer table \ref{resultstable} of ``$\eta$-sureness''
comparing Bayesian results on the example problem with the various
frequentist hypothesis testing solutions considered. The reader is
invited to carefully consider the numbers in this table and ask
themselves which of the various approaches gives the answers that are
correct. We also point out that in 25 of 30 of the cases where the
bullet is found the frequentist confidence that $H_1$ holds is higher
even than the Bayesian posterior probability under the prior that puts
$\frac{3}{4}$ of the probability on $H_1$, despite the frequentist
solution being designed to avoid type I errors where $H_1$ is inferred
but $H_0$ is true.

\begin{table}[hp]
\begin{tabular}{lllllll}
\hspace{1.5 cm}\textbf{Position} & $(-6,0)$ & $(-2,0)$ & $(0,0)$ & $(2,0)$ & $(6,0)$ & not found\\
\textbf{Solution}\\
Bayes (even prior) & 0.118 & 0.014 & 0.500 & 0.986 & 0.882 & 0.500\\
Bayes ($\frac{3}{4}$ on $H_0$) & 0.043 & 0.005 & 0.250 & 0.959 & 0.714 & 0.250\\
Bayes ($\frac{3}{4}$ on $H_1$) & 0.286 & 0.041 & 0.750 & 0.995 & 0.957 & 0.750\\
Min excluded area & 0.776 & 0.500 & 0.776 & 0.879\S & 0.983 & 0.000\\
Discs centre tower 1 & 0.577 & 0.820 & 0.981 & 1.000 & 0.820\S & 0.000\\
$x$-coordinate only & 0.539 & 0.749 & 0.926 & 0.961\S & 0.980 & 0.000\\
Direction from origin (closed sectors) & 0.500 & 0.500 & 1.000\ddag & 1.000 & 1.000 & 0.000\\
Nearly  pseudo-Bayesian\dag & 0.783 & 0.515 & 0.926 & 0.9998 & 0.987 & 0.000\\
Uniform random & $u$ & $u$ & $u$ & $u$ & $u$ & $u$ \\
Full pseudo-Bayesian & 0.283\S & 0.015\S & 0.926 & 0.9998 & 0.987 &
$\frac{u}{2} + \frac{\tan^{-1}2}{2 \pi} + \frac{1}{4}$\\\\
\end{tabular}

\begin{tabular}{lp{13cm}}
$u$ & A random number uniformly distributed on $[0,1]$\\\\
\dag & This is uniformly optimal among \textit{deterministic} frequentist
hypothesis-testing solutions for frequentist confidence levels $\eta \geq 0.5$\\\\
\ddag & But finding the bullet just 1 mm further West would make
confidence only 0.500, as would using open sectors instead of closed ones\\\\
\S & These are the only entries in the frequentist part of the table
when the bullet is found which are less than the Bayesian value using the prior which puts 0.75
of the probability on $H_1$.
\end{tabular}

\caption{
\label{resultstable}
Measures of $\eta$-sureness (whether posterior probability or
frequentist confidence) that the princess is on tower 1 at $(+2,0)$
for various possible places where the bullet may be found or its
absence.}

\end{table}

Third, many users don't appreciate the defects of frequentist methods
because they only use them for 1-dimensional problems, and give little
thought to the properties of the methods they use -- but even for
1-dimensional problems frequentist methods violate the basic
properties required of any inference method, as discussed above in
section \ref{whichsatisfy}. The current 2-dimensional example does,
however, make the difficulties and the dilemma somewhat clearer.

Fourth, the consideration above of ``conclusive'' frequentist
confidence sets for this problem indicates that they can often put
quite high confidence on answers that are actually impossible -- in
this case by saying that e.g. we are $95\%$ frequentist confident that
the princess is on neither tower and $h\in\emptyset$ (section
\ref{conclusive} and figure \ref{freqct1f95}), when we know from the
start that she is on one or the other (and nothing can ever be in the
empty set).

\clearpage

\section{An example where information content can be easily
  calculated}

\label{shells}

\subsection{Introduction}

To assist those who dislike the problem of section \ref{example1}
because of there being only a single data point, we now turn to an
example where in principle we can collect arbitrarily many data
points, though we will stick to $N=20$. This problem has a
2-dimensional discrete sufficient statistic, which we can use to
illustrate how the various inference methods work in a manner similar
to the previous example. In addition, the fact that both the
hypothesis space and the data space are discrete allows us to easily
calculate the True and Apparent Shannon information content (TSI and
ASI respectively) of the various different solutions (more about this
in section \ref{shellsinfointro} below); this allows us to measure the
overall deficits of the various frequentist solutions compared with
both Bayes and the amount of ASI needed to give a definite answer to
the problem (i.e. to answer the User's question of section
\ref{introbackground}).

\subsection{Definition of the ``Shells problem''}

We set this (entirely fictitious) problem in the context of a war
between parties A and B, taking the perspective of somebody supporting
A. We believe that shells used by B are being supplied either by
country 0 or country 1, and we want to determine who is supplying
shells to B in order to put sanctions on them.

We know that two different types of mortar are being used by B. For a
data source we observe that on landing shells exhibit one of three
different fragmentation patterns, and that the probability of each of
these three patterns occurring depends both on the type of mortar used
to fire the shell and on the country of origin of the shell (all
shells otherwise appearing identical).

Prior to collecting live data, we obtain captured mortars of the two
types, and we buy from an arms dealer supplies of shells from each
country. We use these to conduct experiments in a location remote from
the battlefield, and determine the probabilities $P(f|m,c)$ of the
different fragmentation types occurring when each shell type is fired
from each mortar type, where $f\in\{0,1,2\}$ denotes the fragmentation
type, $m\in\{0,1\}$ denotes the type of mortar used, and $c\in\{0,1\}$
denotes the country of origin of the shell.

Going out to collect the real data set, we observe the first $N=20$
shell landings and their fragmentation type, all of which came from
the same mortar (but we don't know which type of mortar was used). We
note that the order in which the shells landed does not affect the
inference, so we simply count the number $n_f$ of each framentation
type that occurs, giving a vector $n=(n_0,n_1,n_2)$ such
that $$\sum_{f=0}^2{n_f}=N.$$ 

Thinking about how to do the inference, we note that the expected
number $k_{f,m,c}$ of each fragmentation type given any particular
combination of mortar type and shell source is given by $k_{f,m,c}=N
P(f|m,c)$, and $k_{m,c}$ will denote the vector $(k_{0,m,c},k_{1,m,c},
k_{2,m,c})$. However, we can't assume that we will always get the
exact expected number of each fragmentation pattern.

\subsection{Formulation of the inference problem}

We therefore formulate our inference problem by the following
equations, which are agreed by both frequentists and Bayesians:

$$\theta = c,$$ $$\Theta=\{0,1\},$$
$$\phi = m,$$ $$\Phi=\{0,1\},$$
$$H = \Theta \times \Phi,$$
$$x = n = (n_0,n_1,n_2),$$ $$X = \left\{n \in \{0,1,...,N\}^3:\sum_{f=0}^2{n_f}=N\right\},$$
$$P(x|\theta,\phi) = P(n|m,c) =
\frac{N!}{\prod_{f=0}^2{(n_f!)}}\prod_{f=0}^2{P(f|m,c)^{n_f}},$$
$$H_\theta = H_c = \{(c,0), (c,1)\}.$$

The actual values of $P(f|m,c)$ determined by our off-battlefield
experiments that we use for this problem in this section are given in
table \ref{pfgivenmc}.

\begin{table}[htp]
\begin{center}

\begin{tabular}{rrrrr}
& \multicolumn{2}{c}{$c=0$} & \multicolumn{2}{c}{$c=1$}\\
& $m=0$ & $m=1$ & $m=0$ & $m=1$ \\
$f=0$ & 0.191 & 0.652 & 0.266 & 0.058 \\
$f=1$ & 0.105 & 0.192 & 0.213 & 0.909 \\
$f=2$ & 0.704 & 0.156 & 0.521 & 0.033 \\
\end{tabular}

\caption{The values of $P(f|m,c)$ used in this section, assumed to
  have been determined by off-battlefield experiments.
\label{pfgivenmc}
}
\end{center}
\end{table}

\begin{figure}[htp]
\begin{center}

\includegraphics[scale=0.5]{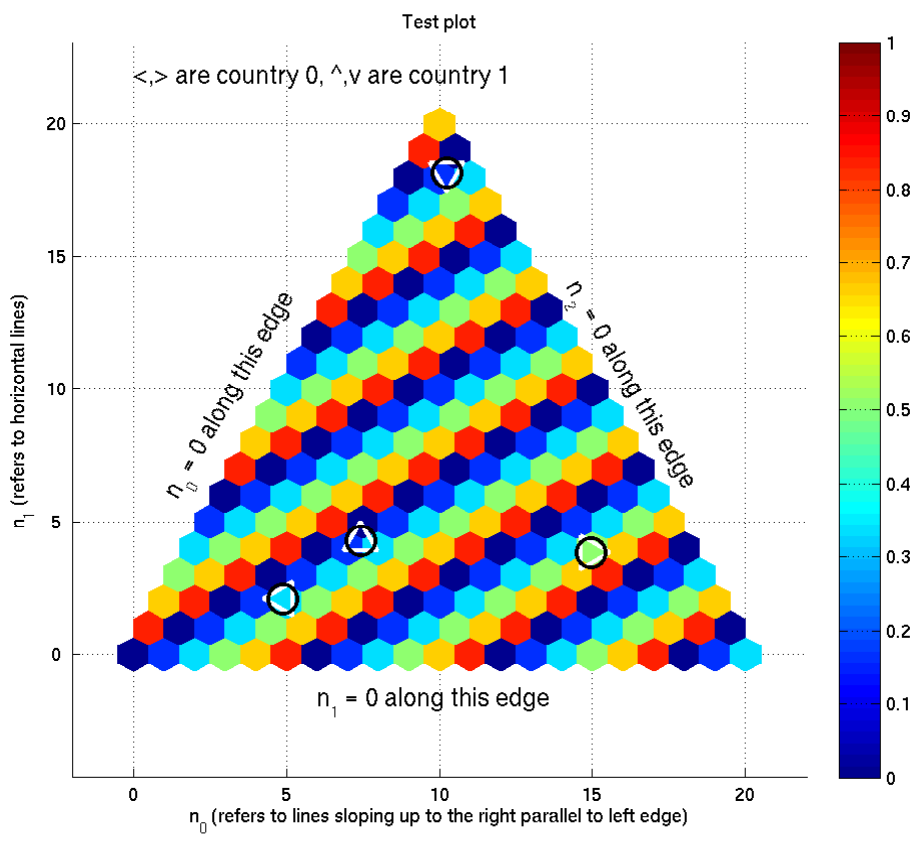}

\caption{The data space of the shells problem of section
  \ref{shells}. The actual colours on this plot mean nothing; they are
  chosen only in order that adjacent hexagons (cells) may be
  distinguished. Each possible data value $n = (n_0,n_1,n_2)$
  corresponds to one cell. The expectation $k_{m,c}$ of $n$ given
  $(m,c)$ is also plotted: triangle pointing left is $k_{0,0}$,
  triangle pointing right is $k_{1,0}$, triangle pointing up is
  $k_{0,1}$, and triangle pointing down is $k_{1,1}$.  The cell at the
  bottom left, black, represents the data value $(0,0,20)$; moving
  upwards along the left edge we find cells of colours red, orange,
  green, light-blue representing successively
  $(0,1,19),(0,2,18),(0,3,17),(0,4,16)$. Moving from the bottom left
  corner horizontally to the right we find cells that are dark-blue,
  light-blue, green, orange representing
  $(1,0,19),(2,0,18),(3,0,17),(4,0,16)$. If we start at the orange
  cell in the middle of the bottom edge representing $(10,0,10)$ and
  move upwards and to the right parallel to the left edge we find
  cells that are green, light-blue, dark-blue, ..., and finally black,
  representing $(10,1,9),(10,2,8),(10,3,7),...$ and finally
  $(10,10,0)$. Starting again at the orange cell in the middle of the
  bottom edge representing $(10,0,10)$ and moving up and to the left
  parallel to the right edge we find cells that are light-blue, black,
  ..., and finally light-blue representing $(9,1,10),(8,2,10)...$ and
  finally $(0,10,10)$.
\label{testplot}
}

\end{center}
\end{figure}

Before plotting the likelihood function $P(n|m,c)$, we first show a
plot of the data (sufficient statistic) space in order to be sure that
later plots are correctly understood. The reader is invited to read
through the caption of figure \ref{testplot} to ensure that he
interprets the later plots correctly.

\begin{figure}[hp]

\begin{tabular}{cc}

\includegraphics[scale=0.5]{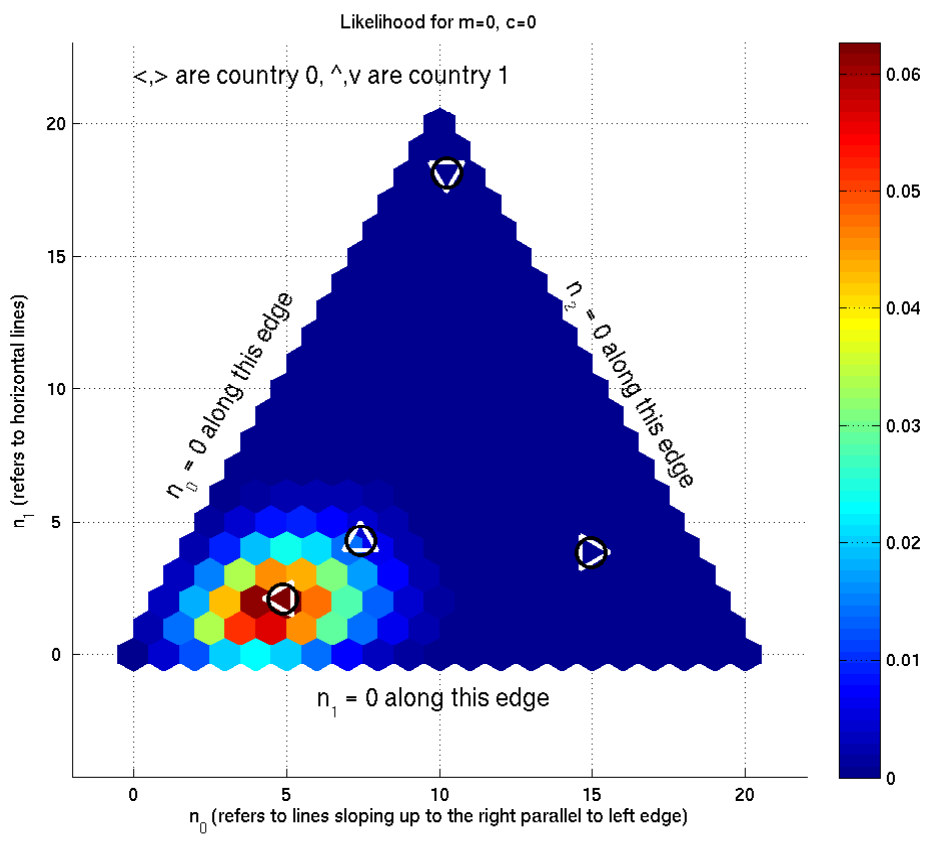} &
\includegraphics[scale=0.5]{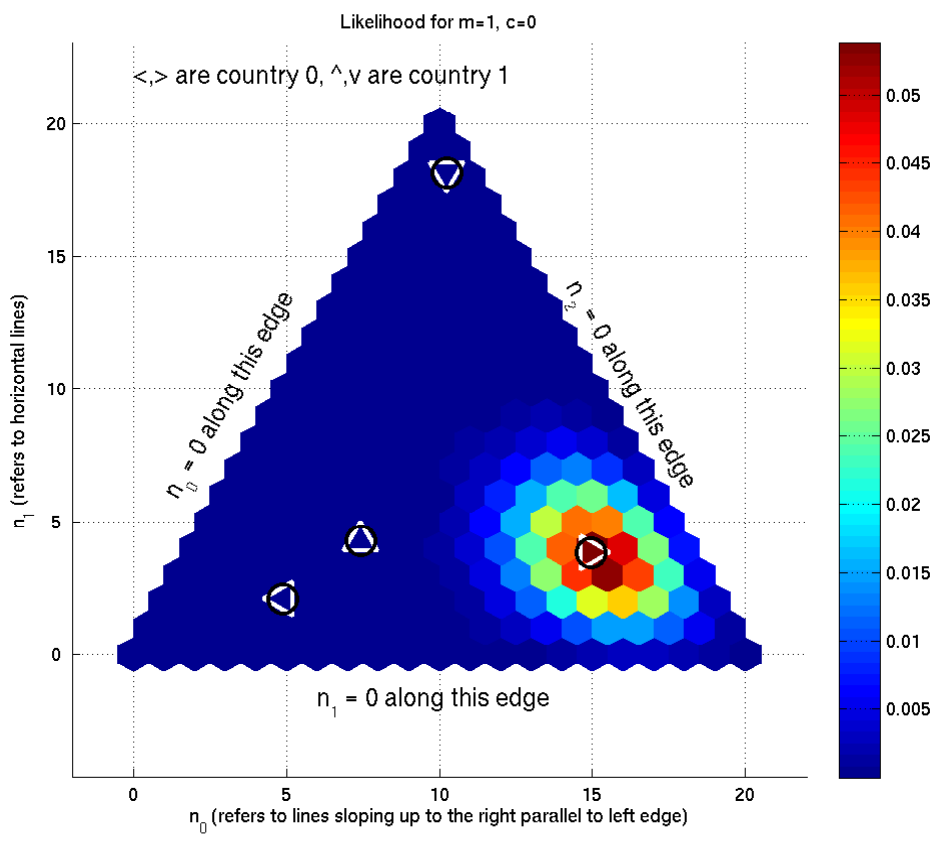} \\
\includegraphics[scale=0.5]{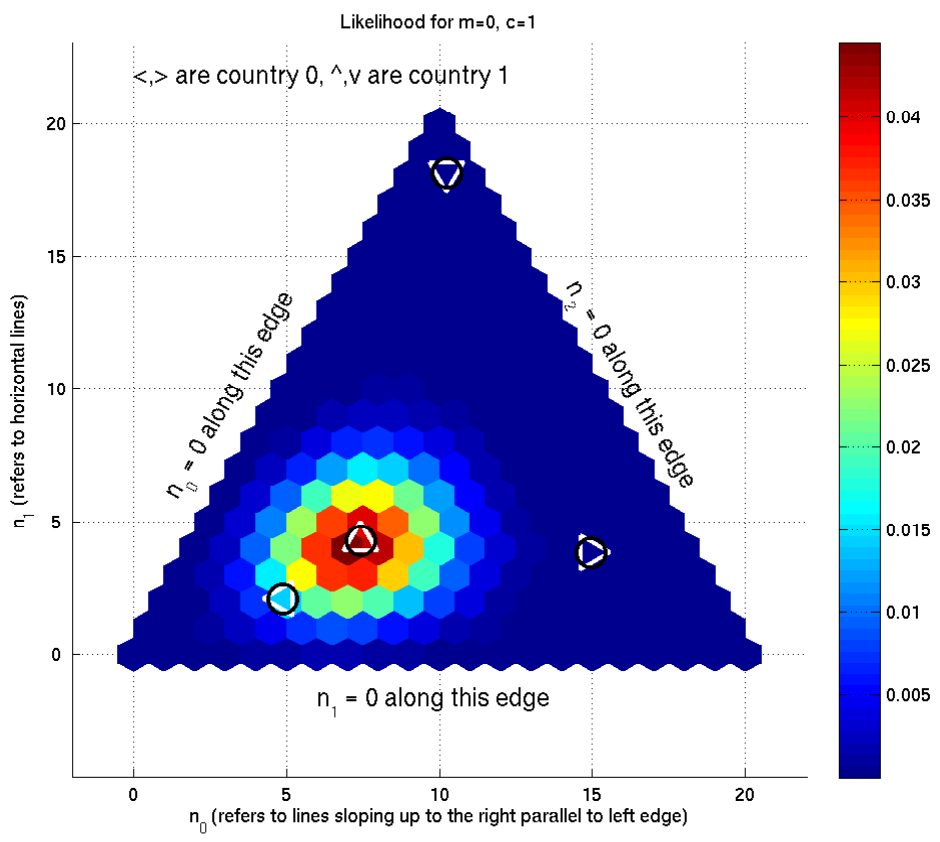} &
\includegraphics[scale=0.5]{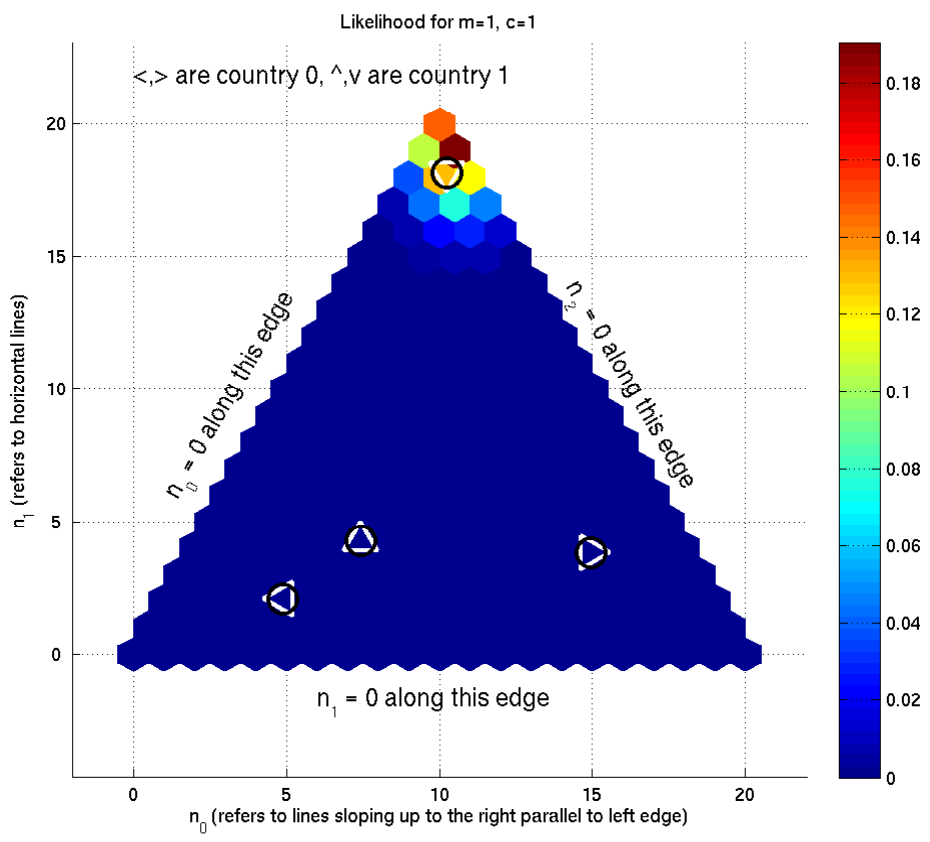}

\end{tabular}

\caption{Plots of the likelihood for each value of $(m,c)$. The
  \textbf{top} two plots are for country 0, and the \textbf{bottom}
  two for country 1. The \textbf{left} two plots are for mortar type
  0, and the \textbf{right} two for mortar type 1. Note that the
  colour scale on each plot is different.
\label{rawdists}
}

\end{figure}

The likelihood function $P(n|m,c)$ is then plotted for each value of
$(m,c)$ in figure \ref{rawdists}. 

\subsection{Information content}
\label{shellsinfointro}

\subsubsection{Intuitive introduction to Shannon information}

For precise details of the definitions of true and apparent Shannon
information (TSI and ASI respectively) see sections \ref{Shannon} and
\ref{ASI}; here we give an short intuitive introduction to help the
reader see why it is of interest to calculate these values.

It is often of interest to know how much one the value of one random
variable $x$ tells about another $\theta$. The posterior probability
distribution $P(\theta|x)$ tells us \textit{what} $x$ tells about
$\theta$; but sometimes we want to know not \textit{what} but
\textit{how much} $x$ tells us about $\theta$, in order that we can
compare this with (say) how much $y$ tells us about $\theta$.

In 1948 Claude Shannon defined in \cite{Shannon} what is now the
generally accepted measure of this, namely the (true) ``Shannon
information'' (TSI) about $\theta$ in $x$. Roughly speaking, this
measures the expectation of the logarithm of the increase in
probability density on the right answer as a result of observing $x$
-- expressed in bits if the logarithm is taken to base 2. 

As an example, suppose before observing $x$ we know that $\theta$ can
take one of two possible values, say 0 and 1, and that each is equally
likely, so that the probability that each is true is $\frac{1}{2}$. If
we then observe data telling us that $\theta$ is definitely 0 (when
indeed $\theta=0$), we have gained one bit of information, because the
posterior probability $P(\theta=0|x)$ is now one, double what it was
before. 

But there is another way this example could be considered, namely that
we have halved the number of possible values of $\theta$. To see that
it is not this that matters, consider the situation when instead of
giving us definite knowledge, the value of $x$ observed instead leads
us to have posterior probabilities $P(\theta=0|x) = 0.99,
P(\theta=1|x=0.01)$ when again in truth $\theta=0$. We have then
\textit{almost} doubled the posterior probability density on the right
answer, without having reduced the possible number of values $\theta$
could take at all -- and it is clearly reasonable to say that we have
gained \textit{almost} one bit of information.

Further, taking the expectation of this increase in probability
density allows us to take suitable note of how often we get how much
advantage from observing $x$. Moreover we find that the definition of
Shannon information resulting leads us to the Data Processing Theorem
(there is a proof in section \ref{basicSi}), which gives us the
intuitively beautiful and appropriate result that we cannot ``magic
up'' information out of thin air, i.e. we cannot increase the amount
of information about $\theta$ that is already present in $x$ by
processing it somehow. For example, if I want to know how much I
weigh, and use as a data source the brightness of the sky, then no
amount of clever signal processing and analysis will allow me to
deduce any more information from this data about my weight than is
present in the data (in this case there is none).

Another relevant quantity is the \textit{apparent} Shannon information
(ASI). Intuitively ASI tells us how much information a particular
solution gives us about what we wanted to know, when it is interpreted
at face value; while TSI tells us how much about what we wanted to
know we can recover by further processing. Suppose I toss a coin that
you can't see, but there are two oracles, one of which tells you what
I actually toss, and the other of which always lies. Both give you 1
bit of TSI about the coin toss, because since you know that the lying
oracle \textit{always} lies you can simply invert its answer. However
the truthful oracle also gives you 1 bit of ASI -- you can take what
it says at face value. The lying oracle gives you $-\infty$ bits of
ASI, because if you believe it you will then have a prior of zero on
the truth for any further analysis, and would then need an infinite
amount of information to make up for this.

\subsubsection{Why the Shannon information in a definite answer varies
  with the prior}

It will be noticed in what follows that the amount of Shannon
information provided by a definite answer varies with the prior on
$\theta$; this subsection aims to give an intuition for why this
should be so.

Suppose I want to know whether or not there has been an earthquake
today in Cambridge in England, where earthquakes are extremely
rare. Then I already know -- before being provided with a definite
answer -- that there probably hasn't. For the same question in San
Francisco the situation is entirely different -- earthquakes there are
common, so I am in more doubt before I start as to whether or not
there might have been one today, and answering the question provides
quantitatively more information than it does in the case of Cambridge.

Of course, if you are told that there actually \textit{has} been an
earthquake today in Cambridge, that provides \textit{more} information
than if I tell that there has been one in San Francisco. But since
\textit{on average} I will be telling you that there has \textit{not}
been one in Cambridge, the expectation of the logarithm of the
increase in probability will be very small.

\subsubsection{Specifics of Shannon information for the shells problem}

In this particular problem the data space is discrete (so contains
only finitely many possible values), while the posterior probability
and each type of frequentist confidence are continuously
distributed. As a result, not only does the posterior probability
retain all the relevant TSI and ASI in the data (as is always the
case), but also all the types of frequentist confidence retain all the
TSI -- simply because the value of frequentist confidence returned for
each of the finitely many possible data values is different, so by
using a lookup table we can recover the data from the frequentist
confidence. We therefore don't report the TSI in these cases.

On the other hand the TSI in the frequentist confidence sets for any
fixed value of the confidence $\eta$ will be smaller, because given
any particular confidence set value there will be several possible
data values that would lead to it. 

\subsubsection{Specifics of ASI calculation for the shells problem}

Turning to ASI, in the notation of section \ref{ASI}, we calculate it
for frequentist confidence by assigning $Q_n(c=1) = \eta$ and
$Q_n(c=0) = 1 - \eta$, where $\eta$ is the value of frequentist
confidence returned. In the case of frequentist confidence sets where
we assume only the value for a particular confidence level $\eta$ to
be delivered, we calculate $Q_n(c)$ for a confidence set function
$g(n,\eta)$ by sharing $\eta$ evenly among those values of $c$ in
$g(n,\eta)$ and sharing $1-\eta$ evenly among those values of $c$ not
in $g(n,\eta)$. Where $g(n,\eta)$ contains (only) a single country,
this gives a probability distribution on $c$; but where either both
countries are in $g(n,\eta)$ or where neither is it gives us a
subprobability distribution (i.e. one that integrates (or sums) to
less than 1). The total probability therein will be $\eta$ where
$g(n,\eta)=\{0,1\}$, and $1-\eta$ where $g(n,\eta)=\emptyset$; this
seems reasonable given that saying that you are confident that
$c\in\emptyset$ is unreasonable anyhow, as is saying that you are less
than 100\% confident that $c\in\{0,1\}$.

We now have the framework for this example in place and can look at
the various possible solutions.

\subsection{The Bayesian solution}

For the Bayesian solution we first choose a
prior $$P(m,c)=P(c)P(m|c)$$ on the unknowns. We then calculate as
follows:
$$P(c|n) = \frac{\sum_{m=0}^1{P(m,c)P(n|m,c)}}{\sum_{c=0}^1{\sum_{m=0}^1{P(m,c)P(n|m,c)}}}.$$

\begin{figure}
\begin{tabular}{cc}

\includegraphics[scale=0.5]{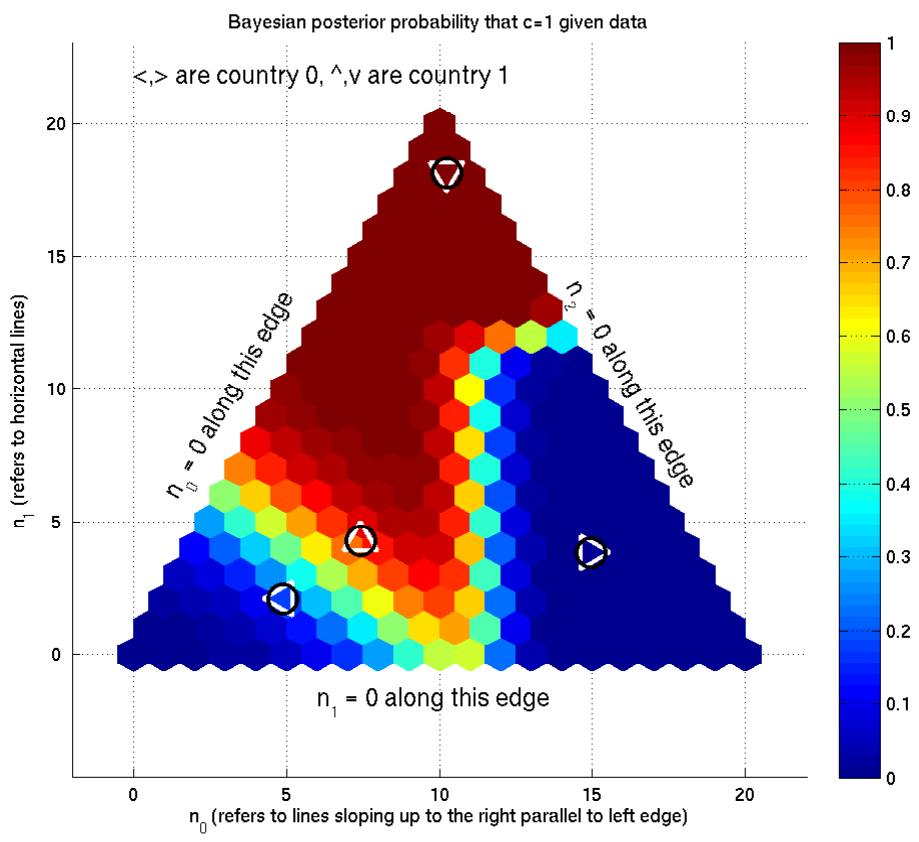} &
\includegraphics[scale=0.5]{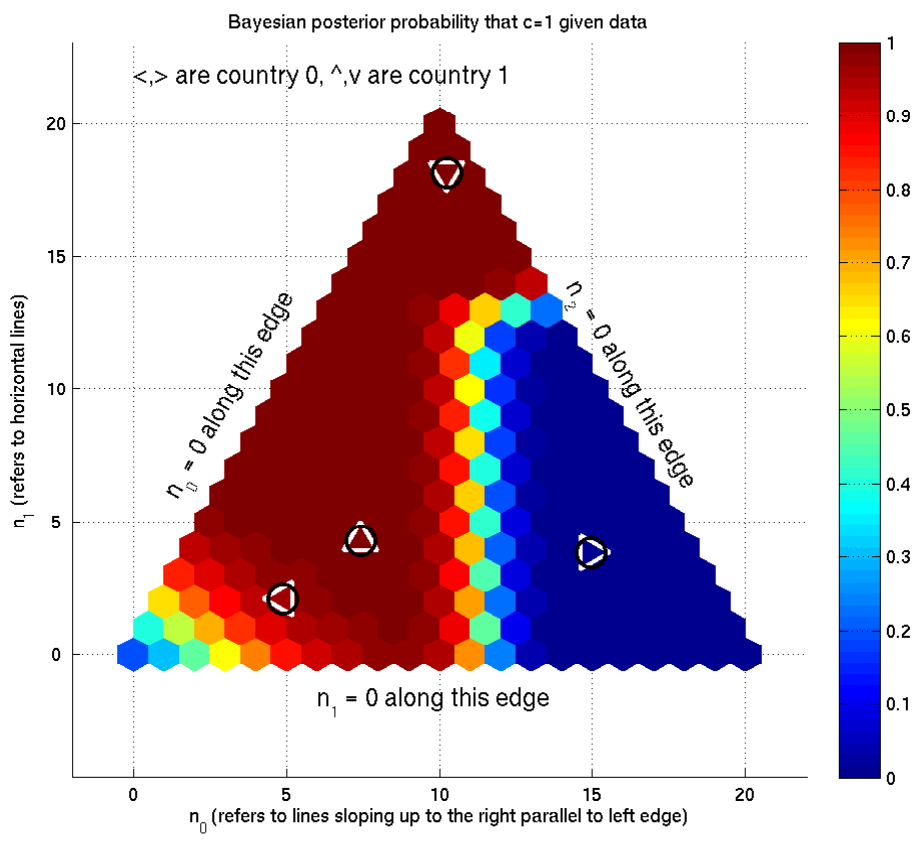} \\
\includegraphics[scale=0.5]{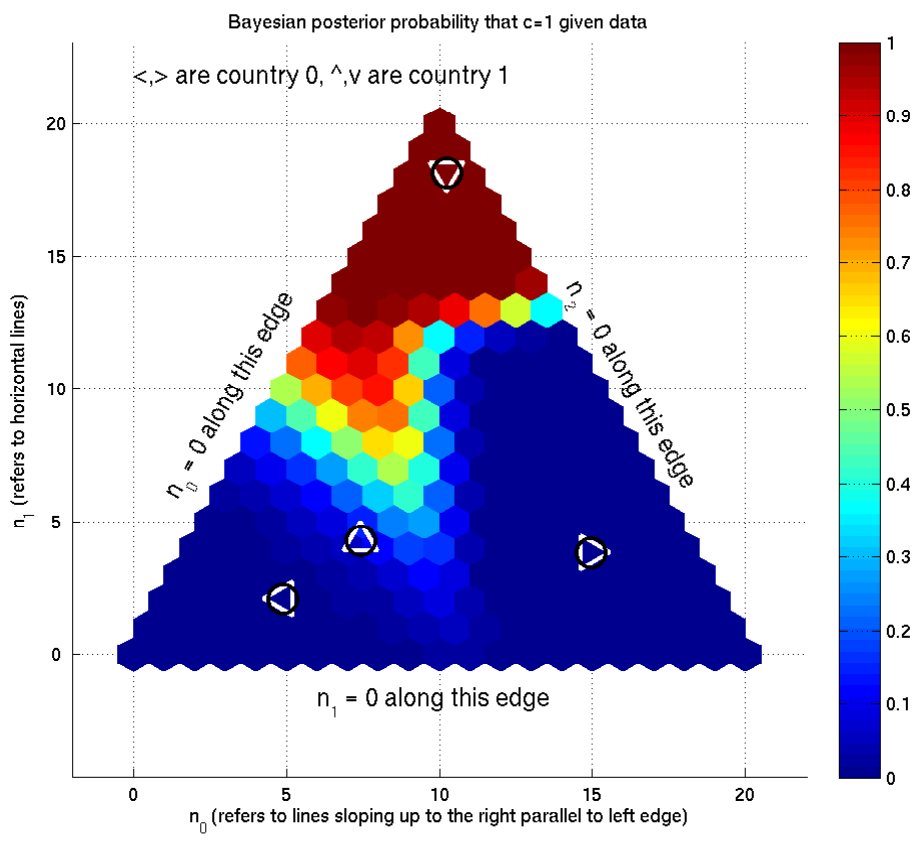} &
\includegraphics[scale=0.5]{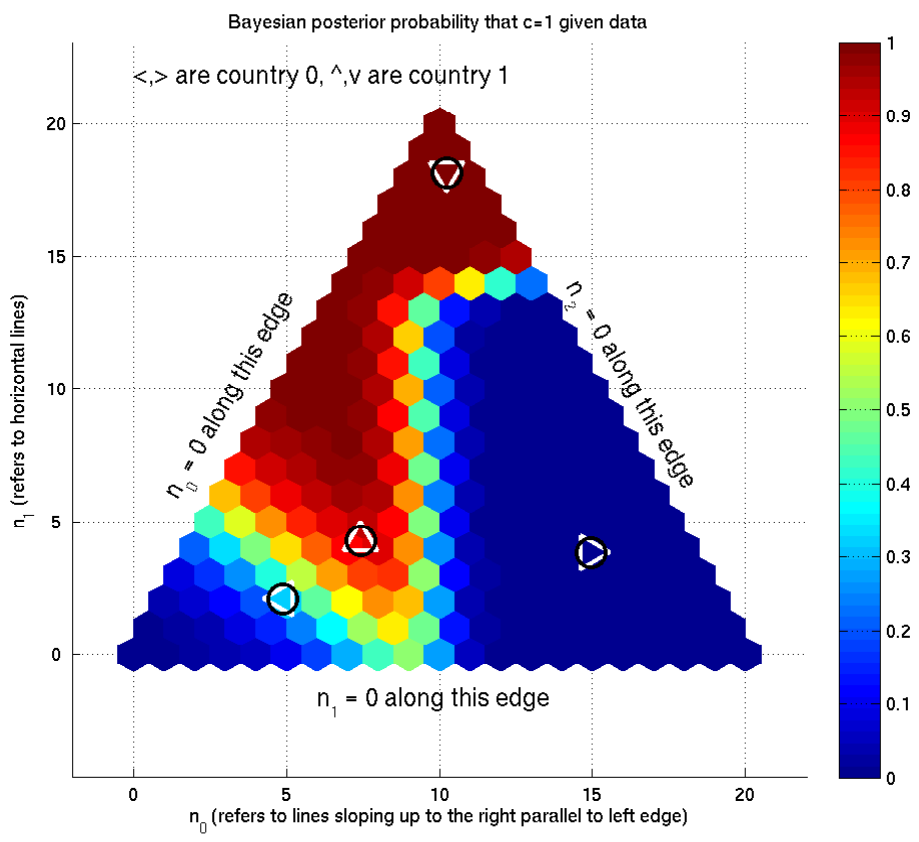} \\

\end{tabular}

\caption{The Bayesian solution for the shells problem of section
  \ref{shells}. The \textbf{top} two plots have $P(c=0) = P(c=1) =
  \frac{1}{2}$. The \textbf{bottom} two plots have $P(c=0) = 0.98,
  P(c=1) = 0.02$. The \textbf{left} two plots have $P(m|c) = 0.5$ for
  all four values of $(m,c)$. The \textbf{right} two plots have
  $P(m|c) = 0.01$ if $m=c$ and $0.99$ otherwise (saying that we think
  type 0 mortars are much more likely to be used with shells from
  country 1, etc.). The ASI given by the Bayesian solution is always
  equal to the TSI (and to the TSI about country in the data); the
  values (reading from top left to top right then bottom left, bottom
  right) are 0.669, 0.888, 0.088, 0.109 bits. This compares with 1 bit
  needed to draw a definite conclusion for the top two plots, and
  0.141 bits for a definite conclusion in the bottom two plots (as in
  the bottom two we already knew more about which answer was likely
  \textit{a priori}).
\label{shellsBayes}
}
\end{figure}

The results of the Bayesian solution are plotted in figure
\ref{shellsBayes} for several possible choices of the prior. The
values of the ASI resulting, in all cases equal to the TSI in the
Bayesian solution and the TSI in the data, and less than the TSI or
ASI required to get a certain correct answer to the User's question,
are shown in the caption.

\subsection{Frequentist hypothesis testing solutions}

\subsubsection{General}

In all cases we treat $H_0$, the hypothesis that country 0
provided the shells, as the null hypothesis, with the alternative
hypothesis being $H_1$. For all these methods, as discussed in section
\ref{shellsinfointro}, the TSI in the solution is equal to that in the
data and in the Bayesian solution, because the data can be deduced
from the frequentist confidence that $H_1$ holds. The ASI, however, is
in general substantially less than the TSI, and is given in the
captions to the various figures.

\subsubsection{Critical regions based on distance from expected values
  under the null hypothesis}

Choosing critical regions to consist of the intersections of the
complements of discs centred on each of the two expected values under
the null hypothesis $H_0$, we get the top left plot in figure
\ref{shellsnear}. Choosing critical regions based on the minimum
distance to an expected value under $H_1$ we get the top right plot in
figure \ref{shellsnear}. In all cases the ASI obtained is
substantially less than that in the Bayesian solution.

\subsubsection{Critical regions based on probability under $H_0$ or
  $H_1$} 

\begin{figure}
\begin{tabular}{cc}

\includegraphics[scale=0.5]{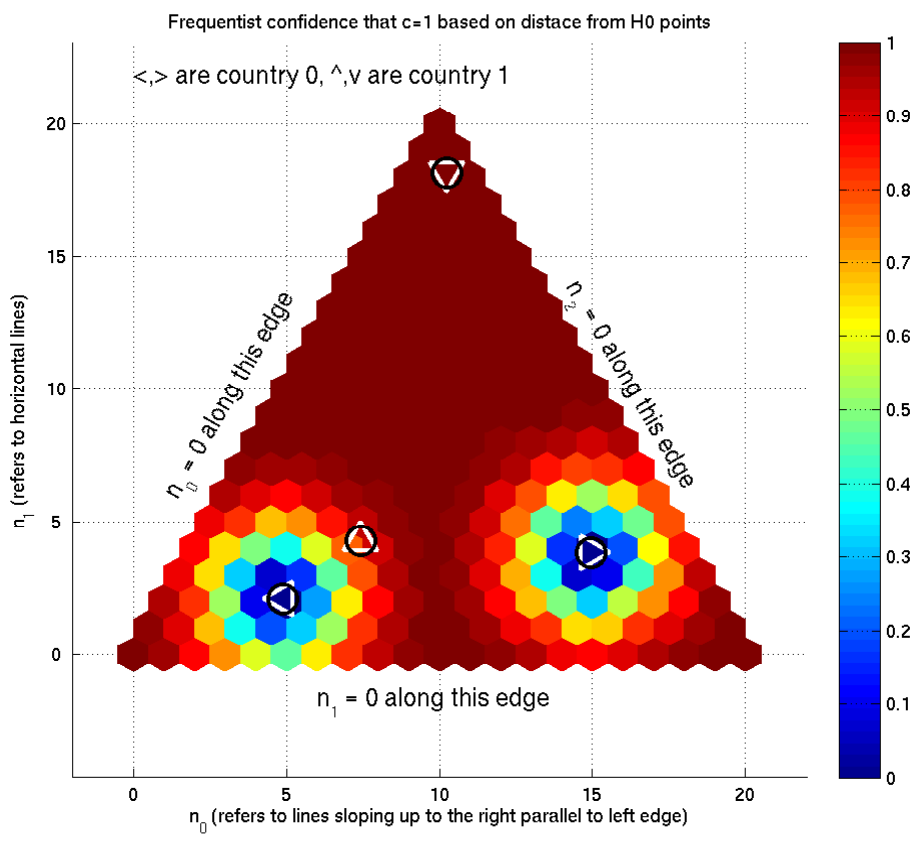} &
\includegraphics[scale=0.5]{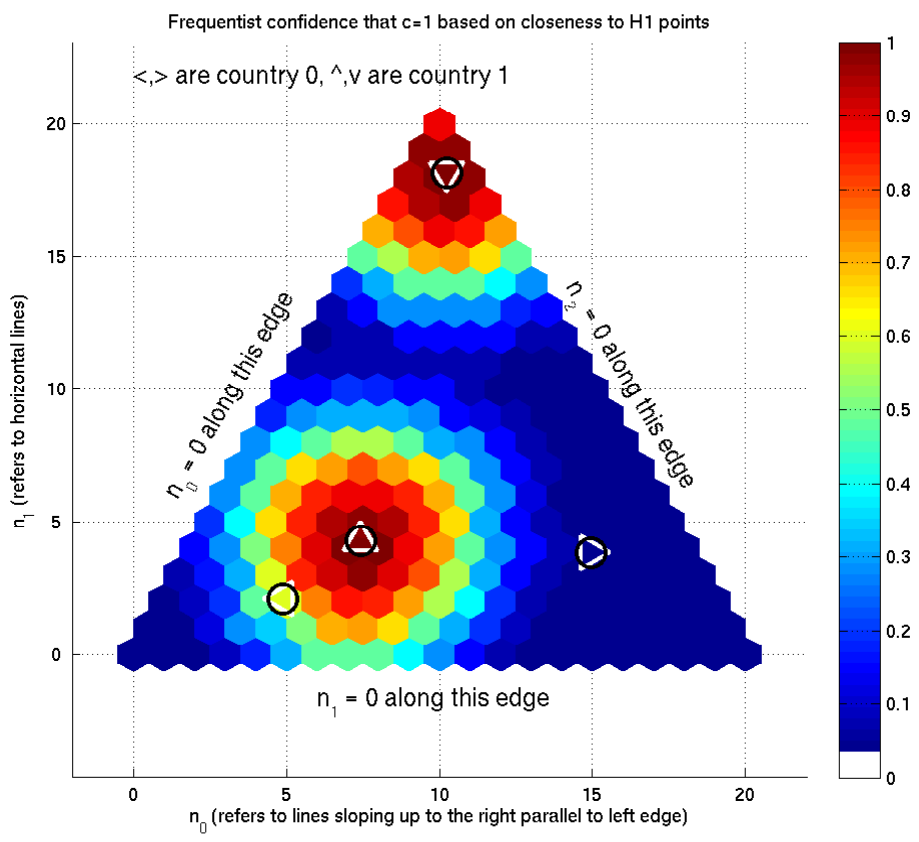} \\
\includegraphics[scale=0.5]{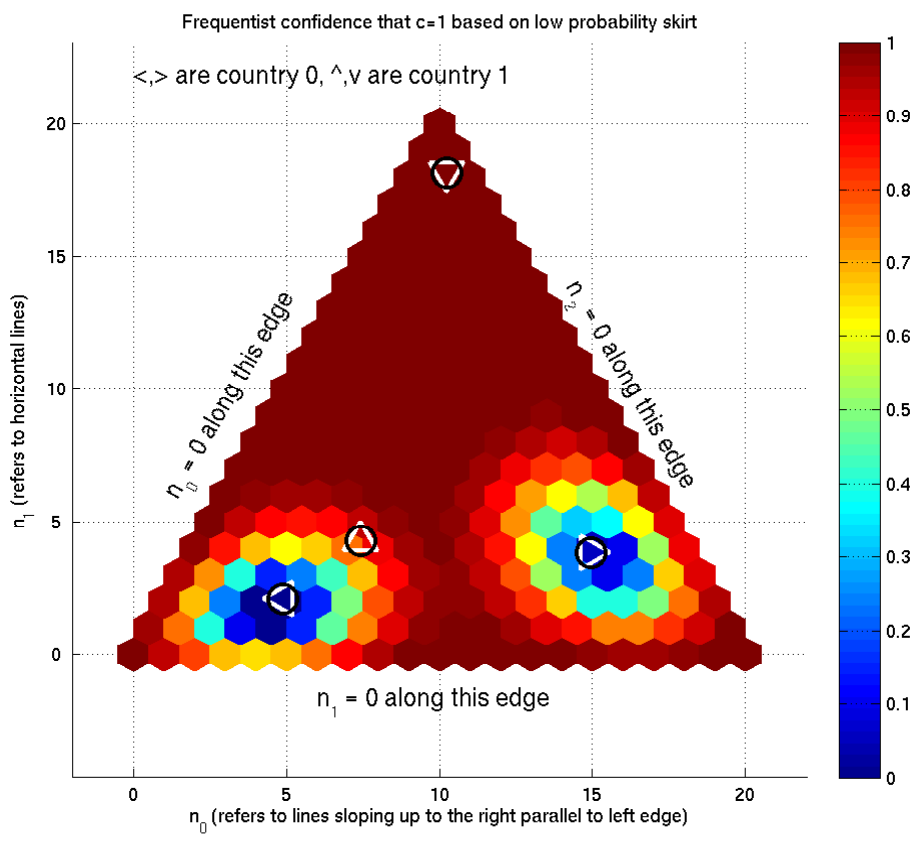} &
\includegraphics[scale=0.5]{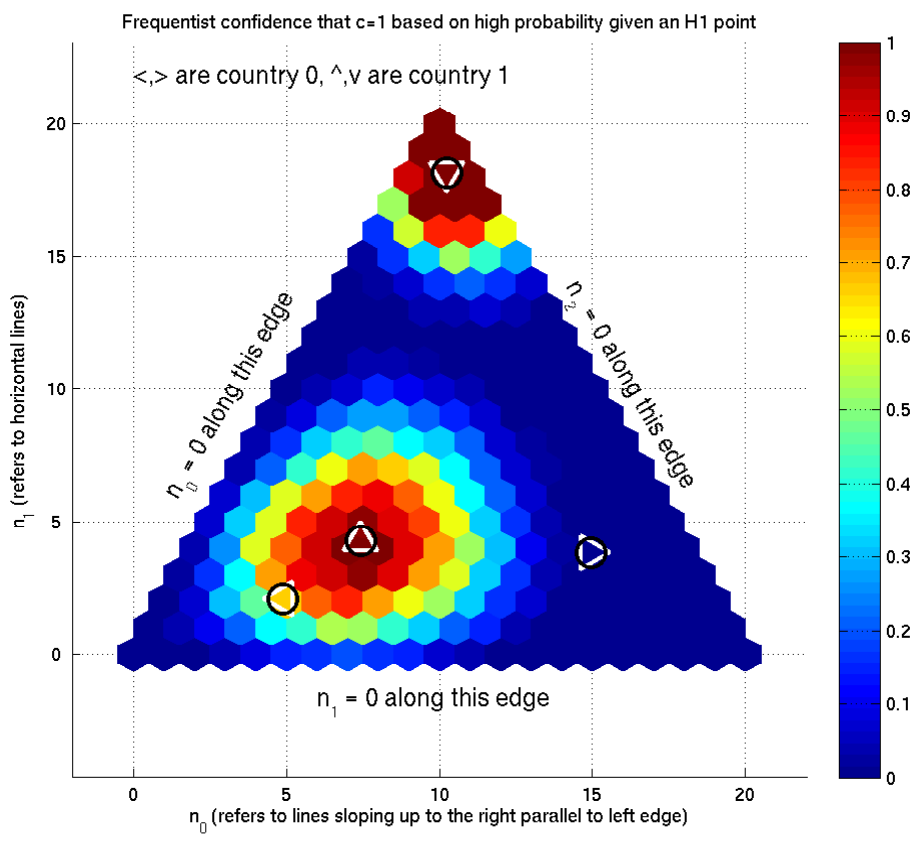} \\

\end{tabular}

\caption{Frequentist confidence that country 1 provided the shells
  given the data, based on: \textbf{Top left:} Critical regions based
  on the minimum distance from an expected value under the null
  hypothesis $H_0$; under the four priors discussed in the caption of
  figure \ref{shellsBayes} the ASI in this solution is 0.118, -0.123,
  -1.196, -1.286 bits respectively, all much lower than the
  corresponding values for the Bayesian solution. \textbf{Top right:}
  Critical regions based on the minimum distance from an expected
  value under the alternative hypothesis $H_1$; under the four priors
  discussed in figure \ref{shellsBayes} the ASI in this solution is
  0.482, 0.718, -0.586, 0.060 bits respectively. \textbf{Bottom left:}
  Critical regions consisting of the lowest likelihood under any
  subhypothesis of $H_0$; under the four priors discussed the ASI in
  this solution is 0.070, -0.269, -1.149, -1.291 bits
  respectively. \textbf{Bottom right:} Critical regions consisting of
  the highest likelihood under any subhypothesis of $H_1$; under the
  four priors discussed the ASI in this solution is 0.451, 0.685,
  -0.619, 0.056 bits respectively.
\label{shellsnear}
}

\end{figure}

Choosing critical regions to include those parts of data space of
lowest likelihood under $H_0$ we get the bottom left plot in figure
\ref{shellsnear}. Choosing them based on highest likelihood under
$H_1$ we get the bottom right plot in that figure. Again, the ASI
obtained is less than that in the Bayesian solution. To be quite
clear, whichever way we choose the \textit{sets} constituting the
various critical regions, we still calculate the frequentist
confidence based on the probability of getting data in those sets
given $H_0$.

\subsubsection{(Basic) pseudo-Bayesian (frequentist) solutions}

\begin{figure}
\begin{tabular}{cc}

\includegraphics[scale=0.5]{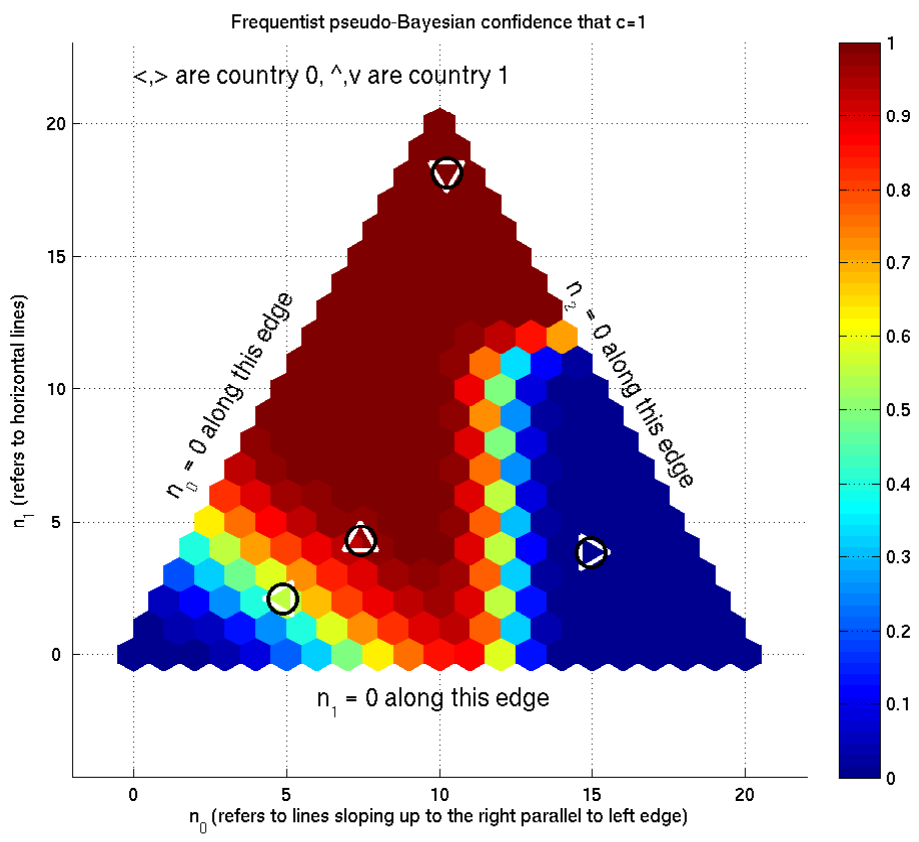} &
\includegraphics[scale=0.5]{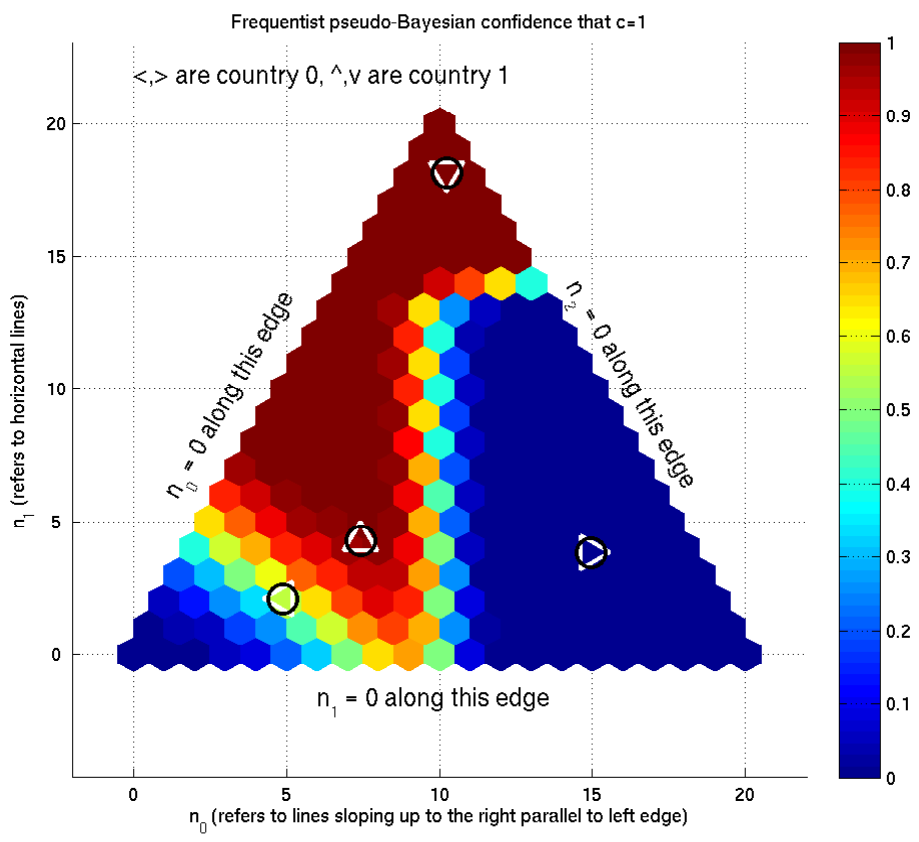} \\
\includegraphics[scale=0.5]{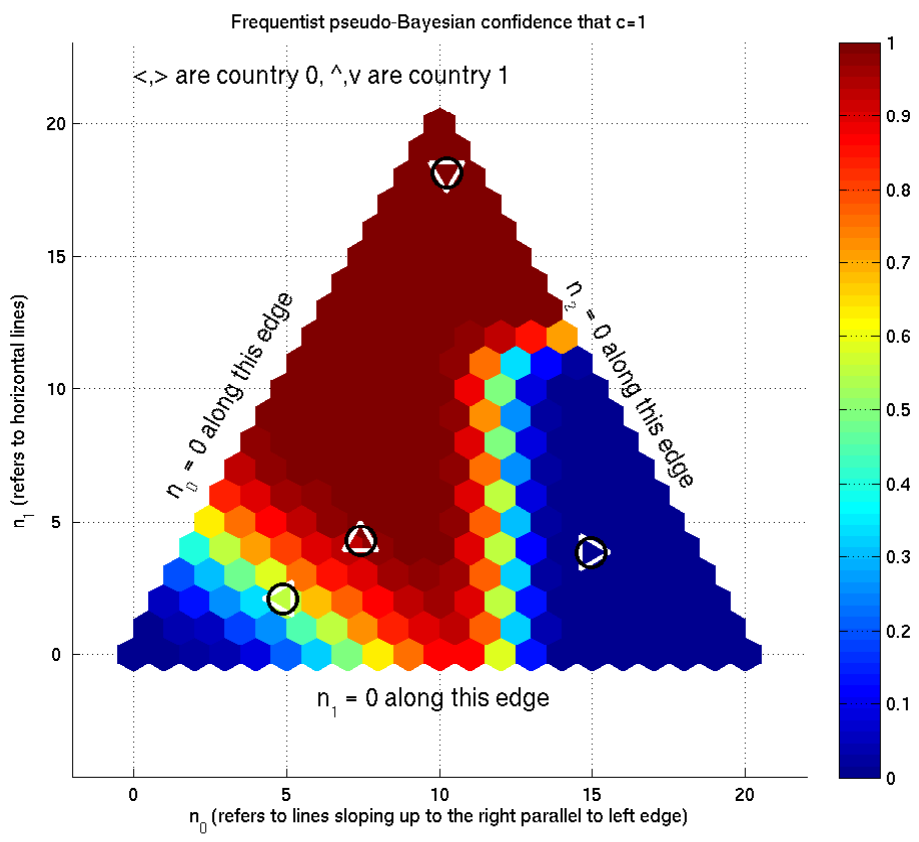} &
\includegraphics[scale=0.5]{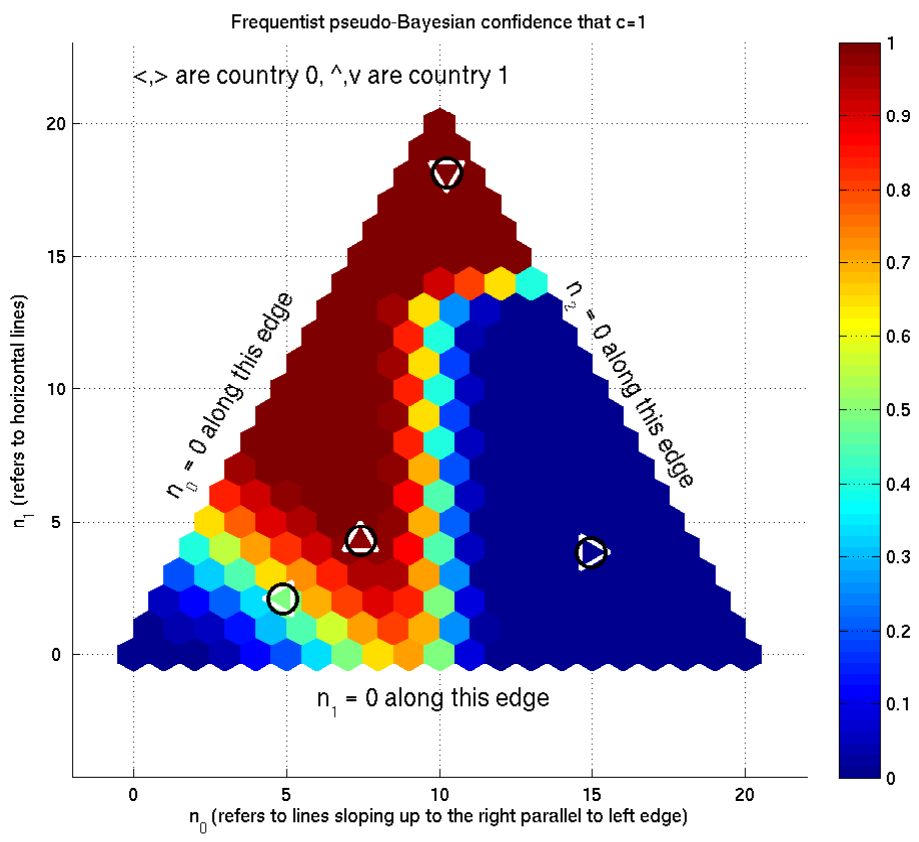} \\

\end{tabular}

\caption{The pseudo-Bayesian frequentist confidence that $H_2$ holds
  under the four priors discussed in the caption to figure
  \ref{shellsBayes}. The \textbf{top} two plots have $P(c=0) = P(c=1)
  = \frac{1}{2}$. The \textbf{bottom} two plots have $P(c=0) = 0.98,
  P(c=1) = 0.02$. The \textbf{left} two plots have $P(m|c) = 0.5$ for
  all four values of $(m,c)$. The \textbf{right} two plots have
  $P(m|c) = 0.01$ if $m=c$ and $0.99$ otherwise (saying that we think
  type 0 mortars are much more likely to be used with shells from
  country 1, etc.). The ASI in the four plots are respectively 0.538,
  0.657, -0.662, 0.104 reading in the same order as in the caption to
  figure \ref{shellsBayes}. Despite the ASIs being different, the
  inference (shown in the plots) in the top row is the same as that in
  the bottom row, but the left-hand plots differ from the right-hand
  plots.
\label{shellspseudoBayes}
}
\end{figure}

Whereas in section \ref{freqpseudoBayes} on the example of section
\ref{example1} there was only one possible pseudo-Bayesian solution
independent of the prior $P(\theta)$, in this example there are many,
depending not on the prior $P(c)$ on country of origin of the shells,
but on the probability $P(m|c)$ on mortar type given country of origin
in the prior $P(m,c)=P(m|c)P(c)$. We illustrate what happens for
different $P(c)$ and different $P(m|c)$ in figure
\ref{shellspseudoBayes}. Despite the fact that the Bayesian prior is
used in constructing the pseudo-Bayesian solutions, the ASI is in
three out of four cases substantially lower than for the Bayesian
solution. 

\subsection{Frequentist confidence sets}

Just as in section \ref{fconfsets} for the example of section
\ref{example1}, we may construct ``inconclusive'' confidence set
functions from any single choice of frequentist critical regions, and
``conclusive'' ones either from one such type used twice (the second
time with the role of null and alternative hypotheses reversed), or
from two different choices of frequentist critical regions. In all
cases we will show only the confidence set functions and the values of
TSI and ASI resulting from restricting the value of the confidence
parameter $\eta$ to the value $0.95$ .

\subsubsection{``Inconclusive'' confidence set functions}

\begin{figure}
\begin{tabular}{cc}

\includegraphics[scale=0.5]{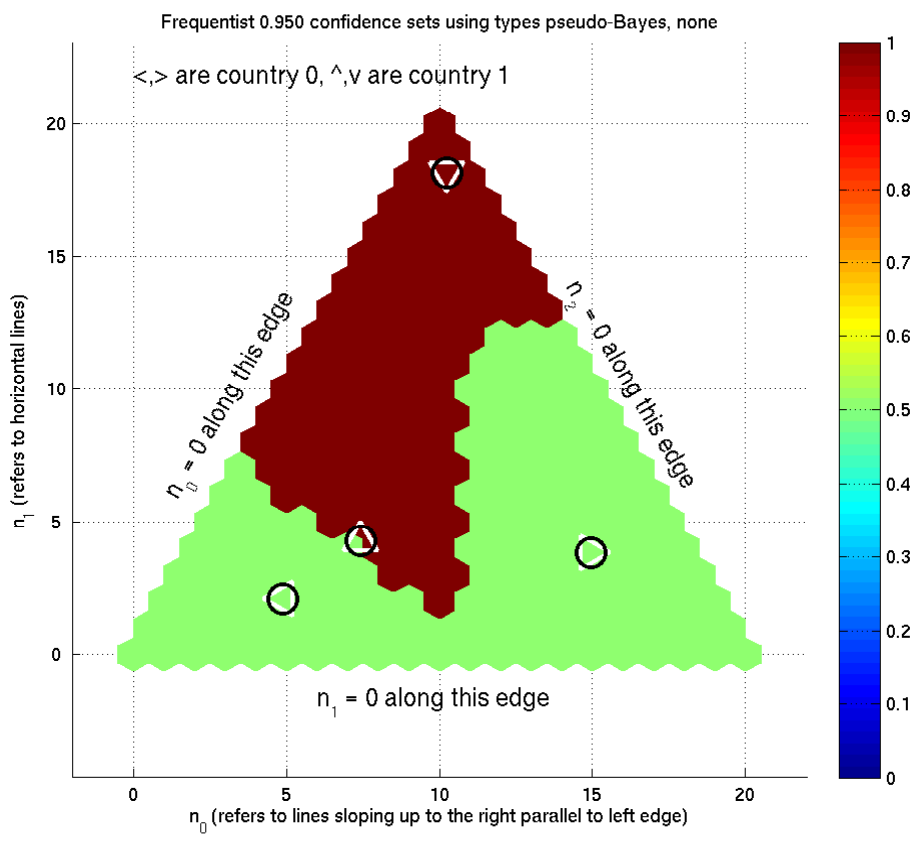} &
\includegraphics[scale=0.5]{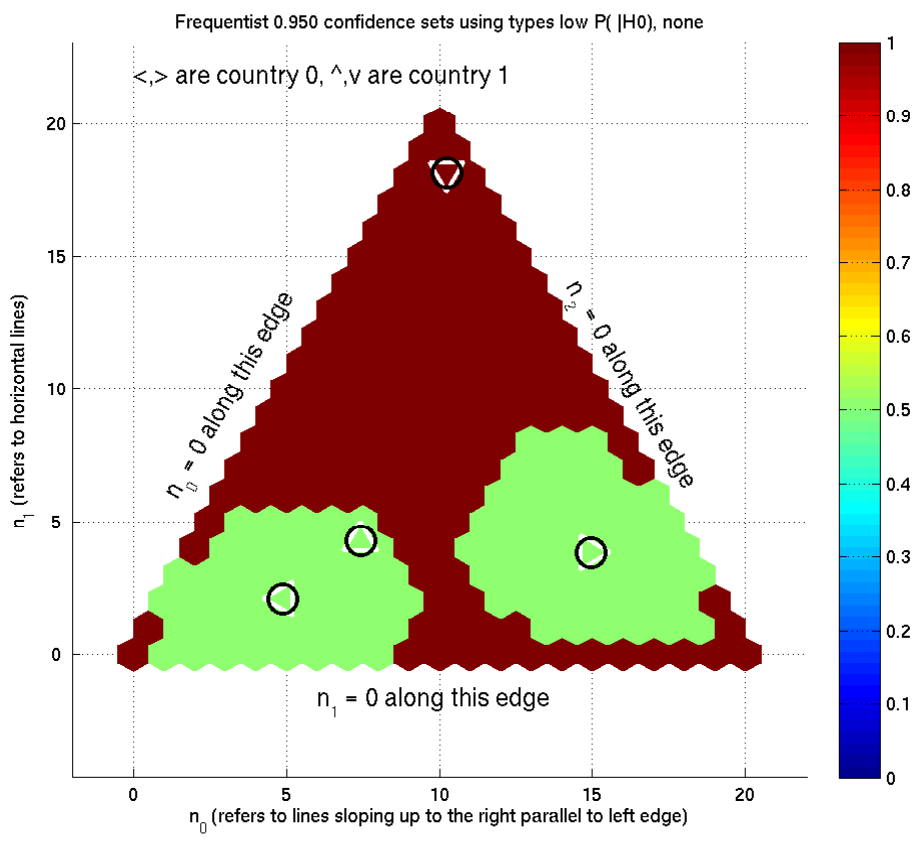} \\
\includegraphics[scale=0.5]{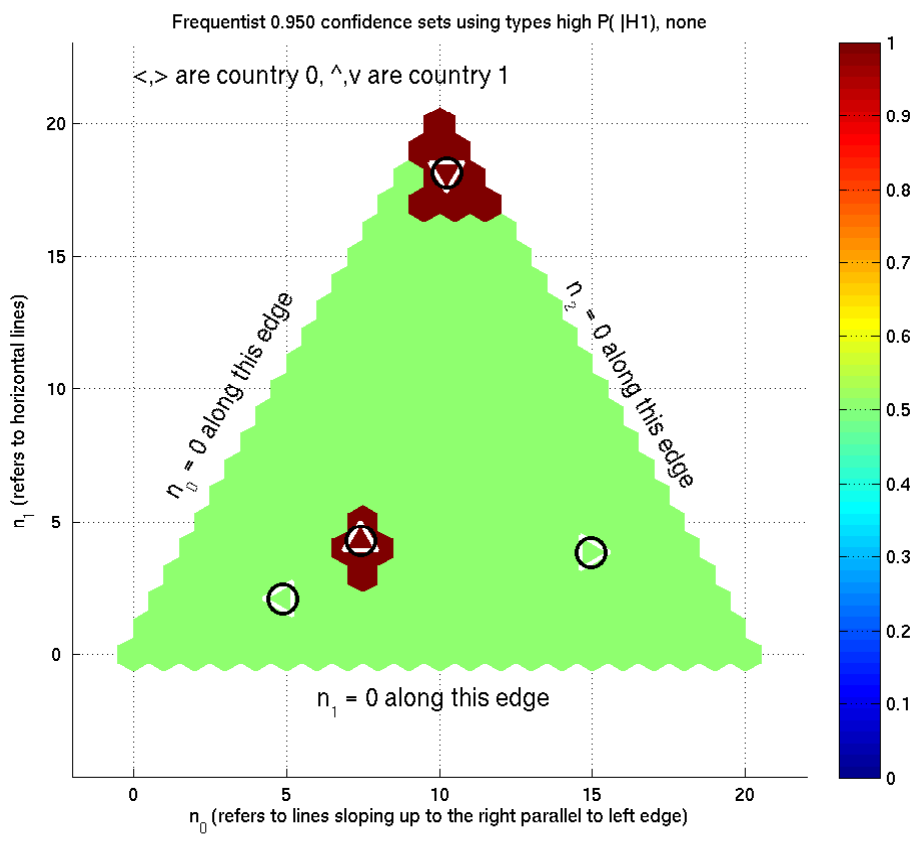} &
\includegraphics[scale=0.5]{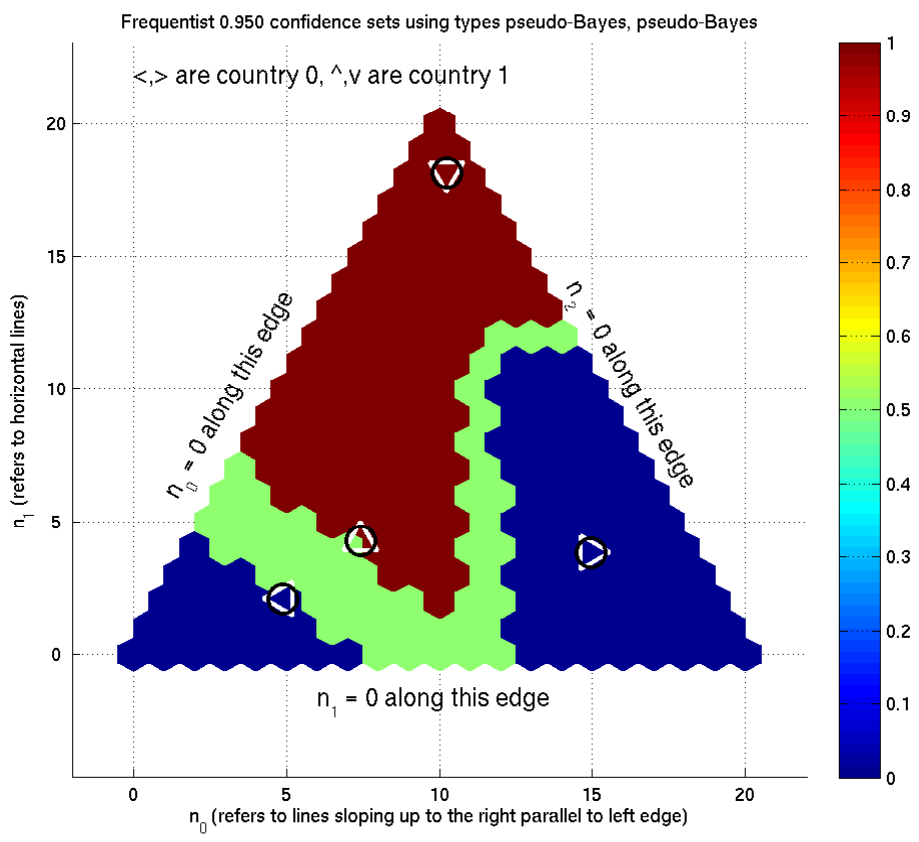} \\

\end{tabular}

\caption{``Inconclusive'' frequentist confidence set functions for
  $\eta=0.95$. In all cases brown means $f(n,\eta)=\{1\}$, dark blue
  means $f(n,\eta)=\{0\}$, green means $f(n,\eta)=\{0,1\}$, and yellow
  means $f(n,\eta)=\emptyset$. The \textbf{top left} plot is for the
  pseudo-Bayesian confidence set function, the \textbf{top right} for
  the one based on low likelihood under the null hypothesis $H_0$, and
  the \textbf{bottom left} for the one based on high likelihood under
  the alternative hypothesis $H_1$. Evaluated only under the first of
  the priors discussed in the caption of figure \ref{shellsBayes}, TSI
  values are respectively 0.489, 0.385, and 0.256 bits, and ASI values
  are respectively 0.260, 0.207, and 0.131 bits.
\label{shellsinconc}
}
\caption{The \textbf{bottom right} plot shows the ``conclusive''
  frequentist confidence set function for $\eta=0.95$ based on
  pseudo-Bayesian critical regions. TSI is 0.601 and ASI is 0.584 bits
  under the first of the priors discussed in the caption to figure
  \ref{shellsBayes}.
\label{shellsconc1}
}
\end{figure}

The top left plot of figure \ref{shellsinconc} shows the inconclusive
frequentist confidence set function based on the pseudo-Bayesian
critical regions; the top right plot that based on low likelihood
under the null hypothesis $H_0$; and the bottom left plot that based
on high likelihood under the alternative hypothesis $H_1$.

\subsubsection{``Conclusive'' confidence set functions}

\begin{figure}
\begin{tabular}{cc}

\includegraphics[scale=0.5]{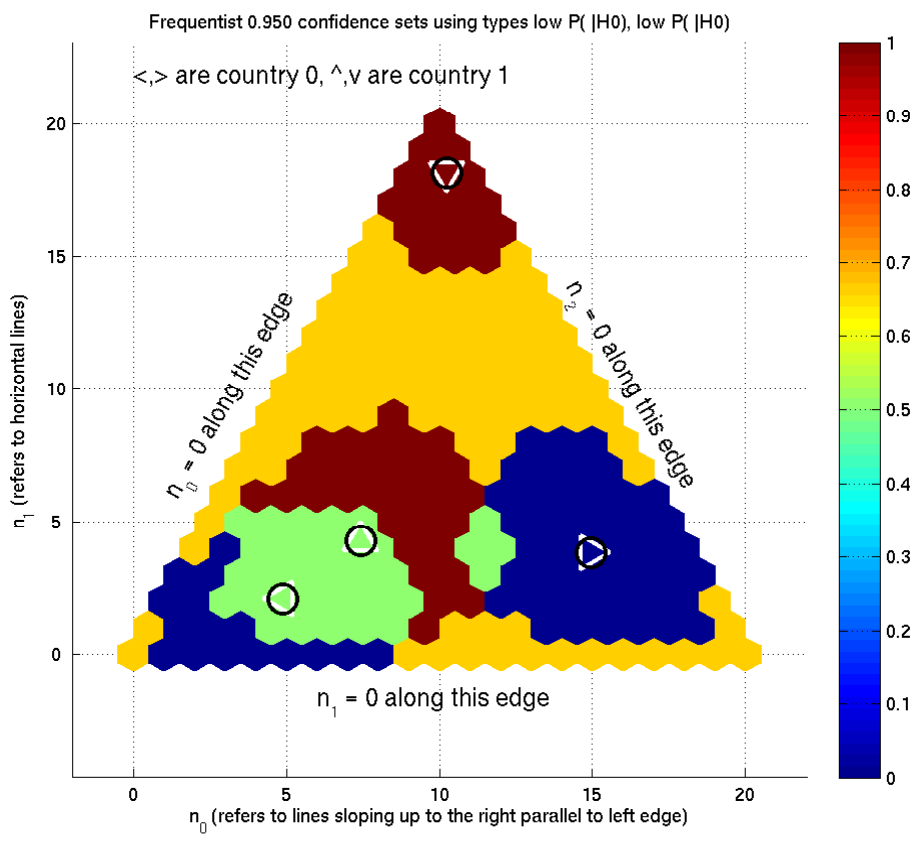} &
\includegraphics[scale=0.5]{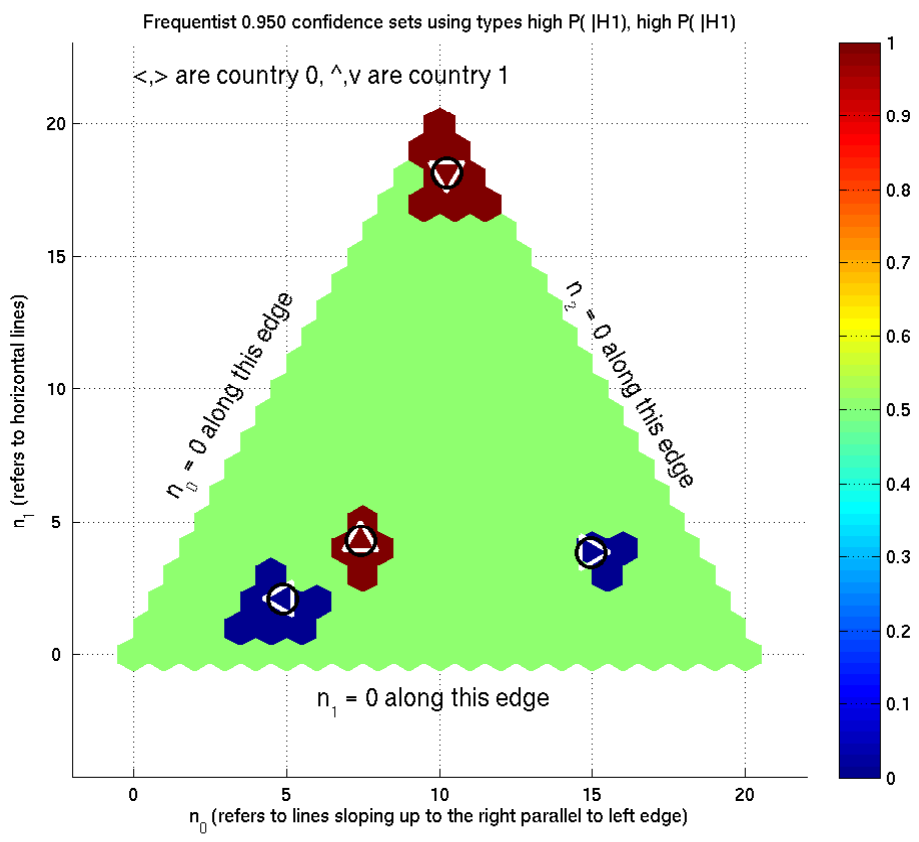} \\
\includegraphics[scale=0.5]{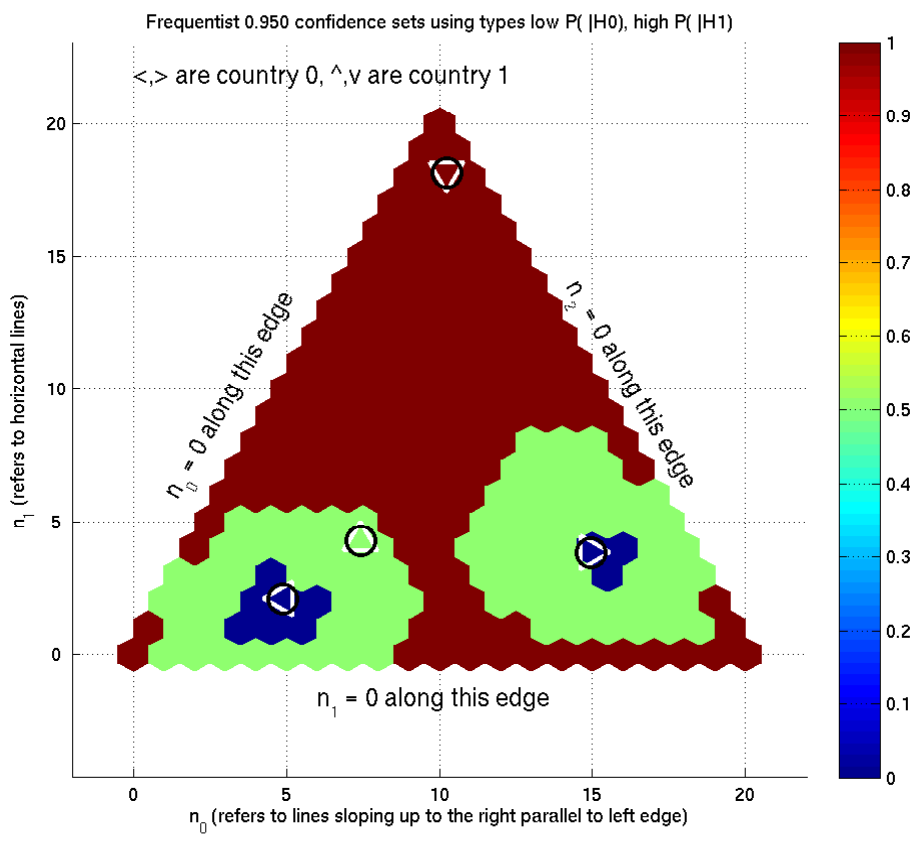} &
\includegraphics[scale=0.5]{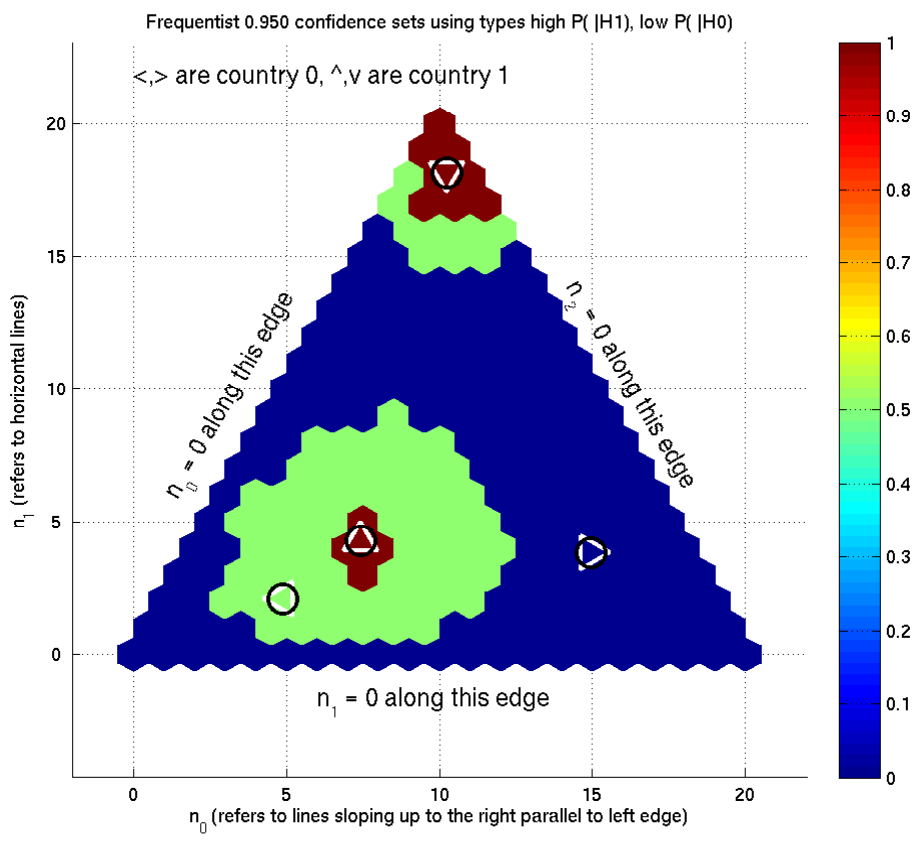} \\

\end{tabular}

\caption{``Conclusive'' frequentist confidence set functions for
  $\eta=0.95$. In all cases brown means $f(n,\eta)=\{1\}$, dark blue
  means $f(n,\eta)=\{0\}$, green means $f(n,\eta)=\{0,1\}$, and yellow
  means $f(n,\eta)=\emptyset$. The \textbf{top left} plot is uses
  critical regions based on low likelihoods under the null hypothesis,
  the \textbf{top right} one on high likelihoods under the alternative
  hypothesis, the \textbf{bottom left} one uses critical regions that
  use low likelihood under $H_0$ to reject $H_0$ but high likelihood
  under $H_0$ to reject $H_1$, while the \textbf{bottom right} plot is
  similar but for $H_1$ and $H_0$ reversed. Evaluated only under the
  first of the priors discussed in the caption of figure
  \ref{shellsBayes}, TSI values are respectively 0.549, 0.282, 0.398,
  and 0.420 bits, while ASI values are respectively 0.420, 0.200,
  0.276, and 0.387 bits.
\label{shellsconc2}
}

\end{figure}

The bottom right plot of figure \ref{shellsconc1} shows the
``conclusive'' frequentist confidence set function for the
pseudo-Bayesian critical regions considered in figure
\ref{shellspseudoBayes}. The top left plot of figure \ref{shellsconc2}
shows that for critical regions based on low likelihood under the null
hypothesis $H_0$, the top right one for those based on high likelihood
under the alternative $H_1$, while the bottom left one uses low
likelihood under $H_0$ to reject $H_0$ but high likelihood under $H_0$
to reject $H_1$, and the bottom right one is similar but for $H_1$ and
$H_0$ reversed. 

Note how for both inconclusive and conclusive confidence set functions
the mapping from data to frequentist confidence set varies enormously
between the different choices of confidence set function one can make.

\clearpage

\subsection{Summary of TSI and ASI values for the specific shells problem}

\label{shellsinfo}

Above we considered only one particular version of the problem, namely
that with the values of $P(f|m,c)$ given in table \ref{pfgivenmc}, for
which we now summarise the ASI values in table \ref{ASItable} and the
TSI values in table \ref{TSItable}.

\begin{table}[H]
\begin{center}

\begin{tabular}{lrrrr}
\textbf{Prior $P(c)$:} & \multicolumn{2}{c}{\textbf{A}} &
  \multicolumn{2}{c}{\textbf{B}} \\  
\textbf{Prior $P(m|c)$:} & \multicolumn{1}{c}{\textbf{C}} & \multicolumn{1}{c}{\textbf{D}} & \multicolumn{1}{c}{\textbf{C}} & \multicolumn{1}{c}{\textbf{D}} \\
\textbf{Method}\\
Needed to be definite & 1.000 & 1.000 & 0.141 & 0.141\\
Bayes & 0.669 & 0.888 & 0.088 & 0.109\\
Pseudo-Bayes & 0.538 & 0.657 & -0.662 & 0.104\\
Far from $H_0$ & 0.118 & -0.123 & -1.196 & -1.286\\
Near $H_1$ & 0.482 & 0.718 & -0.586 & 0.060\\
Unlikely under $H_0$ & 0.070 & -0.269 & -1.149 & -1.291 \\
Likely under $H_1$ & 0.451 & 0.685 & -0.619 & 0.056 \\
IFCS by pseudo-Bayes & 0.260 & 0.119 & -1.031 & -0.928 \\
IFCS by unlikely under $H_0$ & 0.207 & 0.051 & -1.055 & -1.081 \\
IFCS by likely under $H_1$ & 0.131 & 0.046 & -1.071 & -0.931 \\ 
CFCS by pseudo-Bayes & 0.584 & 0.530 & -0.301 & 0.033 \\
CFCS by unlikely under $H_0$ & 0.420 & 0.384 & -0.502 & -0.221 \\
CFCS by likely under $H_1$ & 0.200 & 0.044 & -0.860 & -0.782 \\
CFCS by $H_0$ & 0.276 & 0.049 & -0.844 & -0.932 \\
CFCS by $H_1$ & 0.387 & 0.442 & -0.471 & 0.002 \\
\end{tabular}

\caption{Apparent Shannon information (ASI) values (bits) (section
  \ref{shellsinfointro}) for the different methods on the shells
  problem of section \ref{shells} under various different priors. IFCS
  = inconclusive frequentist confidence set function; CFCS =
  conclusive frequentist confidence set function.  Priors are:
  \textbf{A}: $P(c)=0.5\,\forall c$; \textbf{B}:
  $P(c=0)=0.98,P(c=1)=0.02$; \textbf{C}: $P(m|c)=\frac{1}{2}\,\forall
  m,c$; \textbf{D}: $P(m|c)=\left\{\begin{matrix}0.01 & (m=c)\\ 0.99 &
  (m\neq c)\end{matrix}\right.$. In each case the settings for
  frequentist methods that use prior information matched the prior
  under which ASI was assessed.
\label{ASItable}
}

\end{center}
\end{table}

\begin{table}[H]
\begin{center}

\begin{tabular}{lrrrr}
\textbf{Prior $P(c)$:} & \multicolumn{2}{c}{\textbf{A}} &
  \multicolumn{2}{c}{\textbf{B}} \\  
\textbf{Prior $P(m|c)$:} & \multicolumn{1}{c}{\textbf{C}} & \multicolumn{1}{c}{\textbf{D}} & \multicolumn{1}{c}{\textbf{C}} & \multicolumn{1}{c}{\textbf{D}} \\
\textbf{Method}\\
Needed to be definite & 1.000 & 1.000 & 0.141 & 0.141\\
Bayes & 0.669 & 0.888 & 0.088 & 0.109\\
IFCS by pseudo-Bayes & 0.489 & 0.225 & 0.056 & 0.043 \\
IFCS by unlikely under $H_0$ & 0.385 & 0.147 & 0.044 & 0.016 \\
IFCS by likely under $H_1$ & 0.256 & 0.129 & 0.029 & 0.025 \\
CFCS by pseudo-Bayes & 0.601 & 0.799 & 0.062 & 0.100 \\
CFCS by unlikely under $H_0$ & 0.549 & 0.757 & 0.061 & 0.076 \\
CFCS by likely under $H_1$ & 0.282 & 0.142 & 0.030 & 0.025 \\
CFCS by $H_0$ & 0.398 & 0.155 & 0.044 & 0.017 \\
CFCS by $H_1$ & 0.420 & 0.730 & 0.038 & 0.080 \\
\end{tabular}

\caption{True Shannon information (TSI) values (bits) (section
  \ref{shellsinfointro}) for the different methods on the shells
  problem of section \ref{shells} under various different priors. IFCS
  = inconclusive frequentist confidence set function; CFCS =
  conclusive frequentist confidence set function. TSI for the
  frequentist hypothesis testing methods is the same as for Bayes.
  Priors are: \textbf{A}: $P(c)=0.5\,\forall c$; \textbf{B}:
  $P(c=0)=0.98,P(c=1)=0.02$; \textbf{C}: $P(m|c)=\frac{1}{2}\,\forall
  m,c$; \textbf{D}: $P(m|c)=\left\{\begin{matrix}0.01 & (m=c)\\ 0.99 &
  (m\neq c)\end{matrix}\right.$. In each case the settings for
  frequentist methods that use prior information matched the prior
  under which ASI was assessed.
\label{TSItable}
}

\end{center}
\end{table}

\subsection{Further assessment of other versions of the shells problem}

To see whether the same conclusions hold generally, we also ran 100
different such problems, drawing the values of $P(f|m,c)$ randomly
from a flat Dirichlet distribution on the 3-simplex, and the values of
$P(m,c)$ independently from a flat Dirichlet distribution on the
4-simplex.

\begin{figure}[H]
\begin{center}

\includegraphics[scale=0.7]{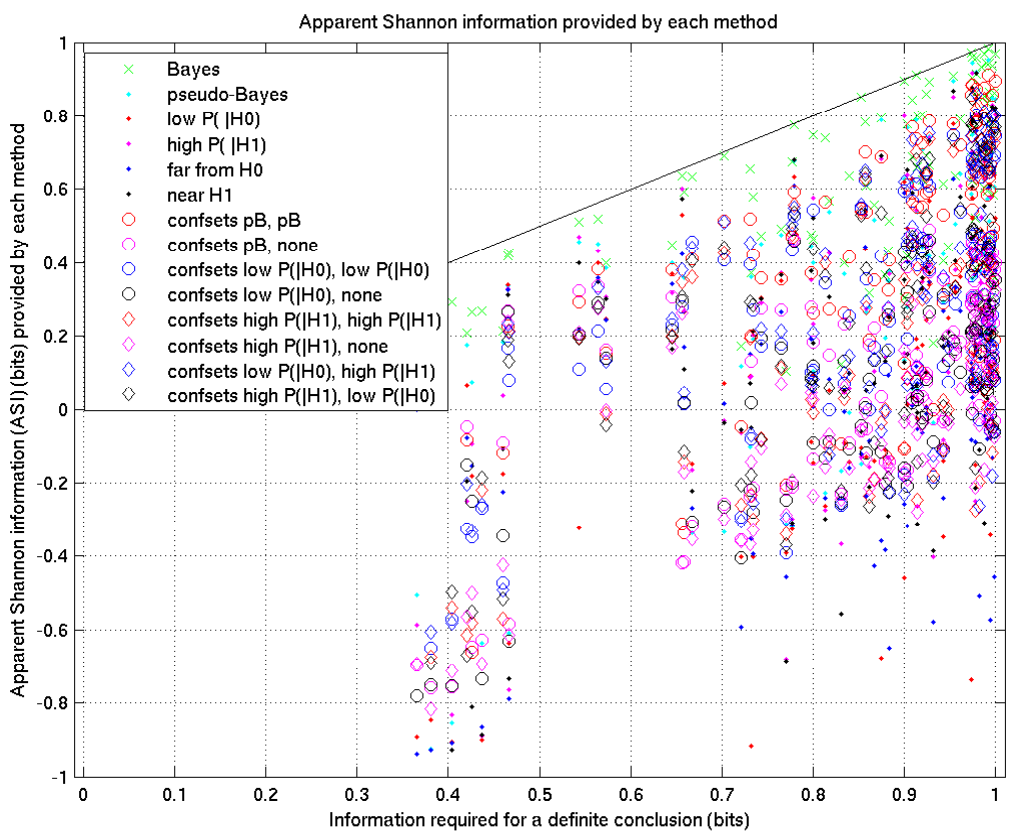}

\caption{The apparent Shannon information (ASI) achieved by each
  method for the Shells problem of section \ref{shells}, plotted
  against the amount of information needed to give a definite answer
  to the User's question on the $x$-axis. Versions of this plot
  without the legend and zoomed in are in figures
  \ref{shellsinfoplot2} and \ref{shellsinfoplot3} below.
\label{shellsinfoplot1}
}

\end{center}
\end{figure}

\begin{figure}[H]
\begin{center}

\includegraphics[scale=0.7]{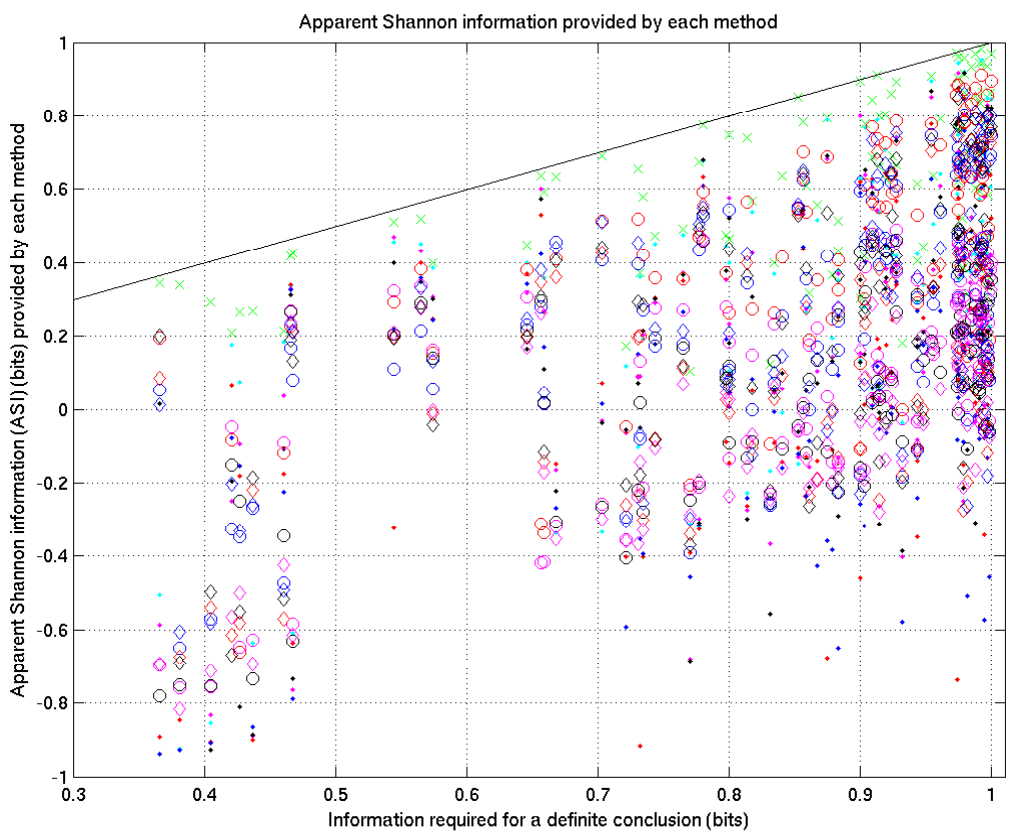}

\caption{A zoomed in version of figure \ref{shellsinfoplot1}. The
  apparent Shannon information (ASI) achieved by each method for the
  Shells problem of section \ref{shells}, plotted against the amount
  of information needed to give a definite answer to the User's
  question on the $x$-axis. The meaning of the various dots is shown
  in the legend to figure \ref{shellsinfoplot1} above.  A still more
  zoomed in version is in figure \ref{shellsinfoplot3} below.
\label{shellsinfoplot2}
}

\end{center}
\end{figure}

\begin{figure}[H]
\begin{center}

\includegraphics[scale=0.7]{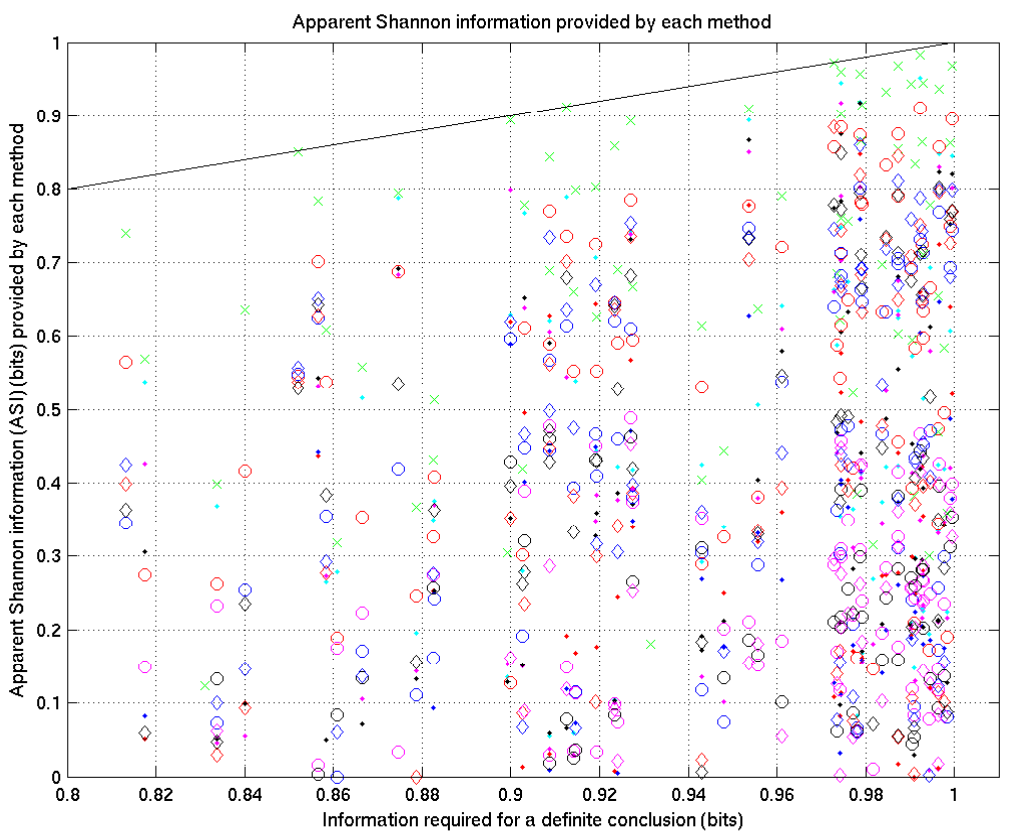}

\caption{A further zoomed in version of figures \ref{shellsinfoplot1}
  and \ref{shellsinfoplot2}. The apparent Shannon information (ASI)
  achieved by each method for the Shells problem of section
  \ref{shells}, plotted against the amount of information needed to
  give a definite answer to the User's question on the $x$-axis. The
  meaning of the various dots is shown in the legend to figure
  \ref{shellsinfoplot1} above. 
\label{shellsinfoplot3}
}

\end{center}
\end{figure}

For each such sample problem, we calculated the ASI resulting from
applying the various methods of section \ref{shells}, which we plot in
figure \ref{shellsinfoplot1} (and again without the legend in
figure \ref{shellsinfoplot2}, and zoomed in to the top right corner
in figure \ref{shellsinfoplot3}). All three plots also show by the
$x$-coordinate and the black diagonal line the amount of TSI and ASI
needed to reach a definite conclusion to the User's problem. In all
cases the green cross marking the Bayesian result comes out above the
column of results of alternative methods, and usually not far below
the black line. In all three plots, many of the frequentist methods
have fallen off the bottom of the plot because of yielding ASI more
negative than -1 bit.

\clearpage

\subsection{Discussion of the shells example}

We have seen that there are a myriad of frequentist methods for
determining the country of origin of the shells in this problem, some
of which even depend on the choice of Bayesian prior that is used in
the calculation. There is no uniformly optimal set of critical regions
for this problem. In contrast there is a single Bayesian method, which
correctly takes prior information into account, and which maximises
both true and apparent Shannon information in the resulting answer,
conserving all the relevant information in the data. All the
frequentist methods shown lose ASI in comparison, while the
95\% frequentist confidence sets also lose TSI. Many of the
frequentist methods are indeed sufficiently misleading that they give
negative ASI, making using the results less good than just going with
the prior and ignoring the data.

Indeed the various frequentist methods give utterly different
results. If we receive the data value $n=(5,10,5)$, for example, the
frequentist methods give the following results in the order of
presentation of the figures: Frequentist confidence that $c=1: 0.999,
0.114, 0.999, 0.066, 0.994, 0.404, 0.994, 0.404$; 95\% confidence
sets: $\{1\}, \{1\}, \{0,1\}, \{1\}, \emptyset, \{0,1\}, \{1\},
\{0\}$. In other words, dependent on purely arbitrary choices of
nested sets of critical regions or confidence set functions, all four
possible frequentist confidence sets may be arrived at with 95\%
frequentist confidence (including that $c$ is in the empty set), while
we may obtain frequentist confidence levels that $c=1$ from 0.066 to
0.999 .

In contrast, depending on the prior in use, the posterior probability
that $c=1$ given by the Bayesian approach with this same data is (for
the four priors considered, in the order of presentation in figure
\ref{shellsBayes}): $0.930, 0.931, 0.213, 0.215\ $. Note that the
prior $P(c)$ for the first two of these values is $0.5$ on each
country, while for the last two it is $0.98$ on $c=0$ and $0.02$ on
$c=1$, a huge difference; it should not be surprising that if we are
almost certain that $c=0$ before collecting the data then we are less
sure that $c=1$ after collecting it than we would be with equal prior
probabilities on the two countries.

The conclusion remains that the Bayesian method retains more ASI than
any other method, and that although it usually doesn't extract
sufficient information to reach a definite answer to the User's
question, it does better than any other possible method.

Even though all the frequentist hypothesis-testing methods retain all
the TSI in their solutions for this particular problem, this is a
distinctly pyrrhic victory; we are just saying that having applied
such a method, and got a poor result with low or even negative ASI, we
needn't give up; rather the best thing to do is to use the method's
output to recover the original data, then apply the Bayesian method.

\clearpage

\section{A more realistic but less dramatic example}
\label{example3}

\subsection{Introduction}

In order to convince readers who are uncomfortable with abstract
problems (such as that presented in section \ref{example1} above) that
frequentist methods are wrong, we now consider a much more realistic
example where we still have a 2-dimensional data space to consider,
and can visualise the posterior probability and frequentist confidence
that results from different observed data. While the results are less
dramatic, they still suffice to make the necessary points in
principle.

We deal with the question of which antibiotic to use to treat a
particular disease. We suppose that we have a set of patients infected
with \textit{Microsoftus gatesii}, and we want to know whether
treatment with jobsucillin or torvaldomycin is better at curing
them. We choose to randomise $N_1+N_2=60$ patients, giving $N_1=40$
patients torvaldomycin and $N_2=20$ patients jobsucillin. Our null
hypothesis $H_0$ is that jobsucillin is better. 

In terms of our standard notation for inference problems from section
\ref{infprobdef} above, we set $$\Theta=\{0,1\},$$ $$\Phi=[0,1]^2$$
(with $\phi_1$, aka $p_1$, the unknown probability of cure with
torvaldomycin and $\phi_2$, aka $p_2$, the unknown probability of cure
with
jobsucillin), $$H=\{(\theta,p_1,p_2)\in\Theta\times\Phi:\theta=[p_1>p_2]\}$$
(where $[\ ]$ denotes the function that takes the value 1 if the
condition inside the brackets is true and 0
otherwise), $$H_0=(\{0\}\times\Phi)\cap
H,$$ $$H_1=(\{1\}\times\Phi)\cap H,$$
$$X=\{0,1,...,N_1\}\times\{0,1,...,N_2\},$$ where $x=(n_1,n_2)$
means that we got $n_1$ cures with torvaldomycin and $n_2$ cures with
jobsucillin, and $$P(n_1,n_2|\theta, p_1, p_2) =
\frac{N_1!N_2!}{n_1!n_2!(N_1-n_1)!(N_2-n_2)!}
p_1^{n_1}(1-p_1)^{N_1-n_1} p_2^{n_2}(1-p_2)^{N_2-n_2}.$$ 

As usual, we consider first the Bayesian answer, then various
frequentist answers. A key point that we want to make is that while
the Bayesian answer depends on the prior -- something that should
indeed influence what we think after collecting the data -- there are
many frequentist answers which differ from one another by making
entirely arbitrary choices, none of which is uniformly optimal.

\subsection{Bayesian solution}
\label{example3Bayes}

Here we set the prior $P(p_1,p_2)=[(p_1,p_2)\in [0,1]^2]$ (the
independent uniform prior) which automatically results in
$P(H_0)=P(H_1)=\frac{1}{2}$. (Obviously we could set different priors,
but for brevity we only explore this one which most people will find
reasonable.)

\begin{figure}[hp]
\begin{center}

\includegraphics[scale=0.5]{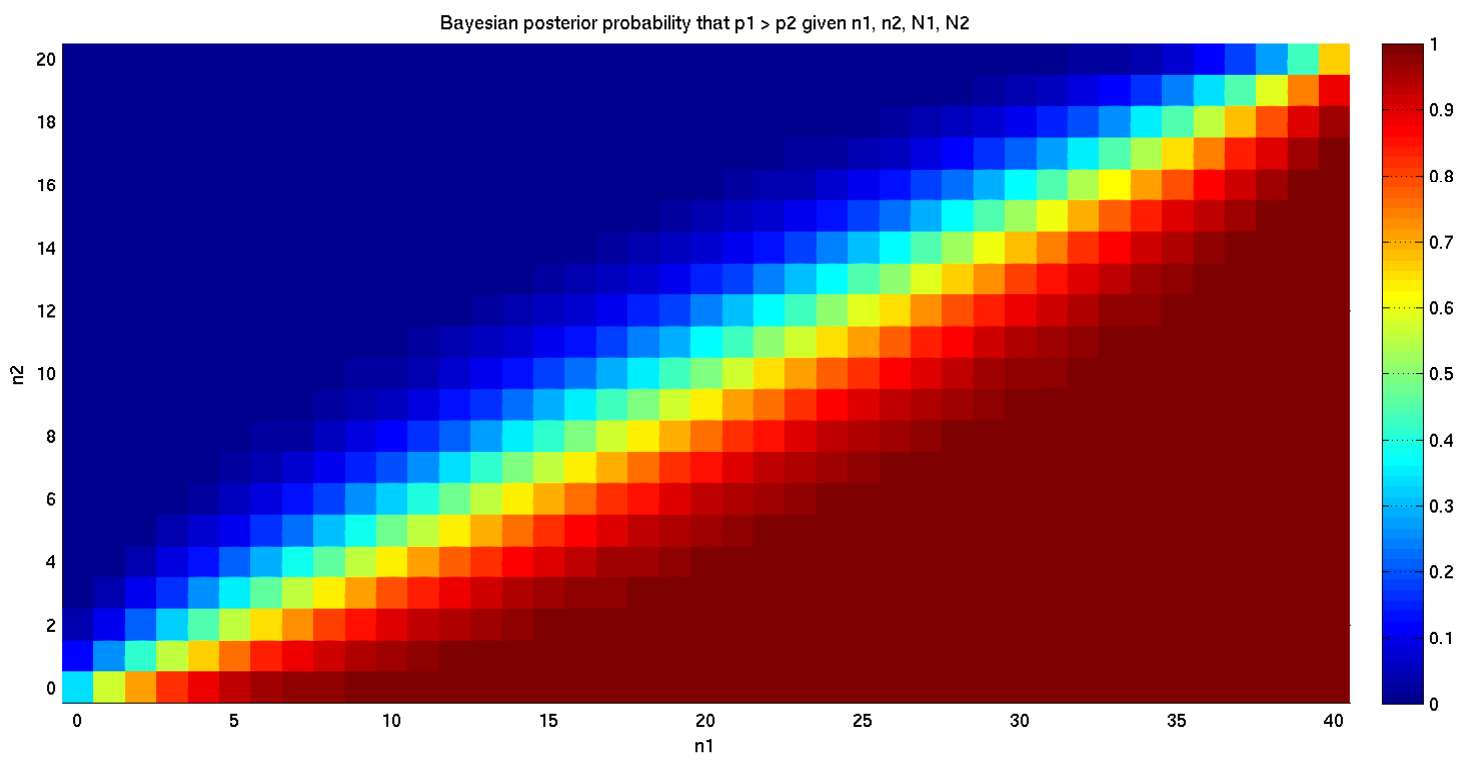}

\caption{Bayesian posterior probability that $H_1$ holds for each
  possible observed value of $(n_1,n_2)$.
\label{bcBayes}
}

\end{center}
\end{figure}

\begin{figure}[hp]
\begin{center}

\includegraphics[scale=0.5]{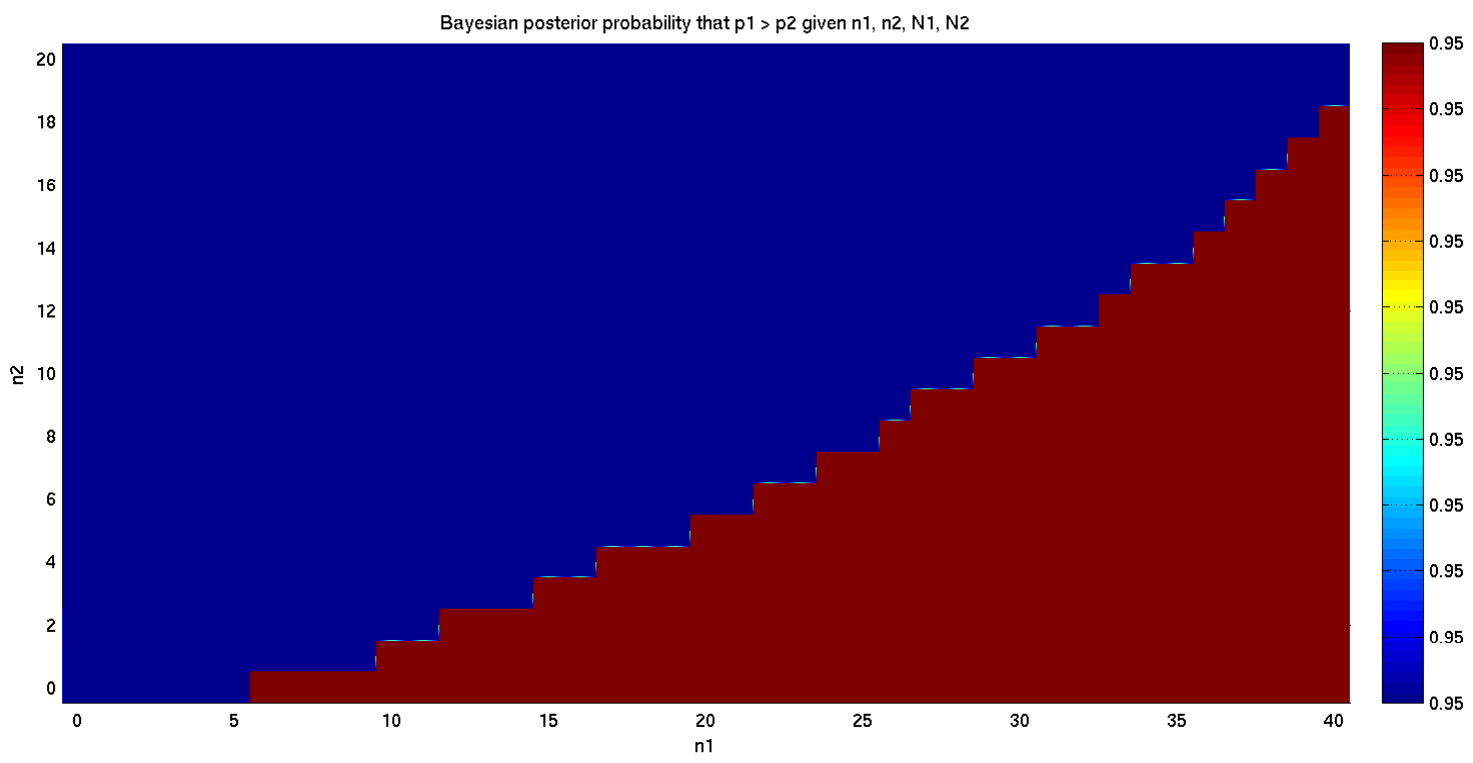}

\caption{Bayesian posterior probability that $H_1$ holds for each
  possible observed value of $(n_1,n_2)$, with each block coloured
  brown or blue according to whether it is above or below 0.95
  respectively.
\label{bcBayes95}
}

\end{center}
\end{figure}

Leaving the reader to do the calculations (or to find them in appendix
\ref{example3Bayescalcs}), the Bayesian solution is shown in figures
\ref{bcBayes} and \ref{bcBayes95}. We observe that it is totally
reasonable, and that for 283 of the possible values of $(n_1,n_2)$
(i.e. of blocks in figure \ref{bcBayes}) the probability that $H_1$
holds is at least $0.95$, and for 283 it is $\leq 0.05$.

\subsection{First exact frequentist solution}
\label{bcfreq1text}

\begin{figure}[hp]
\begin{center}

\includegraphics[scale=0.5]{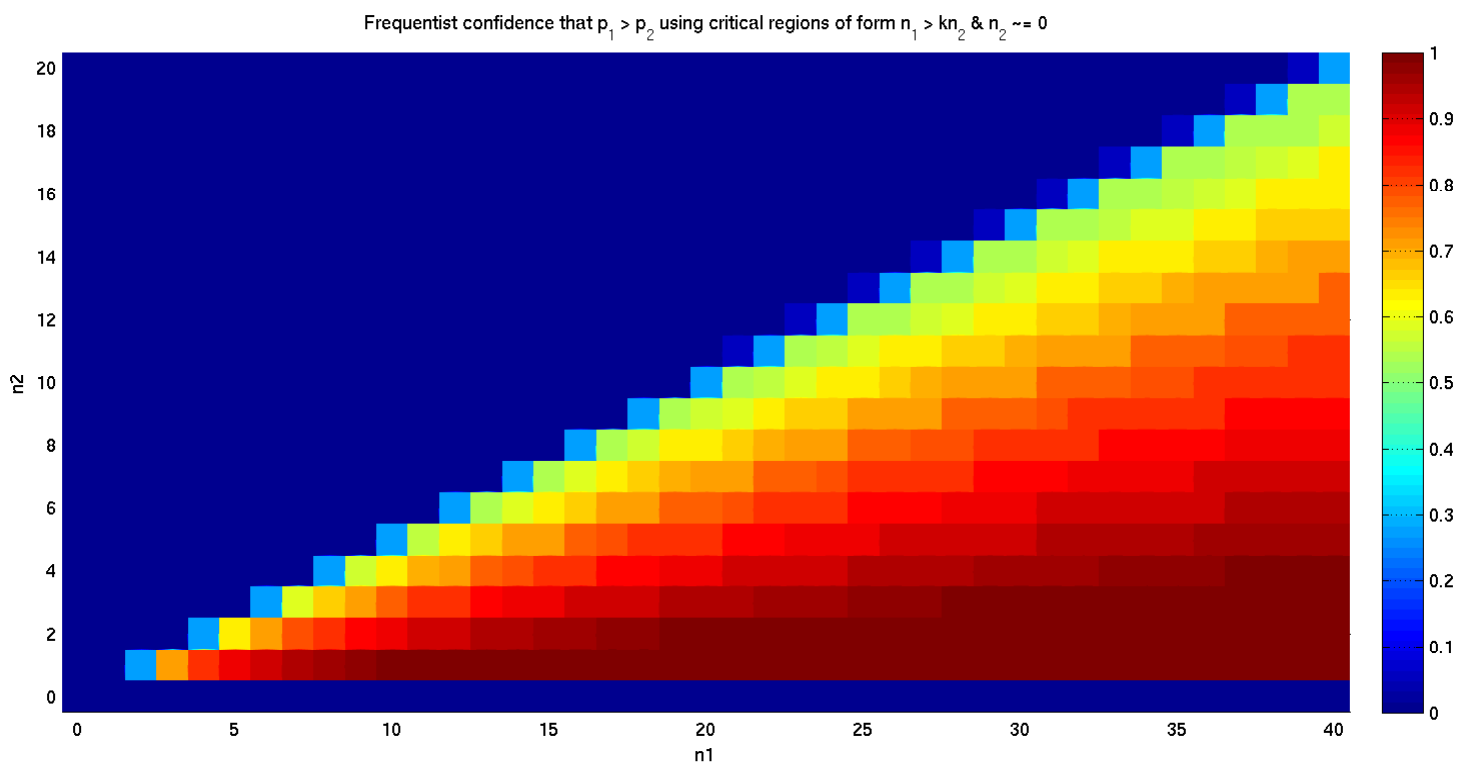}

\caption{First version of exact frequentist confidence that $H_1$
  holds for each possible observed value of $(n_1,n_2)$.
\label{bcfreq1}
}

\end{center}
\end{figure}

\begin{figure}[hp]
\begin{center}

\includegraphics[scale=0.5]{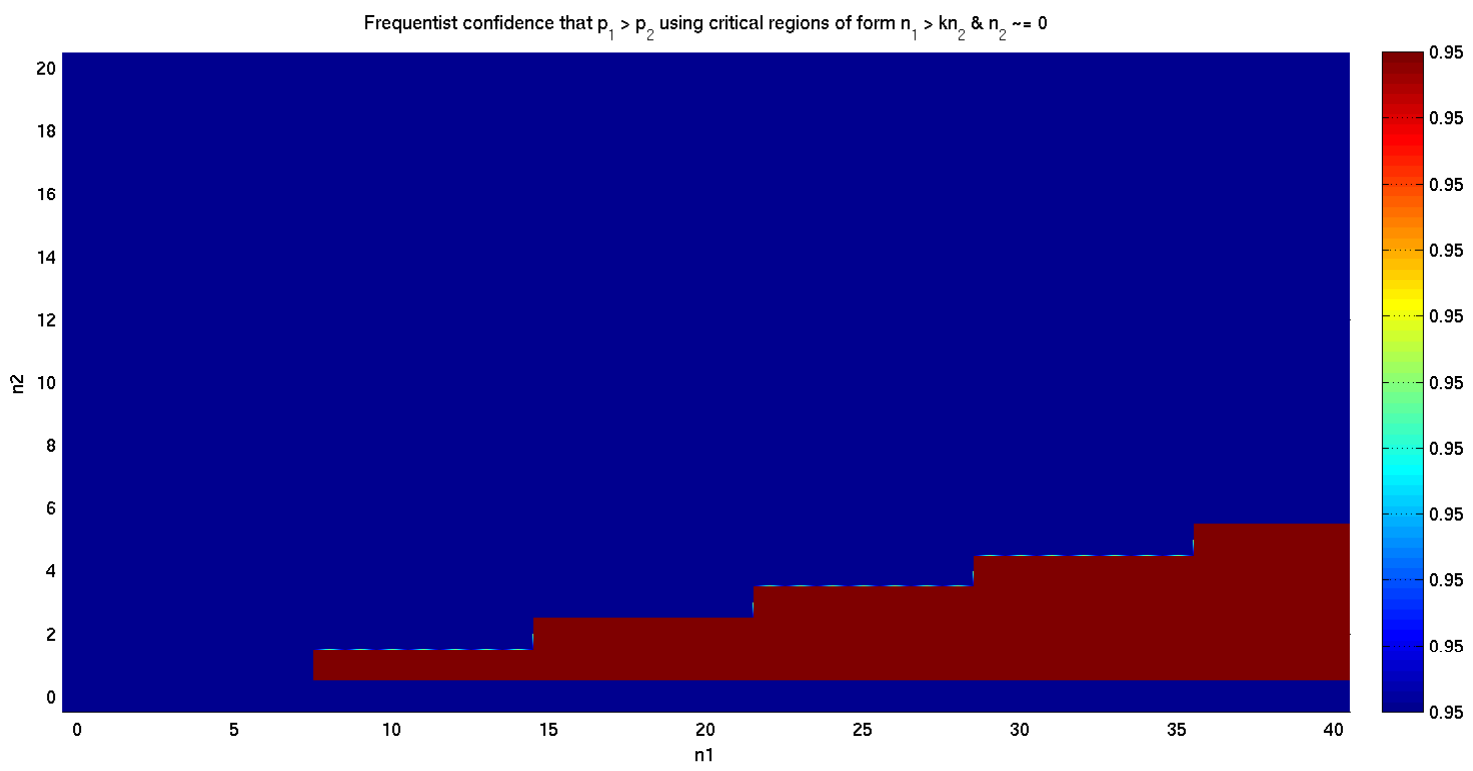}

\caption{First version of exact frequentist confidence that $H_1$
  holds for each possible observed value of $(n_1,n_2)$, with each
  block coloured brown or blue according to whether it is above or
  below 0.95 respectively.
\label{bcfreq195}
}

\end{center}
\end{figure}

Here we use frequentist critical regions of the
form $$C_{\eta(k)}=\{(n_1,n_2):n_2\neq 0, n_1 > kn_2\}$$ for various
values of $k$ and the uniquely determined increasing function $\eta$
to $[0,1]$ that gives the corresponding frequentist confidence. The
result is shown in figues \ref{bcfreq1} and \ref{bcfreq195}. We
observe that it is not really reasonable (e.g. we never conclude that
$H_1$ holds if $n_2=0$, as we have excluded this region of $X$ in
order to get reasonable frequentist confidence elsewhere). There are
95 blocks where frequentist confidence that $H_1$ holds is at least
$0.95$, but 460 out of 861 where it is $\leq 0.05$.

\subsection{Second exact frequentist solution}
\label{bcfreq2text}

\begin{figure}[hp]
\begin{center}

\includegraphics[scale=0.5]{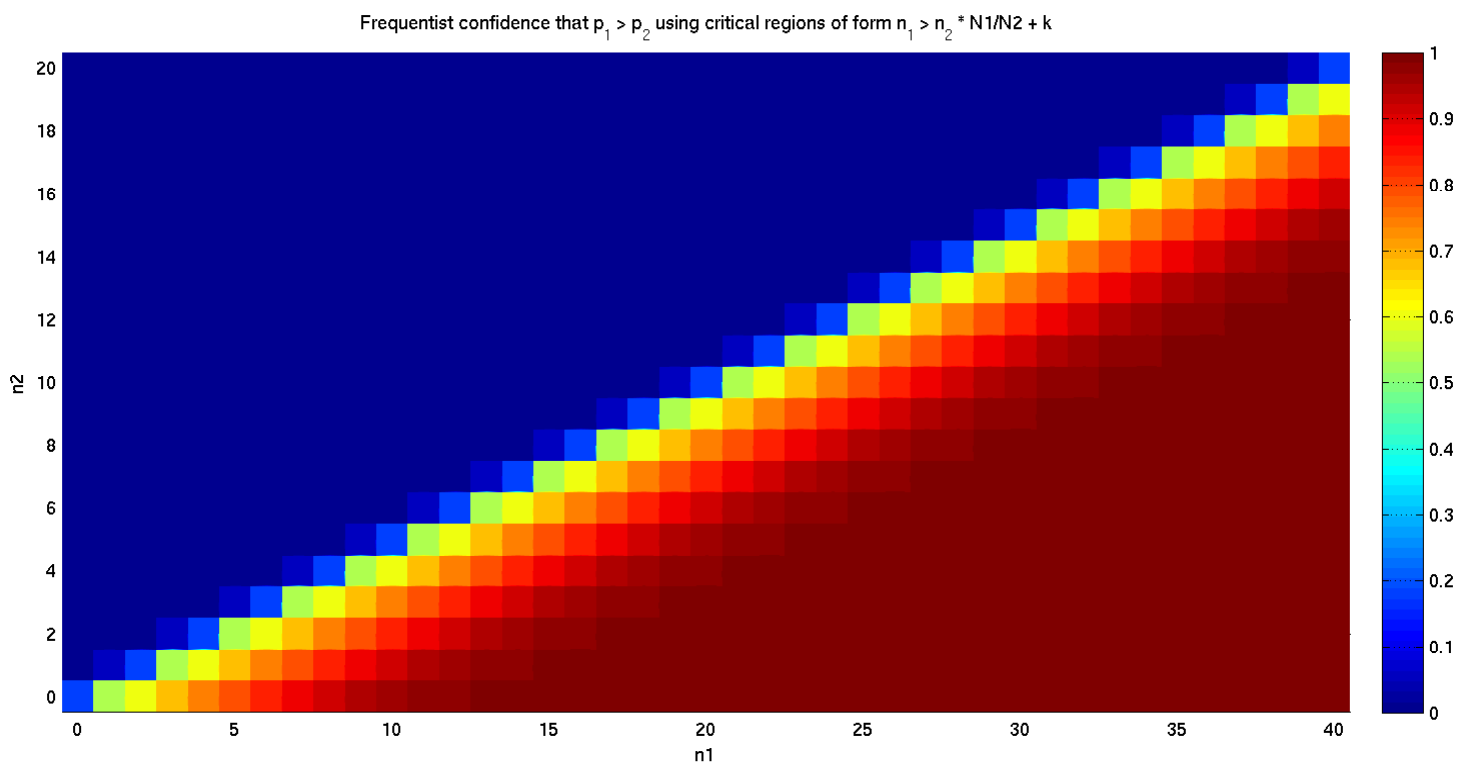}

\caption{Second version of exact frequentist confidence that $H_1$
  holds for each possible observed value of $(n_1,n_2)$.
\label{bcfreq2}
}

\end{center}
\end{figure}

\begin{figure}[hp]
\begin{center}

\includegraphics[scale=0.5]{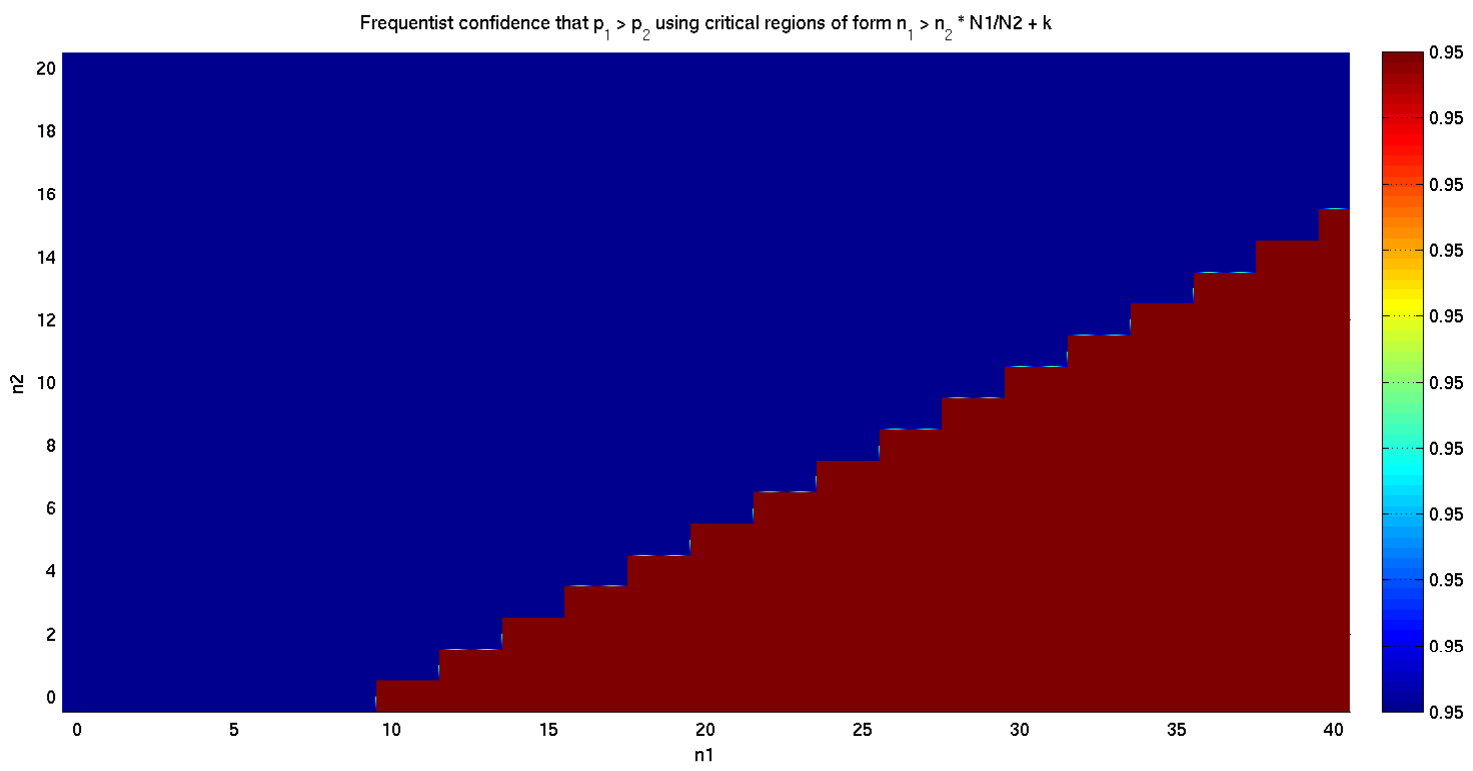}

\caption{Second version of exact frequentist confidence that $H_1$
  holds for each possible observed value of $(n_1,n_2)$, with each
  block coloured brown or blue according to whether it is above or
  below 0.95 respectively.
\label{bcfreq295}
}

\end{center}
\end{figure}

Here we use frequentist critical regions of the
form $$C_{\eta(k)}=\{(n_1,n_2):n_1 > \frac{n_2N_1}{N_2}+k\}$$ for
various values of $k$ and the uniquely determined increasing function
$\eta$ to $[0,1]$ that gives the corresponding frequentist
confidence. The result is shown in figues \ref{bcfreq2} and
\ref{bcfreq295}. We observe that it is not really reasonable (e.g. it
is almost zero on almost half the blocks). There are 256 blocks where
frequentist confidence that $H_1$ holds is at least $0.95$, but 420
out of 861 where it is $\leq 0.05$.

\subsection{$\chi^2$ test, an approximate frequentist solution, and
  the corresponding \mbox{exact} frequentist confidence}
\label{bcchitext}

\begin{figure}[hp]
\begin{center}

\includegraphics[scale=0.5]{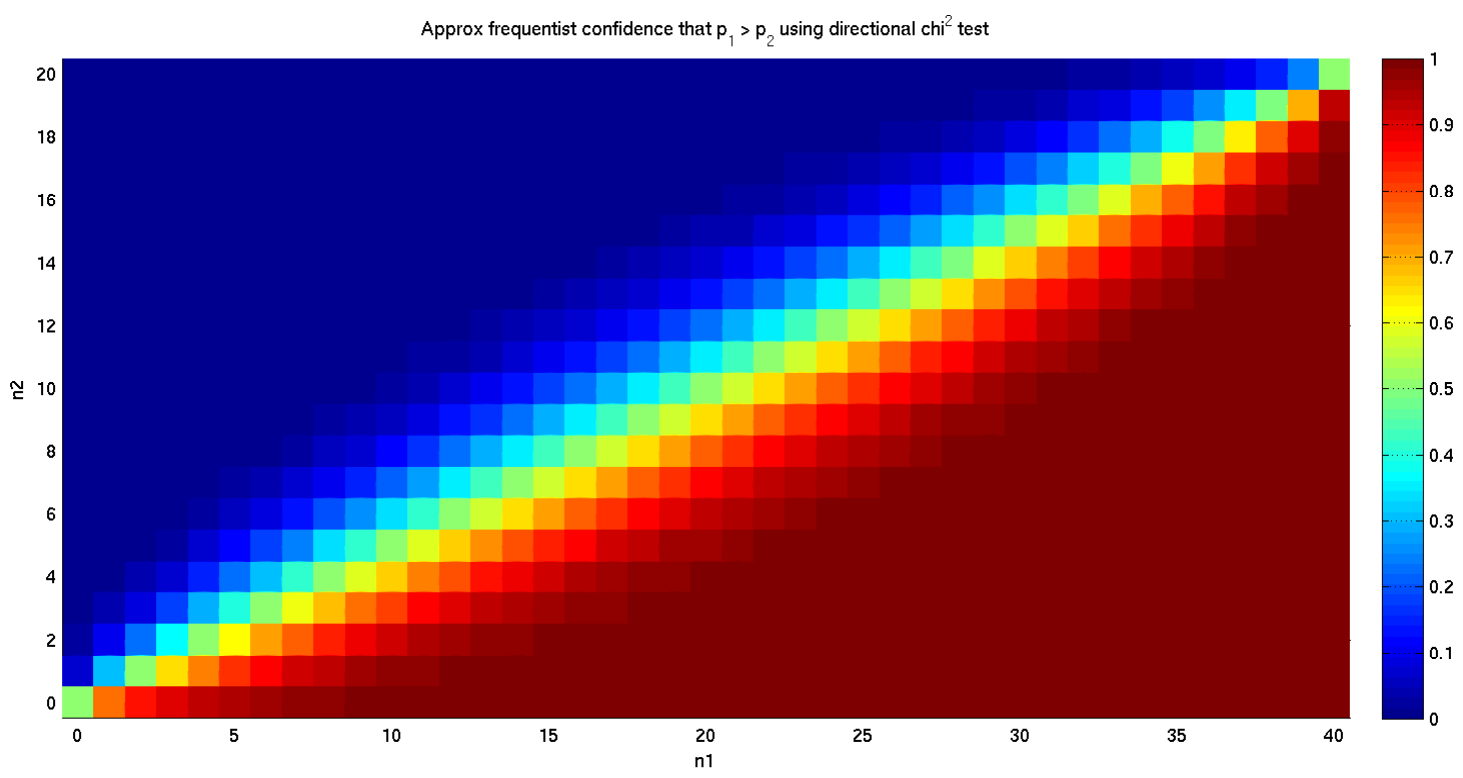}

\caption{Directional $\chi^2$-test version of frequentist confidence
  that $H_1$ holds for each possible observed value of
  $(n_1,n_2)$. Note that these values are an overoptimistic
  approximation (biased in favour of $H_1$). 
\label{bcchiapprox}
}

\end{center}
\end{figure}

\begin{figure}[hp]
\begin{center}

\includegraphics[scale=0.5]{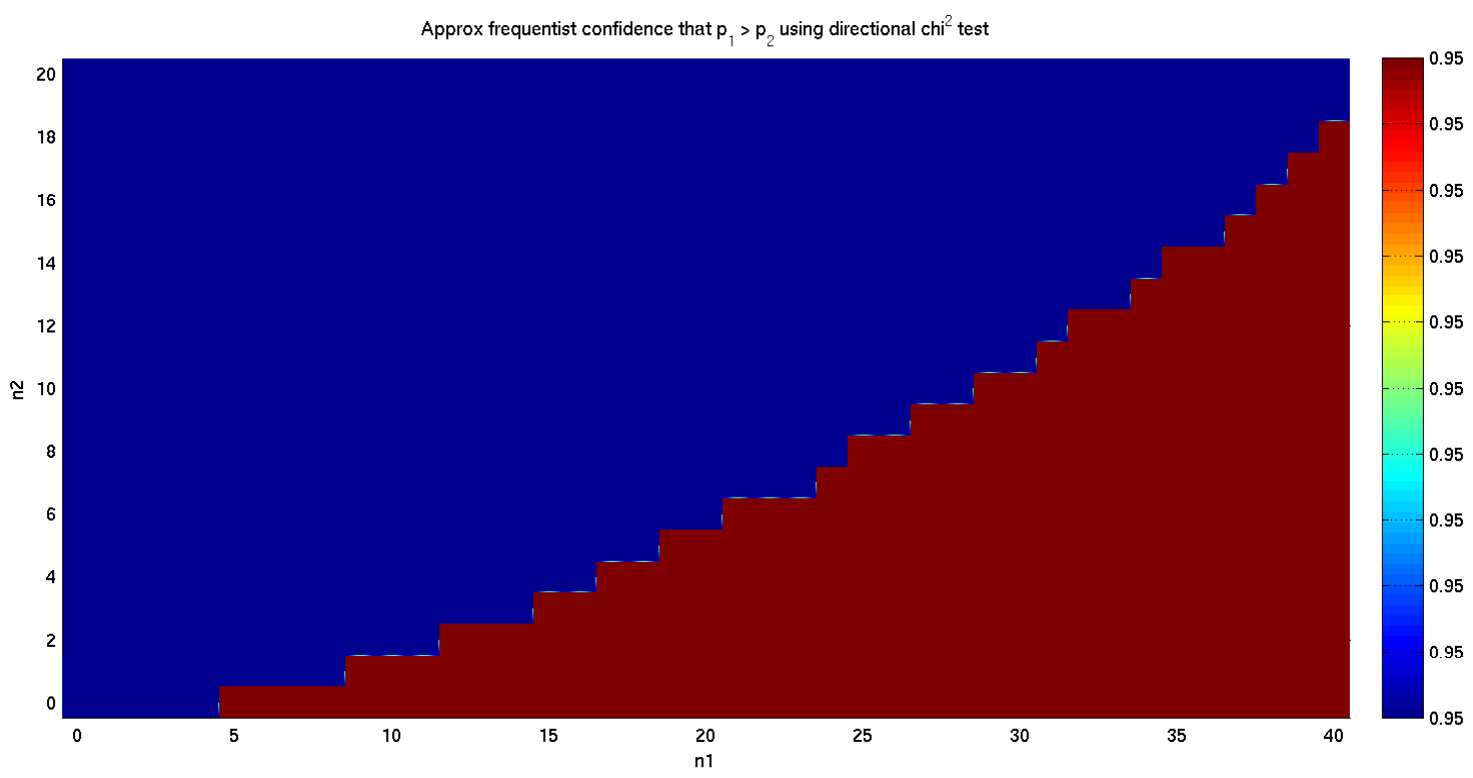}

\caption{Directional $\chi^2$-test version of frequentist confidence
  that $H_1$ holds for each possible observed value of $(n_1,n_2)$,
  with each block coloured brown or blue according to whether it is
  above or below 0.95 respectively. Note that these values are an
  overoptimistic approximation (biased in favour of $H_1$).
\label{bcchiapprox95}
}

\end{center}
\end{figure}

\begin{figure}[hp]
\begin{center}

\includegraphics[scale=0.5]{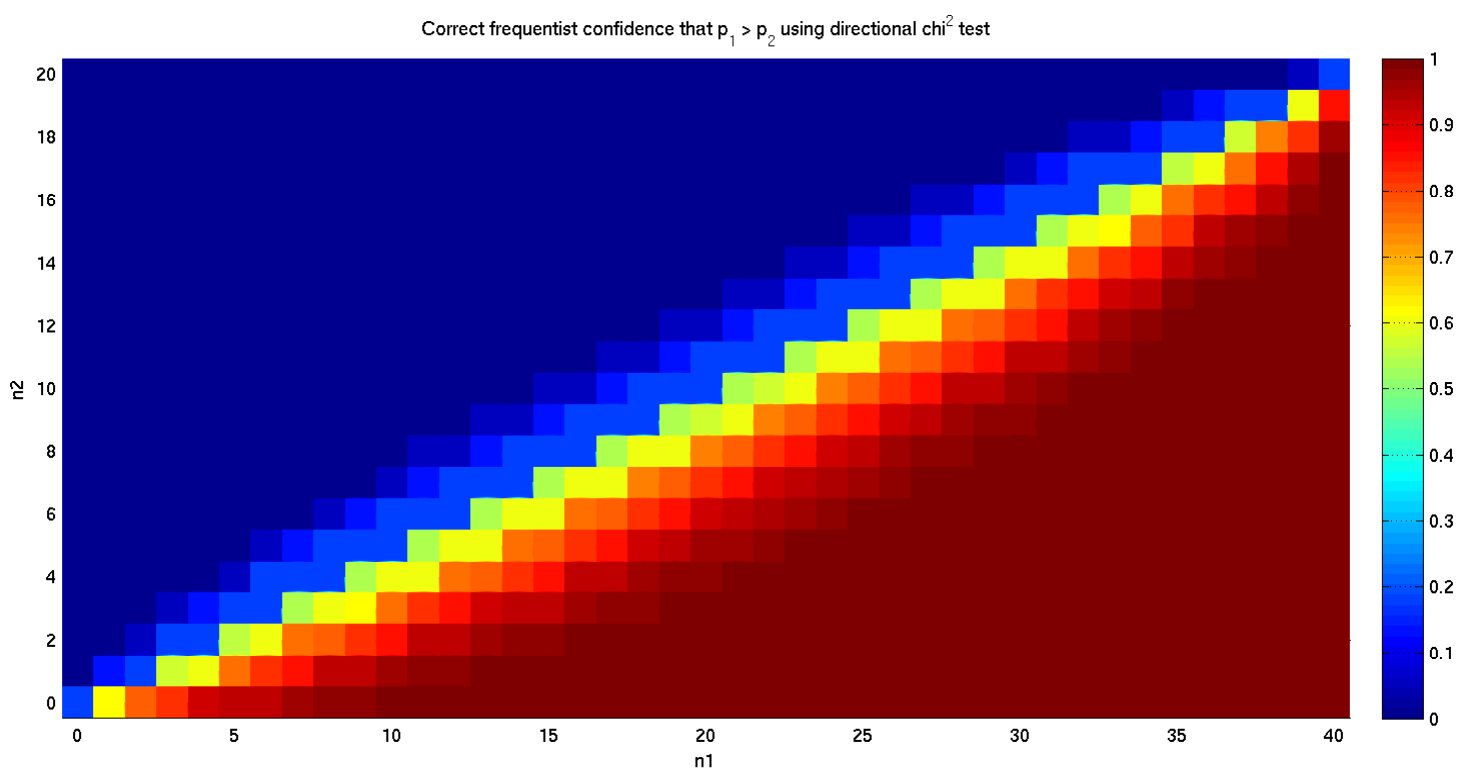}

\caption{Directional $\chi^2$-test version of frequentist confidence
  that $H_1$ holds for each possible observed value of
  $(n_1,n_2)$. These values have been corrected for the errors
  inherent in the $\chi^2$-test.
\label{bcchiexact}
}

\end{center}
\end{figure}

\begin{figure}[hp]
\begin{center}

\includegraphics[scale=0.5]{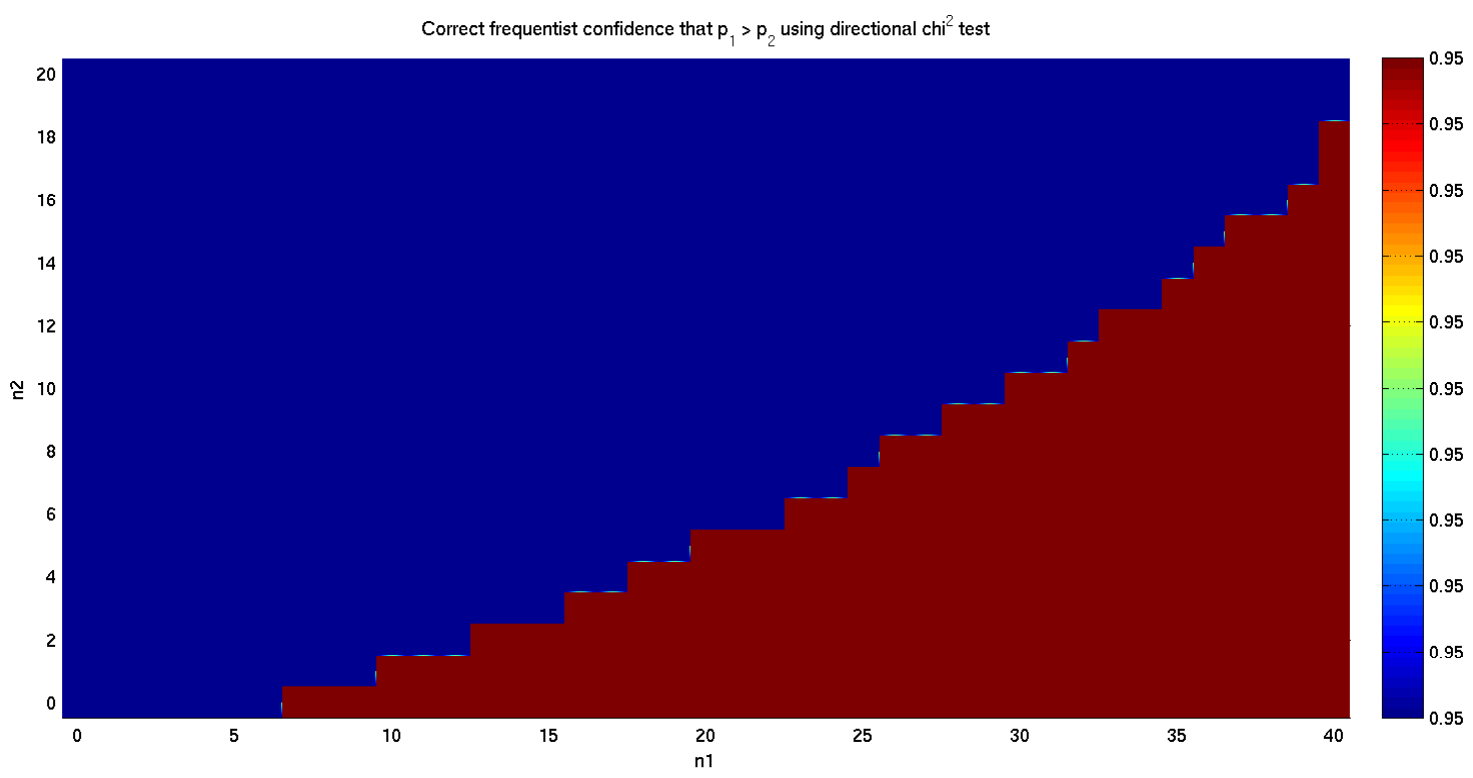}

\caption{Directional $\chi^2$-test version of frequentist confidence
  that $H_1$ holds for each possible observed value of
  $(n_1,n_2)$. These values have been corrected for the errors
  inherent in the $\chi^2$-test. Each block is coloured brown or blue
  according to whether it is above or below 0.95 respectively.
\label{bcchiexact95}
}

\end{center}
\end{figure}

The $\chi^2$ test, here with 1 degree of freedom, is a standard
frequentist test for this problem, which approaches being correct
(from a frequentist point of view) as $N_1$ and $N_2$ approach
infinity (something that in real life never happens). Here we use a
directional version of it (as the standard version only tests whether
$p_1=p_2$ as $H_0$). It returns an over-optimistic estimate of the
frequentist confidence that is biased in favour of $H_1$, as shown in
figures \ref{bcchiapprox} and \ref{bcchiapprox95}; for comparison the
exact frequentist confidence corresponding to the relevant critical
regions is shown in figures \ref{bcchiexact} and \ref{bcchiexact95}.

In the approximation there are 290 blocks with frequentist confidence
above 0.95 and 290 below 0.05; in the exact version there are
respectively 271 and 371. 

\subsection{Critical regions based on Fisher's exact test}
\label{bcFishertext}

``Fisher's exact test'' is a frequentist test for a slightly different
situation, namely that not only is it known in advance how many are in
each treatment group, but also it is known in advance how many get
cured in total. Equivalently, and classically, the situation we have
been addressing so far is equivalent to a lady tasting tea who has
been told that the 60 cups of tea in front of her are some of type A
and some of type B, who is asked to identify which is which; Fisher's
exact test is applicable if she is \textit{also} told in advance how
many there are of type A and how many of type B.

\begin{figure}[p]
\begin{center}

\includegraphics[scale=0.5]{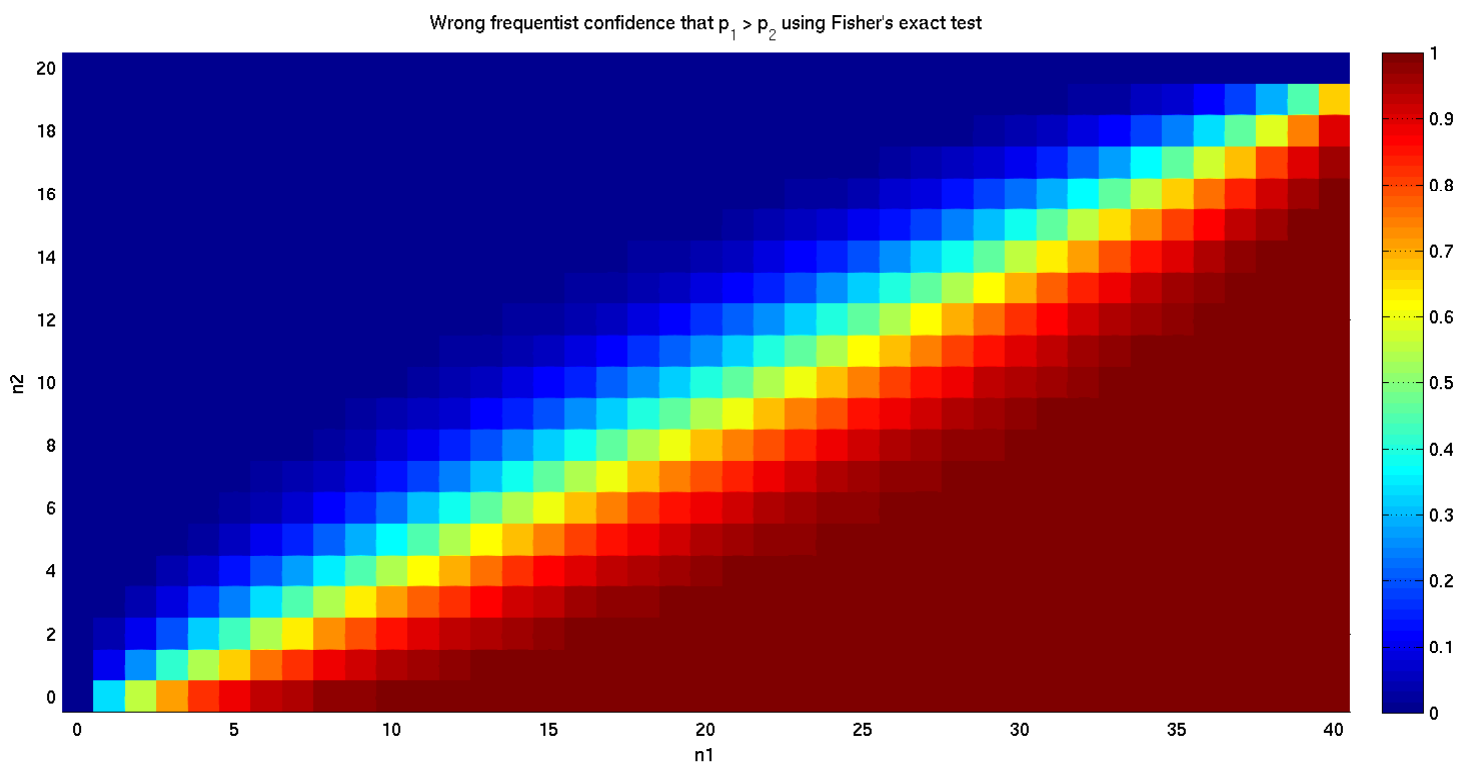}

\caption{Fisher's exact test version of frequentist confidence that
  $H_1$ holds for each possible observed value of $(n_1,n_2)$, as the
  test is often inappropriately used. Note that these values are too
  pessimistic (biased against $H_1$).
\label{bcFisherapprox}
}

\end{center}
\end{figure}

\begin{figure}[p]
\begin{center}

\includegraphics[scale=0.5]{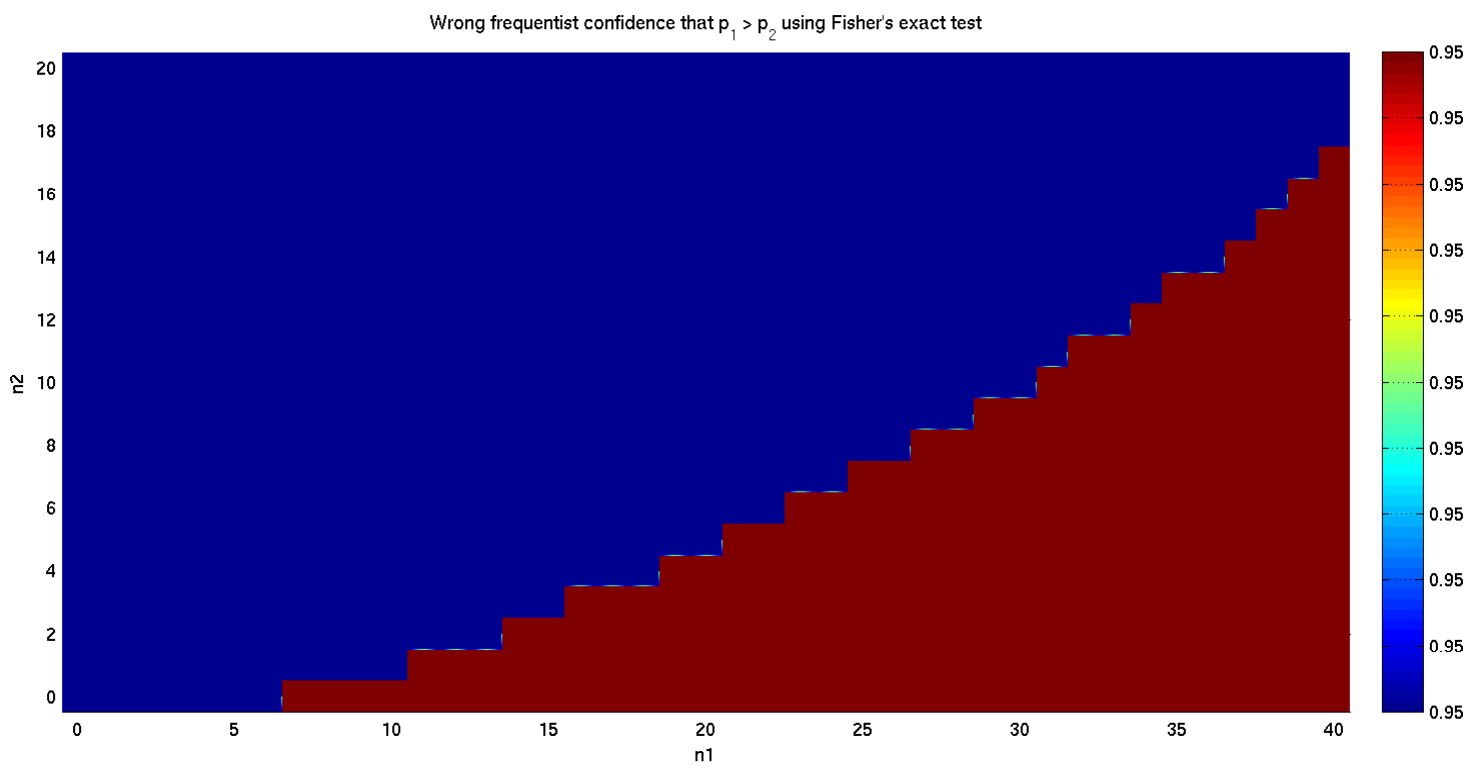}

\caption{Fisher's exact test version of frequentist confidence that
  $H_1$ holds for each possible observed value of $(n_1,n_2)$ as the
  test is often inappropriately used. Note that these values are too
  pessimistic (biased against $H_1$). Each block is coloured brown or blue
  according to whether it is above or below 0.95 respectively.
\label{bcFisherapprox95}
}

\end{center}
\end{figure}

\begin{figure}[p]
\begin{center}

\includegraphics[scale=0.5]{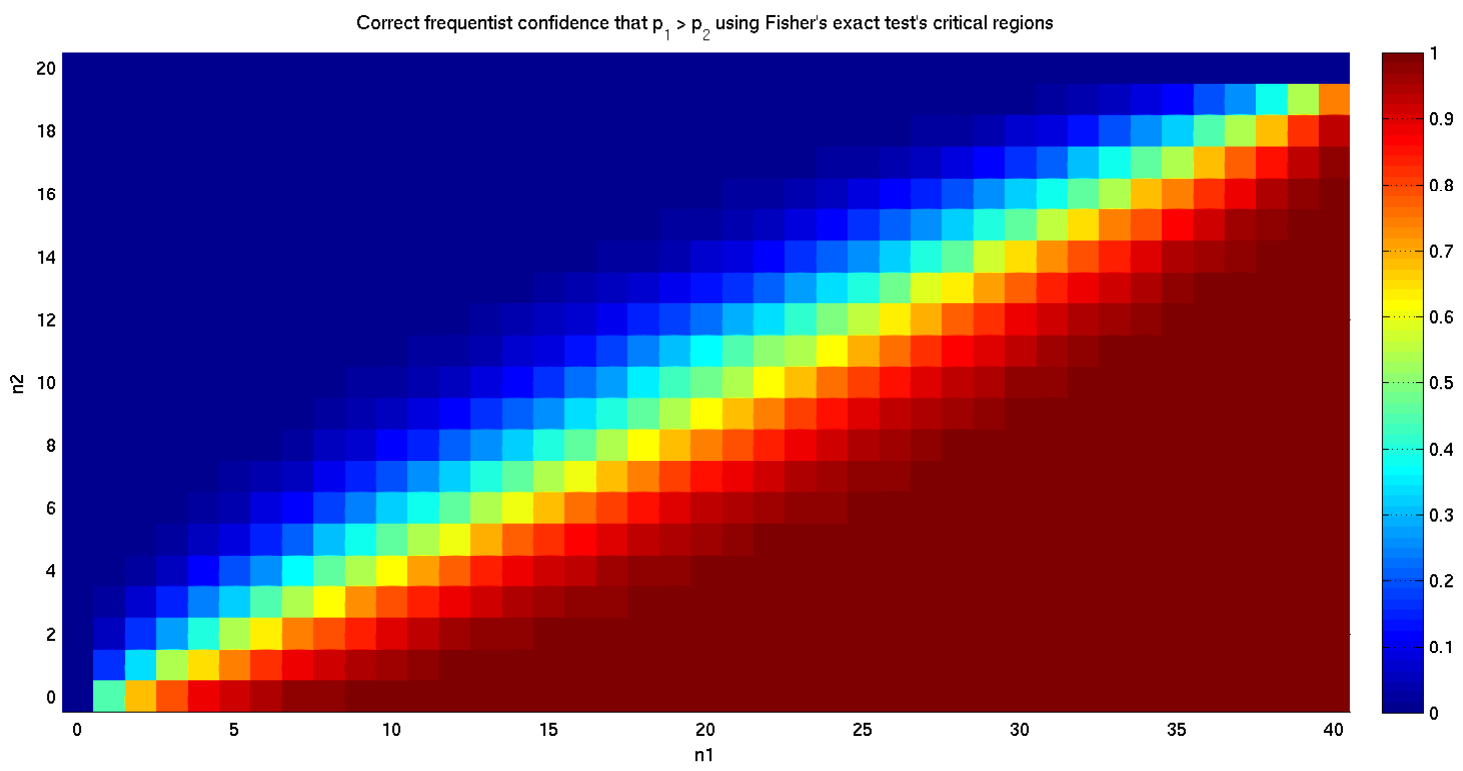}

\caption{Fisher's exact test version of frequentist confidence that
  $H_1$ holds for each possible observed value of $(n_1,n_2)$. These
  values have been corrected for the errors resulting from
  inappropriate use, so are correct for the present problem, and are
  identical with those resulting from using Boschloo's test. 
\label{bcFishercorrect}
}

\end{center}
\end{figure}

\begin{figure}[p]
\begin{center}

\includegraphics[scale=0.5]{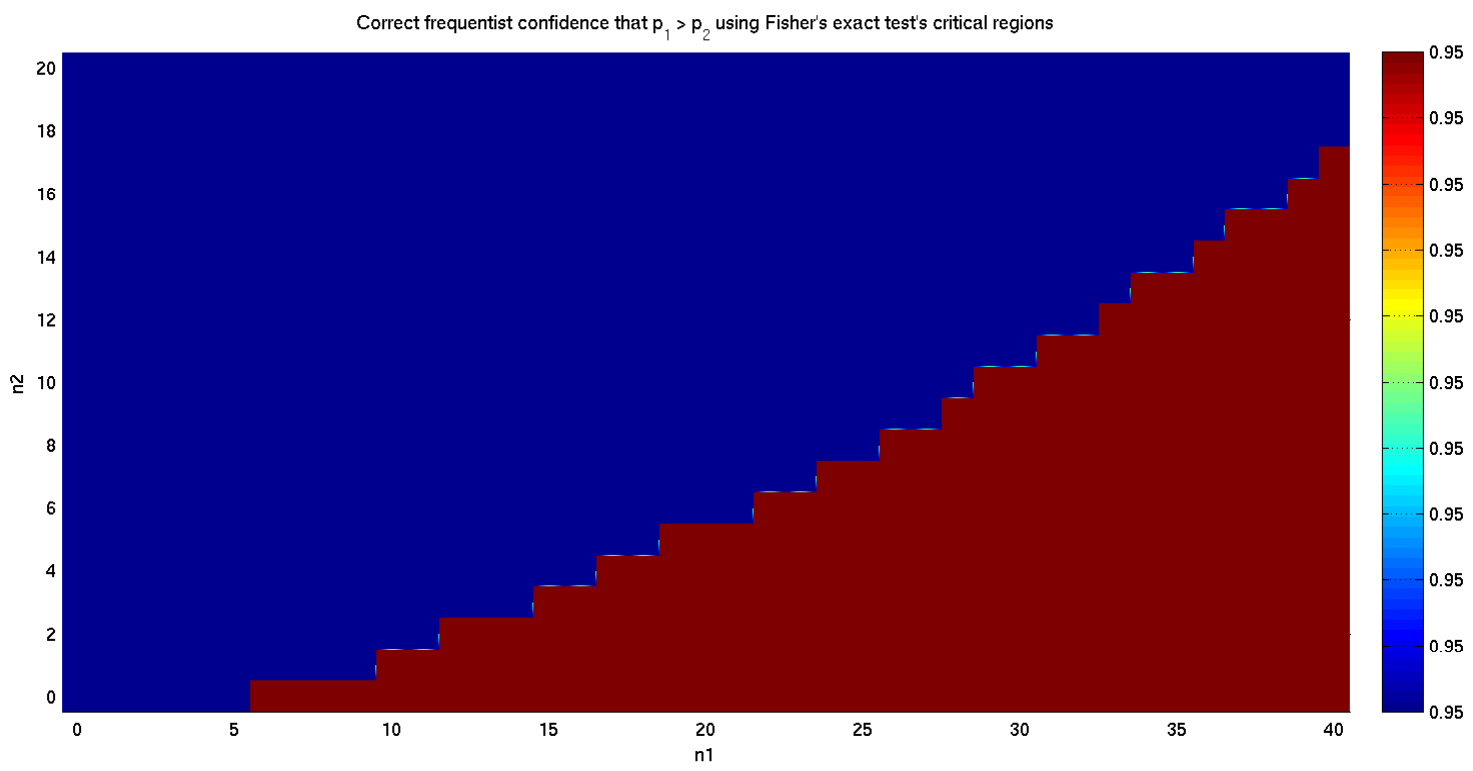}

\caption{Fisher's exact test version of frequentist confidence that
  $H_1$ holds for each possible observed value of $(n_1,n_2)$. These
  values have been corrected for the errors resulting from
  inappropriate use, so are correct for the present problem, and are
  identical with those resulting from using Boschloo's test. Each
  block is coloured brown or blue according to whether it is above or
  below 0.95 respectively.
\label{bcFishercorrect95}
}

\end{center}
\end{figure}

So application of Fisher's exact test to our present situation is
incorrect -- but nonetheless we have often seen it used in this
situation. Moreover, it does define critical regions, and we can
calculate the correct frequentist confidence values for those critical
regions. We present both the naive, incorrect, version of the
frequentist confidence and the corrected version, in figures
\ref{bcFisherapprox} to \ref{bcFishercorrect95}.

In the naive, incorrect, version there are 260 blocks with frequentist
confidence above 0.95 and 321 below 0.05; in the corrected version
there are 280 above 0.95 and 300 below 0.05; the naive version is
therefore too pessimistic (biased against $H_1$).

We also note that Boschloo first proposed using Fisher's exact test
for this present problem with this correction.

\clearpage

\subsection{Effect of choosing which test to use \textit{after}
  collecting the data}
\label{bcchooselatetext}

The experimenter, naturally, is often interested in showing a
``statistically significant'' result, i.e. one with frequentist
confidence $\geq 0.95$. In this situation it is very tempting to not
follow a cardinal rule of frequentist analysis, namely that one must
choose the nested family of critical regions to be used (or
equivalently the frequentist ``test'' to be used) \textit{before}
collecting the data -- and instead to collect the data, then choose
whichever test makes the result have highest frequentist confidence
that $H_1$ holds.

\begin{figure}[p]
\begin{center}

\includegraphics[scale=0.5]{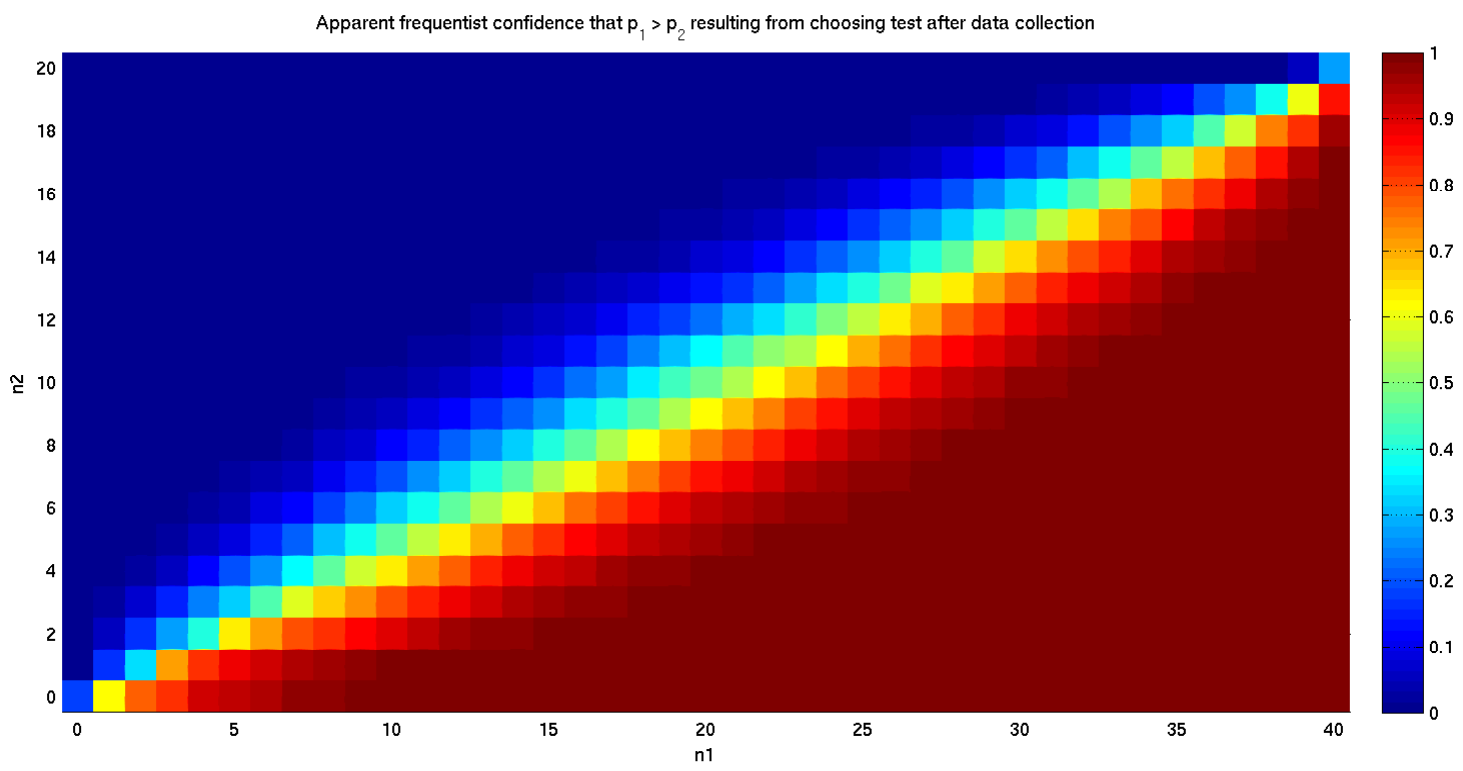}

\caption{The apparent frequentist confidence that $H_1$ holds that
  results from choosing the ``test'' in the light of the data, for
  each possible observed value of $(n_1,n_2)$. These values are
  over-optimistic and biased in favour of $H_1$.
\label{bcchooselate}
}

\end{center}
\end{figure}

\begin{figure}[p]
\begin{center}

\includegraphics[scale=0.5]{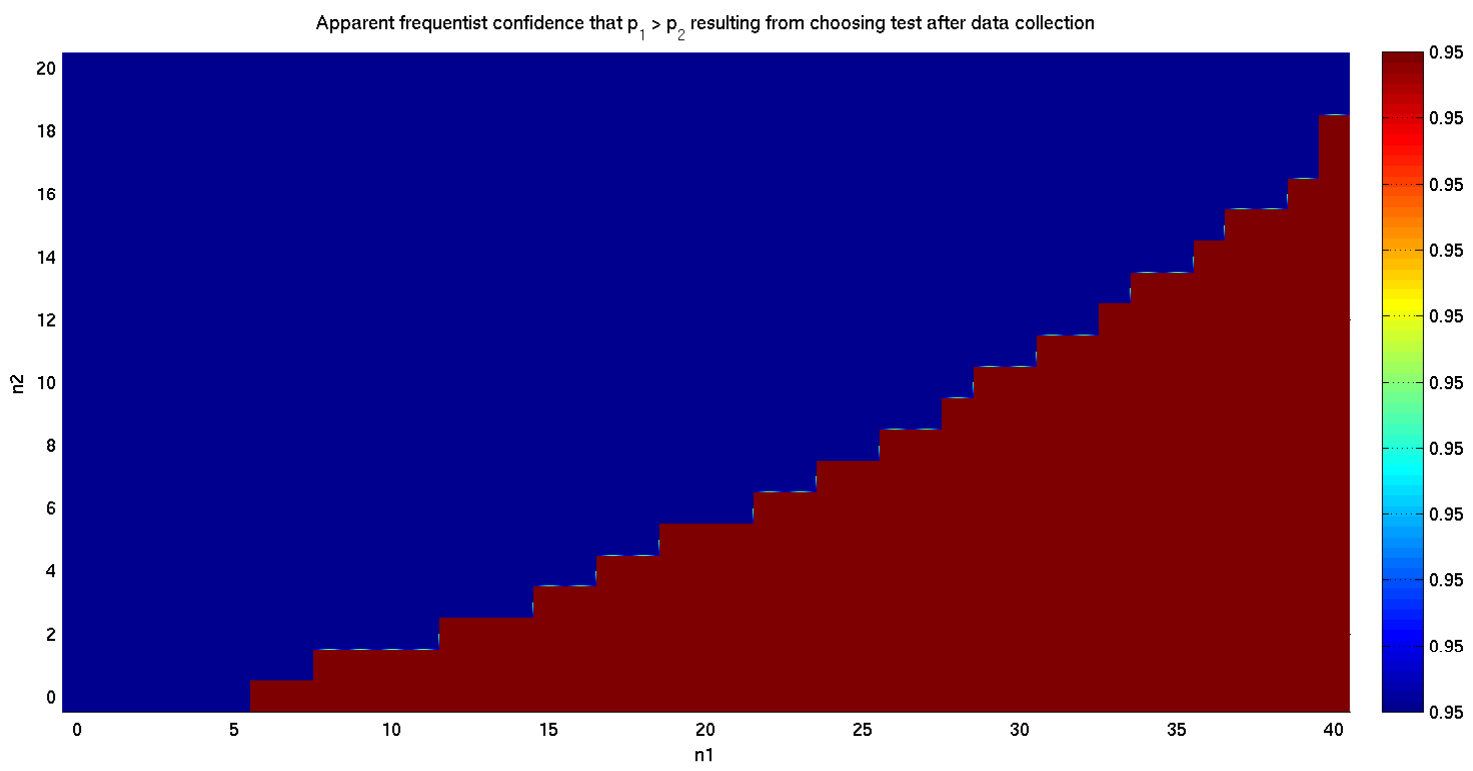}

\caption{The apparent frequentist confidence that $H_1$ holds that
  results from choosing the ``test'' in the light of the data, for
  each possible observed value of $(n_1,n_2)$. These values are
  over-optimistic and biased in favour of $H_1$. Each block is
  coloured brown or blue according to whether it is above or below
  0.95 respectively.
\label{bcchooselate95}
}

\end{center}
\end{figure}

Effectively this results (if we choose from the tests listed in the
subsections above, restricting ourselves to those with correct(ed)
values of frequentist confidence) in the apparent frequentist
confidence levels shown in figures \ref{bcchooselate} and
\ref{bcchooselate95}, in which there are 283 blocks with frequentist
confidence above 0.95 and 297 below 0.05 . In particular, suppose that
we had observed $(n_1,n_2)=(3,1)$: then we would get apparent
frequentist confidence of $0.704$.

\begin{figure}[p]
\begin{center}

\includegraphics[scale=0.5]{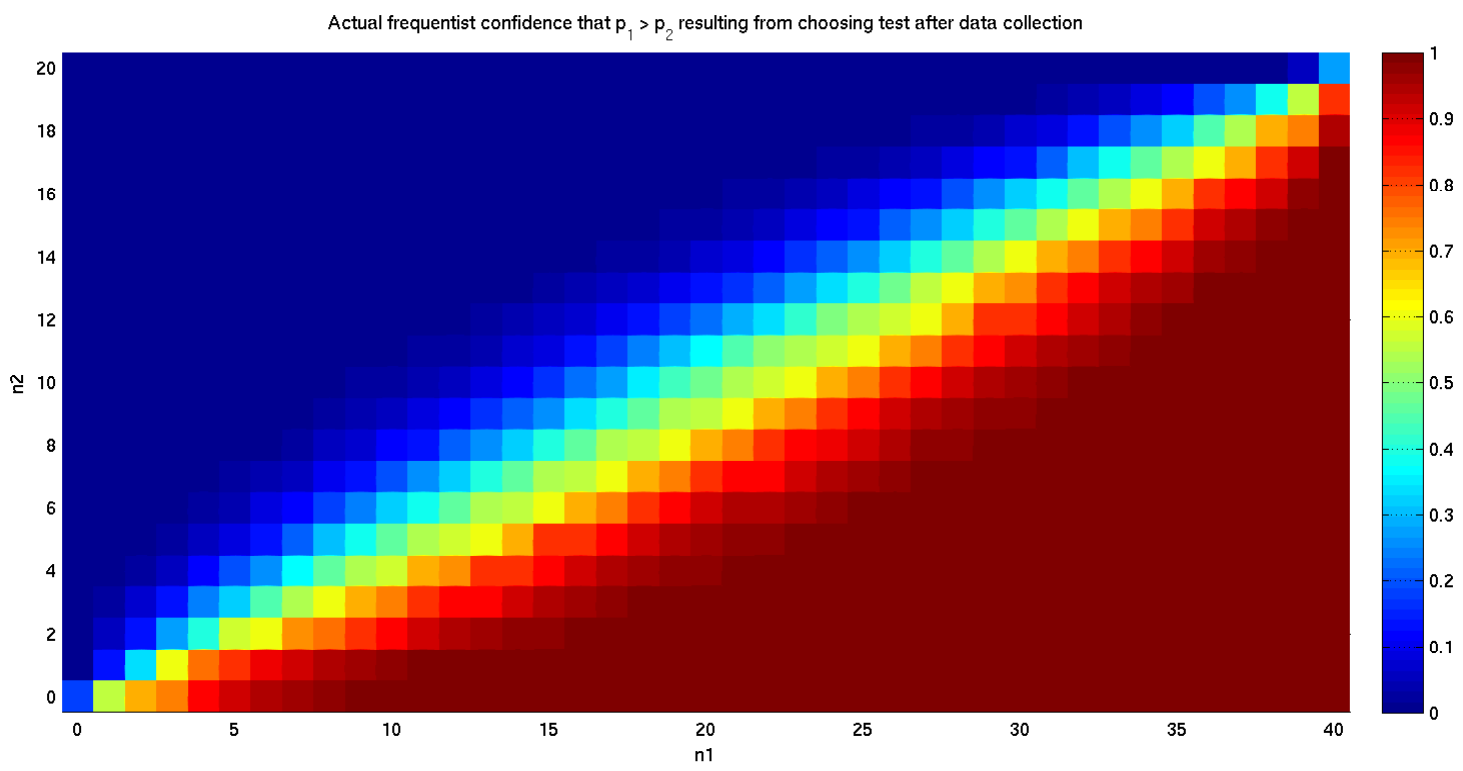}

\caption{The actual frequentist confidence that $H_1$ holds that
  results from choosing the ``test'' in the light of the data, for
  each possible observed value of $(n_1,n_2)$. These values have been
  corrected for the late choice of test. 
\label{bcchooselatefixed}
}

\end{center}
\end{figure}

\begin{figure}[p]
\begin{center}

\includegraphics[scale=0.5]{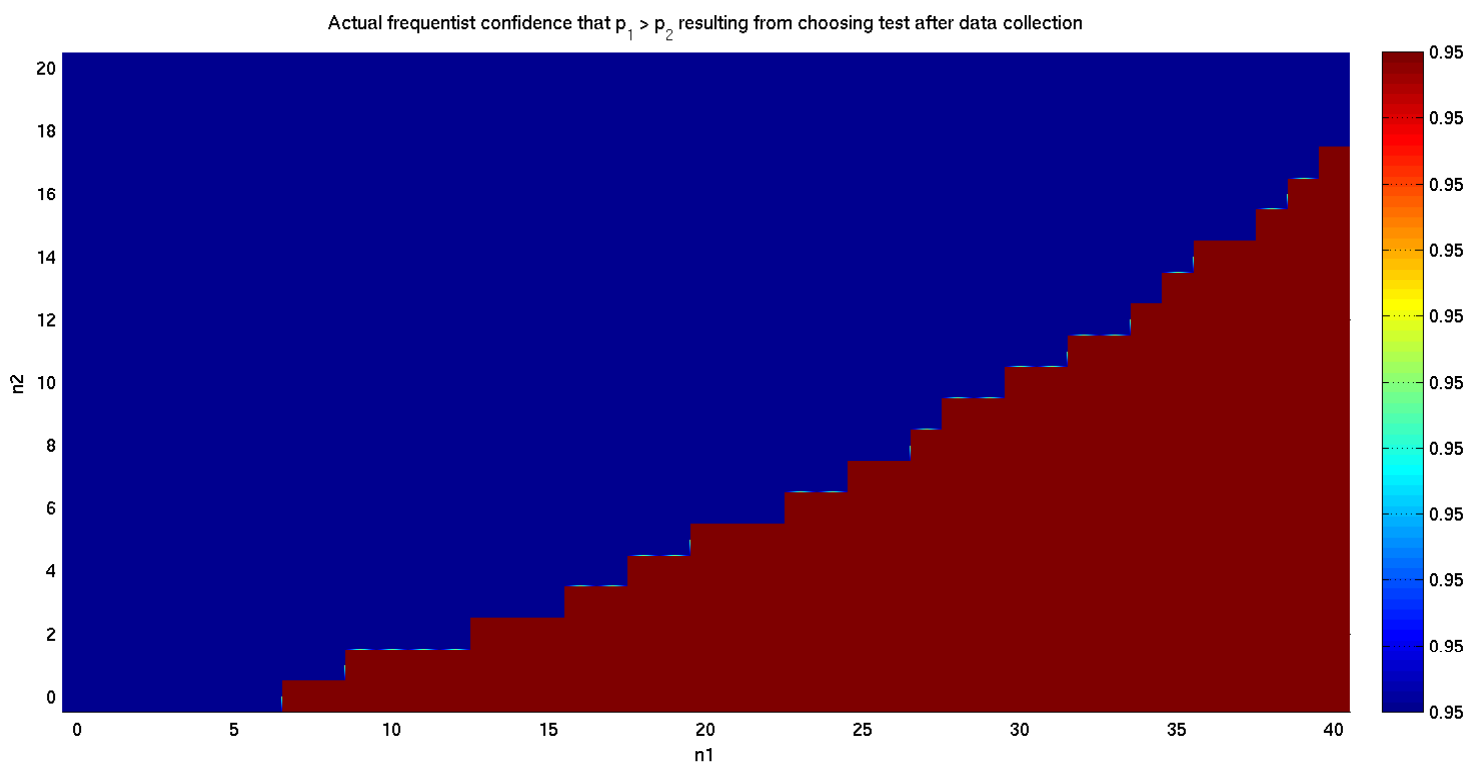}

\caption{The actual frequentist confidence that $H_1$ holds that
  results from choosing the ``test'' in the light of the data, for
  each possible observed value of $(n_1,n_2)$. These values have been
  corrected for the late choice of test. Each block is coloured brown
  or blue according to whether it is above or below 0.95 respectively.
\label{bcchooselatefixed95}
}

\end{center}
\end{figure}

However, if we determine the critical regions that are then in effect,
and calculate their true frequentist confidence, we get the levels
shown in figures \ref{bcchooselatefixed} and
\ref{bcchooselatefixed95}, in which there are 268 blocks with
frequentist confidence above 0.95 and 298 below 0.05 . Data of
$(n_1,n_2)=(3,1)$ now only gives us frequentist confidence of
$0.608$. 

Thus choosing a frequentist ``test'' in the light of the data, without
correcting for this, is cheating. The fact that it is often so
difficult to tell whether or not somebody has chosen a frequentist
test before or after collecting the data is one of the major problems
with frequentism (as shown by its not adhering to \ref{criteria} item
\ref{whetherproceed} above).

But the problem is in reality much worse than this. Suppose in these
circumstances somebody reports that they have used a (corrected)
$\chi^2$-test and found significance at (exactly) the 0.95
level. Unless there is evidence that they had in some way restricted
themselves to using a $\chi^2$-test, and no other test, before
collecting the data, whether or not they actually have found a
statistically-significant result depends on what they \textit{would
  have} done if the $\chi^2$-test had yielded say 90\% rather than
95\% frequentist confidence:

\begin{enumerate}

\item If they would then have tried to publish their paper containing
  only a 90\%-confident result, then they have indeed found 95\%
  frequentist confidence.

\item If, on the other hand, they would then have tried a (corrected)
  Fisher exact test, and if that still didn't get the 95\% level they
  wanted they would then have tried the test of section
  \ref{bcfreq1text}, and if that still didn't get the 95\% level they
  wanted they would then have tried the test of section
  \ref{bcfreq2text}, and would have reported the highest level of
  frequentist confidence of these, then they have not actually found
  95\% frequentist confidence, but some lower level -- even though
  none of these hypothetical events actually happened, and even if one
  (or all !)  of these other tests would in fact have given them
  apparent 95\% frequentist confidence.

\end{enumerate}

This is all akin to frequentism's violation of criterion
\ref{whetherproceed} of section \ref{criteria}, even though those are
not the exact circumstances applying here.

\subsection{Conclusions from this example}

We draw a number of conclusions; since the various plots are
moderately difficult to distinguish from each other by naked eye
assessment of colour shading, the reader is encouraged to consider
specifically the counts of blocks (i.e. data values) which lead to
posterior probabilities or frequentist confidences above 0.95 or below
0.05 (table \ref{numbersabove}):

\begin{table}[hp]
\begin{tabular}{lrr}
\textbf{Test type} & \textbf{Number giving $>0.95$} &
\textbf{Number giving $<0.05$}\\
Bayes & 283 & 283\\
Frequentist section \ref{bcfreq1text} & 95 & 460\\
Frequentist section \ref{bcfreq2text} & 256 & 420\\
Directional $\chi^2$ naive & (290) & (290)\\
Directional $\chi^2$ corrected & 271 & 371\\
Fisher exact naive & (260) & (321)\\
Fisher exact corrected (Boschloo's test) & 280 & 300\\
Late choice (apparent) & (283) & (297)\\
Late choice (corrected) & 268 & 298
\end{tabular}
\caption{
\label{numbersabove}
Performance of various tests for comparison of jobsucillin and
torvaldomycin for treatment of \textit{Microsoftus gatesii}
infection. The numbers given are the numbers of data values resulting
in frequentist confidence or posterior probability $>0.95$ or $<0.05$
that $H_1$ holds as indicated. Numbers in brackets are for tests that
are either only an approximation or are being misused.}
\end{table}

\begin{figure}[hp]
\begin{center}

\includegraphics[scale=0.5]{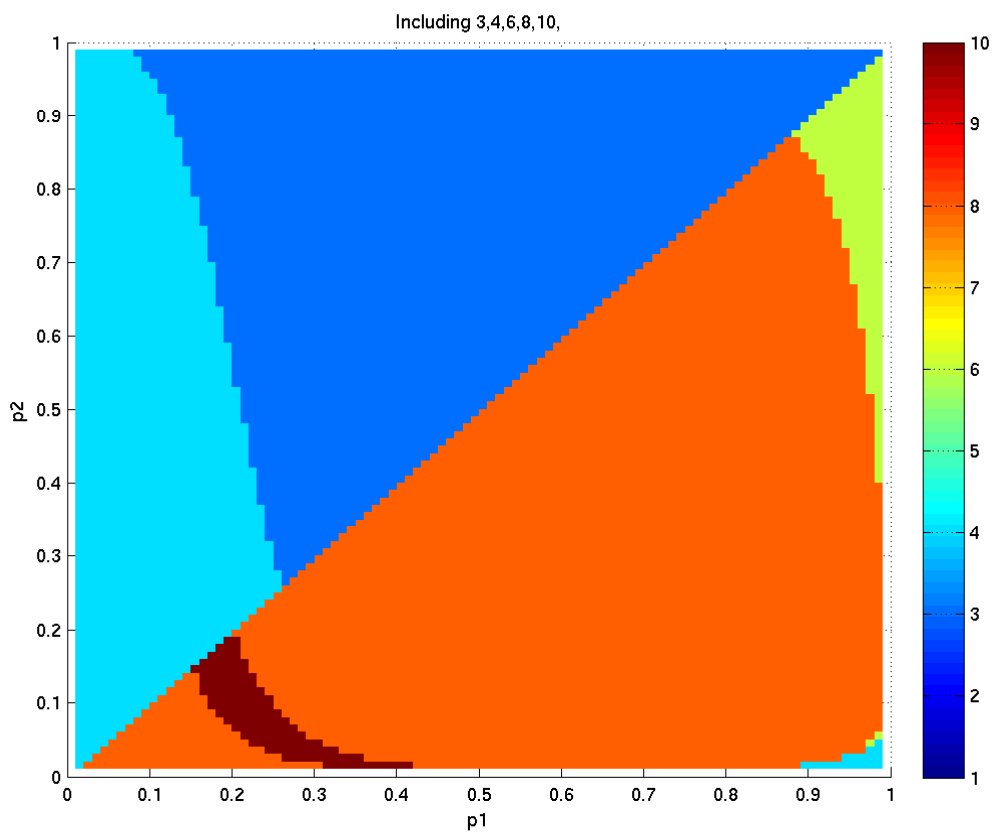}


\caption{The frequentist test described that gives greatest
  resp. least probability of concluding with 95\% frequentist
  confidence that $H_1$ holds when it does resp. doesn't among those
  considered given each value of $(p_1,p_2)$: 3 denotes that of
  section \ref{bcfreq1text}, 4 that of section \ref{bcfreq2text}, 6
  the corrected $\chi^2$-test, 8 the corrected Fisher exact test, and
  10 the corrected frequentist confidence from choosing after
  collecting the data among the others of these tests. (The omitted
  comparisons and numbers correspond to Bayes against thresholds 0.5
  and 0.95 and the inexact frequentist tests, none of which were
  included in this comparison.)
\label{optifreq}
}

\end{center}
\end{figure}

\begin{enumerate}

\item The Bayesian solution gives an entirely reasonable and symmetric
  result; none of the exact solutions in this selection of
  frequentist solutions does.

\item The Bayesian test is more likely to conclude that it is 95\%
  sure that $H_1$ holds than any of the exact frequentist
  tests considered. Moreover the symmetry of the Bayesian test means
  that there are the same number of data values for which it will
  conclude that it is 95\% sure that $H_0$ holds as there are for
  which it will conclude that it is 95\% sure that $H_1$ holds. 

\item The correct(ed) frequentist tests considered are all biased in
  favour of $H_0$ in the sense that there are always more data values
  for which they will conclude that frequentist confidence of $H_1$ is
  below 0.05 than that it is above 0.95, as well as in the sense that
  most of the time the Bayesian posterior probability of $H_1$ is
  greater than the frequentist confidence.

\item The commonly used $\chi^2$ test, in the uncorrected form in
  which it is usually used, is biased in favour of $H_1$ in comparison
  with exact frequentist inference using the same critical regions
  (whether in the directional version we used above or in the standard
  version where $H_1 = \{h\in H : p_1 \neq p_2\}$).

\item The commonly used Fisher's exact test is both inappropriate for
  this problem, and, as usually used, biased against $H_1$ in
  comparison with exact frequentist inference using the same critical
  regions.

\item None of these frequentist tests is uniformly optimal (no such
  test exists for this problem).

  Considering parameter space (as in the usual definition of uniformly
  optimal), figure \ref{optifreq} shows that the best of the
  considered exact frequentist tests to use varies with the actual
  value of $(p_1,p_2)$ applying. Even when $p_1>p_2$, sometimes the
  corrected version of Fisher's exact test does best, sometimes one
  does best using the corrected frequentist confidence resulting from
  choosing among the other exact tests listed after collecting the data.

  Considering data space, if $x=(n_1,n_2)=(8,1)$ then the solution of
  section \ref{bcfreq1text} concludes with $95\%$ frequentist
  confidence that $H_1$ holds, but none of the other exact frequentist
  solutions do. Equally if $x=(n_1,n_2)=(40,18)$ then the (corrected)
  solution of section \ref{bcchitext} so concludes, but none of the
  other exact frequentist solutions do. Or if $x=(n_1,n_2)=(24,7)$ then
  the solutions of sections \ref{bcfreq2text} and (corrected)
  \ref{bcFishertext} so conclude, but the others don't, not even the
  (corrected) one of section \ref{bcchooselatetext}.

\item The Bayesian test reaches its conclusion based on the prior and
  the data, independent of the intentions of the statistical
  analyst. The frequentist one, however, as can be seen from section
  \ref{bcchooselatetext} above, depends on what the intentions of the
  analyst were (unless he can prove he made his choice of test before
  collecting the data).

\item The frequentist test depends on a totally arbitrary choice of
  test, i.e. of nested set of critical regions. The Bayesian test, on
  the other hand, depends on the knowledge possessed before collecting
  the data, and can be reanalysed to evaluate the effect of differing
  assumptions about that knowledge \textit{ad lib}.

\end{enumerate}

In our view the Bayesian test is the only one that makes any sense; it
is also the only one that takes appropriate account of the relevant
prior. 

\section{Examples with long sequences of i.i.d. random variables}
\label{example2}

\subsection{Introduction}
\label{statintro}

For those readers who are more used to seeing inference applied to
problems involving long sequences of independent and identically
distributed random variables, we now consider a different type of
example problem, in two flavours and in two different settings. Here
we concentrate on showing the advantages of the Bayesian approach
rather than showing that the frequentist approach is wrong, as the
latter has been dealt with in sections \ref{critsection} to
\ref{example3} above.

In these examples an experimenter is aiming to show that $h>h_0$ for
some fixed value of $h_0$, so he sets $H_0=\{h\in H : h \leq
h_0\}$. He would like to collect some number $K$ of samples of data
$x_1,x_2,...,x_K$ and reach the conclusion that $h\in H_1$ with some
predetermined level of sureness $\eta$, usually 0.95 (assuming of
course that in truth $h>h_0$). However he does not know the true value
of $h$, and of course the closer $h$ is to $h_0$, the larger $K$ will
have to be to get a high probability that he can conclude that he is
$\eta$-sure that $h\in H_1$.

In the first (Bernoulli) flavour of this example we will assume that
each $x_k$ is independently Bernoulli distributed with parameter $h$,
in other words it is 0 or 1 with probability $h$ that it is 1 (so that
$H=[0,1]$), and take $h_0=0.975$. One might think of this as wanting
to be sure that if we drop a product object on the floor, the
probability of it remaining intact is at least $h_0$. So our data will
consist of a sequence of breakages or survivals of a length somehow
determined. 

In the second (Gaussian) flavour we will assume that each $x_k$ is an
independently Gaussianly distributed real number with unit variance
and mean $h$ (so that $H=\mathbb{R}$) and take $h_0=0$. One might
think of this as wanting to be sure that some measurement gives a
positive result on average. So our data will consist of some number of
measurements and the resulting values. 

Thus in terms of our standard notation for inference problems from
section \ref{infprobdefsub} we are
setting $$\Theta=\{0,1\},$$ $$\Phi=[0,1]\text{ or
}\mathbb{R},$$ $$H=\{(\theta,\phi)\in \Theta\times\Phi: \theta =
[\phi>h_0]\},$$ $$H_k=(\{k\}\times \Phi)\cap H \text{ for
}k\in\{0,1\},$$ $$X=\bigcup_{k=0}^\infty{\{0,1\}^k}\text{ or
}\bigcup_{k=0}^\infty{\mathbb{R}^k},$$
and $$P(x_k|\theta,\phi)=\phi^{x_k}(1-\phi)^{1-x_k}\text{ or
}\frac{1}{\sqrt{2\pi}}e^{-\frac{1}{2}(x-\phi)^2}$$ respectively for
the two flavours, and we can identify $h$ with $\phi$ without
ambiguity. Here saying that we have observed $\mathbf{x}\in
\{0,1\}^K\subset X$ means not only that we have observed the length
$K$ sequence $x_1,x_2,...,x_K$ but also that we have decided to stop
data collection at that point; note that then in the notation of
section \ref{criteria} criterion \ref{unifpow} the likelihood is then
given by $$P(\mathbf{x}|\theta,\phi)=q_K(x_1,...,x_K)
\prod_{k=1}^K{P(x_k|\theta,\phi)(1-q_{k-1}(x_1,...,x_{k-1}))},$$ where
the $q_k$ depend on the data collection plan. Of course, in the
Bayesian solution all the factors involving the $q_k$ are the same in
numerator and denominator of Bayes' theorem, so cancel. Similarly in
the fixed $K$ frequentist approach only $q_K$ is non-zero (indeed it
is 1) and it doesn't depend on the $x_k$, so these factors
disappear. In the general frequentist case, however, the situation is
more complicated.

Note that although the first flavour is reasonably like real-life, the
second is not, as both known variance and Gaussianity are only rarely
encountered in practice; we use this oversimplistic second flavour as
it is then simple enough for the reader to check our calculations
without too much effort.

To give us some lifelike vocabulary we will consider that we are
testing whether or not a particular factory meets some specification
on the products it produces. In the first setting we will assume that
only one aspect of the products produced has to be checked. In the
second setting we will take it that $N=224$ independent similar
aspects of the products have to be checked, and that the factory will
be approved if and only if for all $N$ aspects we are $95\%$-sure
that $h_n > h_0$. Any reader who thinks that one would never need to
check as many as $224$ aspects might be interested to read ISO 20072,
and in particular the example treated in \cite{CTM}, where 16
different aspects of product performance have each to be assessed
under 14 different sets of environmental conditions (although they are
not all of the same flavour).

Whether we are considering a frequentist or a Bayesian approach, we
will consider that the value of $h$ is drawn from a specific input
distribution giving a probability of 0.9 that the factory is good
(i.e. has $h>h_0$ or all $h_n>h_0$ as appropriate) as specified in
table \ref{inputtable}, while for the Bayesian approach we will use a
non-informative broad prior for the actual inference. As a measure of
performance of our inference method we will consider how small $K$ can
be while achieving a suitably high probability (e.g. $0.9$ or $0.99$)
that a good factory drawn from the input distribution will be
approved.

\begin{table}[hp]
\begin{tabular}{llllr}
 & \multicolumn{2}{c}{\textbf{Setting}}\\
\textbf{Flavour} & \textbf{Single aspect} & \textbf{224 aspects} &
\textbf{Prior used for Bayesian inference} & $h_0$\\
Bernoulli & Beta(155,2) & Beta(400,2) & Beta(1,1)=Uniform & 0.975\\
Gaussian & Gaussian(1.3,1) & Gaussian(3.33,1) & Gaussian(0,$10^{-4}$) & 0
\end{tabular}
\caption{
\label{inputtable}
The input distributions used for the $h$ or $h_n$ parameters of
simulated factories to be tested. Parameters of the Gaussian are mean
and scale (=1/variance). All four cases give a probability of
approximately 0.9 that a factory with $h$ or each $h_n$ drawn from the
given distribution will have $h>h_0$ or each $h_n > h_0$ as appropriate.
}
\end{table}

Our aim then will be to be able to approve 90\% (or 99\%) of good
factories drawn from the input distribution. We will consider three
different approaches: the frequentist approach using a predetermined
fixed value of $K$; a  pseudo-Bayesian frequentist approach designed to
be suitable for the target factory pass probability; and the Bayesian
approach using either the uniform prior on $h\in [0,1]$ or a Gaussian
prior of zero mean and standard deviation 100 on $h\in \mathbb{R}$ as
appropriate to the problem flavour. There are of course infinitely
many other possible frequentist approaches, but time and space
restrict us to consider just these here.

We will then examine what average value of $K$ is needed for each of
the three approaches. We first specify the method of calculation for
each method, and then summarise the necessary values of $K$ in a joint
table. All simulations use factories drawn from the input
distribution, discarding any that are not good, until we have 100,000
factories to test. In the case of the Bayesian approach we will also
separately consider the average value of $K$ needed when the bad
factories from the input distribution are also included.

Note that both the frequentist approaches make use of knowledge of the
input distribution (the ``true prior'') in their design (at least as a
worst case that they have to deal with), while the Bayesian approach
\textit{does not} -- it would be exactly the same whatever input
distribution were expected, so long as it is intended to be executed
without knowledge of that distribution, and its parameters depend only
on the nature of the likelihood, on $h_0$, on $N$, and on the desired
probability of passing a good factory. Of course, the Bayesian
approach uses a prior in its execution, namely the flat uniform prior
or a very broad Gaussian prior as appropriate to the problem flavour.

Further, both the frequentist approaches are designed for the purpose
of passing a specific fraction of the given input distribution; the
Bayesian approach, however, requires no such determination in advance
-- it can just be run for as long as it takes to pass any particular
factory (in real life there is only one factory under consideration at
any one time).

\subsection{Frequentist approach with fixed predetermined sample size}

We search for the appropriate value of $K$ by binary chop on a
suitable range of values.

For each potential value of $K$ we calculate the appropriate uniformly
optimal critical region to achieve, for each aspect, 95\% frequentist
confidence that $h\in H_1$. For each good sample factory $f$ drawn
from the input distribution as above we calculate the probability
$p_f$ that the data for a single aspect will lie in the critical
region. The probability of a factory being approved is given
by $$(\mean_f\,{p_f})^N,$$ where $N$ is the number of aspects in the
given setting.

\subsection{Bayesian approach}

We simulate on the case $N=1$ (i.e. on a single aspect), calculating
the corresponding figures for the $N=224$ case afterwards. For each
good sample factory we collect one data point at a time, each time
calculating the posterior probability that $h>h_0$ given the data so
far. When that probability first rises above $p_1$ (resp. falls below
$p_0$) we pass (resp. fail) that aspect, noting the number of
datapoints that have been required. The values of $p_0$ and $p_1$ are
as shown in table \ref{Bayesparmtable}. The manufacturer is left to
decide when to give up testing if the posterior probability passes
neither limit, which time will depend on the cost of the factory, the
cost of testing, and the current posterior probability. In the
simulations, if wanting to pass a good factory with probability
$p=0.9$ say, we will need to pass each aspect with probability
$p^{\frac{1}{N}}$. 

\begin{table}[hp]
\begin{tabular}{lrr|r|r}
& \textbf{Setting} & \textbf{Fraction of}\\
\textbf{Flavour} & \textbf{($N$)} & \textbf{good to pass} & 
$p_0$ & $p_1$\\
\hline
Bernoulli & 1 & 0.9 & $10^{-2}$ & 0.95 \\
& & 0.99 & $10^{-4}$ & 0.95 \\
\hline
Bernoulli & 224 & 0.9 & $10^{-4}$ & 0.95 \\
& & 0.99 & $10^{-5}$ & 0.95 \\
\hline
Gaussian & 1 & 0.9 & $10^{-1}$ & 0.95 \\
& & 0.99 & $10^{-2}$ & 0.95 \\
\hline
Gaussian & 224 & 0.9 & $10^{-4}$ & 0.95 \\
& & 0.99 & $10^{-5}$ & 0.95
\end{tabular}
\caption{
\label{Bayesparmtable}
The settings of the various parameters used for the various Bayesian
tests reported in table \ref{statresults}.}
\end{table}

\subsection{Pseudo-Bayesian frequentist approach}

We set $h_0'>h_0$ such that roughly the desired fraction of good
factories to be passed from the input distribution have $h_n>h_0'$ for
all $n$. We then set $p_0$, $p_1$, $K_{\text{min}}$, and
$K_{\text{max}}$ such that the probability of a factory with $h=h_0$
passing in a single aspect is $\approx 0.05$ under the Bayesian
approach that has exit barred before $K_{\text{min}}$ and after
$K_{\text{max}}$ samples have been collected and otherwise fails if
posterior probability of $H_1$ falls below $p_0$ and passes if it
rises above $p_1$ (as estimated by $10^5$ simulations). We then
proceed as in the Bayesian approach but with the modifications implied
by $h_0'$, $p_0$, $p_1$, $K_{\text{min}}$, and $K_{\text{max}}$.

The specific values used are given in table
\ref{pseudosettings}. Because of the difficulty of achieving exactly
0.95 frequentist confidence, we bracket this figure with two different
sets of values to give a range of expected $K$ required.

Note that many other choices are possible in designing such a 
pseudo-Bayesian test, and we know no way of finding the most efficient
such solution (or whether or not a most efficient one exists), so we are
not claiming that these choices give the most efficient
basic pseudo-Bayesian approach possible. The resulting test is then at least
adaptive to the nature of the actual factory and aspect that it is
currently considering.

\begin{table}[!htp]
\begin{tabular}{lrr|r|r|r|r|r|r}
& \textbf{Setting} & \textbf{Fraction of}\\
\textbf{Flavour} & \textbf{($N$)} & \textbf{good to pass} & $h_0'$ &
$p_0$ & $p_1$ & $K_{\text{min}}$ & $K_{\text{max}}$ & \textbf{conf} \\
\hline
Bernoulli & 1 & 0.9 & 0.979 & 0.001 & 0.95 & 400 & $1 \times 10^5$ & 0.955 \\
& & & 0.979 & 0.001 & 0.95 & 320 & $1 \times 10^5$ & 0.949 \\
\hline
Bernoulli & 1 & 0.99 & 0.9756 & 0.001 & 0.991 & 1500 & $5 \times 10^6$ & 0.956 \\
& & & 0.9756 & 0.001 & 0.991 & 800 & $5 \times 10^6$ & 0.946 \\
\hline
Bernoulli & 224 & 0.9 & 0.97695 & 0.001 & 0.96 & 1500 & $2 \times 10^6$ & 0.953\\
& & & 0.97695 & 0.001 & 0.96 & 1300 & $2 \times 10^6$ & 0.949\\
\hline
Bernoulli & 224 & 0.99 & 0.9752 & 0.0005 & 0.985 & 50000 & $4 \times 10^7$ & 0.959 \\
& & & 0.9752 & 0.0005 & 0.985 & 30000 & $4 \times 10^7$ & 0.948 \\
\hline
Gaussian & 1 & 0.9 & 0.36 & 0 & 0.95 & 0 & 125 & 0.951 \\
& & & 0.35 & 0 & 0.95 & 0 & 125 & 0.949 \\
\hline
Gaussian & 1 & 0.99 & 0.04 & 0.001 & 0.96 & 150 & 4000 & 0.952 \\
 & & & 0.04 & 0.001 & 0.96 & 130 & 10000 & 0.947 \\
\hline
Gaussian & 224 & 0.9 & 0.19 & 0.001 & 0.96 & 4 & 1000 & 0.9505 \\
 &  &  & 0.19 & 0.001 & 0.96 & 3 & 1000 & 0.943 \\
\hline
Gaussian & 224 & 0.99 & 0.02 & 0.001 & 0.96 & 600 & 25000 & 0.9504 \\
 &  &  & 0.02 & 0.001 & 0.96 & 500 & 25000 & 0.944 \\
\end{tabular}
\caption{
\label{pseudosettings}
The values of the various parameters used for the various
 pseudo-Bayesian tests reported in table \ref{statresults}. The column
headed ``conf'' gives the frequentist confidence achieved with these
values, based on $10^5$ simulations. 
}
\end{table}

\subsection{Results}
\label{resultssection}

Table \ref{statresults} gives the values of $K$ resulting for each
approach to each flavour and setting.

\begin{table}[!h]
\begin{tabular}{lrr|r|r|r|r}
&&& \multicolumn{4}{c}{\textbf{Approach}} \\ \cline{4-7}
& \textbf{Setting} & \textbf{Fraction of} & \multicolumn{2}{c|}{\textbf{Frequentist}} &
  \multicolumn{2}{c}{\textbf{Bayesian}} \\ \cline{4-7}
\textbf{Flavour} & \textbf{($N$)} & \textbf{good to pass} & \textbf{Fixed
  $K$} & \textbf{ Pseudo-Bayesian} & \textbf{Good} & \textbf{All}\\
\hline
Bernoulli & 1 & 0.9 & 3000 & 1464 - 1509 & 472 & 568 \\
& & 0.99 & 130000 & 27965 - 29415 & 1145 & 2392 \\
& & 0.996 & & & 2091 & 5573 \\
\hline
Bernoulli & 224 & 0.9 & 5600000 & 323142 - 384169 & 45483 & 46106 \\
& & 0.99 & $>$ 200000000 & 6770414 - 11288141 & 49720 & 57616 \\
\hline
Gaussian & 1 & 0.9 & 19 & 15 - 18 & 3.2 & 3.5\\
& & 0.99 & 1000 & 178 - 187 & 9.8 & 16\\
& & 0.996 & & & 28 & 60\\
\hline
Gaussian & 224 & 0.9 & 16352 & 758 - 972 & 273 & 278\\
& & 0.99 & 672000 & 112092 - 134638 & 285 & 324\\
& & 0.998 & & & 294 & 426\\
\end{tabular}
\caption{
\label{statresults}
The expected number of datapoints needed for all $N$ aspects of a good
factory drawn from the relevant input distribution to be approved,
under the various methods. Values were estimated from 100,000 samples
of good factories from the input distribution, except that the last
column gives the corresponding number for the Bayesian test if the bad
factories are also included in the factories being tested (for
fixed-$K$ this makes no difference, and for  pseudo-Bayesian
these values are not shown but would be higher than those for good
factories). For the  pseudo-Bayesian tests the first figure
resulted from the run for the lower bracketing frequentist confidence
from table \ref{pseudosettings} and the second for the higher. For the
Bernoulli 224-aspect fixed-$K$ case the number of datapoints total
required to get 99\% of the good factories to pass was in the range
$2\times 10^8$ to $2.5\times 10^8$.}
\end{table}

\subsection{Discussion of results}
\label{statresdisc}

\subsubsection{General}
\label{statresgen}

We first remind ourselves of what the different methods are supposed
to achieve:

\begin{itemize}

\item Bayesian approaches that result in a factory being given 0.95
  posterior probability of being good are supposed to result in not
  more than 0.05 of such approved factories actually being bad.

  If the prior matches the actual distribution of factories being
  tested, then that will always be the case; however in the examples
  here the two distributions are very different, and we would like to
  know what fraction of the approved factories are in fact bad.

\item Frequentist approaches that result in a factory getting 95\%
  frequentist confidence that it is good are supposed to have the
  effect that any bad factory gets such approval with probability not
  more than 0.05, irrespective of the input distribution. This will be
  the case with these frequentist solutions in particular.

\end{itemize}

Of course, what actually matters to the users of these factories'
products is whether the ones that are approved meet specification or
not. If most of the factories being tested are in fact bad, then the
frequentist achievement could be that the total number of factories
getting approval is low compared with the number tested, but most of
the approved ones are bad. Moreover, since the frequentist criterion
makes no promises that it won't fail good factories with arbitrarily
high probability, it is even possible in principle that the fraction
of bad factories in the frequentist-approved set is higher than it was
in the original distribution of factories being tested, although this
cannot happen in this specific example.

So we now review from table \ref{statresults} how much data is
required to achieve various fractions of good factories passing, and
then turn to ask what happens to bad factories in this example, which
apart from the oversimplification of the second flavour is very
realistic.

\subsubsection{What happens to good factories with each approach ?}

It is striking how many fewer drops or measurements respectively are
needed to approve any given fraction of good factories by the Bayesian
approach than by either of the other two. The smallest ratio of sample
numbers needed by any of these frequentist solutions to get 90\% passes
of good factories to the corresponding Bayesian number is 2.5 (for the
Bernoulli $N=1$ case, pseudo-Bayes / Bayes), and the second smallest
such number is 2.7 (for the Gaussian $N=224$ case, pseudo-Bayes /
Bayes). The highest similar ratio (for 90\% passes of good factories)
observed here is 121 (for the Bernoulli $N=224$ case, fixed-$K$ /
Bayesian).

For 99\% pass rate of good factories, we observe that the fixed-$K$,
Gaussian $N=224$ version requires over 2000 times as much data as the
Bayesian version does, and the corresponding Bernoulli case requires
over 3000 times as much data as the Bayesian version. Even for the
comparison between  pseudo-Bayesian and Bayesian tests, in the
Gaussian $N=224$ case the  pseudo-Bayesian requires more than 300
times as much data. These comments apply even if the appropriate
fraction of bad factories are included for the Bayesian case, and have
used the lowest of the two bounding figures for the  pseudo-Bayesian
case. (Of course this is not to say that one might not be able to
design a different frequentist solution that did better.)

One should also note that for the Bayesian method, the amount of data
required rises less fast than $N$ (bearing in mind the input
distributions were adjusted to make 90\% of the input factories good
in both cases), while the frequentist versions have data demands which
increase much faster than $N$.

Finally we note that the frequentist test designers need to have some
idea of the input distribution (or prior !) to be able to design their
tests with adequate probability of passing good factories. What
evidence would a frequentist fixed-$K$ designer need to be reasonably
sure that a reasonable prior to design to was as used here, if he
optimally used a flat prior and Bayes (!) to find out ?

For the Bernoulli flavour, for $N=1$, he would need to drop $155$
product samples and observe only one breakage; for $N=224$ he would
need to drop $224\times 400=89600$ devices, observing only one
breakage in each of the 224 sets; in the first case he has already
consumed a quarter as much data as the Bayesian is expected to
need\footnote{Recall that the Bayesian needs no such preparatory data
  as he is using a flat prior.} to pass 90\% of good devices, and in
the second nearly twice as much -- and that's before he's even started
his test proper. For the Gaussian flavour (variance being known) a
single measurement might suffice (depending on its value), but that
single measurement has already used over a quarter of the expected
number of measurements needed by the Bayesian to pass 90\% of good
devices in the $N=1$ case.

\subsubsection{What happens to bad factories with each approach ?}

But we must also consider what happens to bad factories with each
approach, because they are not treated the same, as mentioned in
section \ref{statresgen} above. The frequentist case is easy: no bad
factory has probability greater than 0.05 of passing a test of an
aspect on which it is bad.

In the Bayesian case the situation needs more careful
consideration. First let us consider what the Bayesian method does
from a frequentist point of view. 

\begin{figure}[!htp]
\begin{center}

\includegraphics[scale=0.5]{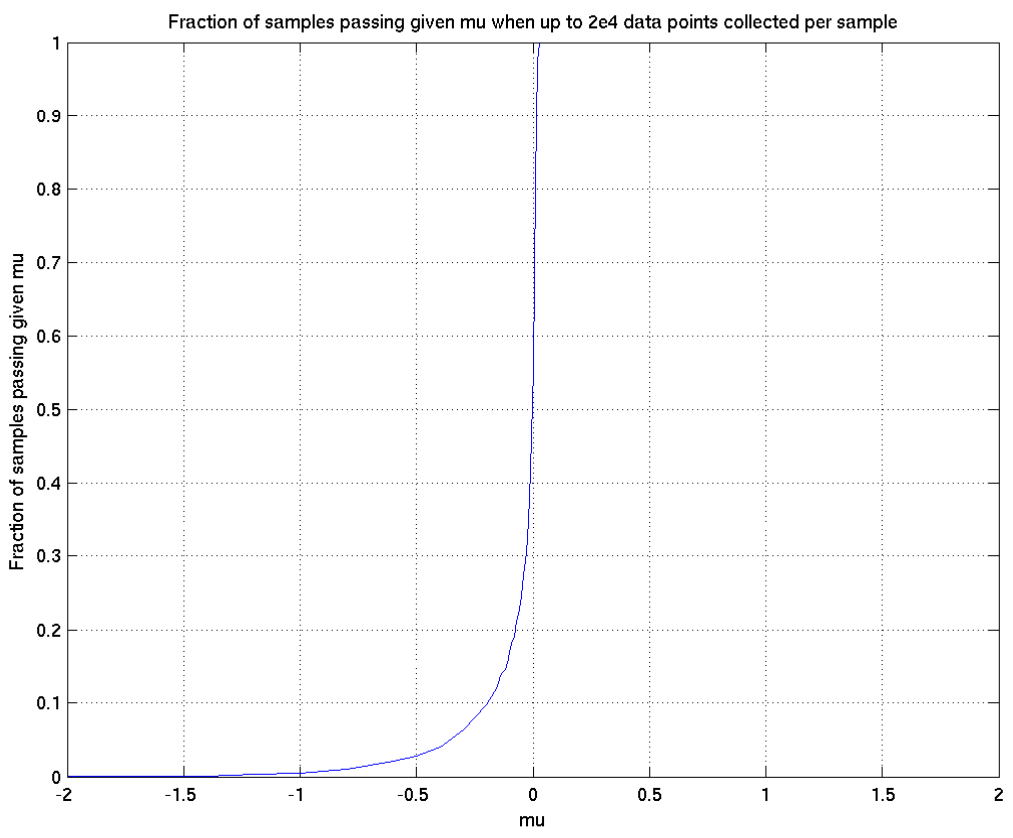}

\caption{The probability of a factory passing a single
  Gaussian-flavour aspect by the Bayesian solution as a function of its
  mean measurement. Here the Bayesian solution was set to fail the
  factory on this aspect if after 20,000 measurements it had not
  passed. (An alternative criterion for giving up might be that the
  posterior probability that the factory was \textit{bad} rose above
  0.95; that would then result in a symmetric curve about zero.)
\label{bayesbads}
}

\end{center}
\end{figure}

Figure \ref{bayesbads} shows the situation in the Gaussian $N=1$
case, where the Bayesian solution was set to abandon the factory (and
fail it) if it had not passed after 20,000 measurements. A consequence
of the property that Bayes will (eventually) pass any good factory is
that as $h$ approaches $h_0$ from below, the probability of passing
approaches 1 (if arbitrarily many measurements are allowed). So one
might reasonably worry that factories that were only just bad might
contaminate the output to the extent that one gets more than 5\% of
the approved factories actually being bad. So does that happen here ?

\begin{table}[h]
\begin{tabular}{lrccc}
\textbf{Flavour} & \textbf{Setting ($N$)} & \textbf{Bads/Total in
  input} & \textbf{Pass rate / Bads} & \textbf{Bads passed / Total passed}\\
Bernoulli & 1 & 0.096 & $0.11$  & $1.02 \times 10^{-2}<0.05$\\
Bernoulli & 224 & $4.4\times 10^{-4}$ & 0.16 & $7 \times 10^{-5} <
2.2\times 10^{-4} =0.05 / 224$\\
Gaussian & 1 & 0.097 & 0.08 & $8.1\times 10^{-3}<0.05$\\
Gaussian & 224 & $4.3\times 10^{-4}$ & 0.11 & $4.8\times
10^{-5}<2.2\times 10^{-4}=0.05/224$ \\
\end{tabular}
\caption{
The probability that an aspect of a factory is bad in the
  given input distribution, the pass rate with the Bayesian solution on
  those bad aspects, and the resulting fraction of the apparent
  passed aspects that are actually bad.
\label{statresbads}
}
\end{table}

The answer is No. Table \ref{statresbads} shows the fraction of the
input distribution that is bad (for a single aspect), the fraction of
those bad factories that pass this single aspect, and the resulting
fraction $f$ of the passed population that are actually bad. In all
cases $Nf$, an overestimate of the population of factories that are
bad, is below 0.05 of the passed population of factories.

For that not to be the case, we would need \textit{both} the fraction
of bads in the input distribution to be higher, \textit{and} the bads
to be concentrated exactly at or just below $h_0$. A factory designer
would not be able to achieve that if he wanted to. Given that of
course in real life we do not have an input population but just a
single factory that we are testing, the corresponding proposition is
that we believe with high probability that our factory has $h$
\textit{just} below $h_0$. In the highly unlikely event that that is
the case, the situation can of course be fixed by using a prior that
matches what we actually think about the factory under test.

\subsubsection{What is the meaning of a bad factory passing in the
  Bayesian case ?}

Now let us think about the meaning of a bad factory passing the test
from a Bayesian point of view. This means that the dataset that has
been collected is at least 0.95 likely to have come from a good
factory, assuming the prior we have used for inference. Given that
that is the only dataset we have to consider on this point, if the
prior we have used reflects our beliefs about how likely different
values of $h$ are, then it is therefore entirely appropriate that it
be approved.

\section{How should we set a prior with a regulatory authority ?}
\label{regulatorprior}

\subsection{Scenario}

If then we are going to use Bayesian methods, we will need to be
prepared to express what we know before collecting the data as a prior
probability distribution on the unknowns. Where only one party at a
time needs to use the deductions or predictions that result, this may
be relatively straightforward.

However, one situation where two parties with different interests may
need to be able to reach a common conclusion about some unknown
parameters is when a regulatory agency needs to be convinced by a
manufacturer that his product satisfies necessary legal
requirements. For example, an asthma inhaler may have a requirement
that the dose delivered is within some particular range for a given
fraction of the doses delivered; or there may be a requirement that
when a device is dropped onto a hard floor from a given height it has
at least a given probability of surviving undamaged. In this situation
a regulator might want a prior on the parameters that was sceptical
that the device meets the given requirements, while the manufacturer
might want one that even before collecting the data considers it very
likely that the device meets the requirements.

It is also possible that regulator and manufacturer even differ on the
type of distribution of the regulated attribute; for example the
manufacturer may think that his inhaler delivers log doses distributed
according to a Gaussian
distribution $$P(x|\mu,\sigma)=\frac{1}{\sigma\sqrt{2
    \pi}}e^{-\frac{1}{2}(\frac{x-\mu}{\sigma})^2},$$ with suitable
prior distributions on the parameters $\mu$ and $\sigma$, while the
regulator thinks the log doses are distributed according to a Student
distribution $$P(x|m,r,\mu)=\frac{1}{\sqrt{2\pi}}\frac{\Gamma(m+\frac{1}{2})}{\Gamma(m)}
\frac{r^m}{(r + \frac{1}{2}(x-\mu)^2)^{m+\frac{1}{2}}},$$ again with
suitable priors on the parameters $m,r,\mu$.

\subsection{Proposed way forward -- ``Let the data decide''}

\label{datadecide}

In setting a prior that is acceptable to both parties, the first step
should obviously be to attempt to reach agreement through discussion
after sharing all available data. If unsuccessful, it may then be
necessary that a formal mechanism is available for the parties to
reach a prior by a mechanism that automatically mediates between their
two desired priors and/or models. Moreover knowing that such a
mechanism exists may also facilitate agreement actually being reached
at the initial discussions.

We therefore propose a mechanism for resolving such disputes about
priors as follows, to which we give the moniker ``Let the data
decide'', following MacKay\cite{DataDecide}.

Let us denote by $M_1$ the model desired by the regulator and by $M_2$
the model desired by the manufacturer, in each case complete with
priors on the parameters ($\theta_1, \theta_2$ respectively) deemed
suitable by the relevant party. Let us denote by $A$ the event that
the actual distribution meets the desired criteria (e.g. has 90\% of
log doses between $a$ and $b$), so that when the data $D$ is
collected, $P(A|M_k,D)$ is the probability that if model $M_k$ holds
the legal requirements are satisfied.

Our proposal, then, is very simple: we make no decision on whether
$M_1$ or $M_2$ is to be used, but instead assume that the test data
$D$ arises as a result of parameters and likelihood taken from one or
other of these models with probability $\frac{1}{2}$ on each.

\subsection{General effect of this proposal}

One might be forgiven for thinking that then we simply
get $$P(A|D)=\frac{1}{2}P(A|M_1,D) + \frac{1}{2}P(A|M_2,D).$$ But this
is \textit{not} what happens. Rather, according to the marginalisation
rule and the chain rule of probability,
\begin{IEEEeqnarray*}{rCl}
P(A|D) &=& P(A,M_1|D) + P(A,M_2|D)\\
&=& P(A|M_1,D)P(M_1|D) + P(A|M_2,D)P(M_2|D).
\end{IEEEeqnarray*}

In other words, each model contributes to the confidence that the
legal requirements are met in proportion to the probability that that
model did indeed give rise to the data.

\subsection{Examples of the use of this proposal}

\subsubsection{Introduction to examples}

For each of these examples, for simplicity, we will consider the
situation that the legal requirement to be satisfied is that when a
device is dropped onto a hard floor a single time from a defined
height the probability $p$ that it survives intact must be at least
$\tau = 0.9$. We further suppose that both parties agree that, given
that no device will be tested more than once, the outcomes of the
various drops will be independent, so that after $N$ drops, the number
of intact devices will be binomially distributed according to
$$P(n|N,p) = \frac{N!}{n!(N-n)!}p^n(1-p)^{N-n}.$$

For all examples we used the optimum Bayesian data collection plan as
in appendix \ref{bayesproofs} criterion \ref{appunifpow}, except that
we stopped each simulation if any of the following occurred:
\begin{itemize}
\item $P(A|D) > 0.95$ (concluding requirements probably met, as in the
  optimum Bayesian data collection plan));
\item $P(A|D) < 0.01$ (concluding requirements probably not met);
\item more than $10,000$ devices had been dropped (leaving the
  decision whether or not to continue up to the manufacturer).
\end{itemize}

We first look at three unrealistic examples simply for the purpose of
showing that ridiculous settings don't lead to ridiculous results. We
will then turn to look at a more realistic scenario, first from a
frequentist point of view then from a Bayesian point of view, before
exploring the effects of mismatch between prior and the actual
distribution of regulatory submissions on a range of priors.

\subsubsection{Both parties pulling in opposite directions, version 1}

\label{bothobstructive}

Suppose the regulator decides to make it as difficult as possible for
the manufacturer to pass the test (but without actually making it
impossible), and sets the prior on $p$ to be given by $$P(p|M_1) =
\left\{\begin{matrix}1-\delta & (p = \tau - \epsilon)\\ \delta & (p =
\tau + \epsilon)\\ 0 & (\text{otherwise})\end{matrix}\right.$$ where
$\delta>0,\epsilon>0$ are small positive numbers perhaps in the region
of $0.001$.

The manufacturer, on the other hand, perhaps misguidedly, attempts to
abuse the system in his own favour by setting the prior $$P(p|M_2) =
\left\{\begin{matrix}1 & (p = 1)\\0 &
(\text{otherwise}).\end{matrix}\right.$$

What then happens as data is collected using the optimal data
collection plan of appendix \ref{bayesproofs} criterion
\ref{appunifpow} ?

First, if every device dropped survives intact, the probability that
the device factory is legally acceptable $P(A|D)$ passes $0.95$ when
$N=\lceil 27.946\rceil = 28$. For comparison, the uniformly optimal
frequentist test for $H_0=(p\in [0,0.9])$ reaches $95\%$ frequentist
confidence when $N=\lceil 28.433\rceil=29$, hardly a large difference.

Second, as soon as a single device fails to survive intact, $P(M_2|D)$
becomes zero, so that the regulator's obstructive prior takes full
effect, and depending on the values of $\epsilon$ and $\delta$ it may
become arbitrarily difficult to pass the test even after dropping many
thousands of devices.

\subsubsection{Both parties pulling in opposite directions, version 2}

If the regulator is still obstructive using the same prior as in
section \ref{bothobstructive}, a manufacturer who is slightly wiser
but still biased in his own favour might choose to set the
prior $$P(p)=\frac{[\tau< p \leq 1]}{1-\tau},$$ i.e. the prior that is
uniform on the interval $(\tau,1]$, again expressing $100\%$ prejudice
  that his factory produces legally acceptable devices. Again using
  the uniformly optimal data collection plan of appendix
  \ref{bayesproofs} criterion \ref{appunifpow}, the number of tests
  needed before $P(A|D)=0.95$ is now a random variable that depends on
  the true value of $p$. If in truth $p=1$, so that no devices break,
  we reach this posterior probability when $N=42$ (the answer to life,
  the universe, and everything according to Douglas Adams\footnote{In
    ``The Hitch-hiker's guide to the Galaxy''}), a number which is
  even larger than the frequentist's 29, while if $p<1$ then it will
  be expected to take longer still. Passing the test is now almost
  synonymous with showing that $M_2$ holds. A full exploration of what
  happens when $p<1$ will be found below in section \ref{examplesims}.

\subsubsection{Both parties pulling in opposite directions, version 3}

Keeping again the same prior for the regulator, the biased
manufacturer might try doing exactly the opposite of what the
regulator does but more so, by setting his prior according to $$P(p) =
[p=\tau + \epsilon],$$ expressing total prejudice that $p$ is just on
the right side of $\tau$. Now it becomes very difficult to get
$P(M_2|D)$ to move upwards from $0.5$, as both models are essentially
saying that $p=\tau$ to a good (but slightly different)
approximation. Given these particular choices, again passing the test
is almost synonymous with showing that $M_2$ holds, which is extremely
difficult. 

\subsubsection{Realistic regulator but manufacturer calls upon
  his past experience to bias prior in his favour}

\label{reasonableregvpast}

In this more realistic situation, the regulator sets the
prior $$P(p)=\frac{\Gamma(25)}{\Gamma(21)\Gamma(4)}p^{20}(1-p)^{3},$$
which is a Beta prior with its $0.025$ and $0.975$ quantiles at
$0.676$ and $0.953$, also implying that there is only a $0.2$
probability that the factory's devices meet requirements. On the other
hand the manufacturer has dropped $100$ similar (but not quite
identical) devices without breakage before, so wants to set the prior
$P(p) = 101p^{100}$, i.e. to say that $p$ is
Beta(101,1)-distributed. We note to start with that under the
manufacturer's prior the probability that the factory is good is
already around $0.99998$, so he would have no need to collect any data
at all were his prior to be accepted by the regulator. As it is, if in
fact $p=1$, then after dropping $17$ devices he will get no breakages,
and the posterior probability that $p>\tau$ will be about $0.956$;
this is because the data ($17$ drops with no breakages) has allowed
the inference that the data arose from $M_2$ rather than $M_1$ with
posterior probability

\begin{IEEEeqnarray*}{rClr}
P(M_2|D) &=& \frac{P(M_2)P(D|M_2)}{P(M_1)P(D|M_1) +
  P(M_2)P(D|M_2)}\\ &=&\frac{P(D|M_2)}{P(D|M_1)+P(D|M_2)} &
(P(M_1)=P(M_2)=\frac{1}{2})\\ &=&\frac{\int{P(p|M_2)P(D|p)\,dp}}{\int{P(p|M_1)P(D|p)\,dp}
  +
  \int{P(p|M_2)P(D|p)\,dp}}\\ &=&\frac{\int{(101p^{100})p^{17}\,dp}}{\int{(\frac{\Gamma(25)}{\Gamma(21)\Gamma(4)}p^{20}(1-p)^{3})p^{17}\,dp}
  + \int{(101 p^{100})p^{17}\,dp}}\\ &=&\frac{\frac{101}{118}}
{\frac{\Gamma(25)}{\Gamma(21)\Gamma(4)}\frac{\Gamma(38)\Gamma(4)}{\Gamma(42)}
  + \frac{101}{118}}\\ &\approx& 0.891
\end{IEEEeqnarray*}
so that
\begin{IEEEeqnarray*}{rCl}
P(A|D) &=& P(M_1|D)P(A|M_1,D) + P(M_2|D)P(A|M_2,D)\\
&\approx& 0.109\times 0.597 + 0.891 \times 0.999994\\
&\approx& 0.956
\end{IEEEeqnarray*}
while in contrast, had the regulator's prior been accepted, one would
have needed $28$ drops without breakage to have reached the same
conclusion. 

Let us also consider the same pair of priors and see what happens if
in fact $p=0.85$, i.e. the devices produced by the factory are
unacceptable ($\tau$ being $0.9$).

In this case the number of drops until we know whether or not the
devices produced by the factory are acceptable is a random variable,
so to find out what happens we conducted simulations. We think of each
simulation as the testing of a factory producing devices. 

Note that the manufacturer was allowed to set a more definite
confidence threshold for failure than the regulator sets for success,
because it is the manufacturer who has to pay for more devices to be
tested -- and in any case, this only affects when he gives up, and
does not affect how sure we have to be that his product is good before
passing it.

Now, of the $100$ runs total, none required more than $10,000$ drops
to reach a decision, $93$ concluded that with probability $>0.99$ the
devices did not meet the requirements after an average of $141$
devices were dropped, and $7$ concluded that with probability $>0.95$
the devices did meet the requirements (although in fact they didn't).

In contrast, if the regulator's prior had been used unmodified, all
$100$ runs would have concluded that with probability $>0.99$ the
devices did not meet the requirements after an average of $145$
devices had been dropped.

We note that it is difficult to compare this with what would happen
with a frequentist test because there are so many possible choices of
critical regions and data collection plans, none of which is in any
sense optimal.

Turning back to devices that ought to pass the test, with $p=0.95$,
with the prior containing both regulator's and manufacturer's
components, of $100$ runs all concluded that the devices met the
requirements after on average $99$ devices had been dropped. In
contrast using only the regulator's prior the same conclusions were
reached, but for those concluded to meet the requirements on average
$162$ devices needed to be dropped.

The astute reader may by now have noticed that we have so far been
testing these priors using a frequentist point of view, contrary to
the spirit of this paper. Let us now consider what happens from a
Bayesian point of view. Here we will look at inference from each
possible point in the space of pairs $(n_{\text{good}},n_\text{bad})$.
In a later section of this paper (\ref{examplesims}) we will visit a
number of priors using factory parameters drawn from the prior in use
and from other priors.

Note that we cannot display the full data space in only $2$
dimensions. Although pairs $(n_{\text{good}},n_\text{bad})$ represent
a sufficient statistic for inference at that point in time (i.e. they
contain all the information for inference at that point without
collecting any more data), the outcome of a test run depends also on
the order in which the results come in -- e.g. $10$ breakages
followed by $90$ survived drops will have a different outcome from the
reverse order.

\begin{figure}[htp]
\begin{center}

\includegraphics[scale=0.5]{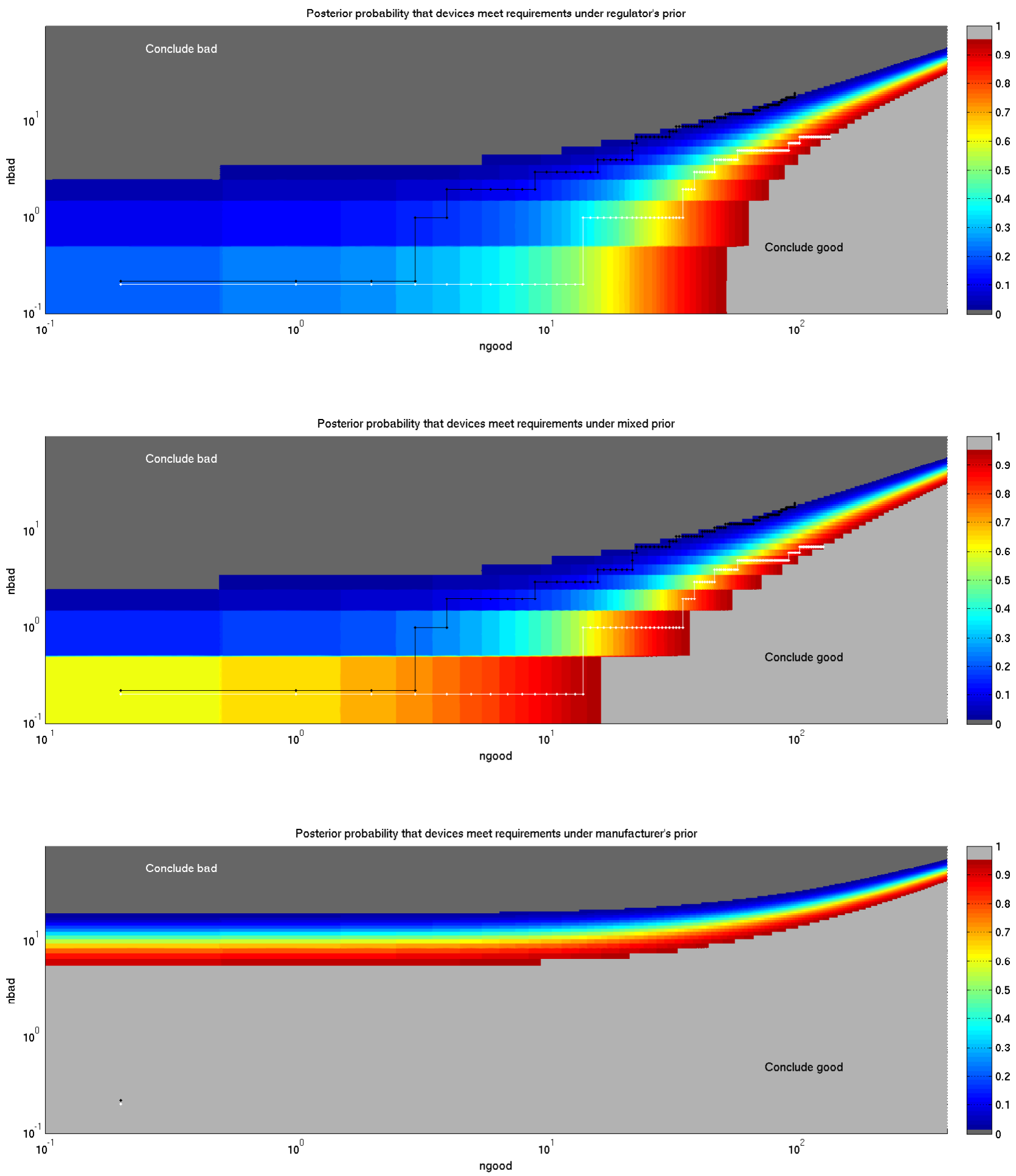}

\caption{Posterior probability that devices meet legal requirements as
  a function of the numbers of devices that are intact $n_\text{good}$
  and broken $n_\text{bad}$ using the regulator's prior that is
  Beta(21,4) and the manufacturer's prior that is Beta(101,1). The top
  plot uses the regulator's prior alone, the bottom one the
  manufacturer's prior alone, and the middle one the proposed mixture
  prior. The starting point at $(0,0)$ has been marked as $(0.2,0.2)$
  on these log scales. The white track shows a device with
  $p=0.95>\tau$ and the black track one with $p=0.85<\tau$. With the
  manufacturer's prior both are deemed good without collecting any
  data. The other two get the same results as each other, though it is
  clear that if there are no or few breakages, the mixed prior will
  conclude sooner than the regulator's prior.
\label{bvreg21c4man101c1}
}

\end{center}
\end{figure}

Figure \ref{bvreg21c4man101c1} shows what happens with two typical
time-courses in data space; the black trace uses devices that do not
meet the requirement as they have $p=0.85<\tau$, while the white one
uses devices that do meet the requirement with $p=0.95>\tau$. Each
time a device is dropped the trace moves either one unit (block) to
the right if it survives, or one block up if it breaks. The posterior
probability $P(A|D)$ is shown in colour, and the testing terminates
when either the time-course enters the light grey lower region, when
it concludes the devices do meet the requirements, or when the
time-course enters the dark grey upper region, when it concludes that
the devices do not meet the requirements.

As can be seen, these two time-courses each reach the correct
conclusion with the regulator's prior or the mixed prior, despite the
manufacturer's prior concluding that both sets of devices meet the
requirements without collecting any data.

\begin{figure}[hpt]
\begin{center}

\includegraphics[scale=0.5]{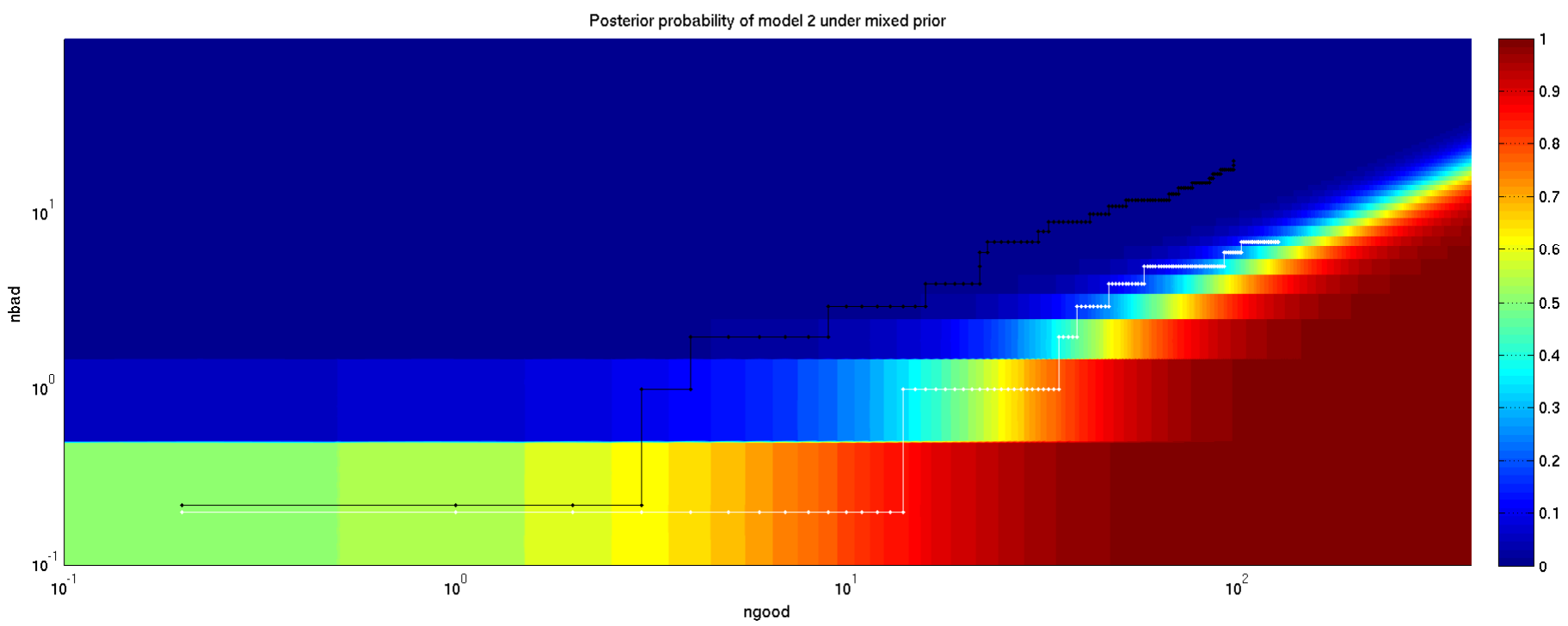}

\caption{Posterior probability that manufacturer's model and priors
  are correct as a function of the numbers of devices that are intact
  $n_\text{good}$ and broken $n_\text{bad}$ using the regulator's
  prior that is Beta(21,4) and the manufacturer's prior that is
  Beta(101,1). The starting point at $(0,0)$ has been marked as
  $(0.2,0.2)$ on these log scales. The white track shows a device with
  $p=0.95>\tau$ and the black track one with $p=0.85<\tau$.
\label{bvmreg21c4man101c1}
}

\end{center}
\end{figure}

Figure \ref{bvmreg21c4man101c1} then shows what happens to $P(M_2|D)$
as the two test runs proceed. 

From a Bayesian point of view the key observation, however, is not the
time courses of these two individual tests, but the underlying colour
plot. As can immediately be seen from the three underlying colour
plots of figure \ref{bvreg21c4man101c1}, the effect of admixing 50\%
of the regulator's prior to the manufacturer's prior is to
dramatically change the inference pattern; in general the less
restrictive (i.e. less informative) prior dominates the inference
pattern in such a mixture.

\subsubsection{Flat prior from regulator, manufacturer calls on past
  experience}

\label{flatregvpast}

Revisiting the results of the previous section
\ref{reasonableregvpast}, we modify the regulator's prior to be flat
on $[0,1]$, while leaving the manufacturer using $$P(p)=101
p^{100}\ [0\leq p \leq 1].$$

If the devices under test are good with in truth $p=0.95$, all of the
$100$ runs drew the correct conclusion after an average of $33$
devices were dropped.

\begin{figure}[hp]
\begin{center}

\includegraphics[scale=0.5]{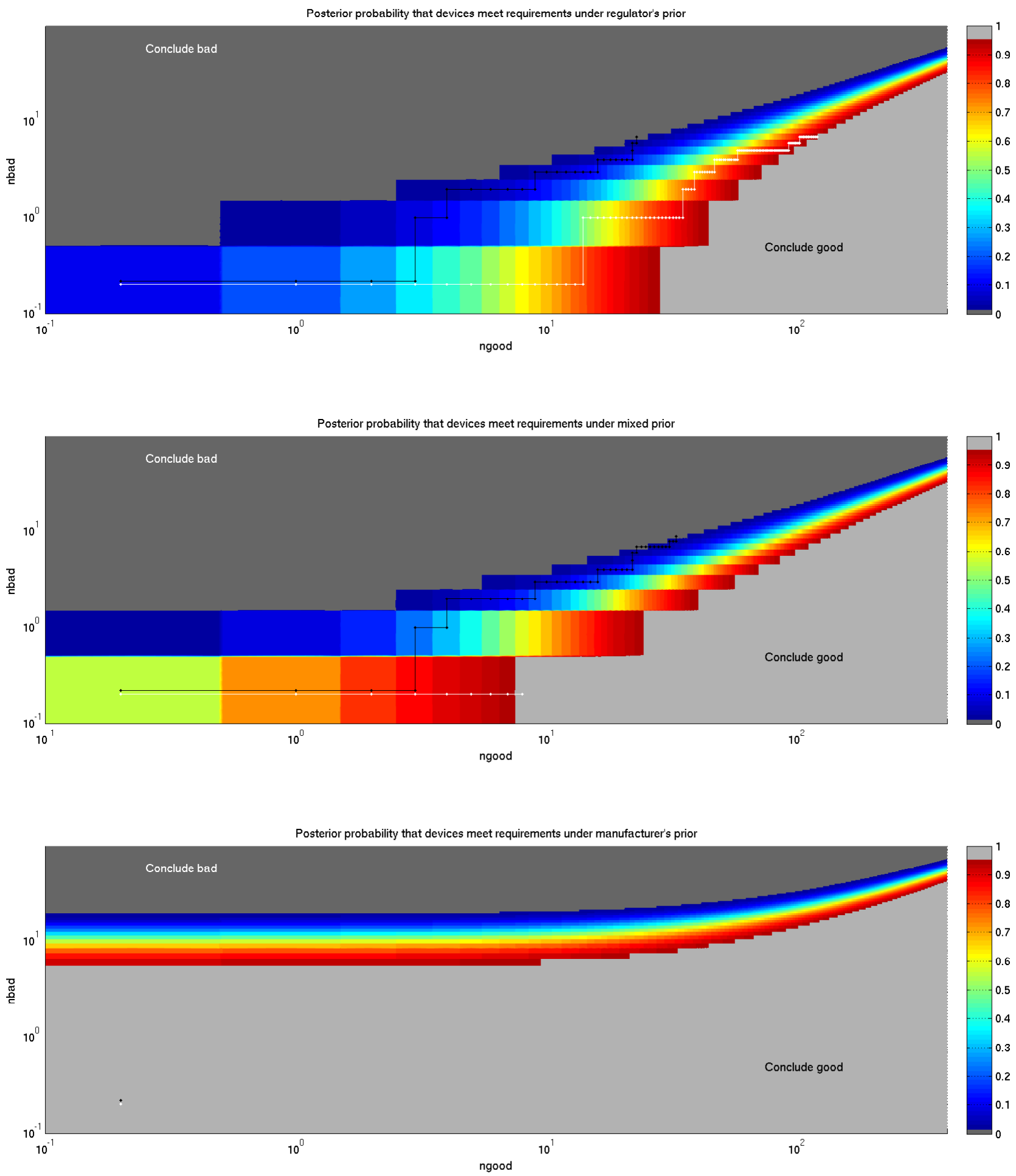}

\caption{Posterior probability that devices meet legal requirements as
  a function of the numbers of devices that are intact $n_\text{good}$
  and broken $n_\text{bad}$ using the regulator's prior that is
  Beta(1,1) (i.e. flat) and the manufacturer's prior that is
  Beta(101,1). The top plot uses the regulator's prior alone, the
  bottom one the manufacturer's prior alone, and the middle one the
  proposed mixture prior. The starting point at $(0,0)$ has been
  marked as $(0.2,0.2)$ on these log scales. The white track shows a
  device with $p=0.95>\tau$ and the black track one with
  $p=0.85<\tau$. With the manufacturer's prior both are deemed good
  without collecting any data. The other two get the same results as
  each other, though it is clear that if there are no or few
  breakages, the mixed prior will conclude sooner than the regulator's
  prior. (In the middle plot, the brief excursion of the black track
  into the dark grey area without terminating is due to inadequate
  resolution in the color scale.)
\label{bvreg1c1man101c1}
}

\end{center}
\end{figure}

\begin{figure}[htp]
\begin{center}

\includegraphics[scale=0.5]{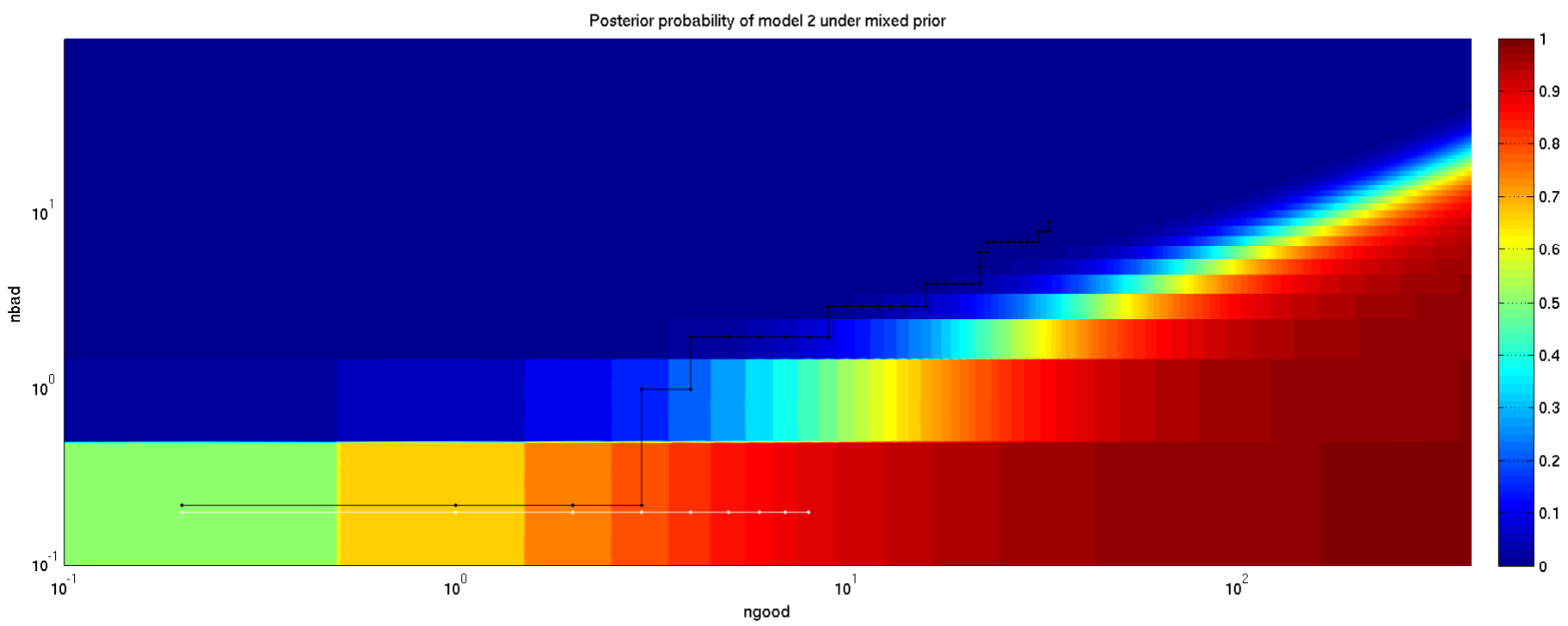}

\caption{Posterior probability that manufacturer's model and priors
  are correct as a function of the numbers of devices that are intact
  $n_\text{good}$ and broken $n_\text{bad}$ using the regulator's
  prior that is Beta(1,1) and the manufacturer's prior that is
  Beta(101,1). The starting point at $(0,0)$ has been marked as
  $(0.2,0.2)$ on these log scales. The white track shows a device with
  $p=0.95>\tau$ and the black track one with $p=0.85<\tau$.
\label{bvmreg1c1man101c1}
}

\end{center}
\end{figure}

The time-course plots for this case are shown in figures
\ref{bvreg1c1man101c1} and \ref{bvmreg1c1man101c1}. As can be seen the
inference pattern is again dominated by the regulator's prior, though
perhaps to a lesser extent than in figure \ref{bvreg21c4man101c1}.

But for a frequentist the point of interest here is what happens to a
factory producing devices with $p=0.85$ i.e. that do not meet the
requirements, and comparing it with the corresponding result in
section \ref{reasonableregvpast}. Of $100$ runs, although $70$ draw
the correct conclusion after an average of $124$ drops, there are now
$30$ that conclude wrongly that the devices do meet the requirements
(with probability $>0.95$) after an average of just $10$ devices are
dropped. Why is this ?

Essentially what is going on here is that the regulator has set a
prior that reckons that a manufacturer aiming to get $90\%$ of devices
to survive a drop unscathed has a probability as high as $0.5$ of
producing devices that only survive $50\%$ of the time. If in truth
the devices have a probability of survival only slightly below the
acceptable value of $0.9$, the regulator's prior is easily
rejected. We will discuss this further below; for now we just note
that the regulator changing his prior to be more in the regulator's
favour has the opposite effect of that desired -- it makes it
\textit{more} likely that a bad device will pass, by making it more
likely that his prior will be rejected by the data.

\subsubsection{Results of runs with true parameters drawn from various
  distributions}

\label{examplesims}

We now investigate what happens using this approach both testing using
a distribution equal to the prior, and testing with various other
distributions. The definitions of the various distributions are in
table \ref{priordefstable}, while the results of the various runs are
in table \ref{examplestable}.

\begin{table}[hp]
\begin{center}
\begin{tabular}{lllrrr}
\textbf{Number} & \textbf{Name} & \textbf{Type} &
\multicolumn{2}{c}{\textbf{Parameters}} & \textbf{Mixture fraction} \\
\hline 
1 & Flat & Beta & 1.0 & 1.0 & 1.00\\
\hline
2 & Somewhat good & Beta & 2.0 & 1.0 & 1.00\\
\hline
3 & Reasonable & Beta & 21.0 & 4.0 & 1.00\\
\hline
4 & Obstructive & Point & 0.89 &  & 0.95\\
& & Point & 0.91 &  & 0.05\\
\hline
5 & Past experience & Beta & 101.0 & 1.0 & 1.00\\
\hline
6 & Uniform prejudiced & Uniform & 0.90 & 1.00 & 1.00\\
\hline
7 & Extremely prejudiced & Uniform & 0.99 & 1.00 & 1.00\\
\hline
8 & Frequentist bad & Point & 0.89 &  & 1.00\\
\hline
9 & Frequentist good & Point & 0.91 &  & 1.00\\
\hline
10 & Perfect & Point & 1.00 &  & 1.00\\
\hline
11 & Catastrophic & Point & 0.00 &  & 1.00\\
\hline

\end{tabular}

\caption{
\label{priordefstable}
Definitions of the various priors and test distributions investigated
in table \ref{examplestable}. The parameters given are $\alpha$ and $\beta$
for Beta distributions, the range for uniform distributions, and the
location of the point mass for point masses. Plots of the various Beta
priors are in figure \ref{showpriors}.}
\end{center}
\end{table}

\begin{figure}[hpt]
\begin{center}

\includegraphics[scale=0.5]{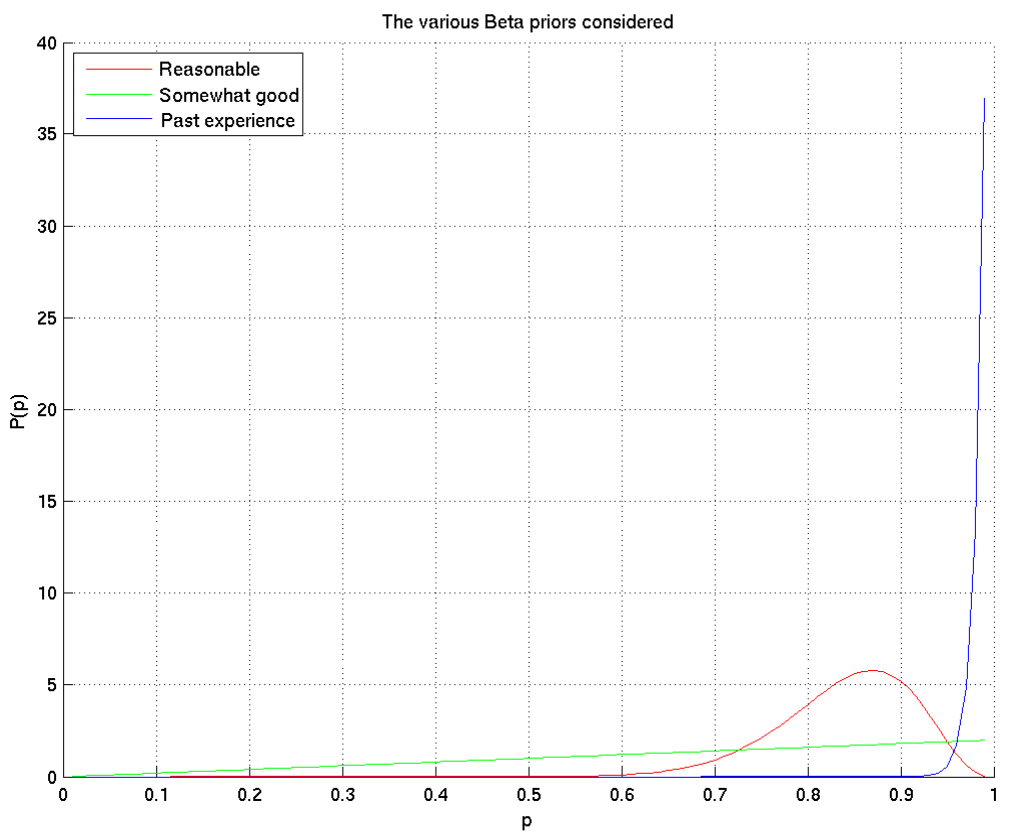}

\caption{
\label{showpriors}
Plots of the various Beta priors defined in table \ref{priordefstable}.
}

\end{center}
\end{figure}

\begin{table}[hp]
\begin{tabular}{llllrrrrr}
\textbf{Test} & \multicolumn{1}{p{2.5 cm}}{\textbf{Prior used}} &
\multicolumn{1}{p{2.5 cm}}{\textbf{Input \mbox{distribution}}} &
\textbf{Truth} & \multicolumn{3}{c}{\textbf{Outcome}} &
\multicolumn{2}{p{2.5 cm}}{\textbf{Mean number of devices dropped}}\\ & Regulator /
Manufacturer & & & \multicolumn{1}{r}{Bad} & \multicolumn{1}{r}{Good}
& \multicolumn{1}{r}{Unkn} & \multicolumn{1}{r}{Bad} &
\multicolumn{1}{r}{Good}\\ \hline 1 & R: Reasonable & Reasonable & Bad: & 339 & 30 & 10 & 307.1 & 139.9 \\
& M: Past experience & Past experience & Good: & 4 & 612 & 5 & 306.8 & 94.2 \\
\hline
2 & R: Reasonable & Reasonable & Bad: & 340 & 33 & 14 & 439.3 & 72.4 \\
& M: Uniform prejudiced & Uniform prejudiced & Good: & 3 & 605 & 5 & 41.7 & 172.2 \\
\hline
3 & R: Reasonable & Reasonable & Bad: & 350 & 26 & 11 & 333.2 & 72.2 \\
& M: Past experience & Uniform prejudiced & Good: & 9 & 594 & 10 & 131.2 & 266.1 \\
\hline
4 & R: Reasonable & Reasonable & Bad: & 344 & 26 & 9 & 403.1 & 94.8 \\
& M: Uniform prejudiced & Past experience & Good: & 1 & 618 & 2 & 8967.0 & 104.3 \\
\hline
5 & R: Flat & Reasonable & Bad: & 284 & 92 & 3 & 198.6 & 14.3 \\
& M: Past experience & Past experience & Good: & 7 & 610 & 4 & 154.7 & 46.7 \\
\hline
6 & R: Reasonable & Flat & Bad: & 431 & 5 & 3 & 97.4 & 17.0 \\
& M: Past experience & Past experience & Good: & 0 & 561 & 0 & NaN & 43.9 \\
\hline
7 & R: Reasonable & Somewhat good & Bad: & 388 & 7 & 1 & 92.9 & 146.6 \\
& M: Past experience & Past experience & Good: & 0 & 601 & 3 & NaN & 49.1 \\
\hline
8 & R: Obstructive & Obstructive & Bad: & 453 & 27 & 9 & 2083.5 & 738.9 \\
& M: Uniform prejudiced & Uniform prejudiced & Good: & 7 & 502 & 2 & 1187.3 & 493.5 \\
\hline
9 & R: Flat & Obstructive & Bad: & 97 & 366 & 26 & 3480.3 & 89.2 \\
& M: Uniform prejudiced & Uniform prejudiced & Good: & 2 & 509 & 0 & 9.5 & 49.1 \\
\hline
10 & R: Obstructive & Obstructive & Bad: & 459 & 30 & 0 & 693.0 & 30.7 \\
& M: Extremely prejudiced & Uniform prejudiced & Good: & 33 & 478 & 0 & 1072.5 & 583.0 \\
\hline
11 & R: Reasonable & Reasonable & Bad: & 362 & 5 & 12 & 337.6 & 832.8 \\
& M: Reasonable & Past experience & Good: & 5 & 610 & 6 & 1665.8 & 154.0 \\
\hline
12 & R: Obstructive & Reasonable & Bad: & 376 & 3 & 0 & 216.1 & 3669.0 \\
& M: Obstructive & Past experience & Good: & 15 & 605 & 1 & 1087.7 & 476.4 \\
\hline
13 & R: Obstructive & Reasonable & Bad: & 370 & 9 & 0 & 214.0 & 779.4 \\
& M: Past experience & Past experience & Good: & 15 & 605 & 1 & 1087.7 & 183.6 \\
\hline
14 & R: Reasonable & Catastrophic & Bad: & 521 & 0 & 0 & 4.0 & NaN \\
& M: Past experience & Perfect & Good: & 0 & 479 & 0 & NaN & 17.0 \\
\hline
15 & R: Reasonable & Frequentist bad & Bad: & 377 & 106 & 38 & 2454.7 & 74.6 \\
& M: Past experience & Frequentist good & Good: & 16 & 459 & 4 & 58.3 & 1353.5 \\
\hline

\end{tabular}

\caption{
\label{examplestable}
A range of simulations done to investigate inference using the
proposed mixtures of regulator's and manufacturer's prior (in equal
parts), and the two listed input distributions (in equal parts). Each
of $1000$ runs were classified using the optimal Bayesian algorithm of
appendix \ref{bayesproofs} item \ref{appunifpow} giving the Outcome,
while the status given by the value of $p$ for the particular factory
tested is as under the Truth column. ``Bad'' denotes $p<\tau$ (if as
an Outcome then with posterior probability $> 0.99$) while ``Good''
denotes $p>\tau$ (and if as an Outcome then with posterior probability
$>0.95$). ``Unknown'' indicates that the trial was stopped after
$10,000$ devices had been dropped with neither of the other possible
conclusions reached. For a full discussion of the results see section
\ref{examplesims}, and for definitions of the priors and input
distributions see table \ref{priordefstable}.}
\end{table}

The first thing to notice is that where the test distribution is the
same as the prior in use (as in tests 1, 2, and 8), the fraction of
runs with Good outcome that should have been classified as Bad is
e.g. $30$ out of $642$ (for test 1), or about $5\%$, in keeping with
the posterior probability required for concluding Good of
$0.95$. Similarly the fraction of runs with Bad outcome that should
have been classified as good is around $1\%$, in keeping with the
$0.01$ posterior probability required before concluding Bad. Note in
particular that we are talking in the first place of $612$ of $642$
Good \textit{outcomes} being correctly classified, not $612$ of $616$
truly Good factories, as a frequentist assessment would report -- we
are interested in what \textit{did} happen in each individual factory
considered alone, not in what \textit{might} have happened when
multiple (e.g. bad) factories considered as a group.

Second, in tests 3 and 4, where the manufacturer mis-guesses either
mildly over-optimistically (test 3) or pessimistically (test 4),
classification performance remains good.

We first notice rather worse performance in test 5 when the regulator
tries too hard to counteract the manufacturer's optimistic prior, and
says that when aiming to design a factory making devices with at least
a $0.9$ probability of surviving being dropped, half of manufacturers
would instead make devices with at least a $0.5$ probability of
breaking. The other way round, though, in tests 6 and 7, where the
prior expects something better than actually turns up, there is no
problem.

Similarly, let us consider tests 8, 9, and 10, in which half of the
input distribution comes from the obstructive distribution and the
other half is uniformly distributed over the good range. In test 8 the
prior is exactly matched to the distribution of factories tested, with
the expected good performance, although the number of devices to be
dropped for each factory becomes rather large (due to the obstructive
nature of this prior, which is designed to make it difficult for the
manufacturer). If the regulator moves his prior yet further to low
values of $p$ as in test 9, the performance he is interested in (what
happens to truly bad factories) gets dramatically worse, while if the
manufacturer tries to boost his own chances by making his own prior
even more prejudiced in his favour (as in test 10) the performance
that he is interested in (what happens to truly good factories) gets
worse, as does his workload.

We think that the moral of the story is not to try to influence the
results by how you set your prior - rather, your prior will be
preferred if it is more realistic, so set it to reflect what you
actually think might be the case. Maybe both parties would then more
easily reach a common mind without resort to this proposed mechanism.

Tests 11 and 12 show what happens without the proposed ``let the data
decide'' mechanism: here we set pure regulator priors that are
unmixed reasonable (test 11) or obstructive (test 12). Comparing test
11 with test 1, the classification performance of test 1 is still up
to specification (5\% and 1\% error rates under good and bad outcomes
respectively), but the manufacturer's testing workload for good
devices is reduced by 39\% compared with test 11. Comparing test 12
with test 13, we see that test 13 with the mixed prior proposed has
just as good classification performance as test 12 with only the
regulator's prior, but the manufacturer's workload is reduced by
62\%. We also note that the regulator's obstructive prior has led to
around a 5\% rate of Bad outcomes that should have been classified as
Good (due to the mismatch with the actual incoming distribution),
whether or not the present proposal is used.

Test 14 shows what happens when either every device survives being
dropped or every device breaks; we get perfect classification, with
respectively 17 or 4 devices needing to be used. The corresponding
numbers for a frequentist test using minimal sample size are 29 and 1,
so on this count also the Bayesian test does better.

But the two odd men out, with spectacularly bad performance, are tests
15 and 9. With test 15, among the 1000 simulated factories being
tested, half \textit{just} failed to meet the requirements with
$p=0.89$ while the other half \textit{just} met the requirements with
$p=0.91$. This is, of course, an extremely unrealistic scenario
(except when doing frequentist testing) -- no manufacturer could
achieve a set of factories that did this if they wanted to. Test 9 has
a very similar issue: here we again have unrealistic spikes in the
input distribution just below and just above $\tau$, and the regulator
has shot himself in the foot by putting in a prior that is far worse
than reality.

\subsubsection{Conclusions from these examples}

We can draw several conclusions from these examples, as follows:

\begin{enumerate}

\item The nature of Bayesian inference is clear from those tests (1,
  2, and 8) where the input distribution exactly matches the prior:
  false positives are $\approx5\%$ of the apparent postives at the
  $95\%$ level and false negatives are $\approx 1\%$ of the apparent
  negatives at the $99\%$ level.

\item When both parties set realistic priors, that good classification
  performance is maintained in the face of mismatch in the form of
  significantly worse input distribution than expected.

\item Important savings of workload for the manufacturer can be
  achieved using the proposed ``let the data decide'' model for
  setting a joint prior for the two parties over simply using the
  regulator's prior.

\item \label{idealaim} That saving is compatible with maintained good
  classification performance if both parties aim to set priors that
  reflect what they really think rather than setting a prior that will
  create prejudice in their own favour. Indeed, as we see from the
  comparisons discussed above, moving a prior away from reality in the
  direction that would be expected to favour the outcome you prefer
  tends to have the opposite effect from that desired.

\item As a result of knowing about point \ref{idealaim} above, there
  is an increased chance of informal discussion reaching an agreed
  prior that is prejudicial to neither side.

\item If a factory does happen to be perfect then that conclusion will
  be reached with $0.95$ posterior probability with fewer devices
  dropped than the frequentist test of minimum sample size.

\item Frequentist testing (with input populations of factories right
  on the borderline of acceptability) does not represent reality, and
  is not an appropriate method of assessing either a Bayesian test
  solution or a real-life test solution.

\end{enumerate}

We note also that thanks to the availability of a uniformly optimal
Bayesian data collection plan a manufacturer:

\begin{itemize}

\item does not need to know exactly how good the devices made by his
  factory are before starting the test; and

\item does not need to adjust the sample size to match the actual
  error rate; and

\item does not run any risk of a truly good factory failing its tests
  through mis-guessing \textit{exactly} how good it is; and

\item does not have to take into account all the other parameters of
  the devices that need to be tested when designing the tests for
  ``this'' parameter. 

\end{itemize}

\section{Can frequentist results be intuitively post-processed to give
  something just as good as Bayesian results ?}
\label{cannotpostfix}

\subsection{Introduction}

We now turn to consider an idea that has been put forward by some
supporters of frequentist methods, who suggest that frequentist
confidence levels are simpler to calculate than Bayesian posterior
probabilities, and can then be intuitively interpreted, perhaps using
intuitively Bayesian thinking, to give results that are just as good
as Bayesian posterior probabilities.

We can easily show that this is not in general the case, although
there is always a particular choice of a nested set of critical
regions for which the resulting frequentist confidence can be so
post-processed\footnote{Such a choice, however, usually requries
  Bayesian calculations to be done first, and does not provide any
  specific algorithm for doing the post-processing.}. Such a proof
rests on two fundamental pillars: the concept of Shannon information,
and the Data Processing Theorem (\cite{DPT}).

\subsection{Intuitive argument}

At one level, it is obvious that we cannot in general recover the
Bayesian posterior probability itself from the frequentist confidence:
compare figures \ref{postprobeven} and \ref{freqxonly}; for any one
value of frequentist confidence in figure \ref{freqxonly} there are
many possible values of posterior probability that may correspond in
figure \ref{postprobeven}. (On the other hand comparing figures
\ref{postprobeven} and \ref{unifopt} we see that sometimes we are able
to recover posterior probability from frequentist confidence, at least
if we exclude some things that happen only with total probability
zero.)

However, something much deeper is true: in general, the frequentist
confidence may tell us strictly less about the unknowns than does the
Bayesian posterior probability, although there are some cases where
both carry the same amount of information about the unknowns.

Intuitively, the argument goes like this (a detailed precise version
follows in subsequent subsections): We want to know $\theta$, and Shannon
has enabled us to define and calculate the amount of information about
$\theta$ contained in the observed data $x$ (or in any other variable). The
posterior probability distribution $P(\theta|x)$ can be shown to contain
exactly the same amount of information about $\theta$ as $x$ does. The Data
Processing Theorem tells us that any other variable that is a function
of $x$ contains no more information about $\theta$ than $x$ does; so
frequentist confidence cannot contain any \textit{more} information
about $\theta$ than the Bayesian posterior probability, and \textit{nor can
  anything deduced from the frequentist confidence}. Moreover, any
function $z(x)$, for which the probability that $P(\theta|z(x))$ is not
equal to $P(\theta|x)$ is greater than zero, contains strictly less
information about $\theta$ than $x$ does. Since this applies to the
frequentist confidence calculated from many choices of nested sets of
critical regions, not only may frequentist confidence contain strictly
less information about $\theta$ than the Bayesian posterior, but there is
nothing one can then post-process out of the frequentist confidence
that does any better -- and in particular one cannot in general
post-process frequentist confidence into Bayesian posterior
probability.

Indeed, the argument can be taken further. If we are given a
probability distribution $Q(\theta|x)$ on $\Theta$, that varies with $x$, we can
define the Apparent Shannon Information (ASI) in $Q$ about $\theta$, which
measures what $Q$ tells us about $\theta$ \textit{without} any further
(intuitive or precise) post-processing. We can show that the ASI never
exceeds the (true) Shannon information in $x$ about $\theta$, and that if
we treat frequentist confidence $c$ as $Q(H_1|x)$ (setting
$Q(H_0|x)=1-c$) then in many cases the ASI is even negative, showing
that frequentist confidence can be downright misleading.

The following subsections make this argument precise.

\subsection{Shannon information}
\label{Shannon}

Suppose we have two random variables $\theta$ and $x$. Then the (true)
Shannon information (TSI) contained in $x$ about $\theta$ is defined
by $$I(\theta;x)=\E_{P(x,\theta)}{\log
  \frac{P(\theta|x)}{P(\theta)}},$$ where $\E_{P(x,\theta)}$ denotes
the expectation using the probability measure $P$ on the random
variables $x$ and $\theta$. If $x$ and $\theta$ are continuous random
variables with density functions (also denoted by $P$) this can also
be written
as $$I(\theta;x)=\int{P(x,\theta)\log\frac{P(\theta|x)}{P(\theta)}\,d(x,\theta)},$$
and if they are discrete random variables then the integral can be
replaced with a summation. In both cases the range of integration or
summation is over all the possible values of $\theta$ and $x$.

If the logarithm is taken to base 2 then $I(\theta;x)$ will be measured in
bits; if taken to base $e$ then $I(\theta;x)$ will be measured in ``nats''
or ``nepers''.

To see that Shannon information is a sensible way of measuring
information content the reader is invited to consider the case that
$\theta$ is a uniform random number from $[0,1]$ that we do not know, and
that $x$ is a statement by a somewhat dishonest politician who does
know $\theta$ telling us whether it is bigger or smaller than
$\frac{1}{2}$, or whether it is bigger or smaller than $\frac{1}{4}$,
etc. 

\subsection{Basic facts relevant to Shannon information}
\label{basicSi}

If $P$ and $Q$ are two probability measures on a random variable $x$,
a simple application of Jensen's inequality to the concave function
$\log$ tells us that\footnote{Here $P(x)$ and $Q(x)$ may be both
  probabilities, both probability densities with respect to the same
  underlying measure, or $\frac{Q(x)}{P(x)}$ may be a Radon-Nikodym
  derivative of the absolutely continuous part of $Q$ with respect to
  $P$ (measurable on the same $\sigma$-algebra as $P$) evaluated at
  $x$. Any set on which $P$ is zero contributes nothing to the
  expectation. If there is any set of positive $P$-probability but
  zero $Q$-probability then the expectation is $-\infty$. Similar
  considerations apply with opposite sign to the following displayed
  equation.} $$\E_{P(x)}{\log\frac{Q(x)}{P(x)}}\leq 0$$ with equality
if and only if $P$ and $Q$ are equal except on a set of
$P$-probability zero. In
consequence $$\E_{P(x)}{\log\frac{P(x)}{Q(x)}}\geq 0.$$
Since $$\E_{P(x,\theta)}=\E_{P(x)}\E_{P(\theta|x)},$$ we can also
deduce that $I(\theta;x)\geq 0$, and by subtraction that for any
probability distribution $Q(\theta|x)$ on $\theta$ that varies with
$x$, $$\E_{P(x,\theta)}{\log\frac{P(\theta|x)}{P(\theta)}}\geq
\E_{P(x,\theta)}{\log\frac{Q(\theta|x)}{P(\theta)}},$$ and hence that
for any random variable $z$ that is a function of
$x$, $$\E_{P(x,\theta)}{\log\frac{P(\theta|x)}{P(\theta)}}\geq
\E_{P(x,\theta)}{\log\frac{P(\theta|z(x))}{P(\theta)}}=\E_{P(z,\theta)}{\log\frac{P(\theta|z)}{P(\theta)}},$$
i.e. $$I(\theta;x)\geq I(\theta; z),$$ which is the Data Processing
Theorem.

Further, since\footnote{To see this, let $A$ be any measurable subset
of $\Theta$, and denote the random variable $P(\theta|x)$ (measurable
on the $\sigma$-field generated by $x$ and taking values in the space
of probability measures on $\Theta$) by $Q$, and $P(\theta\in A|x)$ by
$R$. Then 
\begin{IEEEeqnarray*}{rClr} P(\theta\in A|P(\theta|x)) &=&
  P(\theta\in A|Q)\\
 &=& \int{P(\theta\in A,x|Q)\,dx} & \text{(marginalisation)}\\
 &=& \int{P(x|Q)P(\theta\in A|x,Q)\,dx} & \text{(chain rule)}\\
 &=& \int{P(x|Q)P(\theta\in A|x)\,dx} & \text{($x$ known implies $Q$ known)}\\
 &=& E(P(\theta\in A|x)|Q) & \text{(definition of conditional expectation)}\\
 &=& E(P(\theta\in A|x)|P(\theta|x)) & \text{(definition of $Q$)}\\
 &=& E(P(\theta\in A|x)|P(\theta\in A|x)) & \text{(rest of $Q$ makes no
    difference)}\\
 &=& E(R|R) & \text{(definition of $R$)}\\
 &=& R & \text{(expectation of a r.v. given itself is itself)}\\
 &=& P(\theta\in A|x) & \text{(definition of $R$)}.\\
\end{IEEEeqnarray*}} $P(\theta|P(\theta|x)) = P(\theta|x)$, we note that the Shannon
information about $\theta$ contained in the Bayesian posterior is equal to
that contained in the observed value $x$, i.e. nothing has been lost,
thus proving the ``Information Optimality'' criterion.

We note also that $I(\theta;x)=0$ if $\theta$ and $x$ are independent, since
then $P(\theta|x)=P(\theta)$. By the initial comments about when Jensen's
inequality gives an equality, we note that this only happens if
$P(\theta|x)=P(\theta)$ except perhaps on a set of $P$-probability zero; but
then $x$ and $\theta$ are independent. Thus $x$ and $\theta$ being independent
is equivalent to there being no Shannon information about $\theta$ in $x$.

\subsection{Apparent Shannon information (ASI)}
\label{ASI}

For two variables, such as $\theta$ and $x$, the Shannon information about
$\theta$ in $x$ tells us how much we can learn about $\theta$ by applying
optimal data processing methods to $x$. But if instead we consider
variables such as $\theta$ and some arbitrary probability distribution
$Q(\theta|x)$ on $\theta$ which varies with $x$, we may also want to know how
much $Q$ tells us about $\theta$ if taken at face value without further
processing.

For this we define the Apparent Shannon Information (ASI) by $$J(\theta;Q)
= \E_{P(x,\theta)}{\log\frac{Q(\theta|x)}{P(\theta)}}.$$ (In contrast the (true)
Shannon information in $Q$ about $\theta$ would be defined by $$I(\theta;Q) =
\E_{P(x,\theta)}{\log\frac{P(\theta|Q(\theta|x))}{P(\theta)}}.)$$

Using the same arguments as in section \ref{basicSi} we find that ASI
is maximised by setting $Q(\theta|x)=P(\theta|x)$, the Bayesian posterior.

\subsection{Application: The frequentist confidence does
  not in general contain sufficient information to be able to recover
  the Bayesian posterior from it}

Now consider figure \ref{freqxonly}. It is clear that the frequentist
confidence (in the case of this particular choice of nested critical
regions) is a function of the $x$-coordinate of where the bullet
lands, and that this function is injective. Accordingly the
$x$-coordinate (or the fact that the bullet didn't land) can be
recovered from the frequentist confidence as well as vice versa, and
hence by the data processing theorem the information contained in the
frequentist confidence about $h$ is equal to that contained in the
$x$-coordinate, namely 0.233 bits (as discussed in section
\ref{commonpartsall}). But the information contained in the Bayesian
posterior (making use of both coordinates) is the same as that
contained in the $x$ and $y$ coordinates together, or about 0.258
bits. Thus by the data processing theorem there is no way of
recovering either the Bayesian posterior or anything containing the
same amount of information from the frequentist confidence (for this
particular choice of critical regions, and hence in general).

\subsection{Special cases where the frequentist confidence does
  contain as much information as the Bayesian posterior}

Now consider figure \ref{freqpseudoB} and compare it with figure
\ref{postprobeven}. In this case the Bayesian posterior probability is
a bijective function of the frequentist confidence, and therefore the
Shannon information about $h$ contained in the two are equal. But this
is, of course, a special case - and it is certainly \textit{not}
simpler to use the pseudo-Bayesian critical regions than to use the
Bayesian method.

Moreover, if we interpret the frequentist confidence here as a
probability that $H_1$ holds, we can assess the Apparent Shannon
Information (ASI) in the frequentist confidence and compare it with
that of the Bayesian posterior (the latter being 0.258 bits). The ASI
in the frequentist confidence is defined
as $$\E_{P(x,h)}{\log\frac{c(h|x)}{P(h)}},$$ where for $h=1$, $c(h|x)$
denotes the frequentist confidence that $H_1$ holds given $x$ has been
observed, and for $h=0$ it denotes one minus that value. For these
nearly pseudo-Bayesian critical regions, the ASI in the frequentist
confidence turns out to be $-1.69$ bits -- in other words the
frequentist confidence is downright misleading as it stands (despite
the fact that it can be post-processed to give the Bayesian posterior
probability) -- one would be better off ignoring the frequentist
confidence and sticking with the (even) prior probability.

But if you need to use pseudo-Bayesian critical regions to get a
version of frequentist confidence from which it is possible to recover
the Bayesian posterior, why not just calculate the Bayesian posterior
in the first place ? It would certainly be simpler to do so.

Moreover, we should caution the reader that there are examples in
which a basic pseudo-Bayesian frequentist confidence \textit{cannot}
be post-processed to give the Bayesian posterior, and indeed contains
strictly less Shannon information about $\theta$. Indeed in the next
section we give an example of a problem where deterministic
frequentist solutions only ever give zero information about $\theta$.

\subsection{An inference problem on which deterministic frequentist 
solutions \mbox{provide} no information at all}
\label{zeroinfo}

We return to using the notation of section \ref{infprobdef}.

Let $$\Phi=\{0,1,2\},$$ $$\Theta=\{0,1\},$$ $$H=\{(0,0), (0,1),
(1,2)\}\subset
\Theta\times\Phi,$$ $$X=\{0,1\},$$ $$P(x|h)=P(x|\theta,\phi)=
                [x=[\theta\neq\phi]]$$ (where $[\ ]$ denotes the
                function which takes the value 1 when the expression
                inside the $[\ ]$ is true and 0 otherwise),
                $$H_0=\{(0,0),(0,1)\}=(\{0\}\times \Phi)\cap H =
                \{(0,\phi)\in H\},$$ and $$H_1=\{(1,2)\}.$$

Putting that in words, there are three possible underlying states A =
$(0,0)$, B = $(0,1)$, and C = $(1,2)$, of which A and B are in $H_0$;
the observed data $x$ is always 1 if B or C is the underlying state,
and always 0 if A is the underlying state.

First consider the Bayesian position. We assume that the Bayesian
chooses an arbitrary prior on $H$ which puts non-zero probability on
each of the three possible underlying states. Then if $x=0$ is
observed, we get $P((0,0)|x) = 1$, i.e. state A becomes
certain\footnote{Technically ``almost certain'', i.e. has probability
  1.}, as does $\theta=0$, which they weren't before. Thus $h$ and $x$
are not independent, and therefore $I(h;x)>0$ by the remarks of
section \ref{basicSi} above. Similarly $\theta$ and $x$ are not
independent, and therefore $I(\theta;x)>0$.

In particular, to give a specific example, suppose that the Bayesian
prior is $\frac{1}{3}$ on each of the three possibilities in $H$,
i.e. the Bayesian considers all three underlying states equally likely
\textit{a priori}. Then, considering the possible states
$(0,0),(0,1),(1,2)$ of $h$ and their associated (only possible) $x$
values of $0,1,1$ in order, we have $$I(\theta; x) =
\frac{1}{3}\log\frac{1}{2/3} + \frac{1}{3}\log\frac{1/2}{2/3} +
\frac{1}{3}\log\frac{1/2}{1/3} = \frac{1}{3}\log\frac{27}{16},$$ which
works out at 0.174 nats or 0.252 bits. Moreover if different priors
are in force, the information about $\theta$ contained in the Bayesian
posterior can be anything up to 1 bit (which is the most you can ever
expect about a variable which has only two possible values).

But now consider what nested sets of critical regions the frequentist
could deploy. An empty critical region is no use, as $x$ can never end
up in it. The critical region $\{0\}$ has probability $1$ of $x$ being
in it when $h=(0,0)\in H_0$, so it would only give frequentist
confidence of zero. Similarly $\{1\}$ has probability $1$ of $x$ being
in it when $h=(0,1)\in H_0$, so too would only give frequentist
confidence of zero. The same is true for the critical region
consisting of the whole of $X$ -- and there are no other subsets of
$X$ to consider.

Thus whatever nested set of critical regions the frequentist might
deploy, he will always end up with frequentist confidence $c=0$ that
$H_1$ holds. Then $I(\theta;c)=0$, i.e. the frequentist confidence
contains no information at all about $\theta$, our wanted variable,
whichever choice of nested critical regions we adopt. Moreover
$J(\theta;c)=-\infty$ because $h=(1,2)$ occurs with non-zero
probability, and then the frequentist method always gives zero
confidence that $h\in H_1$ when actually (with the even prior above)
the data allows us to have posterior probability of $\frac{1}{2}$ that
indeed $h=(1,2)$ and $\theta=1$. Nonetheless the Bayesian posterior
does contain a positive amount of information about $\theta$. So for
this problem, no matter what choice of critical regions is made, the
Bayesian posterior can never be deduced from the frequentist
confidence when this is obtained from deterministic critical regions.

\subsection{What about non-deterministic critical regions ?}

However, the reader may wonder why we have not mentioned
non-deterministic critical regions here. It turns out that we can
always make a non-deterministic nested set of critical regions such
that the resulting frequentist confidence contains all the available
information in the Bayesian posterior about whether $H_0$ or $H_1$ is
true. Moreover, we can do so in such a way that the set of critical
regions in question are arbitrarily close to a set of basic pseudo-Bayesian
critical regions in an appropriate sense. Nonetheless, as we shall
see, there is really no point using such critical regions -- it is far
easier to just use the Bayesian method. Moreover designing such
critical regions requires use of the Bayesian prior, thus negating one
of the ``advantages'' that frequentists see in avoiding Bayes.

Take then any hypothesis testing problem, so that $\Theta=\{0,1\}$ and
$H_0=\{(\theta,\phi)\in H:\theta=0\}$, i.e. $\theta$ is either 0 or 1
and tells us whether $H_0$ or $H_1$ holds.

Now choose any $\alpha\in (0,1]$; the smaller the choice of $\alpha$,
  the closer our constructed set of critical regions will be to a
   pseudo-Bayesian set.

If we now consider the basic pseudo-Bayesian critical regions (for
some particular prior and data collection plan) defined
by $$D_p=\{x\in X : P(H_1|x)\geq p\},$$ then we may extend the data
space from $X$ to $X\times [0,1]$ with $P(x,u|h)=P(x|h)[u\in [0,1]]$,
and set $$C_{\eta(p)}=D_p\times [\alpha p,1]$$ for some
$\eta:[0,1]\to[0,1]$ a non-strictly increasing function such that
$\eta(p)$ is the frequentist confidence corresponding to $C_{\eta(p)}$
that $H_1$ holds. Then $D_p$ and $C_{\eta(p)}$ are decreasing as $p$
increases. Moreover it is clear that as $\alpha\to 0$, $C_{\eta(p)}
\to D_p\times (0,1]$ (and, in an appropriate sense, uniformly in $p$).

Now let $p$, except when acting as a dummy variable, denote the random
variable $P(h\in H_1|x)$.

If $\eta(p)$ is not to contain all the information about $\theta$
contained in $p$, we need $\eta$ to not have a left inverse, hence not
be injective, on any set of probability $1$. We will show that there
is in fact some set of zero probability off which $\eta$ is strictly
increasing and injective, whence the information content
$I(\theta;\eta(p))$ is equal to $I(\theta;p)$ as required.

By the definitions of $C$ and $\eta$, we have, for all $p\in
[0,1]$, $$\eta(p) = 1 - \sup_{h\in H_0}{((1-\alpha p)P(x\in D_p|h))} =
1 - (1-\alpha p) \sup_{h\in H_0}{P(x\in D_p|h)},$$ which is strictly
increasing so long as $$\sup_{h\in H_0}{P(x\in D_p|h)>0}$$ for all
$p<1$. Suppose then that this is not the case, and note that this
quantity is a non-strictly decreasing function of $p$. Then
let $$p_0=\sup{( \{p\in [0,1]:\sup_{h\in H_0}{P(x\in D_p|h)>0}\}\cup\{0\})},$$ and
define $$R=\left\{\begin{matrix}(p_0,1) & (\sup_{h\in H_0}{P(x\in
  D_{p_0}|h)>0})\\ [p_0,1) & (\sup_{h\in H_0}{P(x\in
    D_{p_0}|h)=0})\end{matrix}\right.$$ and $$D_R=\bigcup_{p\in
  R}{D_p}.$$ Then $\eta(p)$ is strictly increasing when $p\notin R$,
i.e. off $D_R\setminus D_1$, so it remains only to show that
$P(D_R\setminus D_1)=0$.

To avoid measure-theoretic complications we will assume that $x$ and
$h$ are variables having density functions with respect to an
underlying measure on $H\times X$ (this indeed also addresses the
general case since we may take the underlying measure to be $P$ and
the density function $P(h,x)=1$ everywhere).

Choose then any $p\in R$. Then for all $h\in H_0$, $$0=P(x\in D_p|h)=
\int_{D_p}{P(x|h)dx}.$$ Also, for all $x\in D_p$, $$P(H_1|x)P(x) =
P(H_1,x)$$ and hence $$P(H_1|x)\int_H{P(h)P(x|h)dh} =
\int_{H_1}{P(h)P(x|h)dh}$$ and
therefore $$\int_{D_p}{P(H_1|x)\int_H{P(h)P(x|h)dh}\,dx} =
\int_{D_p}{\int_{H_1}{P(h)P(x|h)dh}\,dx}.$$

Since the integrands are non-negative we may apply Tonelli's theorem
to get $$A:=\int_{H_0\cup H_1}{P(h)\int_{D_p}{P(H_1|x)P(x|h)dx}\,dh} =
\int_{H_1}{P(h)\int_{D_p}{P(x|h)dx}\,dh}=:B.$$

Since $P(H_1|x)\leq 1$, we note that $$A\leq
\int_{H_1}{P(h)\int_{D_p}{P(H_1|x)P(x|h)dx}\,dh} +
\int_{H_0}{P(h)\int_{D_p}{P(x|h)dx}\,dh}.$$ But the inner integral in
the last term is just $P(x\in D_p|h)$ for some $h\in H_0$, which we
know is zero. Therefore $$A\leq
\int_{H_1}{P(h)\int_{D_p}{P(H_1|x)P(x|h)dx}\,dh} \leq B = A,$$ so both
inequalities are equalities,
and $$\int_{H_1}{P(h)\int_{D_p}{P(H_1|x)P(x|h)dx}\,dh} =
\int_{H_1}{P(h)\int_{D_p}{P(x|h)dx}\,dh}.$$ Therefore $$\E_{(x,h)\in
  D_p\times H_1}{P(H_1|x)} = 1.$$ But we already know
that $$P((x,h)\in D_p\times H_0) = 0,$$ so actually $$\E_{x\in
  D_p}{P(H_1|x)}=1,$$ or in other words for almost all $x \in D_p$,
$P(H_1|x)=1$, i.e. $x\in D_1$, so $P(D_p\setminus D_1)=0$. But
$D_R\setminus D_1$ is a countable union of sets of the form $D_p
\setminus D_1$ for various $p\in R$, so also $P(D_R\setminus D_1)=0$
as required, whence frequentist confidence derived from the
non-deterministic nested set of critical regions $C_{\eta(p)}$ carries
all the information in the Bayesian posterior about $\theta$.

Translating this general method onto the specific problem in section
\ref{zeroinfo} above, the solution arrived at (for the case
$\alpha=1$) is equivalent to the following recipe for calculating
frequentist confidence: \begin{quote}If $x=0$ set $c=0$ that $H_1$
  holds, otherwise pick a random number $u$ uniformly from $[0,1]$ and
  set $c=u$.\end{quote} It is then clear that $P(H_1|x)$ can be
recovered from $c$ with probability 1 by the following recipe
(assuming the prior of $\frac{1}{3}$ on each of the three possible
values of $h$): \begin{quote} If $c=0$ set $P(H_1|x)=0$, otherwise
  set $P(H_1|x)=\frac{1}{2}$.\end{quote} Nonetheless, the ASI about
$\theta$ in this non-deterministic $c$
is $$J(\theta;c)=\frac{1}{3}\log\frac{1}{2/3} +
\frac{1}{3}\int_0^1{\log\frac{u}{2/3}\,du} +
\frac{1}{3}\int_0^1{\log\frac{u}{1/3}\,du}$$ which works out at
$-0.030$ nats or $-0.044$ bits, showing that unless one does such
post-processing the frequentist confidence itself, even in this case,
is more misleading than just taking the prior.

But why would one bother to first choose a Bayesian prior, then use it
to calculate the Bayesian posterior in order to set up such a
complicated nested set of non-deterministic critical regions to obtain
a version of frequentist confidence from which the Bayesian posterior
could be recovered by some procedure yet to be determined ? We can see
no reason to -- it is, after all, so much easier just to calculate and
use the Bayesian posterior in the first place.

\section{Discussion of other counter-arguments}
\label{counterarguments}

At this point it seems worth listing and commenting on the
counter-arguments in favour of frequentist hypothesis testing and
confidence sets that we have encountered.

\subsection{The Bayesian method is equivalent to the frequentist
  method}

This one is so obviously false that it is surprising that anybody
would say it. What is true is that for problems satisfying certain
nice conditions (which many do, but not all), as the amount of data
approaches infinity the difference in the conclusions of the two
methods approaches zero in some sense. But we have yet to meet a
real-life inference problem in which the amount of data available
approaches infinity -- and the results reported in table
\ref{statresults} make clear that even in real-life problems for which
the said ``certain nice conditions'' \textit{do} apply, frequentist
and Bayesian methods behave very differently in important ways even
when the amount of data is in the thousands or millions of items.

\subsection{Bayesian results depend on the choice of prior}

This is true -- and so they should, as can clearly be seen from the
case where the quantity of data is zero, or where the data is
independent of the wanted unknowns, when what is known about the
wanted unknowns afterwards is exactly the same as what is known
before.

\subsection{Bayesian methods don't control the type I error rate}

This is also true. However:
\begin{enumerate}
\item as we have seen in the example of section \ref{example1}
  (e.g. in figure \ref{freqct1} or in table \ref{resultstable}),
  control of the type I error rate doesn't actually guarantee that
  type I errors will be uncommon in any real-world sense; and
\item as we have seen in the examples of section \ref{examplesims},
  frequentist testing of type I error rate is highly unrepresentative
  of real-world test scenarios; and
\item controlling the type I error rate, at least in the frequentist
  ways discussed here, leads to the resulting frequentist confidence
  not behaving like what most people understand by ``confidence'' at
  all, as seen in section \ref{criteria}, however much frequentists
  might like us to think that ``confidence'' is a good way to describe
  this concept.
\end{enumerate}

\subsection{Choosing a good critical region involves more than just
  controlling the type I error rate}

Various other desirable attributes of critical regions beyond those
given in section \ref{freqmeth} are sometimes listed, e.g.
\begin{enumerate}
\item \label{morelikely} it should be more likely that $x\in C_\eta$
  if $h\in H_1$ than if $h\in H_0$;
\item \label{highposs} $P(x \in C_\eta|h\in H_1)$ should be as high as possible;
\item \label{wantunifopt} $C_\eta$ should be uniformly optimal.
\end{enumerate}

However of these: 
\begin{itemize}

\item \ref{morelikely} is unclear unless $H_0$ and $H_1$
are singleton sets. Do we mean that $$\forall h_0\in H_0,\forall h_1\in H_1,
P(x\in C_\eta|h_1)>P(x\in C_\eta|h_0),$$ or that $$P(x\in C_\eta|h\in
H_1)>P(x\in C_\eta|h\in H_0) ?$$ The former may not be achievable, and
the latter depends on the prior distribution of $h$, something
frequentists don't want to consider; 

\item moreover \ref{morelikely} is achieved in all of figures
  \ref{freqmaxarea}, \ref{freqct1}, \ref{freqxonly},
  \ref{freqdirection}, \ref{freqpseudoBayesfig}, and it doesn't help
  us choose between them;

\item similarly \ref{highposs} depends on the prior on $h$;

\item and for most problems \ref{wantunifopt} is unachievable, while
  even for those for which it is achievable, the resulting solution is
  often nonsense (e.g. figures \ref{unifopt} and \ref{freqpseudoB}).
\end{itemize}

\subsection{Frequentist methods are easier to understand}

We believe that this depends entirely on whether one was first taught
frequentist statistics or Bayesian inference. The na\"ive novice
learning frequentist methods for the first time usually starts with a
natural intuition that is based around Bayesian posterior
probabilities and has to learn to invert their thinking to understand
frequentist methods. Such a person then understandably finds it hard
to revert to their original way of thinking after many years of
frequentism. On the other hand those who start by learning Bayesian
methods rarely have any difficulty understanding the basic principles.

\subsection{Bayesian solutions are harder to compute}

Or alternatively ``Frequentist solutions are easier to compute'' --
indeed so, but often because they are using approximations that are
only valid as the number of independent and identically distributed
data samples approaches infinity (which never happens in real
life). In most non-Gaussian, non-uniform cases accurate calculation of
frequentist confidence is every bit as difficult or more than
calculation of posterior probability -- consider, for example, the
calculations involved in producing figure \ref{unifopt}.

\subsection{Frequentist methods give more opportunities for research}

Again, this is true -- but we do not believe that such opportunities
are valuable contributions to knowledge of inference methods in the
real world, even though they may be good mental exercise.

As we have seen in the example of section \ref{example1}, frequentist
hypothesis testing and confidence sets require arbitrary choices of
(nested families of) critical regions to be made, and require checks
to be made on the chosen critical regions that they do indeed satisfy
the defining conditions. In 2 dimensions there is a huge choice of
potential critical regions to work on; in hundreds of dimensions there
are even more (in a non-technical sense), so there is an enormous
amount of work that could be done.

But the \textit{right} answer to the actual real-life problem is
already known -- it is the Bayesian answer, not any of the
many possible frequentist wrong answers.

Moreover to see that plenty of work can be done on Bayesian methods,
some of it hugely important, the reader is referred to
\cite{ProppWilson,JamesFill,MurdochandGreen,Green98exactsampling},
later work by these authors, and to the numerous recent preprints by
e.g. Amos Storkey and various co-authors on
\url{https://www.arxiv.org}.

\subsection{Regulators require the use of frequentist methods}

Yes, at present they do -- and in some cases they even try to
obfuscate the issue (see \cite{FDApseudoBayes}). But colluding with
those in power when they are wrong advances science not one whit. This
needs to change.

\subsection{Frequentist methods provide incentives for manufacturers
  to produce better products}

In the context of approvals testing of a factory, it is sometimes
claimed by regulators that insisting on frequentist hypothesis tests
provides an incentive for manufacturers to design better products, as
they are then more likely to pass approvals tests than if they only
just satisfy the approvals criteria.

Let us first point out that Bayesian testing also leads to less
testing when each $h_n$ is substantially greater than $h_0$ than when
it is only just above $h_0$ -- so the incentive is not lost with
Bayes. 

Then, let us consider the Bernoulli $N=224$ scenario of table
\ref{statresults} as discussed in section \ref{statresdisc}. Here a
manufacturer forced to do frequentist testing has done his best to
ensure that his factory satisfies all the requirements and wants his
factory, providing it is good, to pass testing with probability at
least $0.99$. However it is hard for him to be sure that his factory
is good before doing the potentially very costly approvals tests;
indeed his preliminary Bayesian scouting run takes nearly ninety
thousand devices to confirm that of the $224$ aspects less than one in
two thousand are likely to be below the quality threshold of
$h_0=0.975$. But in then designing his frequentist approval test he is
forced to do between six million and 200 million individual drops,
each on a new device, potentially using up more than a year's worth of
production. This is unreasonable, given that Bayesian testing will
achieve the goals of \textit{both} the scouting run and the final
approval test in less than sixty thousand drops, and will
automatically adjust the number of tests on each of the $N$ aspects
according to how far that particular $h_n$ is above $h_0$.

What the manufacturer using frequentist testing loses if he goes for
only a $0.9$ probability of passing a good factory, or if he omits the
scouting run, is the risk of losing the factory altogether -- for
frequentist tests cannot just be repeated once all the type I error
probability has been used up. This is usually an unacceptable risk for
such a manufacturer. So what actually risks happening is not that the
manufacturer executes this excessively demanding frequentist test, but
that he decides to assume that his factory is better than it is (or
even that it is perfect) in order to design a less demanding test,
tacitly planning to resort to cheating\footnote{There are many ways,
  e.g. repeating the test of an aspect that fails hoping it won't next
  time, trying a different frequentist test hoping it will give higher
  confidence, claiming that a particular failed drop wasn't done
  properly so replacing it with one that works, collecting more data
  and pretending that one always planned to, discarding results one
  doesn't like as ``outliers'', or even just testing the whole factory
  again from scratch. An even more inventive way is reported in
  \cite{CTM} section 5.3.1 .}\cite{CTM} if in consequence some aspect
of the test fails. Indeed, the need to cheat risks becoming so
commonplace that it is considered totally normal.

So, far from making it more likely that manufacturers will be
motivated to make near-perfect product, insistence on frequentist
testing is likely to result in undetected cheating and complete
invalidity of the approvals testing done.

\subsection{``I define a `good' frequentist solution to be one that
  agrees with Bayes...''}

``I define a `good' frequentist solution to be one that gives the same
answer as an appropriate Bayesian solution; then I'm happy to use
\textit{good} frequentist solutions.'': 
\begin{itemize}
\item How would you know that the `good' frequentist solution agrees
  with Bayes unless you've run the Bayesian solution ?
\item And if you have, why not then just report the Bayesian solution
  and forget the frequentist one ? 
\item And if you're going to do that, why not skip the obfuscation and
  just say that the frequentist way of designing inference methods is
  wrong ?
\end{itemize}

\subsection{``I refuse to engage in this discussion''}

Very sadly, this is by far the most common response to attempts to
debate these matters. Refusal by those currently in the majority to
discuss where they might be wrong is not helpful, and nor is running a
truce between right and wrong.

\section{Conclusion}
\label{conclusion}

First, it is clear that the pseudo-Bayesian method cannot be better
than the best frequentist methods, as it is itself a frequentist
method. However, whether the Bayesian method is (or is not) better
than the frequentist method is a totally different question.

Failure to adhere to any one of the criteria of section \ref{criteria}
above would imply that the concept of being ``$\eta$-sure'' does not
match the real-world understanding of such an expression; this is
especially so of criterion \ref{containment}. We have already seen
that the Bayesian method adheres to all these criteria, while the
frequentist method violates them all.

In addition we note that what we conclude after collecting some
(inevitably finite amount of) data should depend not only on the data
but \textit{also} on what we knew before collecting the data. In
contrast, what we conclude should \textit{not} depend on arbitrary
choices of family of critical regions $C_\eta$. We further note that
although in a small subset of problems to be solved there is a
uniquely determined ``uniformly optimal'' critical region, in the
majority of problems there is no special critical region and there is
an infinite choice of (nested families of) critical regions satisfying
the necessary conditions.

As we have seen in figure \ref{unifopt} of the example of section
\ref{example1}, even a uniformly optimal nested set of critical
regions (where it exists) does not necessarily give remotely sensible
answers. We conclude that FHT can therefore not be considered a
reliable method of solving inference problems.

Moreover use of frequentist confidence sets in place of frequentist
hypothesis testing solves only one of the problems with the
frequentist approach (namely the lack of symmetry). We saw in section
\ref{conclusive} how confidence sets can be made to treat $H_0$ and
$H_1$ symmetrically, at least when symmetry between $H_0$ and $H_1$
exists to start with. But we also saw how there are many, many
possible choices of confidence set functions, and that not all of them
conserve symmetry (see e.g. figure \ref{freqmixedf2}). We therefore
conclude that FCS also can not be considered a reliable method of
solving inference problems.

Returning to our original concept of an admissible method, it is clear
that any Bayesian method satisfies the conditions for it to be
admissible, and that frequentist methods in general do not. 

Even the issue of the choice of priors to please more than one party
can be addressed, as we have seen in section \ref{regulatorprior}
above, and in such a way that it is the data rather than either
party's particular status, opinion, or interests that determines which
party's prior prevails.

In short, this is a case of \textit{res ipsa loquitur} -- the thing
speaks for itself. It should by now be obvious both that the Bayesian
method is mathematically and logically correct and (as far as we can
tell\footnote{As G\"odel has taught us, we cannot prove it.}) consistent, and
that the frequentist method both requires arbitrary and inappropriate
choices and violates the basic properties required of any inference
method.

We therefore conclude that the Bayesian method \textit{is better} than
the frequentist methods discussed (including both hypothesis testing
and frequentist confidence sets/intervals), despite this conclusion
currently being politically incorrect. We conclude further that
control of type I error rate, characteristic of frequentist
approaches, is not only incompatible with the inference process
matching the real world, but (as seen in the example of section
\ref{example1} and in particular in figure \ref{freqct195} and table
\ref{resultstable}) does not even come close to guaranteeing that a
frequentist solution isn't biased in a real-world sense in favour of
$H_1$. We would even go further and conclude that the Bayesian method
is the \textit{correct} way to do inference, and that the frequentist
methods (including the pseudo-Bayesian method and the overall
uniformly optimal solution if it exists) are \textit{wrong in principle}
as methods for doing inference. Those who employ these frequentist
methods should not pretend to users that they give answers to either
the User's question or the Bayesian question of section
\ref{introbackground}, should not pretend that they are anything other
than a sometimes-convenient short-cut in dealing with such questions,
and should not attempt to force others to do likewise. In particular
regulators should permit the use of (pure) Bayesian solutions without
requiring any control of frequentist type I error rate, i.e. without
requiring them to also be frequentist solutions.

\section{Discussion - What should we do about it ?}
\label{discussion}

If then we are agreed that Bayesian methods are the right way of
solving inference problems and that frequentist methods are wrong,
what should we be doing about it ?

One of the observations we have made in practice is that those who
have never been taught frequentist methods find it easy and intuitive
to understand the Bayesian approach. On the other hand, those who are
first being taught frequentist methods find themselves having to
adjust from thinking intuitively about posterior probabilities to
thinking about likelihoods and critical regions (though often the
huge variety of critical regions possible is neither taught to them
nor realised by them); such individuals then find it far harder to
adjust their thinking to the Bayesian approach, despite that being the
naturally intuitive approach to those who have not been taught
frequentism. 

We therefore think first that it is actually rather important that we
\textbf{\textit{stop} teaching high-school and undergraduate students
  frequentist methods} and teach Bayesian methods instead. Frequentist
methods can then be compared with Bayesian methods at a later stage
and their defects simultaneously pointed out, noting that they are of
historical interest because of the way that they deceived so many
people during the 20th and early 21st centuries.

Second, it is important that \textbf{the public become aware that
  frequentist methods are wrong}. To achieve this we need academics to
be prepared to admit that this is the case -- not an easy thing to do
if ones career has been built around doing research on frequentist
techniques. Clearly vested interests and the reluctance of the general
public to think about anything even vaguely mathematical will make
this a hard goal to reach -- but that doesn't mean we shouldn't try.

Third, it is crucial that \textbf{regulatory authorities stop
  requiring the use of frequentist methods}. Their current insistence
on frequentism causes both cheating\cite{CTM} and an increasing number
of people to think that cheating is not actually cheating but the only
possible way to proceed -- and of course it renders the testing done
by people who cheat invalid and not an assurance of safety at all. It
often also leads to more patients being recruited to clinical trials
than are necessary to reach the desired conclusions, thus putting more
patients at risk of getting inferior treatment than necessary (see
section \ref{example2} and in particular table \ref{statresults}
above).

How to achieve these aims is a question that is largely out of the
realm of mathematics and of our personal expertise; but we do not
believe that running a truce\cite{LeeChu} and pretending that both
methods are equally valid is likely to lead to the desired outcomes.

\appendix

\section{Appendix: A discrete problem with no nuisance variables}
\label{dice}

\subsection{The problem}

This appendix is intended for those who have never encountered the
difference between frequentist and Bayesian methods before. We
consider a very simple problem: we have two fair dice, one with 12
sides (numbered 1 to 12) and one with 20 sides (numbered 1 to
20). Somebody is going to throw one of the dice somewhere where we
can't see it and tell us what number $x$ was thrown. Based on that
number, we are going to try to work out how many sides the die that
was thrown had.

Let $\theta$ be the number of sides the thrown die had. Then $\theta$
takes one of the values in the set $\Theta=\{12,20\}$. 

For the sake of defining notation similar to that of section
\ref{infprobdefsub} we will let $\Phi=\{0\}$ (i.e. $\phi$ always takes
the same value $0$ and can be ignored),
$H=\Theta\cong\Theta\times\Phi$, $H_{12}=\{12\}$, $H_{20}=\{20\}$, and
$x$ will take some value in the set $X=\{1,2,...,20\}$. The likelihood
is given by $$P(x|\theta) = \left\{\begin{matrix} \frac{1}{12}[1 \leq
  x \leq 12] & (\theta=12)\\ \frac{1}{20}[1 \leq x \leq 20] & (\theta
= 20)
\end{matrix}\right.$$
for integers $\theta$ and $x$. (Here $[\ ]$ denotes 1 if the statement
inside the brackets is true and zero if it is false.)

\subsection{The Bayesian solution}

Consider now two events $A$ and $B$, such that $P(B)>0$. If we define
the conditional probability $P(A|B)$ of $A$ given $B$
by $$P(A|B)=\frac{P(A,B)}{P(B)},$$ where $P(A,B)$ means the
probability that both $A$ and $B$ occur, then we
have $$P(A|B)P(B)=P(A,B)=P(B,A)=P(B|A)P(A),$$ and
therefore $$P(A|B)=\frac{P(A)P(B|A)}{P(B)},$$ which is Bayes' theorem
in its most basic form.

Now suppose that $A_n$ is the event that $\theta$ takes a particular
value $n$, and that the only possible values $\theta$ can take are
$1,2,...,N$. Then $$P(B) = \sum_{n=1}^N{P(A_n,B)},$$ so
that $$P(A_n|B)=\frac{P(A_n)P(B|A_n)}{P(B)}=
\frac{P(A_n)P(B|A_n)}{\sum_{n=1}^N{P(A_n,B)}} =
\frac{P(A_n)P(B|A_n)}{\sum_{n=1}^N{P(A_n)P(B|A_n)}},$$ which is the
typical form in which we encounter Bayes' theorem for discrete
variables. For other problems involving continuous variables we may
find the summation in the denominator replaced by an integral, the
probabilities replaced by probability densities, the caveat that the
equation only holds with probability $1$, or combinations of these.

Applying this result to the problem in front of us, Bayes' theorem
tells us that for each of the possible values of $\theta$,
\begin{IEEEeqnarray*}{rCl}
P(\theta|x)&=&\frac{P(\theta)P(x|\theta)}{\sum_{\theta\in\{12,20\}}{P(\theta)P(x|\theta)}},
\end{IEEEeqnarray*}
giving us a formula for the posterior probability $P(\theta|x)$ which
is what we want to know. Now, the problem statement has specified
$P(x|\theta)$, but the Bayesian then next has to specify $P(\theta)$
for each value of $\theta$, i.e. say how likely he thought the two
values of $\theta$ were before the die was thrown. We choose to say
that they are equally likely, so that $P(\theta) =
\frac{1}{2}[\theta\in\{12,20\}].$

Then Bayes' theorem tells us in particular that
\begin{IEEEeqnarray*}{rCl}
P(\theta=12|x)&=&\frac{P(\theta=12)P(x|\theta=12)}{\sum_{\theta\in\{12,20\}}{P(\theta)P(x|\theta)}}\\
&=&\frac{\frac{1}{2}P(x|\theta=12)}{\sum_{\theta\in\{12,20\}}{\frac{1}{2}P(x|\theta)}}\\
&=&\frac{P(x|\theta=12)}{\sum_{\theta\in\{12,20\}}{P(x|\theta)}}.
\end{IEEEeqnarray*}

Now, if $x\in\{1,2,...,12\}$, then that gives us
$$P(\theta=12|x)=\frac{\frac{1}{12}}{\frac{1}{12}+\frac{1}{20}}=\frac{5}{8}$$
while if $x\in\{13,14,...,20\}$ it gives us
$$P(\theta=12|x)=\frac{0}{0+\frac{1}{20}}=0,$$ as we might
hope. Similarly $P(\theta=20|x)=\frac{3}{8}$ if $x\in\{1,2,...,12\}$,
otherwise $P(\theta=20|x)=1$.

In exactly similar manner we find that if we had instead said that
\textit{a priori} $P(H_{12})=\frac{2}{3}$ and $P(H_{20})=\frac{1}{3}$,
i.e. that before knowing $x$ we thought that it was twice as likely
that the 12-sided die would be thrown than the 20-sided one, then we
would still have had $P(H_{12}|x)=0$ for $x>12$, but now if $x\leq 12$
then $P(H_{12}|x)=\frac{10}{13}$.

\subsection{Frequentist solutions in general}

Now, the frequentist's main concern is to ensure that if $H_0$ is true
then the probability is at most $1-\eta$ that we decide that we are
$\eta$ frequentist confident that $H_1$ is true, for any particular
level of frequentist confidence $\eta$; most often frequentists are
interested in $\eta=0.95$, though there is nothing special about this
value.

So we first have to decide which of $H_{12}$ and $H_{20}$ we are going
to call $H_0$; below we will consider each possibility separately (and
call the other one $H_1$). Having made that choice, for each
particular value of $\eta\in[0,1]$, before observing the data $x$, we
have to determine a ``critical region'' $C_\eta$, a subset of the data
range $X$, if $x$ is in which we will draw the conclusion that we are
\textit{at least} $\eta$ frequentist confident that $H_1$ is true.

However, in order that we are consistent, and don't end up with
contradictory conclusions that e.g. we are at least $0.75$ frequentist
confident that $H_1$ is true but we are not even $0.6$ frequentist
confident that $H_1$ is true, we impose the restriction that if
$\eta_1\leq \eta_2$ then $C_{\eta_1}\supseteq C_{\eta_2}$, so that it
is then impossible for $x$ to be in $C_{\eta_2}$ without it also being
in $C_{\eta_1}$. A set of critical regions $(C_\eta)_{\eta\in[0,1]}$
that meets this condition is termed ``nested''.

In order that we meet the frequentist's main concern, we also insist
that\footnote{Here $P(C_\eta|h)$, for a subset $C_\eta$ of $X$, means
  the same as $P(x\in C_\eta|h)$, as we regard the random variable $x$
  as inducing a probability measure on $X$, in the same way that
  $P(H_0)$ means the same as $P(h\in H_0)$.} for any $h\in H_0$,
$P(C_\eta|h)\leq 1-\eta$. In this particular case there will only be
one value of $h$ in $H_0$: if $H_0=H_{12}$ then $h\in H_0$ implies
$h=12$, and similarly if $H_0=H_{20}$ then $h\in H_0$ implies $h=20$;
so this condition can be simplified to $P(C_\eta|H_0)\leq 1-\eta$. But
in general in more complicated problems there may be many possible
values of $h$ in $H_0$, so we will need the more complicated
condition.

There is nothing magic about how to choose such $C_\eta$s; we are
allowed to choose any set of regions satisfying the above conditions,
though usually frequentists try to choose $C_\eta$s that tend to make
$P(C_\eta | H_1)$ be as large as possible, or at least fairly
large. If $C_\eta$ has the property that for any alternative $C'_\eta$
satisfying\footnote{The symbol $\forall$ is read ``for all''.} $\forall h\in H_0,P(C'_\eta|h)\leq 1-\eta$, we have
$\forall h\in H_1,P(C_\eta|h)\geq P(C'_\eta|h)$, then $C_\eta$ is said
to be ``uniformly optimal''; for most problems no uniformly optimal
critical regions exist. 

Finally we can observe the actually occurring value of $x$ and report
the frequentist confidence that $H_1$ is true
as $$c=\sup(\{\eta\in[0,1]: x\in C_\eta\}\cup\{0\}),$$ where $\sup$ denotes the
supremum of the set, i.e. its least upper bound. Since all the sets
for this particular problem are finite, and the inclusion of $0$
ensures that the set is non-empty, the supremum will simply be the
maximum of the set, i.e. the maximum value of $\eta$ for which we are
at least $\eta$ frequentist confident that $H_1$ is true.

\subsection{Frequentist solutions with $H_0=H_{20}$}

The frequentist then first has to choose whether to make $H_0=H_{20}$
or $H_0=H_{12}$. We next consider the first possibility, i.e. that he
chooses to set $H_0=H_{20}$.

The frequentist next has to choose nested critical regions
$(C_\eta)_{\eta\in[0,1]}$ in $X$ such that for each $\eta$ and each
$h\in H_0=H_{20}$ (of which there is only one, namely
$h=\theta=20$), $$P(x\in C_\eta|h)\leq 1 - \eta.$$

Now, for the only $h\in H_0$, $P(x|h)=\frac{1}{20}$ for all $x\in
\{1,2,...,20\}$. 

\subsubsection{First frequentist solution with $H_0=H_{20}$}

So one possibility would be to set $$C_\eta =
\{1,2,...,\lfloor 20(1-\eta)\rfloor\}\subseteq X$$ where $\lfloor a \rfloor$
means the largest integer less than or equal to $a$, and in particular
$C_\eta=\emptyset$ for $\eta > 0.95$. The reader should check for
himself that the conditions on $C_\eta$ given under `Frequentist
question' in section \ref{introbackground} do indeed then hold (and
similarly for later alternative definitions of $C_\eta$).

Then we calculate the frequentist confidence $c$ that $H_1$
(i.e. $H_{12}$) holds, given by $$c=\sup{(\{\eta\in[0,1]:x\in
  C_\eta\}\cup\{0\})}.$$ Thus $c=1-\frac{x}{20}$, so that: if the
number rolled was a $1$ then we become 95\% frequentist confident that
it was the 12-sided die that was thrown; if it was 13, then we become
35\% frequentist confident that it was the 12-sided die that was
thrown; while if we roll a 20 we become 0\% frequentist confident that
it was the 12-sided die that was thrown.

Now some of this is reasonable (e.g. if a 20 was thrown, that we are
0\% frequentist confident that it was the 12-sided die that was
thrown), but other parts of it are not (e.g. that if a 1 was rolled
being 95\% frequentist confident that it was the 12-sided die that was
thrown, or that we are 35\% frequentist confident that it was the
12-sided die that was thrown if a 13 had been rolled).

\subsubsection{Second frequentist solution with $H_0=H_{20}$}

But it would be equally frequentistly valid to instead set $$C_\eta =
\{21-\lfloor 20(1 - \eta) \rfloor,...,19,20\}\subseteq X,$$ and in
particular $C_\eta=\emptyset$ for $\eta>0.95$ (and of course there are
many other possible settings which we haven't time to consider). Then
$c=\frac{x-1}{20}$, so that if the number rolled was 20 we become 95\%
frequentist confident that the die rolled was the 12-sided one (er
what ?!?), while if it was a 1 that was thrown we are only 0\%
frequentist confident that it was the 12-sided die that was thrown.

\subsubsection{Pseudo-Bayesian frequentist solution with $H_0=H_{20}$}

Alternatively we could set $$C_\eta = \left\{\begin{matrix}
X & (\eta = 0)\\
\{1,2,...,12\} & (0 < \eta \leq \frac{2}{5})\\
\emptyset & (\eta > \frac{2}{5}),
\end{matrix}\right.$$ so that if $x\leq 12$ then $c=\frac{2}{5}$,
while if $x>12$ then $c=0$. It is then entirely reasonable, if we have
rolled a 13 or bigger, that we then have no confidence that it was the
12-sided die that was thrown, but if a 12 or less was rolled it is not
reasonable that we still think it more likely that the 20-sided die
was thrown -- we can never even get 50\% frequentist confidence
that the 12-sided die was thrown.

We note that this solution corresponds to the (only) basic
pseudo-Bayesian solution with $H_0=H_{20}$ as defined in section
\ref{concisedescription}.

\subsection{Frequentist solutions with $H_0=H_{12}$}

But perhaps the problem is that we have picked the ``wrong'' $H_0$ ?
(But we are free to pick either, so that can't be it.) Let's see what
happens if we set $H_0=H_{12}$. 

Then for $x>12$ and $h\in H_0$ we have $P(x|h)=0$, while for $x\leq
12$ we have $P(x|h)=\frac{1}{12}$.

\subsubsection{First frequentist solution with $H_0=H_{12}$}

So we could e.g. set $$C_\eta = \{1,2,...,\lfloor
12(1-\eta)\rfloor\}\subseteq X$$ (with $C_\eta=\emptyset$ for $\eta >
\frac{11}{12}$), in which case $c=0\lor(1-\frac{x}{12})$ (where $a\lor
b$ is the maximum of $a$ and $b$), so that we can never be as much as
95\% frequentist confident that the 20-sided die was thrown (even if a
20 had been thrown).

\subsubsection{Second frequentist solution with $H_0=H_{12}$}

Or we could set $$C_\eta = \{13-\lfloor 12(1-\eta) \rfloor,...,19,20\}$$
which gives us $c=1\land \frac{x-1}{12}$ (where $a\land b$ is the
minimum of $a$ and $b$), so that for $x\geq 13$ we are 100\%
frequentist confident that the 20-sided die was thrown (Hurrah !), for
$x=12$ we are $\frac{11}{12}$ frequentist confident that the 20-sided die was
thrown, and for $x=1$ we are 0\% frequentist confident that the
20-sided die was thrown (the latter two points not being so good).

\subsubsection{Pseudo-Bayesian solution with $H_0=H_{12}$}

Alternatively we could set $$C_\eta = \left\{\begin{matrix}
X & (\eta = 0)\\
\{13,14,...,20\} & (\eta > 0),
\end{matrix}\right.$$ so that if $x>12$ then $c=1$, otherwise
$c=0$. It is then obviously good that if a 13 or bigger was rolled, we
are 100\% frequentist confident that the 20-sided die was thrown, but
it is not good that we have zero frequentist confidence that the
20-sided die was thrown if the number rolled was 12 or less.

Again, this corresponds to the (only) basic pseudo-Bayesian solution
with $H_0=H_{12}$ as defined in section \ref{concisedescription}.

\subsection{Conclusion on the dice example}

As we will see so often in this paper, the Bayesian solution makes
sense, and none of the frequentist ones really do (though some are
better than others). Intuitively it should only matter whether the
number rolled is $\leq 12$ or $>12$, and it should not matter whether
it is 4 or 5, and it should not matter whether it is 19 or 20. Only
the Bayesian one and the two  pseudo-Bayesian frequentist
solutions (of those examined) satisfy these points, and neither of the
 pseudo-Bayesian ones give sensible answers for $x\leq 12$ (but
there are numerous other frequentist possibilities -- how do we choose
which frequentist solution to use ?).

\section{Appendix: Probability measures, integration, and density functions}
\label{probmeas}

Here we give the briefest of introductions to this subject without any
proofs etc. For the reader who would like to know more, we recommend
\cite{Royden} (or \cite{Chung}).

\subsection{Definitions}
\label{probmeasdef}

Suppose we have a set $\Omega$ and a set of subsets $M$ of $\Omega$ to
each of which we want to associate a number called its ``probability''
by a function $R$. Doing so in a consistent and sensible way turns out
to need the following conditions to hold:

\begin{enumerate}

\item \label{sigmafield} About $M$:

\begin{enumerate}

\item $\Omega\in M$;

\item If $A\in M$ then $\Omega\setminus A\in M$ (here $\Omega\setminus A$
  denotes the set of everything in $\Omega$ that is not in $A$);

\item If $A_1,A_2,...\in M$ then $\bigcup_{n=1}^\infty{A_n}\in M$.

\end{enumerate}

\item \label{meas} About $R$:

\begin{enumerate}

\item $R$ is a function from $M$ to $\mathbb{R}\cup\{\infty\}$, the set of real
  numbers with +infinity;

\item For all $A\in M$, $R(A)\geq 0$.

\item \label{totalprob} $R(\Omega) = 1$;

\item If $A_1,A_2,...\in M$ and for all $n_1\neq n_2$, $A_{n_1}\cap
  A_{n_2}= \emptyset$,
  then $$R(\bigcup_{n=1}^\infty{A_n})=\sum_{n=1}^\infty{R(A_n)}.$$

\end{enumerate}

\end{enumerate}

If conditions \ref{sigmafield} hold then we say that $M$ is a
\textbf{$\sigma$-field} (or \textbf{$\sigma$-algebra}) of subsets of
$\Omega$. If conditions \ref{sigmafield} and \ref{meas} hold then we
say that $R$ is a \textbf{probability measure} on the $\sigma$-field
$M$ of measurable sets of $\Omega$, or more briefly that $R$ is a
probability measure on $\Omega$. If all except \ref{totalprob} hold,
then we say that $R$ is a \textbf{measure} on $M$ (or on $\Omega$).

Intuitively the conditions can be interpreted as follows (see also the
examples of probability measures in section \ref{probmeasexamples}
below).

In the case of \ref{sigmafield}, we need to be able to talk about the
probability that \textit{some} outcome occurs; if we can talk about
the probability that $A$ occurs, then we also need to be able to talk
about the probability that $A$ does not occur; and if we can talk
about each of the probabilities that $A_n$ occurs, then we need to be
able to talk about the probability that at least one of the $A_n$
occurs.

In the case of \ref{meas}, probabilities cannot be negative; the
probability that some outcome happens (anything at all) must be $1$;
and if $A_1,A_2,...$ are mutually exclusive events, then the
probability that at least one of them happens must be the sum of their
individual probabilities.

All of these conditions are eminently reasonable; and as far as this
definition and the subsequent examples in section
\ref{probmeasexamples} go, frequentists and Bayesians agree on them.

In particular if $\Omega$ is the set $\Theta$ of possible values of
the unknown in which we are interested, then $M$ is the set of those
subsets of $\Theta$ to which we can attach a probability using
$R$. Then some brief thought should be sufficient to convince one that
all these conditions should hold for $R$ to be a sensible answer to
the associated inference problem given some observed $x\in X$ (on
this, frequentists might not agree). The only point that may not be
obvious is that $M$ need not be the set $\mathbb{S}(\Theta)$ of
\textit{all} subsets of $\Theta$; unfortunately it turns out that
requiring $M=\mathbb{S}(\Theta)$ is unduly restrictive on what
$\Theta$ can be, hence our not insisting on this equality, and only
requiring that $M\subseteq \mathbb{S}(\Theta)$.

\subsection{Examples}
\label{probmeasexamples}

\subsubsection{Unit interval}

For a simple example let $\Omega=[0,1]$, $M$ be the smallest subset of
$\mathbb{S}(\Omega)$ containing the half-open intervals
$[0,b)\subset[0,1]$ that satisfies the conditions of \ref{sigmafield}
  of \ref{probmeasdef}, and $R$ be the only function satisfying the
  conditions of \ref{meas} such that $R([0,b))=b$. (That $M$ and $R$
    exist and are unique will have to be taken on trust here (or see
    \cite{Royden}).)

Then $R$ could be considered to be the probability measure that
describes the output of a uniform random number generator, and $R(A)$
tells us how likely such a generator's output is to be in the set $A$.

\subsubsection{Unit square}

Alternatively let $\Omega=[0,1]\times [0,1]$, the unit square, $M$ be
the smallest subset of $\mathbb{S}(\Omega)$ containing the half-open
rectangles $[0,a) \times [0,b)\subset[0,1]^2$ that satisfies the
    conditions of \ref{sigmafield} of \ref{probmeasdef}, and $R$ be
    the only function satisfying the conditions of \ref{meas} such
    that $R([0,a)\times[0,b))=ab$.

Then similarly $R(A)$ tells us how likely it is that the values
$(x,y)$ from two independent uniform random number generators are to
lie in the 2-d set $A$.

\subsubsection{Gaussian on the real line}
\label{Gaussmeas}

Alternatively let $\Omega=\mathbb{R}$, $M$ be the smallest subset of
$\mathbb{S}(\Omega)$ containing the open intervals $(a,b)$ that
satisfies the conditions of \ref{sigmafield} of \ref{probmeasdef}, and
$R$ be the only function satisfying the conditions of \ref{meas} such
that for $a\leq b$, $R((a,b)) = G(b)-G(a)$, where $G$ here is the
cumulative distribution function of the standard unit Gaussian given
by $$G(a)=\int_{-\infty}^a{\frac{1}{\sqrt{2\pi}}e^{-\frac{1}{2}x^2}\,dx}.$$
Then $R(A)$ tells us how likely a standard normal random variable is
to have a value in the set $A$.

\subsubsection{Lebesgue measure on the real line}

Alternatively let $\Omega=\mathbb{R}$, $M$ be the smallest subset of
$\mathbb{S}(\Omega)$ containing the open intervals $(a,b)$ that
satisfies the conditions of \ref{sigmafield} of \ref{probmeasdef}, and
$R$ be the only function satisfying the conditions of \ref{meas}
except for \ref{totalprob} such that for $a\leq b$, $R((a,b)) =
b-a$. Then $R(A)$ tells us the length of the set $A$, and $R$ is a
measure that is not a probability measure. This is the measure that is
assumed when doing ordinary integration of functions on the real line
unless stated otherwise.

\subsubsection{Finite set}
\label{finsetmeas}

Alternatively let $\Omega=\{\omega_1,\omega_2,...,\omega_N\}$ be a
finite set and $p_1,p_2,...,p_N$ be non-negative real numbers such
that $\sum_{n=1}^N{p_n}=1$, and let $M=\mathbb{S}(\Omega)$. Then for
any subset $A$ of $\Omega$ we can write $$R(A) = \sum_{n:\omega_n\in
  A}{p_n},$$ giving us a probability measure on $\Omega$. In this case
(see section \ref{integration} below) if $g:\Omega\to\mathbb{R}$,
then $$\int_\Omega{g(\omega)\,dR(\omega)}=\sum_{n=1}^N{g(\omega_n)p_n}.$$

\subsection{Integration, density functions, and the Radon-Nikodym theorem}
\label{integration}

Just as we can do integration on the real line using Lebesgue measure,
so also we can integrate\footnote{How to do this occupies the best
  part of a term's course in the third year of a maths degree,
  although it is in a sense intuitively obvious.} with respect to any
other measure $R$, and denote the resulting integral of a function
$f:\Omega\to\mathbb{R}$ by $$\int{f(\omega )\,dR(\omega )}.$$ If two measures
$R_1$ and $R_2$ on the same set of measurable sets $M$ and an
$M$-measurable\footnote{$f:\Omega\to \mathbb{R}$ is
  \textbf{$M$-measurable} if for all $r\in\mathbb{R}$,
  $\{\omega\in\Omega: f(\omega)<r\}\in M$.} function
$f:\Omega\to\mathbb{R}$ have the relationship that for any $A\in
M$, $$R_2(A)=\int_A{f(\omega )\,dR_1(\omega )},$$ then we say that $f$ is a
density function for $R_2$ with respect to $R_1$ on $M$. Then for any
other integrable function $g:\Omega\to\mathbb{R}$ and $A\in M$ we have
the relationship $$\int_A{g(\omega )\,dR_2(\omega )}=\int_A{g(\omega )f(\omega )\,dR_1(\omega )}.$$

Thus for example if $R_1$ is Lebesgue measure on the real line and
$R_2$ is the Gaussian measure of section \ref{Gaussmeas} above,
then $$f(x)=\frac{1}{\sqrt{2\pi}}e^{-\frac{1}{2}x^2}$$ is a density
function for $R_2$ with respect to $R_1$. Similarly if $R_2$ is the
probability measure on the finite set $\Omega$ of section
\ref{finsetmeas} above, and $R_1$ is the measure on the same set such
that $R_1(A)=|A|$, the number of elements of $A$, then
$f(\omega_n)=p_n$ defines a density function for $R_2$ with respect to
$R_1$.

\textit{A fortiori} if $P$ is a probability measure on
$\Omega=\mathbb{R}$ with a density function $f$ with respect to
Lebesgue measure then we can write the expectation of $g(x)$
as $$\E_x{
  g(x)}=\int_\Omega{g(x)\,dP(x)}=\int_{\mathbb{R}}{g(x)f(x)\,dx};$$ in
this case we often (by abuse of notation) also denote the density
function by $P$, in which case we
get $$\E_x{g(x)}=\int_{\mathbb{R}}{g(x)P(x)\,dx}.$$

Construction of density functions can be achieved with the
Radon-Nikodym theorem, which states that if $R_1$ and $R_2$ are
measures on a $\sigma$-field $M$ of subsets of $\Omega$, and for all
$A\in M$ such that $R_1(A)=0$ we also have $R_2(A) = 0$, then $R_2$
has a density function with respect to $R_1$ on $M$. Note that there
may be many such density functions, any two of which differ only on a
set of probability zero under $R_1$. Note also that if $M'$ is a
sub-$\sigma$-field of $M$, then $R_2$ also has density functions with
respect to $R_1$ on $M'$, and that in general they are \textit{not}
the same as the density functions on $M$. Such a density function is
also known as a ``Radon-Nikodym derivative'' and may be denoted
$\frac{{\rm d}R_2}{{\rm d}R_1}$.

\section{Appendix: Induced measures}
\label{induced}

Here we define what is meant by the probability measures ``induced''
by a given probability measure $P(a,b,c)$. To avoid getting into
measure-theoretic technicalities, we simplify by assuming that all
variables are continuous with values in some finite-dimensional real
space $\mathbb{R}^n$; the reader may wish to assume that $n=1$ for
further simplicity\footnote{The general case is essentially the same
  but taking the underlying measure to be $P$ rather than Lebesgue
  measure, so that the density $P(a,b,c)=1$ everywhere.}. Then the
measure $P$ may be expressed as a joint density function $P(a,b,c)$ on
$\mathbb{R}^{3n}$.

Then the probability measure $P(a,b)$ induced by $P$ is given
by $$P(a,b)=\int_C{P(a,b,c)\,dc},$$ where the range of integration is
over the set $C$ of all possible values of $c$.

\textit{A fortiori} the probability measure $P(a)$ is given
by $$P(a)=\int_B{P(a,b)\,db} = \int_C{P(a,c)\,dc},$$ (the second
equality being guaranteed by Tonelli's theorem) and similarly for
$P(b), P(b,c), P(c)$, etc.

Further $P(a,b|c)$ is given by $$P(a,b|c)=\frac{P(a,b,c)}{P(c)},$$ at
least when $P(c)\neq 0$. Alternatively, for those familiar with the
Radon-Nikodym theorem, we may use it to define $P(a,b|c)$ in such a
way that for all $c$ not in a set of probability zero $P(a,b|c)$ is a
probability measure on $A\times B$, the product of the ranges of $a$
and $b$, that satisfies $$P(a,b,c)=P(a,b|c)P(c),$$ and such that for
any two probability measures $P_1(a,b|c)$ and $P_2(a,b|c)$ satisfying
this definition we have $P_1(a,b|c)=P_2(a,b|c)$ except on some set of
probability zero.

\textit{A fortiori} the probability measure $P(a|c)$ is given
by $$P(a|c) = \int_B{P(a,b|c)\,db} = \frac{P(a,c)}{P(c)}.$$

\section{Appendix: Basics of Bayesian calculation}
\label{basicBayes}

If $a$ and $b$ are random variables having density functions with
respect to some underlying measure, then Bayes' theorem tells us
that $$P(a|b)=\frac{P(a)P(b|a)}{\int_A{P(a)P(b|a)\,da}}$$ except
possibly on an event which has probability zero, where $A$ is the
range of possible values of $a$. Here each of $a$ and $b$ may be
scalar random variables or vectors of many random variables,
i.e. many-dimensional random variables.

Given a prior $P(h)=P(\theta,\phi)$ on $H$ (which is assumed to be
zero for $(\theta,\phi)\notin H$) and a likelihood $P(x|\theta,\phi)$
we calculate the Bayesian posterior distribution on $\Theta$ using
Bayes' theorem by $$P(\theta,\phi|x) =
\frac{P(\theta,\phi)P(x|\theta,\phi)}{\int_H{P(\theta,\phi)P(x|\theta,\phi)\,d(\theta,\phi)}}$$
and hence by integrating with respect to $\phi$ $$P(\theta|x) =
\frac{\int_{\{\phi\in\Phi:(\theta,\phi)\in
    H\}}{P(\theta,\phi)P(x|\theta,\phi)\,d\phi}}{\int_H{P(\theta,\phi)P(x|\theta,\phi)\,d(\theta,\phi)}},$$
and in particular
$$P(H_1|x) =
\frac{\int_{H_1}{P(\theta,\phi)P(x|\theta,\phi)\,d(\theta,\phi)}}{\int_H{P(\theta,\phi)P(x|\theta,\phi)\,d(\theta,\phi)}}.$$

(While the above calculation will always work, sometimes it is easier
to do things slightly differently to simplify the working (e.g. in
appendix \ref{simpleexample} section \ref{simpleBayes}).)

\section{Appendix: Proofs of claims in section \ref{concisedescription}}
\label{conciseproofs}

Suppose we have chosen and fixed a prior $P(h)$, $H_0\subseteq
H_0'\subseteq H$, $H_1'=H\setminus H_0'$, and a data collection
plan. We define $f,f_1,g:[0,1]\to [0,1]$ and $B,B_1,C:[0,1]\to
\mathbb{S}(X)$ by $$B(p) = \{x\in X:P(h\in H_1'|x)> p\}$$ $$B_1(p) =
\{x\in X:P(h\in H_1'|x)\geq p\}$$ $$f(p) = 1 - \sup_{h\in H_0}P(x\in
B(p)|h)$$ $$f_1(p) = 1 - \sup_{h\in H_0}P(x\in B_1(p)|h)$$ $$ g(\eta)
= \inf\{p\in [0,1]:f(p)\geq \eta\}$$
$$C_\eta = \left\{\begin{matrix}B_1(g(\eta)) & (f_1(g(\eta))\geq
\eta)\\ B(g(\eta)) & (f_1(g(\eta)) < \eta).\end{matrix}\right.$$ We
will need to prove that $f,f_1,g$ are (non-strictly) increasing, that
for $p_1<p_2$, $f(p_1)\leq f_1(p_2)\leq f(p_2)$, that $f$ is
right-continuous, that for $p,\eta\in [0,1]$, $$g(f(p))\leq
p$$ $$f(g(\eta))\geq \eta,$$ and that $(C_\eta)_{\eta\in [0,1]}$ are a
valid nested set of critical regions.

The event-valued function $B(p)$ is decreasing and right-continous in
$p$, therefore so also is $P(x\in B(p)|h)$ by continuity of
probability. Then a supremum of decreasing right-continuous functions
is right-continuous and decreasing; thanks to the minus sign, $f$ is
therefore increasing and right-continuous. Moreover for $p_1 < p_2$,
$B(p_1)\supseteq B_1(p_2) \supseteq B(p_2)$, so $f(p_1)\leq
f_1(p_2)\leq f(p_2)$, whence $f_1$ is also increasing. Note also that
$f(1) = 1$.

Let $$A(\eta) = \{p\in [0,1]:f(p)\geq \eta\}.$$ Then $A(\eta)$ is
decreasing in $\eta$, therefore its infimum is increasing, as is $g$.

Now $$g(f(p)) = \inf\{q\in [0,1]:f(q)\geq f(p)\}$$ and $p$ belongs to
this set; therefore $$g(f(p))\leq p.$$

Considering the other composition of $f$ and $g$, for
$\eta\in[0,1]$, $$f(g(\eta)) = f(\inf\{p \in [0,1]:f(p)\geq
\eta\}).$$ This set is non-empty as $f(1)=1\geq\eta$, so since $f$ is
increasing and right-continuous, $f(g(\eta))\geq \eta$.

It remains to show that $(C_\eta)_{\eta \in [0,1]}$ are a valid nested
set of critical regions. For $p_1<p_2$ we have both $B(p_1)\supseteq
B_1(p_2)$ and $B_1(p_1) \supseteq B(p_2)$, so since $f_1$ and $g$ are
increasing $C_\eta$ is nested and decreasing. Then for any $h\in H_0$
and any $\eta\in [0,1]$ such that $f_1(g(\eta))<\eta$,
\begin{IEEEeqnarray*}{rCl}
P(x\in C_\eta|h) &=& P(x\in B(g(\eta))|h)\\
&\leq& \sup_{h\in H_0}P(x\in B(g(\eta))|h)\\
&=& 1 - f(g(\eta))\\
&\leq& 1 - \eta,
\end{IEEEeqnarray*}
while for any $\eta$ such that $f_1(g(\eta))\geq \eta$, we are guaranteed
that 
\begin{IEEEeqnarray*}{rCl}
P(x\in C_\eta|h) &=& P(x\in B_1(g(\eta))|h)\\
&\leq& \sup_{h\in H_0}P(x\in B_1(g(\eta))|h)\\
&=& 1 - f_1(g(\eta))\\
&\leq& 1 - \eta
\end{IEEEeqnarray*}
completing the proof.

\section{Appendix: Proofs of claims in section \ref{fullpseudoBayes}}
\label{fullpseudoBayesproofs}

Building on the notation of appendix \ref{conciseproofs}, we recall
that:
 $$S(\eta) = \{\eta' < \eta : C_{\eta'} \supsetneq C_\eta\}$$ $$D_\eta
= \left(X\cap \bigcap_{\eta'\in S(\eta)}{C_{\eta'}}\right) \supseteq
C_\eta$$ $$\zeta_1(\eta) = 1 - \sup_{h\in H_0}{P(x\in
  C_\eta|h)}$$ $$\zeta_0(\eta) = \sup_{\eta'\in S(\eta)}{\left(1 -
  \sup_{h\in H_0}{P(x\in C_{\eta'}|h)}\right)} = 1 - \inf_{\eta'\in
  S(\eta)}{\sup_{h\in H_0}{P(x\in C_{\eta'}|h)}}\leq
\zeta_1(\eta)$$ $$C'_\eta=\left\{\begin{matrix} (C_\eta\times [0,1])
\cup \left(D_\eta\times\left[0,\frac{
    \zeta_1(\eta)-\eta}{\zeta_1(\eta)-\zeta_0(\eta)}\right]\right) &
(\zeta_0(\eta)<\zeta_1(\eta))\\ C_\eta\times [0,1] &
(\zeta_0(\eta)=\zeta_1(\eta)).
\end{matrix}\right.$$ 

Now, we already know that $\zeta_1(\eta)\geq \eta$. We will also need
that $\zeta_0(\eta)\leq \eta$: so suppose otherwise. Then there is an
$\eta'<\eta$ such that $C_{\eta'}\supsetneq C_\eta$ and $$1 -
\sup_{h\in H_0}{P(x\in C_{\eta'}|h)}>\eta.$$ But either
$C_{\eta'}=B(g(\eta'))$ or $C_{\eta'}=B_1(g(\eta'))$, so either way
$C_{\eta'}\supseteq B(g(\eta'))$, and therefore $$f(g(\eta')) = 1 -
\sup_{h\in H_0}{P(x\in B(g(\eta'))|h)} \geq 1 - \sup_{h\in H_0}{P(x\in
  C_{\eta'}|h)} > \eta.$$ But $$g(\eta) = \inf\{p\in[0,1]:f(p)\geq
\eta\}$$ and $f$ and $g$ are increasing, so $g(\eta')=g(\eta).$ But
then since $C_{\eta'}\supsetneq C_\eta$,
$$C_{\eta'}=B_1(g(\eta))=B_1(g(\eta'))\text{ and }C_\eta=B(g(\eta)).$$
This implies that $$f_1(g(\eta')) = f_1(g(\eta))<\eta$$ by
definition of $C_\eta$. But then since $C_{\eta'}=B_1(g(\eta))$, from
the choice of $\eta'$ we have $$f_1(g(\eta)) = 1 - \sup_{h\in
  H_0}{P(x\in C_{\eta'}|h)} > \eta,$$ a contradiction. 

Therefore indeed $\zeta_0(\eta)\leq \eta \leq \zeta_1(\eta)$.

Now, we next need to prove that the $C'_\eta$ are nested
decreasing. So take any $\eta' < \eta$. If $C_{\eta'}\supsetneq
C_{\eta}$ then $C_{\eta'} \supseteq D_\eta$ and the necessary
inclusion is obvious. If $C_{\eta'}=C_\eta$ then $S(\eta')=S(\eta)$,
$D_{\eta'}=D_\eta$, $\zeta_1(\eta')=\zeta_1(\eta)$, and
$\zeta_0(\eta') = \zeta_0(\eta)$, given which we again have
$C'_{\eta'}\supseteq C'_\eta$.

Finally we need to show that $\sup_{h\in H_0}{P((x,u)\in
  C'_\eta|h)}\leq 1-\eta$. For this, fix a particular
$\eta\in[0,1]$. The case $\zeta_0(\eta)=\zeta_1(\eta)$ is easy because
then $C'_\eta=C_\eta\times [0,1]$. So suppose that
$\zeta_0(\eta)<\zeta_1(\eta)$, and let us abbreviate $\zeta_0(\eta)$
by $\zeta_0$ and similarly for $\zeta_1$.

Then for any $h\in H_0$, 
\begin{IEEEeqnarray*}{rCl}
P(x'\in C'_\eta|h) &=& P(x\in C_\eta|h) + \left(P(x\in D_\eta|h)-P(x\in
C_\eta|h)\right)\frac{\zeta_1-\eta}{\zeta_1-\zeta_0}\\ &=&P(x\in
C_\eta|h)\frac{\eta-\zeta_0}{\zeta_1-\zeta_0} + P(x\in
D_\eta|h)\frac{\zeta_1-\eta}{\zeta_1-\zeta_0}.
\end{IEEEeqnarray*}
Then for all $\eta'\in S(\eta)$, $$P(x\in D_\eta|h)\leq P(x\in
C_{\eta'}|h),$$ therefore $$\sup_{h\in H_0}{P(x\in D_\eta|h)}\leq
\sup_{h\in H_0}{P(x\in C_{\eta'}|h)},$$ and so 
\begin{IEEEeqnarray*}{rCl}
\sup_{h\in H_0}{P(x\in D_\eta|h)} &\leq&
\inf_{\eta'\in S(\eta)}{\sup_{h\in H_0}{P(x\in C_{\eta'}|h)}}\\
&=&1 - \zeta_0(\eta).
\end{IEEEeqnarray*}
Therefore by non-negativity of the two fractions involving the
$\zeta$s,
\begin{IEEEeqnarray*}{rCl}
\sup_{h\in H_0}{P(x'\in C'_\eta|h)} &\leq &
(1-\zeta_1)\frac{\eta-\zeta_0}{\zeta_1-\zeta_0} + (1 -
\zeta_0)\frac{\zeta_1-\eta}{\zeta_1-\zeta_0}\\
&=& 1-\eta
\end{IEEEeqnarray*}
as required.

\section{Appendix: Proofs that the pure Bayesian method satisfies the criteria
  of section \ref{criteria}}
\label{bayesproofs}

The statements being proved can be found at the given criterion numbers.

\begin{enumerate}

\item Criterion \ref{complementarity} (``Complementarity''): Since
  $P(h|x)$ is a probability measure on $H$, and since $$H_0\cup H_1 =
  H_0'\cup H_1' = H,$$ we have $$\eta' = P(H_1'|x) = 1 - P(H_0'|x) = 1
  - P(H_1|x) = 1-\eta$$ as required.

\item Criterion \ref{containment} (``Inclusion''): Since $P(h\in
  H_1|x) = \int_{H_1}{P(h|x)\,dx}$ and density functions are
  non-negative, $P(h\in H_1'|x) \geq P(h\in H_1|x)$. Note that the
  prior on $h$ does not depend on whether we are using $H_0$ or
  $H_0'$.

\item Criterion \ref{whetherproceed}
  (``Intention''): \label{stopindep} Let $S$ denote the event that we
  stop data collection before collecting $x_2$. Then since $S$ depends
  only on the value of $x_1$, $S$ is conditionally independent of $h$
  given $x_1$, i.e. $P(S|x_1) = P(S|x_1,h)$. Therefore

\begin{IEEEeqnarray*}{rCll}P(h|S,x_1) &=&
  \frac{P(h)P(S,x_1|h)}{\int{P(h)P(S,x_1|h)\,dh}}&\text{ (Bayes)}\\
  &=&
  \frac{P(h)P(x_1|h)P(S|x_1,h)}{\int{P(h)P(x_1|h)P(S|x_1,h)\,dh}}&\text{
    ( chain rule of probability)}\\
  &=&
  \frac{P(h)P(x_1|h)P(S|x_1)}{\int{P(h)P(x_1|h)P(S|x_1)\,dh}}&\text{
    (conditional independence)} \\
  &=& \frac{P(h)P(x_1|h)}{\int{P(h)P(x_1|h)\,dh}}&\text{ (cancellation)} \\
  &=& P(h|x_1)&\text{ (Bayes)}.
\end{IEEEeqnarray*}
and hence by integrating with respect to $h$ over
$H_1$, $$P(H_1|S,x_1)=P(H_1|x_1).$$

\item Criterion \ref{productprob} (``Conjunction''): If, for all $n$,
  $P(h_n\in H_{1,n}|x_n)=\eta\in(0,1)$, then since the $N$ systems are
  independent, we have $P(\forall n\in \{1,...,N\},\,h_n\in
  H_{1,n}|x_1,...,x_N)=\eta^N<\eta.$

\item Criterion \ref{sumprob} (``Disjunction''): If, for all $n$,
  $P(h_n\in H_{1,n}|x_n)=\eta\in (0,1)$, then since the $N$ systems
  are independent, we have
\begin{IEEEeqnarray*}{rCl}
P(\exists n\in \{1,...,N\}:\,h_n\in H_{1,n}|x_1,...,x_N) &=& 1 -
P(\forall n\in \{1,...,N\}:\,h_n\notin H_{1,n}|x_1,...,x_N) \\
&=& 1 - (1 - \eta)^N\\ & > & 1 - (1-\eta) \\ &=& \eta.
\end{IEEEeqnarray*}

\item Criterion \ref{multisystem} (``Multiplicity''): The formula for
  $P(h_1\in H_{1,1}|x_1)$ contains no reference to either $N$ or to
  any $x_n$ for $n>1$ and the result is therefore unaffected by these
  variables.  Note also that in the Bayesian paradigm, for any fixed
  posterior probability level that constitutes ``apparent positive'',
  the expected number of false positive results is proportional to the
  number of apparent positive results, and is unaffected by any vast
  number of accompanying apparent negative results (in contrast to the
  frequentist paradigm wherein the bound on the number of false
  positives is proportional to the number of actual negatives).

\item Criterion \ref{unifpow}\label{appunifpow} (``Sequential
  Optimality''): Let $$q_k(x_1,...,x_k) = \left\{\begin{matrix}1 &
  (P(h\in H_1|x_1,...,x_k) \geq \eta)\\ 0 & (P(h\in H_1|x_1,...,x_k) <
  \eta),\end{matrix}\right.$$ so that in particular $q_0=[P(h\in
    H_1)\geq\eta]$, i.e. $1$ if the prior probability that $h\in H_1$
  is greater than or equal to $\eta$ and $0$ otherwise. Then since we
  are taking, with the maximum possible probability, every permitted
  opportunity of concluding that $P(h\in H_1|x_1,...,x_k) \geq \eta$,
  any other data collection plan will lead to lower probabilities of
  thus concluding by any particular $k$. (Note that by the same
  argument as in \ref{stopindep} above, since the decision on when to
  conclude depends only on the data, the stopping decision does not
  bias the conclusion.)

\end{enumerate}

\section{Appendix: Examples where the frequentist method \mbox{violates} the
  criteria of section \ref{criteria}}
\label{freqexamples}

We give examples here showing that frequentist hypothesis testing
violates the criteria of section \ref{criteria}. However we also show
in appendix \ref{hypoconf} how to translate these examples into ones
involving frequentist confidence sets, showing that they are no better.

The following are all based on the same basic problem. In terms of our
standard inference problem notation from section \ref{infprobdefsub}
we set $$\Theta=\{0,1\},$$ $$\Phi=\mathbb{R},$$ $$H=\{(\theta,\phi)\in
\Theta\times\Phi: \theta=[\phi >0]\},$$
and $$H_k=(\{k\}\times\Phi)\cap H\ \ \ \ (\text{for }k=0,1),$$ and
then identify $H$ with $\Phi$ without introducing any ambiguity (so
that then $H_0=\{h\in H:h\leq 0\}$ and $H_1=\{h\in H:h>0\}$) and set
the likelihood
by $$P(x|h)=\frac{1}{\sqrt{2\pi}}e^{-\frac{1}{2}(x-h)^2},$$ the
Gaussian of mean $h$ and known unit variance.

In other words we will be testing whether the mean of a unit-variance
Gaussian is $\leq 0$ ($H_0$) or $>0$ ($H_1$).

Where $N$ independent systems are called for, they will all be taken
to be thus but with potentially different values of $h_n$. $[a<b]$ will
denote the indicator function of the set on which $a<b$, i.e. it is 1
if $a<b$ otherwise it is zero.

Throughout the following proofs we denote the cdf of the standard unit
Gaussian by $G$, and we set $\alpha =
0.975,\beta=G^{-1}(\alpha)\approx 1.96,\gamma=G^{-1}(1 -
\frac{1-\alpha}{2})\approx 2.2414$.

The criterion for which each is an example of violation by frequentist
methods is as shown.

\begin{enumerate}

\item Criterion \ref{complementarity} (``Complementarity''): In this
  case we will vary the problem slightly by instead taking $$\hat{H}=H
  \cap (\Theta \times \{-\beta,+\beta\}),$$ $$\hat{H}_0=\hat{H}\cap
  H_0,$$ $$\hat{H}_1=\hat{H}\cap H_1,$$ and using the $\hat{H}$
  variables in place of the $H$ variables; in other words we will
  assume that we know to start with that $h$ is either $-\beta$ or
  $+\beta$, with $\hat{H}_0$ being that $h=-\beta$.

We set the standard (and uniformly optimal) critical regions
$$C_\eta=[G^{-1}(\eta) - \beta, +\infty).$$ Suppose that $x=0$. Then
  we are $\alpha$-frequentist-confident that $h>0$, indeed that
  $h=+\beta$. But now
  setting $$\hat{H}_0'=\hat{H}_1,$$ $$\hat{H}_1'=\hat{H}_0,$$ so that
  $\hat{H}_0'$ is now that $h=+\beta$, and again forming the uniformly
  optimal critical regions $$C_\eta'=(-\infty,G^{-1}(1-\eta)+\beta],$$
and again observing $x=0$, we find that we are now
$\alpha$-frequentist-confident that $h\leq 0$, indeed that
$h=-\beta$. 

But $\alpha=0.975$, which is near 1, and nowhere near zero, and
nowhere near $1-\alpha=0.025$, violating the criterion.

\item Criterion \ref{containment}
  (``Inclusion''): \label{monotonicity} Let $H_0'=\{0\}$ so that
  $H_1'=\{h\in \mathbb{R}:h\neq 0\}\supset H_1$. Take the standard
  critical regions, namely
  $$C_\eta = [G^{-1}(\eta), +\infty)$$
    and $$C_\eta'=(-\infty,G^{-1}(\frac{1-\eta}{2})] \cup [G^{-1}(1 -
    \frac{1-\eta}{2}),+\infty),$$ noting that frequentists design the
    nested set of critical regions for a particular $H_0$, and that
    changing $H_0$ standardly leads to a different design of the
    critical regions. Suppose that $x=\beta$. Then $c=\sup\{\eta\in
    [0,1]: x\in C_\eta\}=\alpha$ but $c'=\sup\{\eta\in [0,1]: x\in
    C_\eta'\} = 0.95 < \alpha = 0.975$, so $c'<c$ even though
    $H_1'\supset H_1$.

\item Criterion \ref{whetherproceed} (``Intention''): Consider the
  following two data collection plans:
\begin{enumerate}

\item $D$: We collect exactly 1 sample, $x_1$, and
  conclude if $x_1 \geq \beta$, otherwise we never
  conclude. (I.e. $q_0=0,q_1(x_1)=[x_1 \geq \beta],q_2=q_3=...=0$.)

\item $D'$: We collect a single sample $x_1$. If
  $x_1\geq \beta$ then we stop data collection. Otherwise we collect a
  second sample and exit if $x_2\geq \beta$; otherwise we never
  conclude. (I.e. $q_0=0,q_1(x_1)=[x_1\geq \beta],q_2(x_1,x_2)=[x_2
    \geq \beta], q_3=q_4=...=0.$)

\end{enumerate}

Then for $D$ we have $$C_\eta=\{\mathbf{x}:x_1 \geq \beta\},$$ so that
$\eta=\alpha$. For $D'$ we have $$C_\eta' = \{\mathbf{x}: (x_1 \geq
\beta) \lor \left((x_1 < \beta) \land (x_2 \geq \beta)\right)\},$$ so
that $\eta=1 - (1 - \alpha + \alpha \times (1 - \alpha)) \approx
0.9506 < \alpha = 0.975$.

Now suppose that we obtain $x_1=\beta$. Whether we were using $D$ or
$D'$ we stop data collection. If we were using $D$ we conclude that we
are $\alpha$-confident that $h\in H_1$. If we were using $D'$ we
conclude that we are 0.9506-confident that $h\in H_1$. So our
conclusion depends on which data collection plan we were intending to
use, even though the data we have collected (namely $x_1$) is the same
in both cases, violating the criterion.

\item Criterion \ref{productprob} (``Conjunction''): Let $N=2$ and
  suppose we observe $x_1=\beta$ on system $1$ and $x_2=\beta$ on
  system $2$, so that we are $\alpha$-confident that $h_1>0$ and
  $\alpha$-confident that $h_2>0$. Let us consider how confident we
  are that \textit{both} $h_1$ and $h_2$ are $>0$.

Let $H_0=\{(h_1,h_2):(h_1\leq 0) \lor (h_2\leq 0)\}$ so that
$H_1=\{(h_1,h_2):h_1>0,h_2>0\}$. The corresponding uniformly optimal
critical region of maximal frequentist confidence for this data is
$C_\eta=\{(x_1,x_2):x_1\geq \beta,x_2\geq \beta\}$. Then for $h\in
H_0$ the supremum of $P((x_1,x_2)\in C_\eta|h)$ is approached as
$h\to(+\infty,0)$ giving $\sup_{h\in H_0}{P((x_1,x_2)\in C_\eta|h)}=1
- \alpha$ so that $\eta=\alpha$ -- i.e. we are equally sure that both
the $h_n$ are in their respective $H_{1,n}$ as we are that either one
alone is, despite the fact that $\alpha < 1$, violating the criterion.

\item Criterion \ref{sumprob} (``Disjunction'') (also
  \ref{containment}): Again let $N=2$ and suppose we observe
  $x_n=\beta$ on both systems, so that on each system we are
  $\alpha$-confident that $h_n>0$ using the standard critical regions
  $C_{n,\eta}=(G^{-1}(\eta), +\infty)$. Let us consider how
  frequentist confident we are that \textit{at least one} of $h_1$ and
  $h_2$ are $>0$.

Let $H_0=\{(h_1,h_2):h_1\leq 0, h_2\leq 0\}$ so that
$H_1=\{(h_1,h_2):(h_1 > 0) \lor (h_2 > 0)\}$. Then the standard set of
critical regions are $$C_{\zeta(\eta)}=\left((G^{-1}(\eta),+\infty)
\times \mathbb{R}\right) \cup \left(\mathbb{R}\times
(G^{-1}(\eta),+\infty)\right)$$ for some appropriate function $\zeta$,
so that $$(x_1,x_2)\in C_{\zeta(\eta)} \iff \left((x_1 \in
C_{1,\eta}) \lor (x_2 \in C_{2,\eta})\right).$$ To find $\zeta(\eta)$,
note that for $(h_1,h_2)\in H_0$,
\begin{IEEEeqnarray*}{rCl}
P((x_1,x_2)\in C_{\zeta(\eta)}|h_1,h_2) &\leq& P((x_1,x_2)\in C_{\zeta(\eta)}|h_1=0,h_2=0)\\
&=& 1 - P((x_1,x_2)\notin C_{\zeta(\eta)}|h_1=0,h_2=0)\\
&=& 1 - P(x_1\notin C_{1,\eta}|h_1=0)P(x_2\notin C_{2,\eta}|h_2=0)\\
&=& 1 - (1 - (1 - \eta))^2\\
&=& 1 - \eta^2,
\end{IEEEeqnarray*}
and that $(0,0)\in H_0$, so that $\zeta(\eta)=\eta^2$, and so we are
only $\alpha^2$-sure, i.e. strictly less than $\alpha$-sure, that at
least one of $h_1\in H_{1,1}$ and $h_2\in H_{1,2}$ holds. But the
criterion required that we should be strictly more than $\alpha$-sure
of this, so not only is the criterion violated, but for the example in
question and many other similar problems we actually have an
inequality in the opposite direction.

\item Criterion \ref{multisystem} (``Multiplicity''): Here we offer
  not a proof but a description of standard frequentist practice. This
  requires a penalty to compensate for multiple hypothesis tests
  \cite{MultipleComparisons}. If $N$ independent frequentist
  hypothesis tests are done, each giving confidence $c$, and $H_{0,n}$
  holds for all $n$, then each test has a probability of up to $1-c$
  of concluding that its relevant $H_{1,n}$ is true, making the bound
  on the overall expected number of false positives $N(1-c)$. In order
  to ensure that the overall expected number of false positives is
  only $(1-c)$, frequentists usually apply a correction or penalty
  making it harder to conclude that $H_{1,n}$ holds, for each $n$. The
  most often used correction is the Bonferroni correction
  \cite{Bonferronioriginal, Bonferroni}, but there is controversy
  about which (if any) correction should be used when. The magnitude
  of the penalty increases as $N$ does.

  Or to illustrate this by anecdote: Adam compares height of
  Icelanders with that of Swedes, and finds a significant difference
  so publishes the result; Beth compares height of Swedes with that of
  Norwegians and finds a significant difference, so publishes the
  result; Charles compares height of Latvians with that of Lithuanians
  and finds a significant difference, so publishes the result. But
  Laura is more hard-working than the others, so she does all three
  comparisons, but finds no significant result because of corrections
  for multiplicity, so can't get her work published. Poor Laura !

\item Criterion \ref{unifpow} (``Sequential Optimality''): Let $y_k$
  denote the mean of $x_1,x_2,...,x_k$.

Suppose that we have such a uniformly most powerful data collection
plan $D$ for the specific desired confidence $\eta=\alpha$. The idea
will be to reach a contradiction by finding an alternative plan $D'$,
time $k$, and value of $h\in H_1$, such that $p_k(D',h)>p_k(D,h)$, while
nonetheless only exiting the plan if we can conclude that we are still
$\alpha$-frequentist-confident that $h\in H_1$.

Suppose first that $q_0=1-\alpha=0.025$, i.e. that without collecting
any data we conclude with probability $0.025$ that we are
$\alpha$-confident that $h\in H_1$, while with the remaining probability
we proceed. Then for the plan to be valid (except in the trivial case
that $H_0=\emptyset$) we are forced to have $q_1=q_2=...=0$ a.e. . (Note in
passing that this results in a data collection plan $D_0$ that while
perverse and almost useless is nonetheless frequentistly valid, even
though no frequentist would use it in practice.)

Then for $h=2\beta$ we can improve on the probability of $0.025$ for
concluding before collecting $x_2$ by setting $q_1(x_1) = [x_1
  \geq \beta]$ and reducing $q_0$ to $0$, getting us a probability of
$\alpha = 0.975 > 0.025$ that we conclude before collecting $x_2$,
contrary to the uniformly most powerful nature of $D$.

Otherwise $q_0 < 1 - \alpha$, in which case we can consider the
previous (perverse and almost useless) data collection plan $D_0$,
which gives an increased probability of concluding without collecting
any data. Therefore $D$ was not a uniformly most powerful data
collection plan, completing the proof that such a $D$ does not exist
for this example.

Now, the reader may consider that allowing $q_0$ to be non-zero is
itself perverse, and we therefore provide a further proof under the
restriction that we must always have $q_0=0$, as follows.

Suppose that we have such a uniformly most powerful data collection
plan $D$ with $q_0=0$ for the specific desired confidence
$\alpha=0.975$. 

Then for $D$ to be valid it is necessary that $\E(q_1(x_1)|h=0)\leq 1
- \alpha$, and unless we have equality here we can increase $q_1$ by
some constant and set $q_k=0$ for $k > 1$, increasing
$p_1(D,h=\epsilon)$ for some small $\epsilon>0$ and contradicting $D$
being uniformly most powerful. Therefore \mbox{$\E(q_1(x_1)|h=0) = 1 -
  \alpha$,} and $q_k=0$ for $k\geq 2$.

Now, similarly to the previous proof, suppose first that $q_1(x_1) =
[x_1 \geq G^{-1}(\alpha)] = [x_1 \geq \beta]$ with $q_k=0$ for
$k\geq 2$; let us call this specific plan $D_1$. Then with $D_1$, the
probability given $h=2\beta$ of concluding before collecting $x_3$
is the same as that of concluding before collecting $x_2$, i.e. it
is $1 - G(G^{-1}(\alpha) - h) = \alpha$.

But if instead we set $q_1(x_1) = [x_1 \geq G^{-1}(1 - \frac{1 -
    \alpha}{2})] = [x_1 \geq \gamma]$ and $q_2(x_1,x_2) = [y_2 \geq
  \frac{\gamma}{\sqrt{2}}]$ (with $q_k=0$ for $k > 2$), then the
probability given $h=2\beta$ of concluding before collecting $x_2$ is
$\Phi(2\beta - \gamma)\approx 0.95338$ and the probability of
concluding just after collecting $x_2$ but not before that is about
$0.046191$ (based on $10^8$ simulations), giving a total probability
of concluding before collecting $x_3$ of $0.99957 > \alpha$, so $D_1$
was not uniformly most powerful.

Trivially, if for all $k$, $q_k=q_k'\text{ a.e.}$, then also
$p_k(D,h)=p_k(D',h)$. So the only remaining case to consider is that
$D$ is not a.e. equal to $D_1$; in this case we will show that for
$h=2\beta$, $p_1(D,2\beta)<p_1(D_1,2\beta)$, contradicting the
uniformly most powerful nature of $D$:

Taking $q_k$ to denote the quit-probability functions of the plan $D$,
we have shown above
that $$\int_0^\infty{P(x_1|h=0)q_1(x_1)\,dx_1}=\E(q_1(x_1)|h=0)
=1-\alpha$$ and that for $k\geq 2$, $q_k=0$; so for the remaining
argument we may drop the subscript 1s from $q_1$ and $x_1$ for
simplicity. Looking for a contradiction, suppose
$$\int_0^\beta{P(x|h=0)q(x)\,dx}=0.$$ Then also
$$\int_\beta^\infty{P(x|h=0)q(x)\,dx}=1-\alpha.$$ But
$$\int_\beta^\infty{P(x|h=0)\,dx}=1-\alpha$$ and $0\leq q(x)\leq
1$. Therefore $q(x)=[x\geq \beta]$ a.e., and by the previous remarks
$D$ is effectively the same as $D_1$ so is not uniformly most
powerful.

Therefore $\int_0^\beta{P(x|h=0)q(x)\,dx}>0$.

Now let $q'(x)=[x\geq \beta]$, so
that $$\int_0^\infty{P(x|h=0)q'(x)\,dx}=1-\alpha$$ also. Let
$f(x)=q'(x)-q(x)$, so that $$\int_0^\infty{P(x|h=0)f(x)\,dx}=0.$$ Then
on $x<\beta$, $f(x)\leq 0$, while on $x\geq\beta$, $f(x)\geq 0$. 

Moreover on $0\leq x<\beta$, $P(x|h=0)> P(x|h=2\beta)>0$, while
on $x> \beta$, $0 < P(x|h=0)< P(x|h=2\beta)$.

Therefore
\begin{IEEEeqnarray*}{rCcCcCl} 
0 & \geq & \int_0^\beta{P(x|h=2\beta)f(x)\,dx} & \geq &
\int_0^\beta{P(x|h=0)f(x)\,dx}\\ \text{and}&&\int_\beta^\infty{P(x|h=2\beta)f(x)\,dx}&\geq&
\int_\beta^\infty{P(x|h=0)f(x)\,dx} & \geq 0.
\end{IEEEeqnarray*}

Summing these inequalities we deduce
that $$\int_0^\infty{P(x|h=2\beta)f(x)\,dx}\geq
\int_0^\infty{P(x|h=0)f(x)\,dx},$$ with equality only if both summand
inequalities are equalities, and we already know that the RHS is
zero. But since on $[0,\beta)$, $P(x|h=2\beta) - P(x|h=0)$ is negative
  and continuous (so that on any closed subinterval it is bounded away
  from zero) and $f$ does not change sign, equality of the first
  summand implies that on $0\leq x \leq \beta$, $f(x)=0$ a.e.; but
  here $f(x)=-q(x)$, so $q(x)=0$ a.e., which has been excluded since
  $\int_0^\beta{P(x|h=0)q(x)\,dx}>0$.

Therefore $$\int_0^\infty{P(x|h=2\beta)f(x)\,dx}>0,$$
so that $$\int_0^\infty{q'(x)P(x|h=2\beta)\,dx} >
\int_0^\infty{q(x)P(x|h=2\beta)\,dx},$$
i.e. $p_1(D',2\beta)>p_1(D,2\beta)$, contradicting the uniformly most
powerful nature of $D$.

Thus $D$ is not uniformly most powerful, contradicting the assumption,
and completing the proof for the case that $q_0$ is required to be
$0$. 

We note further that the example used is one that, for a
\textit{fixed} amount of data collected, does have a uniformly most
powerful critical region -- and that even this property therefore
does not imply the existence of a uniformly most powerful data
collection and analysis plan.

\end{enumerate}

\section{Appendix: A very simple worked example calculation by each method}
\label{simpleexample}

We here give an extremely simple example problem and solutions by each
method for the benefit of anybody who has not encountered inference
problems in abstract before. No mathematics harder than integrating
$1/x$ will be needed to follow the working - one just needs to read it
slowly enough to make sure that each line has been understood. Indeed
for any reader unfamiliar with calculus, we note that for
$b>a$, $$\int_a^b{f(x)\,dx}$$ simply means the area under the curve
$y=f(x)$ between $x=a$ and $x=b$ (counting it negative if it is below
the $x$-axis), and that for
$b>a>0$, $$\int_a^b{\frac{1}{x}\,dx}=\log\frac{b}{a},$$ where $\log$
denotes the natural logarithm to the base $e$.

\subsection{The problem}

Two real numbers $\phi_1$ and $\phi_2$ are transmitted over a noisy
channel. The channel adds noise $n_1$ and $n_2$ respectively, where
$(n_1,n_2)$ is uniformly distributed over the unit square centred on
the origin, so that the probability density of the noise is given
by $$P(n_1,n_2)=\left[-\frac{1}{2}\leq n_1 \leq
  +\frac{1}{2}\right]\land \left[-\frac{1}{2}\leq n_2 \leq
  +\frac{1}{2}\right],$$ where $[\ ]$ denotes the function that takes
the value $1$ when the statement inside the brackets is true and $0$
when it is false, and where $a \land b$ denotes the minimum of $a$ and
$b$ (and $a \lor b$ will denote the maximum)\footnote{Note that
  $\land$ is also used to mean ``and'', and $\lor$ is also used to
  mean ``or''; if ``true'' is represented by $1$ and ``false'' by $0$,
  the two meanings each of $\land$ and $\lor$ are compatible.}.

We then observe $$x=(x_1,x_2)=(\phi_1+n_1,\phi_2+n_2),$$ and want to
know whether $H_0$ or $H_1$ is true, where $H_0$ is that both
$\phi_1\leq 0$ and $\phi_2\leq 0$, while $H_1$ is that
$(\phi_1>0)\lor(\phi_2>0)$. We are told that our solution need only
deal with $x$ in the square of side 3 centred on the origin (both for
simplifying calculation and so that we can display our solutions
graphically without losing detail near the origin).

Putting all this into our standard notation for inference problems
from section \ref{infprobdefsub}, we are
setting $$\Theta=\{0,1\},$$ $$\Phi=\mathbb{R}^2,$$ $$H=\{(\theta,\phi_1,\phi_2)\in
\Theta\times\Phi:\theta=\left[(\phi_1>0)\lor(\phi_2>0)\right]\},$$ $$H_0=\{(0,\phi_1,\phi_2)\in
\Theta\times\Phi:(\phi_1\leq 0) \land (\phi_2\leq
0)\},$$ $$H_1=\{(1,\phi_1,\phi_2)\in \Theta\times\Phi:(\phi_1> 0) \lor
(\phi_2> 0)\},$$
$$X=\mathbb{R}^2,$$ and $$P(x_1,x_2|\theta,\phi_1,\phi_2) =
P(x_1,x_2|\phi_1,\phi_2) = \left[|x_1-\phi_1|\leq \frac{1}{2}\right]
\left[|x_2-\phi_2|\leq \frac{1}{2}\right].$$

Note that we can afford to identify $h$ with $(\phi_1,\phi_2)$
(ignoring $\theta$) as $\theta$ is uniquely determined by whether
$(\phi_1,\phi_2)$ is in the bottom left quadrant of the plane or not. 

We now turn to calculate the Bayesian solution and two different
frequentist solutions; then we will consider what happens if we decide
to be frequentists but don't choose which of these two sets of
critical regions we are going to use until after seeing the data; and
finally we will look at a simple pseudo-Bayesian (and therefore
frequentist) solution.

\subsection{The Bayesian solution}
\label{simpleBayes}

As Bayesians we first have to choose a prior on $H$, i.e. on
$(\phi_1,\phi_2)$. We choose the uniform distribution on the square of
side $K$ centred on the origin, where $K$ is some number greater than
or equal to 4 -- it will turn out that so long as $K\geq 4$ it will
make no difference at all to the solution for values of $(x_1,x_2)$ in
the square of interest $$\left[-\frac{3}{2},+\frac{3}{2}\right]^2$$
(where here $[a,b]$ denotes the closed interval from $a$ to $b$, and
the $^2$ denotes the Cartesian product with itself). In other words we
are choosing to
set $$P(\phi_1,\phi_2)=\frac{1}{K^2}\left[-\frac{K}{2}\leq \phi_1\leq
  +\frac{K}{2}\right]\left[-\frac{K}{2}\leq \phi_2\leq
  +\frac{K}{2}\right].$$

Now we can ``turn the handle of Bayesian inference'' (as David MacKay
would have put it). We apply Bayes' theorem to deduce
that $$P(\phi_1,\phi_2|x_1,x_2)=\frac{P(\phi_1,\phi_2)P(x_1,x_2|\phi_1,\phi_2)}
{\int{P(\phi_1,\phi_2)P(x_1,x_2|\phi_1,\phi_2)\,d(\phi_1,\phi_2)}},$$
where the range of integration is over all possible values of
$(\phi_1,\phi_2)$, i.e. over all of $\mathbb{R}^2$.

We then notice that the second factor in the integrand is zero except
for $(\phi_1,\phi_2)$ in the square of side $1$ centred on
$(x_1,x_2)$, so
that $$P(\phi_1,\phi_2|x_1,x_2)=\frac{P(\phi_1,\phi_2)P(x_1,x_2|\phi_1,\phi_2)}
{\int_{x_1-\frac{1}{2}}^{x_1+\frac{1}{2}}{\int_{x_2-\frac{1}{2}}^{x_2+\frac{1}{2}}{P(\phi_1,\phi_2)P(x_1,x_2|\phi_1,\phi_2)\,d\phi_2}\,d\phi_1}},$$
and similarly the second factor in the numerator is only non-zero on
the same square.

Now, bearing in mind that we are told we are only interested in
$x_1,x_2$ of absolute value less than $\frac{3}{2}$, we notice that for all
$(\phi_1,\phi_2)$ of interest, $P(\phi_1,\phi_2)=\frac{1}{K^2}$, so
that these factors cancel and we
have $$P(\phi_1,\phi_2|x_1,x_2)=\frac{P(x_1,x_2|\phi_1,\phi_2)}
{\int_{x_1-\frac{1}{2}}^{x_1+\frac{1}{2}}{\int_{x_2-\frac{1}{2}}^{x_2+\frac{1}{2}}
    {P(x_1,x_2|\phi_1,\phi_2)\,d\phi_2}\,d\phi_1}}.$$ Further, for all
$(\phi_1,\phi_2)$ in the range of integration,
$P(x_1,x_2|\phi_1,\phi_2)=1$, so that the integral in the denominator
is $1$ and we are left
with $$P(\phi_1,\phi_2|x_1,x_2)=P(x_1,x_2|\phi_1,\phi_2).$$ (The
reader should most definitely not assume that such a simple
relationship always holds in every situation !)

\begin{figure}[hpt]
\begin{center}

\includegraphics[scale=0.4]{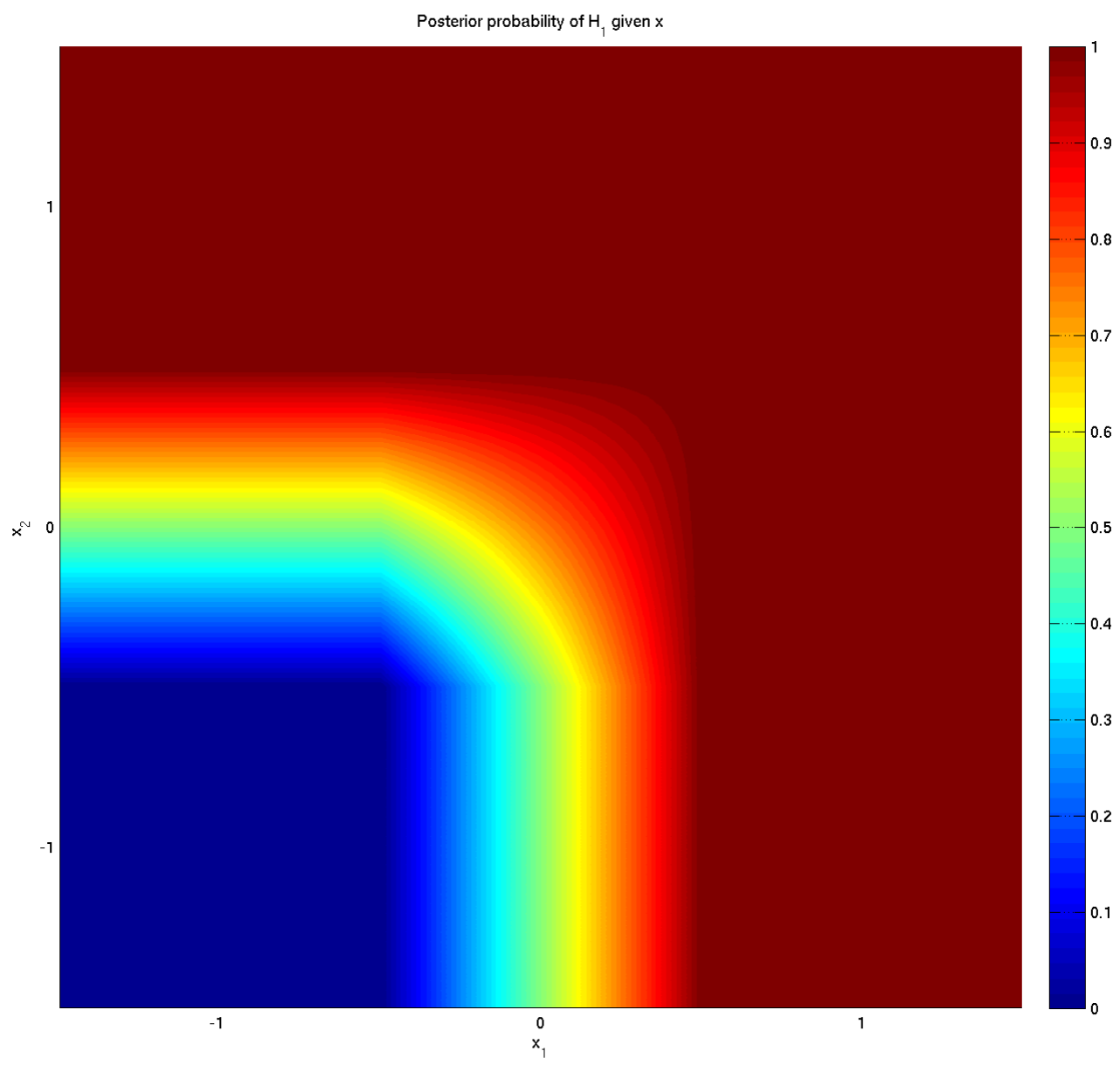}

\caption{
\label{simpleBayesfig}
Bayesian solution to the simple problem of appendix
\ref{simpleexample}. The posterior probability of $H_1$ is plotted as
a function of the observed data $(x_1,x_2)$.
}

\end{center}
\end{figure}

We can now integrate over $(\phi_1,\phi_2)\in
H_0=(-\infty,0]\times(-\infty,0]$ to
get $$P(H_0|x_1,x_2)=\int_{-\infty}^0{\int_{-\infty}^0{P(x_1,x_2|\phi_1,\phi_2)\,d\phi_2}\,d\phi_1}.$$
Taking note of the special form of the
integrand $$P(x_1,x_2|\phi_1,\phi_2) = \left[|x_1-\phi_1|\leq
  \frac{1}{2}\right] \left[|x_2-\phi_2|\leq \frac{1}{2}\right],$$ we
find that $P(H_0|x_1,x_2)$ is zero for $x_1>\frac{1}{2}$ or
$x_2>\frac{1}{2}$, is one if both $x_1\leq-\frac{1}{2}$ and
$x_2\leq-\frac{1}{2}$, and otherwise is given
by $$P(H_0|x_1,x_2)=\int_{0\land(x_1-\frac{1}{2})}^{0\land (x_1+\frac{1}{2})}{
  \int_{0\land(x_2-\frac{1}{2})}^{0\land(x_2+\frac{1}{2})}{1\,d\phi_2}\,d\phi_1} =
\left(0\lor\left(\frac{1}{2}-x_1\right)\land
1\right)\left(0\lor\left(\frac{1}{2}-x_2\right)\land 1\right).$$ Putting
these three possibilities together we find that $$P(H_0|x_1,x_2) =
\left(0\lor\left(\frac{1}{2}-x_1\right)\land
1\right)\left(0\lor\left(\frac{1}{2}-x_2\right)\land 1\right),$$ and
hence $$P(H_1|x_1,x_2) = 1-
\left(0\lor\left(\frac{1}{2}-x_1\right)\land
1\right)\left(0\lor\left(\frac{1}{2}-x_2\right)\land 1\right),$$
completing the Bayesian solution, which is then shown in figure
\ref{simpleBayesfig}.

\clearpage

\subsection{First frequentist solution}
\label{simplefreq1}

\begin{figure}[p]
\begin{center}

\includegraphics[scale=0.4]{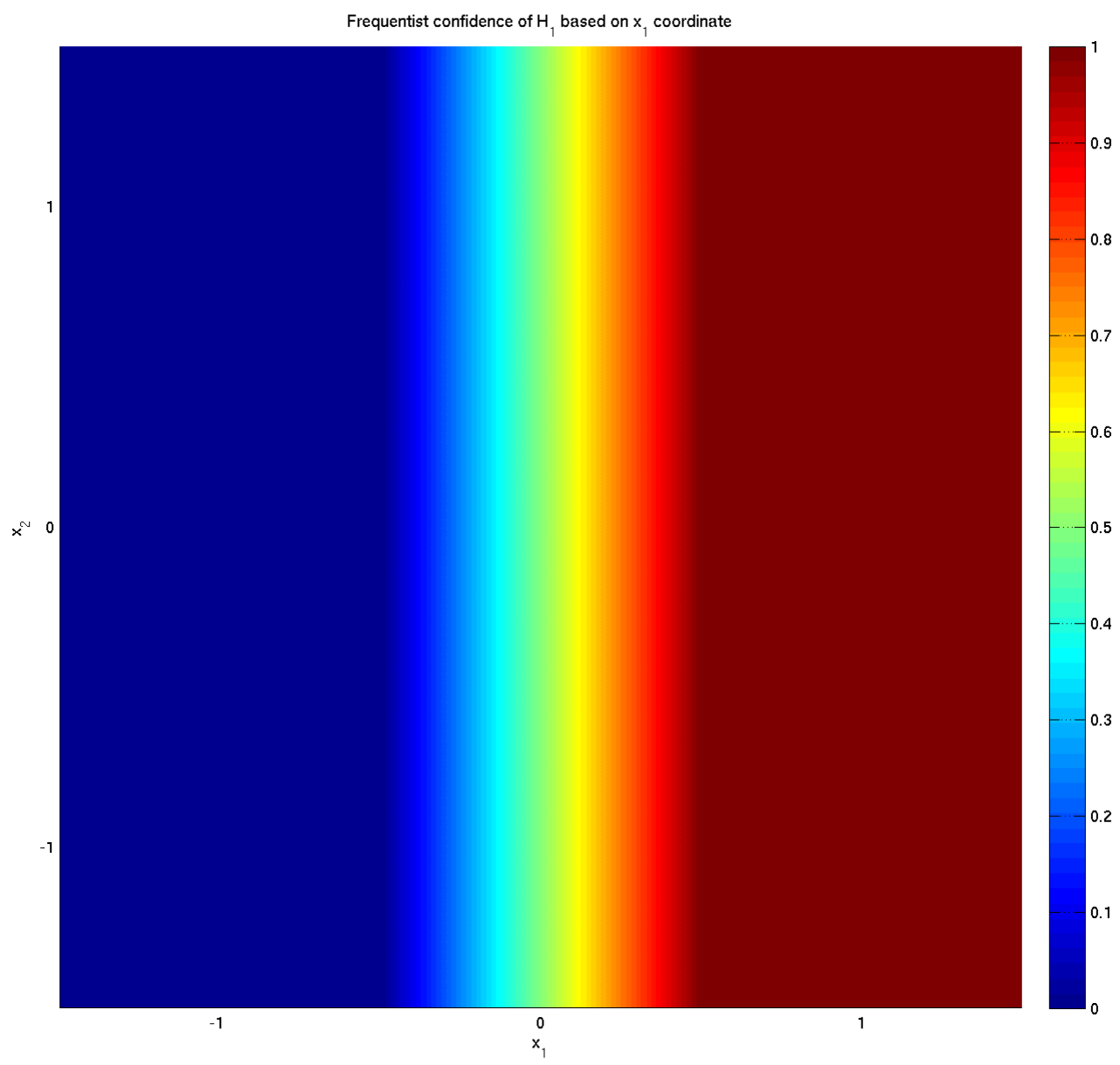}

\caption{
\label{simplefreq1fig}
First frequentist solution to the simple problem of appendix
\ref{simpleexample}. The frequentist confidence in $H_1$ is plotted as
a function of the observed data $(x_1,x_2)$.  }

\end{center}
\end{figure}

For our first frequentist solution we decide (arbitrarily) to use
critical regions of the form $$B_t=\{(x_1,x_2)\in X:x_1>t\},$$ which
form a nested set decreasing as $t$ increases. We then have to
calculate the value of the frequentist confidence $\eta(t)$
corresponding to each value of $t$. This is given by $$\eta(t)=1 -
\sup_{h\in H_0}P((x_1,x_2)\in B_t|h).$$ Since $h=(\phi_1,\phi_2)$ this
is saying that $$\eta(t)=1 - \sup_{(\phi_1,\phi_2)\in
  H_0}P((x_1,x_2)\in B_t|\phi_1,\phi_2).$$ Since whether or not
$(x_1,x_2)\in B_t$ depends only on $x_1$ which in turn depends only on
$\phi_1$, this simplifies to $$\eta(t)=1 - \sup_{\phi_1\leq
  0}P(x_1>t|\phi_1).$$ The probability on the RHS is maximised by
setting $\phi_1=0$, and then
$$P(x_1>t|\phi_1=0)=\left(0\lor\left(\frac{1}{2}-t\right)\land
1\right),$$ so that $$\eta(t)=1 -
\left(0\lor\left(\frac{1}{2}-t\right)\land 1\right) =
\left(0\lor\left(t+\frac{1}{2}\right)\land 1\right),$$ which function
is shown in figure \ref{simplefreq1fig}.

\subsection{Second frequentist solution}
\label{simplefreq2}

\begin{figure}[hpt]
\begin{center}

\includegraphics[scale=0.4]{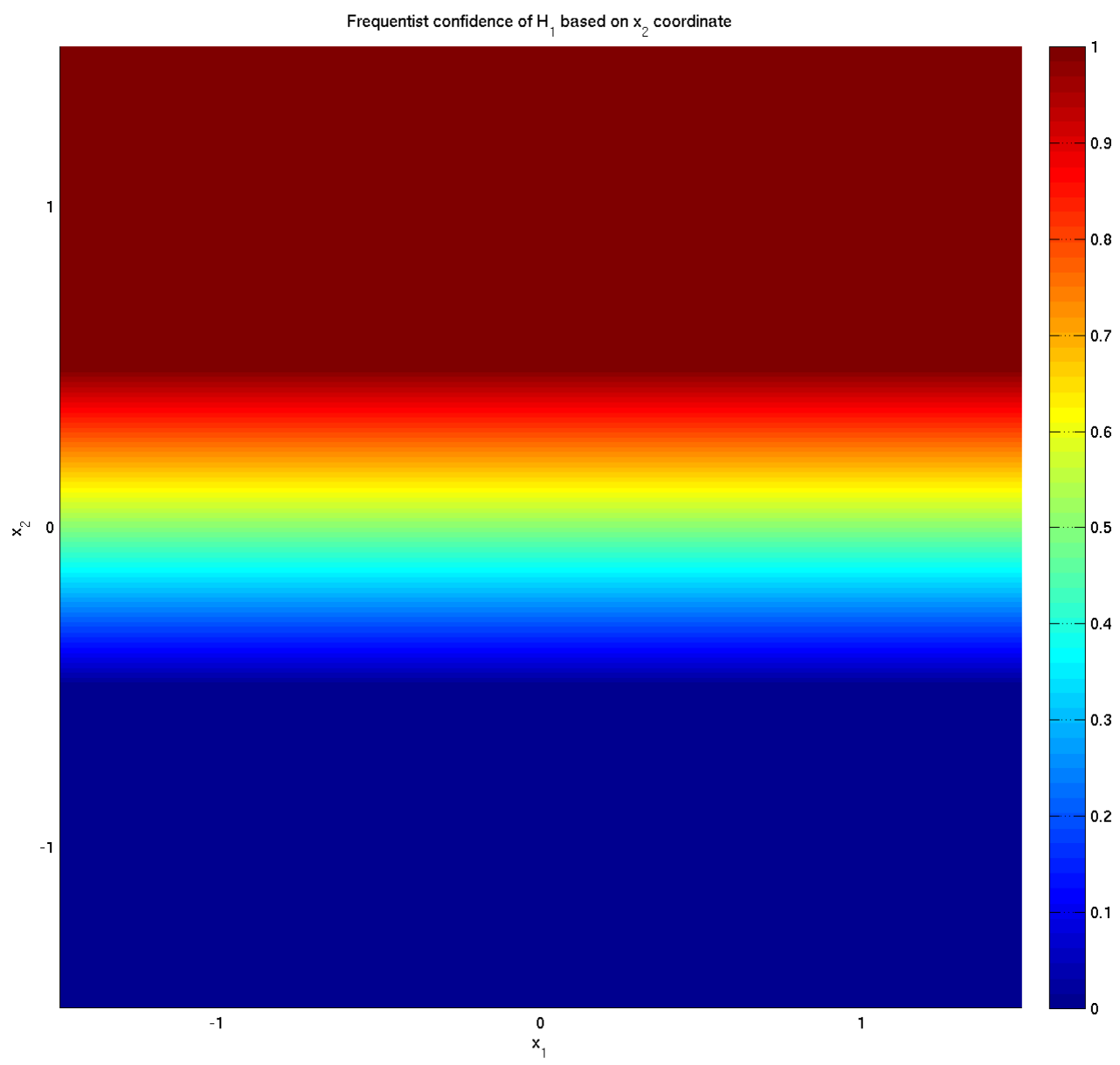}

\caption{
\label{simplefreq2fig}
Second frequentist solution to the simple problem of appendix
\ref{simpleexample}. The frequentist confidence in $H_1$ is plotted as
a function of the observed data $(x_1,x_2)$.  }

\end{center}
\end{figure}

For our second frequentist solution we decide arbitrarily to use
critical regions depending only on $x_2$ in the same way that those of
section \ref{simplefreq1} depended only on $x_1$. That gives the
frequentist confidence shown in figure \ref{simplefreq2fig}.

\subsection{What if we break the fundamental frequentist rule ?}

Suppose now that we have an interest in getting the answer to be that
$H_1$ holds with a high degree of frequentist confidence, and that we
have only (so far) thought of the solutions in sections
\ref{simplefreq1} and \ref{simplefreq2} above. But we don't know
whether our data might be $(2,-1)$ (in which case the solution of
\ref{simplefreq1} would be best) or $(-1,2)$ (in which case the solution
of \ref{simplefreq2} would be best).

Something that might occur to us is to decide which to use not before
but after collecting the data, in such a way as to choose whichever
gives us greatest frequentist confidence in $H_1$ (which of course
breaks the frequentist rules laid down in section \ref{freqmeth}).

\begin{figure}[p]
\begin{center}

\includegraphics[scale=0.4]{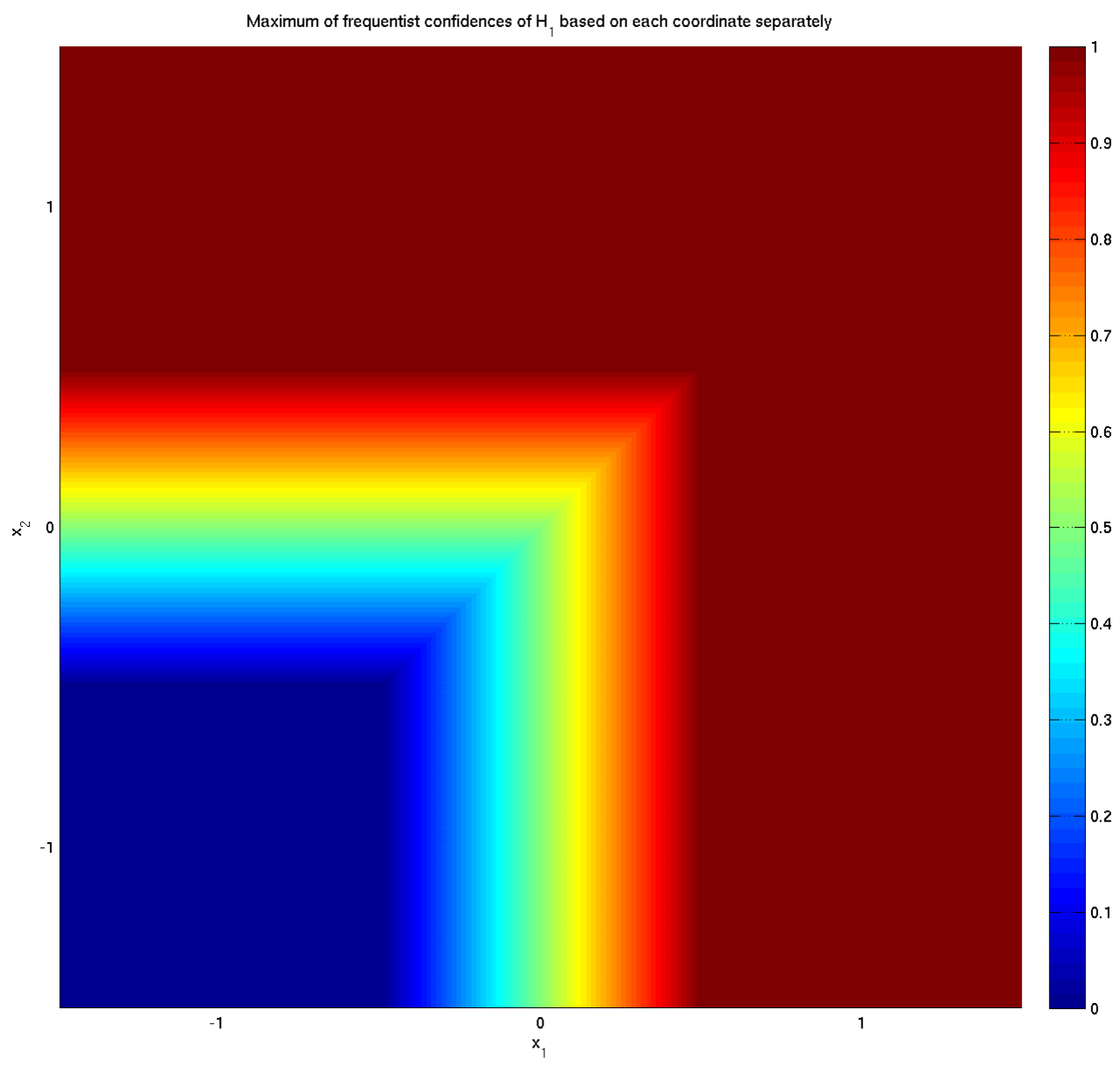}

\caption{
\label{simplefreqlatefig}
The consequence of choosing whether to use the solution of section
\ref{simplefreq1} or \ref{simplefreq2} \textit{after} collecting the
data. The apparent (but wrong) frequentist confidence in $H_1$ is
plotted as a function of the observed data $(x_1,x_2)$.  }

\end{center}
\end{figure}

That would result in the apparent frequentist confidence, as a
function of $(x_1,x_2)$ shown in figure \ref{simplefreqlatefig}. 

But the problem is that the condition in section \ref{freqmeth} that
for all $h\in H_0$ and all $\eta\in [0,1]$, $P(x\in C_\eta|h)\leq
1-\eta$ no longer holds. For example, for $h=(\phi_1,\phi_2)=(0,0)$
and $\eta=\frac{1}{2}$ we have $P(x\in C_\eta|h)=\frac{3}{4}\not\leq 1
- \frac{1}{2}$ as then $C_\eta=\{(x_1,x_2)\in X:(x_1>0)\lor(x_2>0)\}$.

However, we can still use this particular nested set of critical
regions $$B_t=\{(x_1,x_2)\in X:(x_1>t)\lor(x_2>t)\},$$ we just have to
downgrade the frequentist confidence associated with each
one. Checking the conditions needed, we proceed as follows.

\begin{figure}[p]
\begin{center}

\includegraphics[scale=0.4]{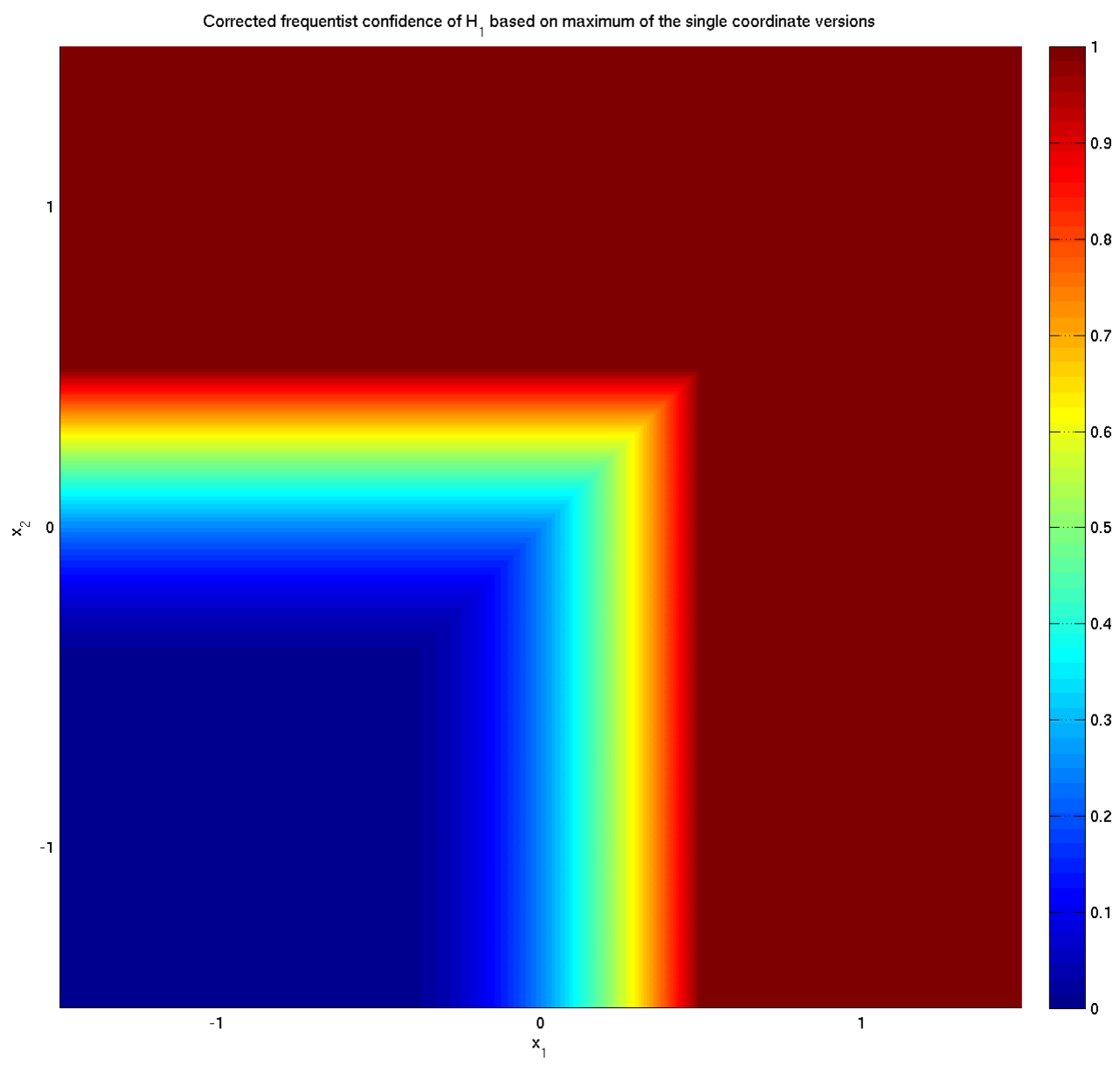}

\caption{
\label{simplefreqlatefixedfig}
The consequence of choosing whether to use the solution of section
\ref{simplefreq1} or \ref{simplefreq2} \textit{after} collecting the
data, corrected for the then necessary adjustment of frequentist
confidence that results. The resulting frequentist confidence in $H_1$
is plotted as a function of the observed data $(x_1,x_2)$.  }

\end{center}
\end{figure}

We want to know, for each $t$, the value of $$1 - \sup_{h\in
  H_0}{P((x_1,x_2)\in B_t|h)} = 1 - \sup_{(\phi_1,\phi_2)\in
  H_0}{P((x_1,x_2)\in B_t|(\phi_1,\phi_2))}.$$ Inspection of figure
\ref{simplefreqlatefig} shows us that the supremum is achieved for
$(\phi_1,\phi_2)=(0,0).$ Then $$P((x_1,x_2)\in B_t|\phi_1=0,\phi_2=0)
= 2 \zeta - \zeta^2,$$ where $$\zeta=\left(0\lor
\left(\frac{1}{2}-t\right)\land 1\right).$$ Therefore $$\eta(t) = 1 -
2\zeta + \zeta^2 = (1 - \zeta)^2 =
\left(0\lor\left(t+\frac{1}{2}\right)\land 1\right)^2.$$ This gives us
the corrected plot in figure \ref{simplefreqlatefixedfig}, but as can
be seen this substantially reduces our chance of getting a high value
of frequentist significance compared with figure
\ref{simplefreqlatefig}. Note that if we had straight off decided to
use these particular $B_t$ then we would have arrived at the same set
of critical regions $C_{\eta(t)}=B_t$ that we have just obtained.

\subsection{Pseudo-Bayesian solution}

Alternatively we could consider a pseudo-Bayesian solution (in this
case the basic and full pseudo-Bayesian solutions coincide). Keeping
the original $H_0$ and $H_1$, we simply calculate the frequentist
confidence that goes with each subset $$B_p=\{(x_1,x_2)\in
X:P(H_1|x_1,x_2)\geq p\},$$ noting that as always these form a nested
family, decreasing as $p$ increases.

Letting $p=P(H_1|x_1,x_2)$, from section \ref{simpleBayes}, we
have $$p=1-\left(0\lor\left(\frac{1}{2}-x_1\right)\land
1\right)\left(0\lor\left(\frac{1}{2}-x_2\right)\land 1\right).$$ By
inspection of figure \ref{simpleBayesfig}, the probability
$P((x_1,x_2)\in B_p|h)$ is maximised at $h=(\phi_1,\phi_2)=(0,0)$, so
we are then interested in $P((x_1,x_2)\in B_p|\phi_1=0,\phi_2=0)$ for
$|x_1|,|x_2|\leq \frac{1}{2}$. The above relationship then simplifies
to $$p=1-\left(\frac{1}{2}-x_1\right)\left(\frac{1}{2}-x_2\right),$$
so by rearrangement for fixed $p$ we
find $$x_2=\frac{1}{2}-\frac{1-p}{\frac{1}{2}-x_1},$$ giving us a
curve bounding this part of $B_p$.

Now, $B_p$ is the region above and to the right of this curve, so
$P(x\in B_p|h=(0,0))$ is one minus the area that is under this curve
and above the line $x_2=-\frac{1}{2}$ between $x_1=-\frac{1}{2}$ and
$x_1=+\frac{1}{2}$. Thus the frequentist confidence $\eta(p)$ that we
can give to $H_1$ is equal to that area under the curve, since it is
one minus the supremum of that probability.

\begin{figure}[hpt]
\begin{center}

\includegraphics[scale=0.4]{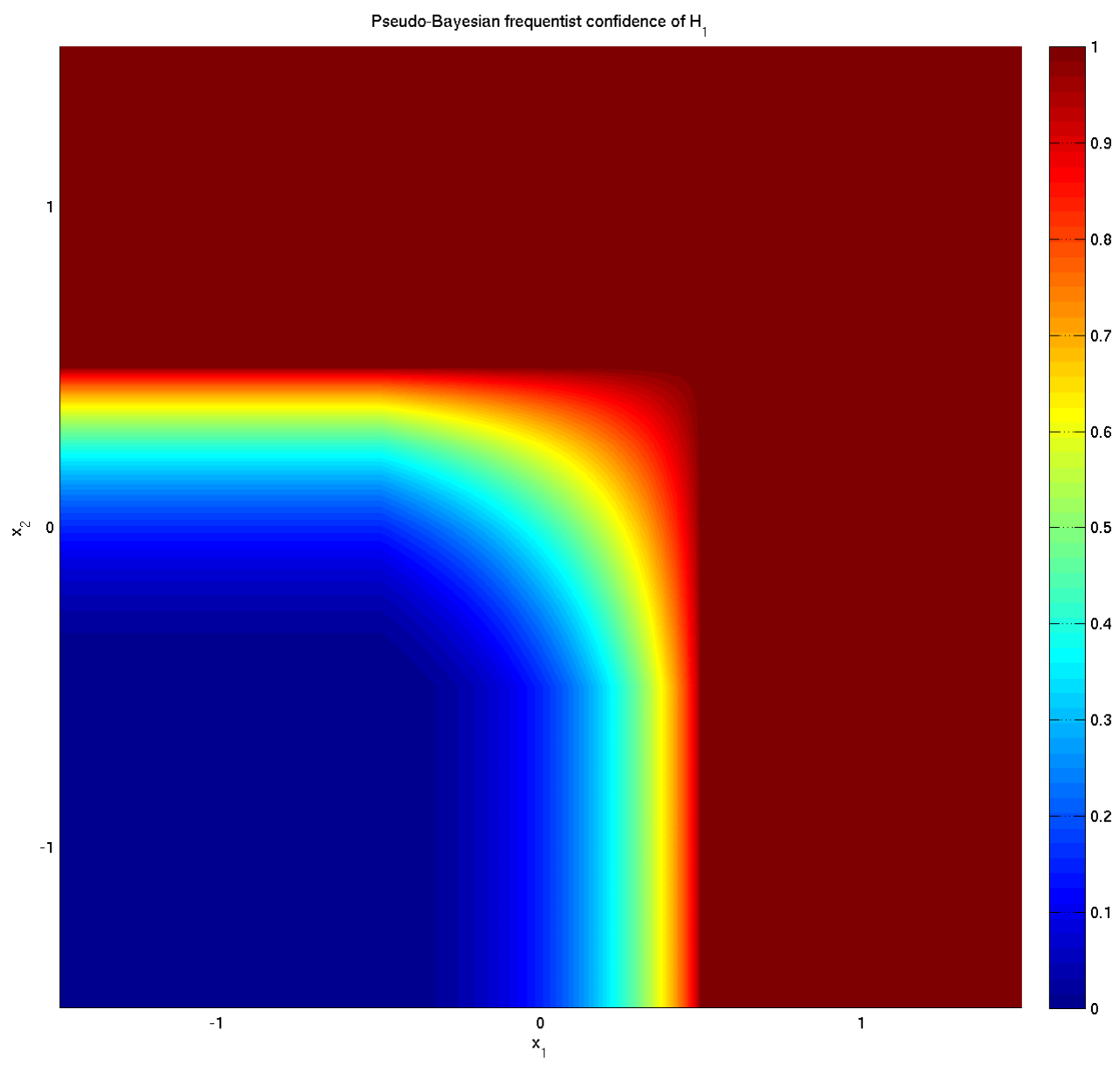}

\caption{
\label{simplepseudoBayesfig}
The frequentist confidence in $H_1$ resulting from this
pseudo-Bayesian solution for the problem of appendix
\ref{simpleexample} is plotted as a function of the observed data
$(x_1,x_2)$.  }

\end{center}
\end{figure}

Thus $$\eta(p)=\int_{-\frac{1}{2}}^{+\frac{1}{2}}{0\lor(\frac{1}{2}+\frac{1}{2}-\frac{1-p}{\frac{1}{2}-x_1})\,dx_1}.$$
Simplifying we get 
\begin{IEEEeqnarray*}{rCl}
\eta(p)&=&\int_{-\frac{1}{2}}^{p-\frac{1}{2}}{\left(1 -
  \frac{1-p}{\frac{1}{2}-x_1}\right)\,dx_1}\\
&=& p -
\int_{-\frac{1}{2}}^{p-\frac{1}{2}}{\frac{1-p}{\frac{1}{2}-x_1}\,dx_1}\\
&=& p -
(1-p)\left[-\log\left(\frac{1}{2}-x_1\right)\right]_{-\frac{1}{2}}^{p-\frac{1}{2}}\\
&=& p + (1-p)\log(1-p),
\end{IEEEeqnarray*}
where for $p=1$ we also have $\eta(p)=1$. This is shown in figure
\ref{simplepseudoBayesfig}. Since $\log(1-p)<0$ for $p>0$, we note
that this frequentist confidence is always lower than the Bayesian
posterior probability (for this solution of this problem).

\subsection{Discussion}

We have seen the Bayesian solution (figure \ref{simpleBayesfig})
corresponding to any of a range of uniform priors on the region of
interest, three correct frequentist solutions (figures
\ref{simplefreq1fig}, \ref{simplefreq2fig},
\ref{simplefreqlatefixedfig}), and one pseudo-Bayesian solution
(figure \ref{simplepseudoBayesfig}) (which is of course a frequentist
solution, not a Bayesian solution).

Of these, the frequentist ones are essentially arbitrary by choice of
their critical regions (though one could argue that there is something
special about the pseudo-Bayesian solution, there are in fact many
other possible pseudo-Bayesian solutions to choose from). 

Moreover, none of the frequentist solutions is uniformly optimal
(i.e. none is most likely to get any given degree of frequentist
confidence that $H_1$ holds when it does irrespective of the actual
value of $(\phi_1,\phi_2)\in H_1$). E.g. if
$(\phi_1,\phi_2)=(0.001,-1)\in H_1$ then with the first frequentist
solution the probability of getting frequentist confidence $>0.5$ is
$0.501$, greater than the zero probability achieved by the second
solution, the approx $1-\frac{1}{\sqrt{2}}\approx 0.293$ for the solution of figure
\ref{simplefreqlatefixedfig}, or the approx $0.187$ for the
pseudo-Bayesian solution. On the other hand for
$(\phi_1,\phi_2)=(-1,0.001)\in H_1$ it is the second frequentist
solution that has most chance of getting that level of frequentist
confidence. 

We would suggest that the Bayesian solution not only makes most sense,
but also gives symmetric posterior probability about the lines
separating $H_0$ and $H_1$ except near the corner at the origin, and
always gives greater or equal certainty for $H_1$ than any of the
frequentist solutions, all of which are biased against $H_1$ from the
Bayesian point of view. Moreover it avoids having to choose (before
collecting the data) which frequentist solution (among these and many
others) to use.

\section{Appendix: Density calculation for bullet arrival point}
\label{densitycalcs}

The base of the tower from which the bullet is fired is at
$(x,y)=(a,b)$ and the height of the tower is $z$.

We change the origin of our coordinate system to be the point from
which the bullet is fired, keeping the axes parallel to their original
directions, and labelling the new axes $u,v,w$; then the base of the
tower is at $(u,v,w)=(0,0,-z)$, and the $x,y$-coordinates of a point
$(u,v,-z)$ on the plane are given by $$x=u+a,$$ $$y=v+b.$$

First, it is evident that the bullet will hit the plane if and only if
the direction in which it is fired is below the horizontal, and that
therefore the probability of it not hitting the plane at all is
$\frac{1}{2}$. 

Now, a nice way of generating a random direction uniformly around the
sphere is draw a random sample from the 3-dimensional unit Gaussian
whose probability density is given by $$P(u,v,w) =
(2\pi)^{-\frac{3}{2}}e^{-\frac{1}{2}(u^2+v^2+w^2)}$$ and then divide
by the distance of the sample from the firing point, as this density
can also be written $$P(u,v,w) =
(2\pi)^{-\frac{3}{2}}e^{-\frac{1}{2}r^2}$$ where
$r=\sqrt{u^2+v^2+w^2}$, the distance from the origin, and the
probability of the sample being exactly at the origin is zero.

Therefore, if $(u,v,w)$ gives us the unnormalised direction in which
the bullet was fired, the point at which it hits the plane will be
at $$(-\frac{zu}{w}, -\frac{zv}{w}, -z),$$ i.e. at $$(x,y) =
(a-\frac{zu}{w}, b-\frac{zv}{w}),$$ so long as $w<0$.

We now aim to apply the standard formula for transforming a density
$P(u,v,w)$ to a density $P(x,y,w)$ (keeping the third coordinate in
its original form). To do that we first express the unnormalised
direction $(u,v,w)$ in terms of the point of impact and $w$, i.e. of
$(x,y,w)$, by the relationships
\begin{IEEEeqnarray*}{rCl}
u &=& -\frac{w}{z}(x-a)\\
v &=& -\frac{w}{z}(y-b)\\
w &=& w
\end{IEEEeqnarray*}
then calculate the matrix of partial derivatives
$$\left(\begin{matrix}\frac{\partial u}{\partial x} & \frac{\partial
    u}{\partial y} & \frac{\partial u}{\partial w}\\ \frac{\partial
    v}{\partial x} & \frac{\partial v}{\partial y} & \frac{\partial
    v}{\partial w}\\ \frac{\partial w}{\partial x} & \frac{\partial
    w}{\partial y} & \frac{\partial w}{\partial w}\end{matrix}\right)
= \left(\begin{matrix}-\frac{w}{z} & 0 & -\frac{x-a}{z}\\ 0 &
  -\frac{w}{z} & -\frac{y-b}{z}\\ 0 & 0 & 1
\end{matrix}\right)$$ which as a triangular matrix has determinant $\frac{w^2}{z^2}$. We then
deduce that $$P(x,y,w) =
\frac{w^2}{z^2}(2\pi)^{-\frac{3}{2}}e^{-\frac{1}{2}(\frac{w^2}{z^2}(x-a)^2
  + \frac{w^2}{z^2}(y-b)^2 + w^2)}.$$

All that remains now is to use the marginalisation rule of probability
to integrate out $w$, getting
\begin{IEEEeqnarray*}{rCl}
P(x,y) &=& \int_{-\infty}^0{P(x,y,w)\,dw}\\
&=& (2\pi)^{-\frac{3}{2}} z^{-2} \int_{-\infty}^0{w^2
  e^{-\frac{w^2}{2z^2}\left((x-a)^2+(y-b)^2+z^2\right)}\,dw}.
\end{IEEEeqnarray*}
We substitute $$s=\frac{w^2}{2z^2}\left((x-a)^2+(y-b)^2+z^2\right)$$
so that 
$$\frac{ds}{dw} = w \frac{(x-a)^2+(y-b)^2+z^2}{z^2}$$
and 
$$\frac{dw}{ds} = w^{-1}\frac{z^2}{(x-a)^2+(y-b)^2+z^2},$$ giving us
\begin{IEEEeqnarray*}{rCl}
P(x,y) &=&
\frac{1}{(2\pi)^{\frac{3}{2}}\left((x-a)^2+(y-b)^2+z^2\right)}\int_{\infty}^0{we^{-s}\,ds}\\
&=&
\frac{z\sqrt{2}}{(2\pi)^{\frac{3}{2}}\left((x-a)^2+(y-b)^2+z^2\right)^{\frac{3}{2}}}\int_0^{\infty}{s^{\frac{1}{2}}e^{-s}\,ds}
\end{IEEEeqnarray*}
since $w$ is the \textit{negative} sqare root of a positive multiple of
$s$. But $$\int_0^\infty{s^{\frac{1}{2}}e^{-s}\,ds} =
\int_0^\infty{s^{\frac{3}{2}-1}e^{-s}\,ds} =
\Gamma\left(\frac{3}{2}\right),$$ by definition of the Gamma function,
and $\Gamma\left(\frac{3}{2}\right)=\frac{\sqrt{\pi}}{2}$, so $$P(x,y)
= \frac{z}{4\pi\left((x-a)^2+(y-b)^2+z^2\right)^{\frac{3}{2}}}$$ as desired.

\section{Appendix: Bayesian calculations to go with section
  \ref{example3Bayes}}
\label{example3Bayescalcs}

We have treated respectively $N_1,N_2$ patients with torvaldomycin and
jobsucillin and observed $n_1,n_2$ cures.

We set an independent conjugate Beta prior on $(p_1,p_2)$ with
parameters $((\alpha_1,\beta_1),(\alpha_2,\beta_2))$ (with the special
case of uniform independent priors being
$((\alpha_1,\beta_1),(\alpha_2,\beta_2))=((1,1),(1,1))$ ), so
that $$P(p_1,p_2)=\frac{\Gamma(\alpha_1+\beta_1)}{\Gamma(\alpha_1)\Gamma(\beta_1)}p_1^{\alpha_1-1}(1-p_1)^{\beta_1-1}\frac{\Gamma(\alpha_2+\beta_2)}{\Gamma(\alpha_2)\Gamma(\beta_2)}p_2^{\alpha_2-1}(1-p_2)^{\beta_2-1}.$$

Then applying Bayes' theorem we
get 
\begin{IEEEeqnarray*}{rCl}
P(p_1,p_2|N_1,N_2,n_1,n_2)&\propto&
\frac{\Gamma(\alpha_1+\beta_1)}{\Gamma(\alpha_1)\Gamma(\beta_1)}p_1^{\alpha_1-1}(1-p_1)^{\beta_1-1}\frac{\Gamma(\alpha_2+\beta_2)}{\Gamma(\alpha_2)\Gamma(\beta_2)}p_2^{\alpha_2-1}(1-p_2)^{\beta_2-1}\\&&\times
\frac{N_1!N_2!}{n_1!n_2!(N_1-n_1)!(N_2-n_2)!}p_1^{n_1}p_2^{n_2}(1-p_1)^{N_1-n_1}(1-p_2)^{N_2-n_2}\\
&\propto&
p_1^{\alpha_1-1}(1-p_1)^{\beta_1-1}p_2^{\alpha_2-1}(1-p_2)^{\beta_2-1}p_1^{n_1}p_2^{n_2}(1-p_1)^{N_1-n_1}(1-p_2)^{N_2-n_2}\\
&=&p_1^{\alpha_1+n_1-1}p_2^{\alpha_2+n_2-1}(1-p_1)^{\beta_1+N_1-n_1-1}(1-p_2)^{\beta_2+N_2-n_2-1}
\end{IEEEeqnarray*}
where the constants of proportionality may depend on $N_1,N_2,n_1,n_2$
but do not vary with $p_1$ or $p_2$. Therefore 
\begin{IEEEeqnarray*}{rCl}
P(p_1>p_2|N_1,N_2,n_1,n_2)&=&\frac{\int_{p_1>p_2}{p_1^{\alpha_1+n_1-1}p_2^{\alpha_2+n_2-1}(1-p_1)^{\beta_1+N_1-n_1-1}(1-p_2)^{\beta_2+N_2-n_2-1}\,d(p_1,p_2)}}{\int_{[0,1]^2}{p_1^{\alpha_1+n_1-1}p_2^{\alpha_2+n_2-1}(1-p_1)^{\beta_1+N_1-n_1-1}(1-p_2)^{\beta_2+N_2-n_2-1}\,d(p_1,p_2)}}\\
&=&\frac{\Gamma(\alpha_1+\beta_1+N_1)\Gamma(\alpha_2+\beta_2+N_2)}{\Gamma(\alpha_1+n_1)\Gamma(\alpha_2+n_2)\Gamma(\beta_1+N_1-n_1)\Gamma(\beta_2+N_2-n_2)}\\
&&\times\int_{p_1>p_2}{p_1^{\alpha_1+n_1-1}p_2^{\alpha_2+n_2-1}(1-p_1)^{\beta_1+N_1-n_1-1}(1-p_2)^{\beta_2+N_2-n_2-1}\,d(p_1,p_2)}\\
&=&\frac{\Gamma(\alpha_1+\beta_1+N_1)\Gamma(\alpha_2+\beta_2+N_2)}{\Gamma(\alpha_1+n_1)\Gamma(\alpha_2+n_2)\Gamma(\beta_1+N_1-n_1)\Gamma(\beta_2+N_2-n_2)}\\
&&\times\int_0^1{\int_0^{p_1}{p_1^{\alpha_1+n_1-1}p_2^{\alpha_2+n_2-1}(1-p_1)^{\beta_1+N_1-n_1-1}(1-p_2)^{\beta_2+N_2-n_2-1}\,dp_2}\,dp_1}\\
&=&\frac{\Gamma(\alpha_1+\beta_1+N_1)}{\Gamma(\alpha_1+n_1)\Gamma(\beta_1+N_1-n_1)}\int_0^1{\beta_{\alpha_2+n_2,\beta_2+N_2-n_2}(p_1)\,dp_1},
\end{IEEEeqnarray*}
where $\beta_{a,b}(p)$ denotes the incomplete Beta function. This
final integral is best evaluated numerically, for example in Matlab,
although one can also evaluate the one in the line above numerically
instead.

However the reader is warned that direct numerical integration, or
even calculation of the integrand of the penultimate line above, will
not usually work due to numerical underflow and/or overflow. Instead
the calculations should be done by keeping the logarithm of each
value, and using the function $f(x)=\log(1+e^x)$ to do addition. Thus
if $x>y$ and $a=\log(x),b=\log(y),c=\log(x+y),d=\log(xy)$,
calculate $$c = a + f(b-a),$$ and obviously $$d=a+b.$$ Moreover
$\log\Gamma(x)$ should be evaluated without first evaluating
$\Gamma(x)$, e.g. by using the Matlab function gammaln().

\section{Appendix: Relationship between confidence sets and complete
  families of critical regions}
\label{hypoconf}

In order to see more generally that all the problems with frequentist
hypothesis testing using critical regions apply equally to using
frequentist confidence sets, we show here a bijective relationship
between \textit{complete families} of critical regions and confidence
set functions.

We must first make this definition:

A family $(C_{\eta, \theta})_{\eta\in [0,1],\theta\in \Theta}$ of
subsets of $X$ is a \textbf{complete family of critical regions} if
for all $(\theta,\phi)\in H$ and all $\eta,\eta_1,\eta_2\in[0,1]$ we
have:

\begin{itemize}

\item $\eta_1 < \eta_2 \implies C_{\eta_2,\theta}\subseteq C_{\eta_1,\theta}$; and

\item $P(x\in C_{\eta,\theta} | \theta,\phi) \leq 1-\eta$.

\end{itemize}

Intuitively, for each $\theta_0\in \Theta$, the family
$(C_{\eta,\theta_0})_{\eta\in[0,1]}$ is an ordinary family of critical
regions for the null hypothesis that $\theta=\theta_0$.

Note that if we have only been given an ordinary family $(C_\eta)$ of
critical regions (for the null hypothesis $\theta=\theta_0$), we can
immediately construct a complete family containing it by setting, for
each $\theta\neq \theta_0$ and all $\eta$,
$C_{\eta,\theta}=\emptyset$. There may of course be many other such
complete families. In particular this applies in the case that
$\Theta=\{0,1\}$ and $\theta_0=0$.

Then for any complete family of critical regions $(C_{\eta,\theta})$
we can define a frequentist confidence set function $f$ by $$f:X\times
[0,1] \to \mathbb{S}(\Theta): (x,\eta) \mapsto
\{\theta\in\Theta:x\notin C_{\eta,\theta}\}.$$

It is trivial to show that the mapping from the set of complete
families of critical regions to the set of frequentist confidence set
functions is bijective, with the inverse giving $(C_{\eta,\theta})$ in
terms of $f$ by $$C_{\eta,\theta} = \{x\in X:\theta\notin f(x,\eta)\}.$$

By way of example of translating one of the counter-examples of
appendix \ref{freqexamples} into one involving confidence sets, take
the example given there to criterion \ref{complementarity}
(``Complementarity''). The family of critical regions for $\hat{H}_0$
gives us half of a complete family, namely the $(C_{\eta,0})$, while
that for $\hat{H}'_0$ gives us the $(C_{\eta,1})$. Translating the
resulting $(C_{\eta,\theta})$ to a frequentist confidence set function
gives us a function $f$ such that $f(0,0.95)=\emptyset$, thus saying
that we are 95\% sure that $\theta\in\emptyset$, i.e. that both
$\theta\neq 0$ and $\theta\neq 1$, i.e. that both $\theta=1$ and
$\theta=0$, which is obviously nonsense and specifically violates
criterion \ref{complementarity}.

\bibliography{ms}
\bibliographystyle{ieeetr}

\end{document}